\documentclass[11pt]{amsart}

\usepackage{amsmath}
\usepackage{amssymb}
\usepackage{amscd}
\usepackage{latexsym}
\usepackage[latin1]{inputenc}

\usepackage{graphicx}
\usepackage{changebar,color}

    \def\P{{\mathbb P}}
    \def\N{{\mathbb N}}
    \def\Z{{\mathbb Z}}
    \def\Q{{\mathbb Q}}
    \def\R{{\mathbb R}}
    \def\C{{\mathbb C}}

    \def\H{{\mathbb H}}
    \def\A{{\mathbb A}}
    \def\bp{{\mathbf p}}
    
    \def\OK{{{\mathcal O}_K}}
    
    \renewcommand{\Re}{\,{\rm Re}\,}
    \renewcommand{\Im}{\,{\rm Im}\,}

    \newcommand{\Spec}{{\rm Spec}\,}
    
    \newcommand{\Hom}{{\rm Hom}}
    \newcommand{\ihom}{{\mathcal Hom}_{{\mathcal O}_X}}
    \newcommand{\iho}{\ihom (F,G)}
     \newcommand{\ihob}{F^\lor \otimes G}

    \newcommand{\Coker}{{\rm Coker}}
    \newcommand{\Oker}{{\rm Coker}\,}
    \newcommand{\Ker}{{\rm Ker}}
    \newcommand{\Ext}{\mbox{\rm Ext}}
    \newcommand{\Exthat}{\widehat{\rm Ext}}
    \newcommand{\id}{\mbox{id}}
    \newcommand{\ol}[1]{\overline{#1}}
    \newcommand{\refg}[1]{{(\ref{#1})}}
    \newcommand{\XR}{{X_\R}}
    \newcommand{\an}{{\rm hol}}
    \newcommand{\CinfX}{\mathcal{C}^\infty_{X_\Sigma}}
    \newcommand{\OanX}{{\mathcal O}^\an_{X_\Sigma}}
    \newcommand{\dolb}{{\overline{\partial}}}
    \newcommand{\ubar}{{\underline{\,\,\,}}}
    
    \newcommand{\size}{{\mathfrak s}}
    
    \newcommand{\rk}{{\rm{rk}}\,}
    \newcommand{\cf}{\emph{cf.}\,}
    \newcommand{\ie}{\emph{i.e.}\,}
    \newcommand{\eg}{\emph{e.g.}\,}
    \newcommand{\dega}{\widehat{\rm deg}\,}
    \newcommand{\degan}{\widehat{\rm deg}_{n}\,}
    \newcommand{\udega}{{\rm u}\widehat{\rm deg}_{n}\,}
     \newcommand{\udegan}{{\rm u}\widehat{\rm deg}_{n}\,}
    \newcommand{\muamax}{\widehat{\mu}_{\rm max}}
    \newcommand{\muamin}{\widehat{\mu}_{\rm min}}
    \newcommand{\mua}{\widehat{\mu}}
    \newcommand{\phit}{\tilde{\phi}}
    \newcommand{\at}{\tilde{a}}
    \newcommand{\pr}{\text{\upshape{pr}}\,}
    \newcommand{\condp}[3]{${{\mathbf P}(#2/#1,\,#3})$}
    \newcommand{\Vs}{$\mathcal V$-small }
    \newcommand{\cV}{{\mathcal V}}
    \newcommand{\FS}{\mathfrak{S}}
    \newcommand{\FP}{\mathfrak{P}}
    \newcommand{\pp}{\mathbf{p}}
      \newcommand{\qq}{\mathbf{q}}

    \newcommand{\Co}{{\dot{C}}}
\newcommand{\Eo}{{\dot{E}}}
\newcommand{\Fo}{{\dot{F}}}
\newcommand{\Go}{{\dot{G}}}

\newcommand{\DOX}{D({\mathcal O}_X\mbox{\scriptsize\rm -mod})}

\def\Oc{{\hat{\mathcal O}}}
\def\Kc{{\hat{K}}}
\def\m{{\mathfrak{m}}}

\newcommand{\cX}{{\mathcal X}}
\newcommand{\cA}{{\mathcal A}}
\newcommand{\cC}{{\mathcal C}}
\newcommand{\cE}{{\mathcal E}}
\newcommand{\cI}{{\mathcal I}}
\newcommand{\cAmod}{{\mathcal A}-\mathbf{mod}}
\newcommand{\cXmod}{{\mathcal O}_{\mathcal X}-\mathbf{mod}}
\newcommand{\cXqc}{{\mathcal O}_{\mathcal X}-\mathbf{qc}}
\newcommand{\HD}{H_{\rm{Dolb}}}
\newcommand{\cl}{{\rm cl}}
\newcommand{\clar}{\widehat{{\rm cl}}}

\newcommand{\Triv}{{\rm{Triv}}}

\newtheorem{theorem}[subsubsection]{Theorem}
\newtheorem{corollary}[subsubsection]{Corollary}
\newtheorem{lemma}[subsubsection]{Lemma}

\newtheorem{proposition}[subsubsection]{Proposition}

\newtheorem{example}[subsubsection]{Example}
\newtheorem{question}[subsubsection]{Question}

\numberwithin{equation}{section}

\begin{document}

\date{\today}

\title[Extensions on arithmetic schemes I]
{Hermitian vector bundles and extension groups on arithmetic schemes.
I.
Geometry of numbers}

\author{Jean-Beno{î}t Bost}
\address{J.-B. Bost, D{é}partement de Math{é}matiques, Universit{é}
Paris-Sud,
B{â}timent 425, 91405 Orsay cedex, France}
\email{jean-benoit.bost@math.u-psud.fr}
\author{Klaus K{ü}nnemann}
\address{K. K{ü}nnemann, Mathematik, 
Universit{ä}t Regensburg, 93040 Regensburg,
Germany}
\email{klaus.kuennemann@mathematik.uni-regensburg.de}

\begin{abstract}
We define and investigate extension groups in the context of Arakelov
geometry.
The ``arithmetic extension groups" $\widehat{\rm Ext}^i_{X}(F,G)$ we
introduce are extensions by groups of analytic types of the usual
extension groups ${\rm Ext}^i_{X}(F,G)$
attached to ${\mathcal O}_X$-modules $F$ and $G$ over an arithmetic scheme
$X$.
In this paper, we focus on the first
arithmetic extension group $\widehat{\rm Ext}^1_{X}(F,G)$ --- the
elements of which may be described in terms of admissible short exact
sequences of hermitian vector bundles over $X$ --- and we especially
consider the case 
when 
$X$ is an ``arithmetic curve", namely the spectrum $\Spec {\mathcal O}_K$ of the
ring of integers in some number
field $K$. Then the study of arithmetic extensions over $X$ is
related to
old and new problems concerning lattices and the geometry of numbers.

Namely, for any two hermitian vector bundles $\ol{F}$ and
$\ol{G}$ over $X:=\Spec \OK$, we attach a logarithmic \emph{size}
$\size_{\ol{F}, \ol{G}}(\alpha)$ to any element $\alpha$ of
$\widehat{\rm
Ext}^1_{X}(F,G)$, and we give an upper bound on $\size_{\ol{F},
\ol{G}}(\alpha)$
 in terms of
slope invariants of $\ol{F}$ and $\ol{G}$. We further illustrate this
notion
by relating the sizes of restrictions to points in $\P^1(\Z)$ of the
universal extension over $\P^1_{\Z}$ to the geometry of $PSL_{2}(\Z)$ 
acting on Poincar\'e's upper half-plane, and by deducing some
quantitative results in reduction theory from our previous
upper bound 
on sizes. Finally, we investigate the behaviour of size by base change
(\ie,  under extension of the ground field $K$ to a larger number field
$K'$): when the base field
$K$ is $\Q$, we establish that the size, which cannot increase under
base
change, is actually invariant when the field $K'$ is an abelian
extension of $K$, or
when $\ol{F}^\lor \otimes \ol{G}$ is a direct sum of root lattices
and of lattices of Voronoi's first kind.

The appendices contain  results concerning extensions in categories
of
sheaves on ringed spaces, and lattices of Voronoi's first kind
which might also be of independent interest.

\bigskip

\end{abstract}

\date{January 11,2007}
\maketitle

\setcounter{tocdepth}{2}
{
\tableofcontents
}

\setcounter{section}{-1}
\section{Introduction}
The aim of this paper is to introduce and to study \emph{arithmetic
extensions} and the \emph{extension groups} they define in the
framework
of Arakelov geometry.

\medskip

\noindent {\bf 0.1.}
Arithmetic extensions are objects which arise naturally at
various  places in arithmetic geometry.
Let $X$ be an arithmetic scheme -- namely a separated scheme of
finite type over $\Z$ -- such that $X_{\C}$ is smooth, and let
$X(\C)$ be
the complex manifold of its complex points.
By definition, for any two locally free coherent ${\mathcal O}_{X}$-modules $F$
and $G$,
an arithmetic extension $(\mathcal{E},s)$ of $F$ by $G$ is given by an
extension of ${\mathcal O}_X$-modules
\[
\mathcal{E}:\,0\longrightarrow  G\longrightarrow  E\longrightarrow
F\longrightarrow  0
\]
together with a $\mathcal{C}^\infty$-splitting over $X(\C)$
$$s: F_{\C} \longrightarrow E_{\C},$$
invariant under
complex conjugation, of the extension of complex vector bundles over
$X(\C)$
\[
\mathcal{E}_{\C}:\,0\longrightarrow  G_{\C}\longrightarrow  E_{\C}
\longrightarrow  F_{\C}\longrightarrow  0
\]
deduced from $\mathcal{E}$ by extending the scalars from $\Z$ to
$\C$. 

Recall that an \emph{hermitian vector bundle }
$\ol{V}:=(V,\|.\|)$
over $X$ is the data of a locally free coherent sheaf $V$ over $X$,
together with a $\mathcal{C}^\infty$-hermitian metric $\|.\|$ on the
attached vector bundle $V_{\C}$ on $X(\C)$ that is invariant under
complex conjugation. Arithmetic extensions arise for instance from
\emph{admissible extensions}
\begin{equation}\label{admissex}
\ol{\mathcal{E}}:\,0\longrightarrow  \ol{G}\longrightarrow  \ol{E}
\longrightarrow
\ol{F}\longrightarrow  0,
\end{equation}
of hermitian vector
bundles over $X,$ namely from the
data of an extension
\[
\mathcal{E}:\,0\longrightarrow  G\longrightarrow  E\longrightarrow  F
\longrightarrow  0
\]
of the underlying ${\mathcal O}_{X}$-modules such that the hermitian metrics
$\|.\|_{\ol{F}}$ and $\|.\|_{\ol{G}}$ on $F_{\C}$ and $G_{\C}$ are
induced (by restriction and quotients) by the metric $\|.\|_{\ol{E}}$
on
$E_{\C}.$ In this case, orthogonal projection determines a
$\mathcal{C}^\infty$-splitting $s^\perp:F_{\C}\to E_{\C}$ of
$\mathcal{E}_{\C},$ and $(\mathcal{E},s^\perp)$ is an
arithmetic extensions of $F$ by $G.$

It turns out that, by means of the Baer sum construction, one may
define an addition law on the set $\widehat{\rm Ext}^1_{X}(F,G)$ of
isomorphism classes of arithmetic extensions of $F$ by $G$, which in
this way
is endowed with a natural structure of an abelian group.
Moreover, in analogy to the arithmetic Chow groups, the arithmetic
extension
group
$
\widehat{\rm Ext}_X^1(F,G)
$
is an extension of the ``classical" extension
group ${\rm Ext}_{{\mathcal O}_X}^1(F,G)$, defined in the context of sheaves of
${\mathcal O}_X$-modules, by a group of analytic type.
More precisely, it fits into an exact sequence
\begin{equation}\label{basicexact}
{\rm Hom}_{{\mathcal O}_X}(F,G)
{\longrightarrow }
{\rm Hom}_{\mathcal{C}^\infty_{X}}(F_\C,G_\C)^{F_\infty}
\stackrel{b}{\longrightarrow }
\widehat{\rm Ext}^1_{X}(F,G)
\stackrel{\nu}{\longrightarrow}  {\rm Ext}_{{\mathcal O}_X}^1(F,G)
\longrightarrow 0,
\end{equation}
where $F_\infty$ acts on $X(\C),$ $F_{\C},$ and 
$G_{\C}$ by complex conjugation. 
We  may also define an homomorphism
\[\Psi:\widehat{\rm Ext}_{X}^1(F,G)\longrightarrow
Z_{\overline{\partial}}^{0,1}(X_\R,\ihob) \]
to the group $Z_{\overline{\partial}}^{0,1}(X_\R,\ihob)$ of
$F_{\infty}$-invariant
$\ol{\partial}$-closed forms of type $(0,1)$ on $X(\C)$ with
coefficients in $F^\lor_\C\otimes G_\C$, by sending the class of an
arithmetic
extension $(\mathcal{E},s)$ to its ``second 
fundamental form"  $\ol{\partial}s.$

The arithmetic extension group $\widehat{\rm Ext}^1(F,G)$ actually
admits an interpretation
in terms of homological algebra, in the spirit of the well-known
identification of the 
``classical" extension group ${\rm Ext}^1_X(F,G)$, originally defined
by classes of 1-extensions
equipped with the Baer sum, with the ``cohomological" extension group
$\Hom_{D({{\mathcal O}_{X}-\mathbf{mod}})}(F, G[1]),$ defined  as a
group
 of morphisms  
in the derived category of (sheaves of) ${\mathcal O}_X$-modules over $X$ .
Indeed, if $(X_\R, {\mathcal C}_\R^\infty)$ denotes the ringed space
quotient of $(X(\C), {\mathcal C}^\infty_{X(\C)})$ by the action of complex
conjugation (acting both on $X(\C)$ and on values of 
$C^\infty$-functions), and if $$\rho: (X_\R, {\mathcal C}_\R^\infty)
\longrightarrow (X, {\mathcal O}_X)$$ is the natural map of ringed spaces,
then, for any ${\mathcal O}_X$-module $G$ on $X$, we may consider the
adjunction map
\begin{equation}\label{adG}
{\rm ad}_G: G \longrightarrow \rho_\ast \rho^\ast G
\end{equation}
--- it maps any local section $g$ of $G$ to the  section $g_\C$, seen 
as a $\mathcal{C}^\infty$-section of $G_\C$, invariant under the complex 
conjugation $F_\infty$ --- 
 and its cone $C({\rm ad}_G),$ 
namely (\ref{adG}) seen as
 complex of length $2$, with $G$ (resp. $\rho_\ast \rho^\ast G$)
sitting in degree $-1$ 
 (resp. $0$). Then, for any two locally free coherent sheaves 
 $F$ and $G$ on $X,$ we have a natural isomorphism of abelian groups
:
 $$\widehat{\rm Ext}^1_X(F,G) \stackrel{\sim}{\longrightarrow}
 \Hom_{D({{\mathcal O}_{X}-\mathbf{mod}})}\bigl(F,C({\rm ad}_G)\bigr)$$
 between our arithmetic extension group and the group of
morphisms from $F$ to $C({\rm ad}_G)$ in the derived category 
$D({{\mathcal O}_{X}-\mathbf{mod}})$ 
of the abelian category of sheaves of ${\mathcal O}_{X}$-modules over $X$
(see 
  \ref{derived} \emph{infra}).

In a forthcoming part of this work, the above cohomological
interpretation of $\widehat{\rm Ext}^1(F,G)$ will motivate us to
consider higher arithmetic extension groups $\widehat{\rm
Ext}^i(F,G)$, defined for any integer $i \geq 1$ by means of the
Dolbeault complex $(A^{0,.}_{X(\C)},\ol{\partial})$ on $X(\C)$ and its
subcomplex $(A^{0,.}_{X_{\R}},\ol{\partial})$ of conjugation
invariant forms, which defines a complex of sheaves of modules on the 
ringed space $(X_{\R}, \mathcal{C}^\infty_{\R})$. 

For any $\mathcal{C}^\infty_{\R}$-module $F$ on $X_\R,$ we 
get the ``Dolbeault resolution"  ${\mathcal Dolb}(F)$ of $F$ by applying 
the functor $F\otimes_{ \mathcal{C}^\infty_{\R}}.$ to this complex.
In particular, for any sheaf $G$ of ${\mathcal O}_X$-modules, we may consider 
the associated sheaf $\rho^\ast G$ of $\mathcal{C}^\infty_{\R}$-modules 
over $X_\R$, and the naive truncation ${\mathcal Dolb}(\rho^\ast
G)_{\leq i-1}$ of its Dolbeault resolution.
The adjunction map (\ref{adG}) extends to a morphism of complexes
$${\rm ad}_{G}^{i-1}: G \rightarrow \rho_{\ast}({\mathcal Dolb}(\rho^\ast
G)_{\leq i-1}),$$
and its cone $C({\rm
ad}_G^{i-1})$ is a complex of (sheaves of) ${\mathcal O}_X$-modules. 

For any two  ${\mathcal O}_X$-modules $F$ and $G$, 
 we shall define
$$\widehat{\rm Ext}^i_X(F,G) := 
\Hom_{D({{\mathcal O}_{X}-\mathbf{mod}})}\bigl(F,C({\rm
ad}_G^{i-1})[i-1]\bigr).$$
This group may be interpreted as an ``hyper-extension group":
$$\Hom_{D({{\mathcal O}_{X}-\mathbf{mod}})}\bigl(F,C({\rm
ad}_G^{i-1})[i-1]\bigr)
\simeq
{\rm Ext}^i_X\bigl(F, C({\rm
ad}_G^{i-1})[-1]\bigr),$$
where, by the very definitions of the Dolbeault resolution and of a cone, 
the ``shifted cone" $C({\rm ad}_G^{i-1})[-1]$ is the following complex of
length $i+1$ of sheaves of ${\mathcal O}_X$-modules, with $G$ sitting 
in degree $0$:
\[
0\longrightarrow G\stackrel{-{\rm ad}_G}{ \longrightarrow} \rho_{\ast}\rho^{\ast} G
\stackrel{-\ol{\partial}_G}{\longrightarrow}
\rho_{\ast}(\rho^{\ast} G \otimes_{\mathcal{C}^\infty_{\R}} A^{0,1}_{X_{\R}} )
\stackrel{-\ol{\partial}_G}{\longrightarrow} \cdots
\stackrel{-\ol{\partial}_G}{\longrightarrow} \rho_{\ast}(\rho^{\ast} G \otimes_{\mathcal{C}^\infty_{\R}}
A^{0,i-1}_{X_{\R}}) \longrightarrow 0.
\]

\medskip

\noindent {\bf 0.2.}
Classical constructions in algebraic and differential
geometry provide natural instances of admissible and arithmetic
extensions. In the second part of this paper
\cite{bostkuennemann2}, we shall discuss three of these constructions,
which give rise  to the
  \emph{arithmetic Atiyah extension}, the
\emph{arithmetic Hodge extension}, and the \emph{arithmetic Schwarz
extension}. To advocate the investigation of arithmetic extensions, we
want to indicate briefly their constructions:

(i) Let $\overline{E}$ be an hermitian vector bundle on an arithmetic
scheme
$X$.
The bundle of $1$-jets of $E$ induces an extension of
${\mathcal O}_{X}$-modules
\[
0\longrightarrow  \Omega^1_{X/\Z}\otimes E\longrightarrow
\mathcal{J}_{X/\Z}^1(E)\longrightarrow  E\longrightarrow  0,
\]
the Atiyah extension of $E$. 
The holomorphic vector bundle $E_\C$ carries a unique
$\mathcal{C}^\infty$-connection  which is compatible with the metric
and
the complex structure, its so-called Chern 
connection, which induces a 
$\mathcal{C}^\infty$-splitting $s$ of the Atiyah
extension and yields a canonical arithmetic extension class
\[
\widehat{\rm at}(\ol{E})\in\widehat{\rm
   Ext}^1_X(E,\Omega_{X/\Z}^1\otimes E).
\]
It is a refinement both of the algebraic Atiyah class ${\rm at}(E) =
\nu(\widehat{\rm at}(\ol{E}))$ in ${\rm
   Ext}^1(E,\Omega_{X/\Z}^1\otimes E)$ and of the curvature form of
the
   Chern connection of $\ol{E}_\C$, which coincides with
$\Psi(\widehat{\rm
   at}(\ol{E}))$ (up to some normalization factor). Applying a trace map 
   to  $\widehat{\rm at}(\ol{E})$,
   we get an arithmetic first Chern class in ``arithmetic Hodge
   cohomology":
   $$\hat{c}^H_1(\ol{E}) \in  \widehat{\rm
   Ext}^1_{X}({\mathcal O}_X,\Omega_{X/\Z}^1).$$

(ii) Let $f:X\to Y$ be a smooth
proper
morphism of arithmetic schemes such
that the Hodge to de Rham spectral sequence
$E_1^{p,q}=R^pf_*\Omega_{X/Y}^q\Rightarrow
R^{p+q}f_*\Omega_{X/Y}^\cdot$
degenerates at $E_1$.
The spectral sequence defines the so-called Hodge extension
\begin{equation}\label{Hodge}
{\rm Hdg(X/Y)}:
0\longrightarrow  f_*\Omega_{X/Y}^1\longrightarrow
R^1f_*\Omega_{X/Y}^\cdot\longrightarrow  R^1f_*{\mathcal O}_X\longrightarrow  0
\end{equation}
whose interest was already advocated by Grothendieck in
\cite{grothendieck66}.
Complex Hodge theory 
equips Hdg$(X/Y)$ with a canonical structure of an arithmetic
extension\footnote{Namely, the  vector bundle
over $Y(\C)$ defined by the
relative algebraic de Rham cohomology $R^1f_*\Omega_{X/Y}^\bullet$ 
may be identified with the relative first Betti cohomology with complex
coefficients of $X(\C)/Y(\C)$; the complex conjugation on coefficients acts on
Betti cohomology and maps the $\C$-analytic sub-vector bundle
$(f_*\Omega_{X/Y}^1)_{\C}$ of 
$(R^1f_*\Omega_{X/Y}^\bullet)_{\C}$
onto a ${\mathcal C}^\infty$ direct summand of
$(f_*\Omega_{X/Y}^1)_{\C}$,  which
provides a ${\mathcal C}^\infty$-splitting of the extension of $\C$-analytic
vector bundles over $Y(\C)$ defined by Hdg$(X/Y)$.}.
We thus obtain the class of the arithmetic Hodge extension
\[
\widehat{\rm Hdg}(X/Y)\in \widehat{\rm
Ext}_{Y}^1(R^1f_*{\mathcal O}_X,f_*\Omega_{X/Y}^1).
\]

(iii) Let $f:C\to X$ be a smooth, projective curve of genus $g\geq 2$
over an
arithmetic scheme $X$.
Using Deligne's definition of the torsor of projective connections
on relative curves in \cite{deligne70}, I.5\footnote{Strictly
speaking, the definition in \emph{loc. cit.} is stated in the
framework of complex analytic spaces. However, it is formulated in a
general geometric language, which makes it meaningful in the context
of smooth relative curves over an arbitrary scheme.}, one obtains
a canonical extension of ${\mathcal O}_{X}$-modules
\[
\mathcal{S}_{C/X}\,:\,0\longrightarrow  f_*\Omega_{C/X}^{\otimes 2}
\longrightarrow  S_{C/X}
\longrightarrow  {\mathcal O}_X\longrightarrow  0,
\]
the splittings of which correspond to projective connections on $C/X$.
Complex uniformization by the upper half-plane $\H$ induces a
${\mathcal C}^\infty$ projective
connection on $C(\C)/X(\C)$ that is holomorphic along the fibers --- 
hence a ${\mathcal C}^\infty$-splitting of
$\mathcal{S}_{C/X}$ over $Y(\C)$ --- and allows one to define from
$\mathcal{S}_{C/X}$
the arithmetic Schwarz extension and its class
\[
\widehat{\mathcal{S}}_{C/X}\in
\widehat{\rm Ext}^1_{X}({\mathcal O}_X,f_*\Omega_{C/X}^{\otimes 2}).
\]

The non-vanishing of each of the above classes $\hat{c}^H_1(\ol{E}),$
$\widehat{\rm Hdg}(X/Y),$ or $\widehat{\mathcal{S}}_{C/X}$ is an
intriguing issue, related to deep problems in Diophantine geometry
and transcendence theory.

\medskip

\noindent {\bf 0.3.}
In this paper, after introducing the arithmetic extension
groups $\widehat{\rm Ext}^1_{X}(F,G)$ and discussing their basic
properties in Section 2, we concentrate on the case where $X$
is an ``arithmetic curve", namely the spectrum $\Spec {\mathcal O}_K$ of the
ring of integers in some number
field $K$. It turns out that the study of arithmetic extensions over
$X$ is related to
old and new problems concerning lattices and the geometry of numbers.

Namely, if $F$ and $G$ are vector bundles over $X:=\Spec {\mathcal O}_{K}$ (\ie,
projective
${\mathcal O}_K$-modules), we obtain from the basic exact sequence
(\ref{basicexact}) 
a canonical isomorphism
\begin{equation}
    \widehat{\rm Ext}^1_{X}(F,G) \simeq
    \frac{\Hom_{{\mathcal O}_{K}}(F,G)\otimes_\Z\R}{\Hom_{{\mathcal O}_K}{(F,G)}}.
    \label{eq:exttor}
\end{equation}
Consequently the arithmetic extension group $\widehat{\rm
Ext}_X^1(F,G)$
carries a canonical structure of a real torus.
Moreover, if $F$ and $G$ are equipped with hermitian metrics, which
makes them hermitian vector bundles $\ol{F}$ and  $\ol{G}$, we get an
induced
Riemannian metric on this real torus.
In Section 3, we define the \emph{size} 
$\size_{\ol{F},\ol{G}}(\mathcal{E},s)$ of an
arithmetic extension of
  $F$ by $G$ as the logarithm (in $[-\infty, + 
\infty[$) of the distance to zero of the
corresponding point in the torus (\ref{eq:exttor}).
Let $\ol{\mathcal{E}}$ be an admissible extension \refg{admissex}
with associated arithmetic extension $(\mathcal{E},s^\perp)$ as above,
and let 
\[
\varphi\,:E\stackrel{\sim}{\longrightarrow} G\oplus F
\]
be an
isomorphism of ${\mathcal O}_K$-modules
compatible with the extension $\mathcal{E}$ (that is, such that
$\varphi^{-1}\circ (Id_{G},0) : G \rightarrow E$ and ${\rm pr}_{2}\circ
\varphi : E\rightarrow F$ coincide with the morphisms defining
$\mathcal{E}$). Then 
\[
{\frac{1}{[K:\Q]}
\sum\limits_{\sigma:K\hookrightarrow \C}\|\varphi_{\sigma}
\|_{\ol{E}^\lor\otimes (\ol{G}\oplus \ol{F}),\sigma}^2} =
{\frac{1}{[K:\Q]}
\sum\limits_{\sigma:K\hookrightarrow \C}\|\varphi^{-1}_{\sigma}
\|_{(\ol{G}\oplus \ol{F})^\lor\otimes\ol{E},\sigma}^2}
\geq {\rk_{{\mathcal O}_K}E},
\]
and 
the minimum value achieved by the left-hand side when $\varphi$ runs
over all the isomorphisms of $\OK$-modules as above is precisely
$$\rk_{{\mathcal O}_K}E + \exp(\,2\size_{\ol{F},\ol{G}}(\mathcal{E},s^\perp))$$
(see Proposition \ref{proposition.splittingsize} and 
Corollary \ref{cor.splittingsize}
\emph{infra}).

Motivated by analogous results concerning vector bundles on
projective curves over a field, we show that the size of arithmetic
extensions
satisfies the following upper bound:
\begin{equation}\label{fuest}
\size_{\ol{F},\ol{G}}(\mathcal{E},s)\leq
\muamax (\ol{F})-\muamin (\ol{G})
+\frac{\log |\Delta_{K}|}{[K:\Q]}
+ \log\,\frac{\rk_K F_K\cdot\rk_K G_K}{2},
\end{equation}
where $\widehat{\mu}_{\rm max}(\ol{F})$ and $\widehat{\mu}_{\rm
min}(\ol{G})$  denote the maximal and minimal normalized
slopes
of $\ol{F}$ and $\ol{G}$ (see \ref{Adeg},
\emph{infra}), and  $\Delta_{K}$ the discriminant of the number field
$K$.
To establish (\ref{fuest}), we rely on (i) some upper bound on
the Arakelov degree of a sub-line bundle in
the tensor product of two hermitian vector bundles over $\Spec {\mathcal O}_K$,
and (ii)  some ``transference
theorem" from the geometry of numbers, which relates the inhomogeneous
minimum (also called the covering radius) of a lattice in a euclidean
vector space to the first of the successive minima of the dual
lattice.

Section 4 is devoted to further examples and applications of the
notion of size. In particular, using the 
inequality (\ref{fuest}), we derive an avatar, in 
the
framework of Arakelov geometry over arithmetic curves, of the main
result of the classical reduction theory of positive quadratic forms.
It claims the existence of some
``almost-splitting" for any hermitian vector bundle $\ol{E}$ over
$\Spec {\mathcal O}_{K}$, namely the existence of $n:= \rk E$ hermitian lines
bundles
$\ol{L}_{1},$ \ldots, $\ol{L}_{n}$ over $\Spec {\mathcal O}_K$, and of an
isomorphism
of ${\mathcal O}_{K}$-modules
\[
\phi: E \stackrel{\sim}{\longrightarrow} \bigoplus_{i=1}^{n}L_{i}
\]
such that the archimedean norms of $\phi$ and $\phi^{-1}$, computed by
using the hermitian structures on $\ol{E}$ and on the orthogonal
direct sum $\bigoplus_{i=1}^{n}\ol{L}_{i},$
are bounded in terms of $K$ and $n$ only (Theorem \ref{thm.reducArak}
\emph{infra}).

Besides, for any rational point $P\in \P^1(\Q)=\P^1(\Z)$, we
calculate the
size of the inverse image
$P^*\ol{\mathcal{E}}$ of the universal extension
\[
\ol{\mathcal{E}}:0\longrightarrow 
\ol{S}\longrightarrow  \ol{{\mathcal O}_X^{\oplus 
2}}\longrightarrow
\ol{{\mathcal O}_X(1)}\longrightarrow  0
\]
over the projective line $X=\P^1_\Z$ equipped with its natural
structure of
an admissible extension.
The extension class of $P^*\ol{\mathcal{E}}$ is trivial iff
$P\in\{0,\infty\}$.
For $P\in\P^1(\Q)\setminus \{0,\infty\}$, we show that the
size of $P^*\ol{\mathcal{E}}$
is related to the usual height $h(P)$ of $P$ by the inequalities
\[
-\frac{1}{2}\log\,2+h(P)\leq \size(P^*\ol{\mathcal{E}})
\leq -\log\,2+2\,h(P).
\]
We also give a geometric description of
the size $\size(P^*\ol{\mathcal{E}})$ by means of so-called Ford
circles (namely the images under elements in $SL_2(\Z)$ of the
horocycles $\{ \Im z = 1 \}$ in Poincar\'{e}'s upper half-plane).

The final section of this paper is devoted to the intriguing question
of the invariance of size under base change. Recall  that  an
extension of number fields $K'/K$ defines a morphism
$g:\Spec{\mathcal O}_{K'}\to \Spec{\mathcal O}_K$
of ``arithmetic curves".
For hermitian vector bundles $\ol{F}$ and 
$\ol{G}$ over $S$, there is an induced
morphism
\[
g^*:\widehat{\rm Ext}_S^1(F,G)
\longrightarrow  \widehat{\rm Ext}^1_{S'}(g^*F,g^*G).
\]
It is easy to see that the inequality
\begin{equation}\label{size.ineq}
\size_{g^*\ol{F},g^*\ol{G}}(g^*e)\leq \size_{\ol{F},\ol{G}}(e)
\end{equation}
holds for every extension class $e\in \widehat{\rm Ext}_S^1(F,G)$.
Motivated again by geometric considerations, we ask -- at least if $K$
is the field $\Q$ -- wether
the size of extensions of $\ol{F}$ by $\ol{G}$
is invariant under the base change $g$, namely
whether the
inequality \refg{size.ineq} is indeed an equality  for any extension
class
$e\in \widehat{\rm Ext}_S^1(F,G)$.

Let $\ol{E}$ denote the hermitian vector bundle $\ol{F}^\vee\otimes
\ol{G}$
over Spec ${\mathcal O}_K$.
The extension of scalars ${\mathcal O}_K \hookrightarrow {\mathcal O}_{K'}$ defines a
natural
$\R$-linear map
\[
\Delta:E_\R=E\otimes_\Z\R\longrightarrow
(g^*E)_\R=(E\otimes_{{\mathcal O}_K}{\mathcal O}_{K'})\otimes_\Z\R.
\]
Let $\mathcal{V}(\ol{E})\subseteq E_\R$ denote the Voronoi cell
of the euclidean lattice $E\subseteq E_\R$ underlying $\ol{E}$.
Then the size of extensions of $\ol{F}$ by $\ol{G}$
is invariant under the base change $g$ if and only if
\begin{equation}\label{vertcond}
\Delta\bigl(\mathcal{V}(\ol{E})\bigr)\in \mathcal{V}(g^*\ol{E}).
\end{equation}
Clearly \refg{vertcond} holds iff $\Delta$ maps
the set of vertices of the polytope 
$\mathcal{V}(\ol{E})$ to $\mathcal{V}(g^*\ol{E})$.

Here are some results which point towards a positive answer to our
question
in the case where the base field is $\Q$.
Hence assume $K=\Q$, put $L=K'$, and define $\ol{E}$ as above.
Then we show that \refg{size.ineq} is an equality
for any extension class $e\in \widehat{\rm Ext}_S^1(F,G)$ if either

(i) $L/\Q$ is an abelian extension, or

(ii) $\ol{E}$ is an orthogonal direct sum of hermitian line bundles,
or

(iii) $\ol{E}$ is a root lattice, or

(iv) $\ol{E}$ is a lattice of Voronoi's first kind
(hence in particular if $\rk_\Z E\leq 3$).

We use condition \refg{vertcond} to prove these results.
For abelian extensions, we reduce to the cyclotomic case and use some
auxiliary results of Kitaoka, which he established when investigating
minimal vectors in tensor products of euclidean lattices.
Using the elementary inequality
\[
\sum\limits_{\sigma:L\hookrightarrow \C}|\sigma(\alpha)|^2
-\sum_{\sigma:L\hookrightarrow \C}\Re \sigma(\alpha)\geq
\sum\limits_{\sigma:L\hookrightarrow \C}|\sigma(\alpha)|^2
-\sum_{\sigma:L\hookrightarrow \C}|\sigma(\alpha)|\geq
0
\]
satisfied by any integral element $\alpha\in{\mathcal O}_L$, we show
that \refg{vertcond} holds when  $\ol{E}$ has rank one,
and consequently when it splits as a direct sum of
hermitian line bundles.
Our proof for root lattices relies on the computation of the vertices
of the
Voronoi cells of the irreducible root lattices $A_n$, $D_n$, $E_6$,
$E_7$,
and $E_8$ by Conway and Sloane (\cite{conwaysloane99}, Chapter 21).
Our treatment of lattices of Voronoi's first kind uses the
description of
the Voronoi cell of an euclidean lattice with strictly obtuse
superbase
which is given in  Appendix B.

Finally, 
as a consequence of our ``reduction theorem" and of 
case (ii), we show, in the case
where the
base field $K$ is $\Q$, that equality holds in (\ref{size.ineq})
 ``up to some constant". Namely, we derive the existence of
a non-negative real constant $c(\rk F, \rk G)$ --- depending on the
ranks
of $F$ and $G$ only --- such that the inequality 
\[
\size_{\ol{F},\ol{G}}(e)\leq \size_{g^*\ol{F},g^*\ol{G}}(g^*e) +
c(\rk F, \rk G)
\]
holds for any class $e\in \widehat{\rm Ext}_S^1(F,G)$.

 Appendix A gathers ``well known" facts concerning extension
groups of sheaves of modules. In particular, it specifies 
sign conventions which enter in the construction of canonical
isomorphisms between variously defined extension groups.

Appendix B contains a self-contained 
presentation of lattices of Voronoi's
first kind, a description of their Voronoi cells, and various facts
concerning these lattices  which might be of independent interest.

\medskip

\noindent {\bf 0.4.}
The starting points of this paper have been, in 1998, (i) the
observation  that, for any two hermitian vector bundles $\ol{F}$ and
$\ol{G}$ over an arithmetic curve $X,$ the set of isomorphism classes
of admissible extensions of $\ol{F}$ by $\ol{G}$ becomes an abelian
group when the Baer sum of two admissible extensions is equipped with
the hermitian structure defined by formula (\ref{drole de somme})
\emph{infra}, and (ii) Grothendieck's remark  in
\cite{grothendieck66} on the non-trivial information encoded in the
extension class of the Hodge extension (\ref{Hodge}). 

Related ideas have been investigated in Mochizuki's preprints
\cite{mochizuki99p}. Let us emphasize a major difference between his
approach and ours: Mochizuki thinks of the Hodge extension in the
context of Arakelov geometry as some kind of \emph{non-linear}
geometric object, while we see it as an element of some naturally
defined \emph{abelian} extension \emph{group}.
Moreover, Mochizuki's earlier work \cite{mochizuki96} has been an
inspiration for considering the arithmetic Schwarz extension.

Let us finally indicate that, in  \cite{chamberttschinkel01},
Chambert-Loir and Tschinkel have defined and
investigated ``arithmetic torsors'' under some group scheme $G$ on an
arithmetic
scheme $X$, at least when $G$ is deduced by base change from a group
scheme over an ``arithmetic curve". Their definition easily extends
to the case of general smooth affine group schemes over $X$, and
specialized to vector groups  of the form $\check{E}\otimes F,$ where
$E$ and $F$ are vector bundles over $X$, is equivalent to our
definition of arithmetic  extensions of $E$ by $F$ (see
\ref{arithmetic torsors}, \emph{infra}).

It is a pleasure to thank E. Bayer for helpful remarks on euclidean lattices, and R. Bost
for his help in the preparation of the figure.
We are grateful to the TMR network `Arithmetic geometry'
and the DFG-Forschergruppe `Algebraische Zykel und $L$-Funktionen'
for their support and to the universities of Regensburg and
Paris (6,7,11,13) for their hospitality.

\section{Preliminaries}

\subsection{Arithmetic schemes}
We work over an {\it arithmetic ring} $R=(R,\Sigma,F_\infty)$ in the
sense of Gillet and Soul\'e,
\cite[3.1.1]{gilletsoule90}.
Recall that this means that $R$ is an excellent regular
noetherian integral domain, $\Sigma$ is a
finite nonempty set of monomorphisms from $R$ to $\C$, and $F_\infty$
is a conjugate-linear
involution of $\C$-algebras $F_\infty:\C^\Sigma\to \C^\Sigma$ such
that
$F_\infty\circ\delta=\delta$ for the canonical map
$\delta:R\to \C^\Sigma=\prod_{\sigma\in\Sigma}\C$.

Let $S$ be the spectrum of an arithmetic ring $R$, and $K$ its field
of fractions.
An {\it arithmetic scheme}\footnote{We use the terminology
\emph{arithmetic scheme} for what is called an \emph{arithmetic
variety} in \cite{gilletsoule90} and subsequent papers by Gillet and
Soul\'e, in order to avoid confusion with quotients of symmetric domains
by the action of arithmetic groups.}
$X$ over $R$ is a
separated $S$-scheme $X$ of finite type
such that each base change $X_\sigma =X×_{R,\sigma} \C$, $\sigma\in
\Sigma$,
is smooth over ${\rm Spec}\,\C$ (or equivalently such that $X_{K}$ is 
smooth over $K$).
For $\sigma$ in $\Sigma$, we write $X_\sigma=X\otimes_{R,\sigma}\C$.
We obtain a scheme
\[
X_\Sigma=X\otimes_{R,\delta}\C^\Sigma=\coprod_{\sigma\in
\Sigma}X_\sigma
\]
and a complex manifold
\[
X_\Sigma(\C)=\coprod_{\sigma\in \Sigma}X_\sigma(\C).
\]
We write $X(\C)$ instead of $X_\Sigma(\C)$ if
$\Sigma=\{\sigma:R\hookrightarrow \C\}$.

\subsubsection{}\label{aring}
The most prominent example of an arithmetic ring is
${\mathcal O}_K=({\mathcal O}_K,\Sigma, F_\infty)$ where ${\mathcal O}_K$ is the ring
of integers in a number field $K$, $\Sigma$ is the set of complex
embeddings $\sigma:K\hookrightarrow \C$,
and $F_\infty$ is given by
\[
F_\infty:\C^\Sigma\longrightarrow  \C^\Sigma\,\,,\,\,
(z_\sigma)_{\sigma\in\Sigma}
\mapsto (\overline{z_{\overline{\sigma}}})_{\sigma\in\Sigma}.
\]
Then an arithmetic scheme over ${\mathcal O}_{K}$ is precisely a separated
$\Z$-scheme $X$ of finite type such that $X_{\Q}$ is smooth, equipped 
with a scheme morphism to $\Spec {\mathcal O}_{K},$  and $X_{\Sigma}(\C)$ is the
complex manifold $X(\C)$ of all complex points of $X.$

\subsubsection{}\label{arbach}
There are natural morphisms of locally ringed spaces
\[
j:(X_\Sigma(\C),{\mathcal O}_{X_\Sigma}^{\,\an}) {\longrightarrow } (X,{\mathcal O}_X)
\]
where $\OanX$ denotes the sheaf of holomorphic functions on
the complex manifold $X_\Sigma(\C)$ and
\[
\kappa:(X_\Sigma(\C),\CinfX){\longrightarrow } (X_\Sigma(\C),\OanX)
\]
where $\CinfX$ denotes the sheaf of complex valued smooth functions.
The morphism $j$ is flat by \cite[Exp. XII]{SGA1}.
To any ${\mathcal O}_X$-module $F$ on $X$ is associated an $\OanX$-module
$F_\C^{\an}=j^*F$ on $X_\Sigma(\C)$ and an $\CinfX$-module
$F_\C=\kappa^*j^*F$. The so-defined functor $F \longmapsto F_\C$ is exact, 
as a consequence of following Lemma:

\begin{lemma}
    
The morphism $\kappa$ is flat, i.e. $\mathcal{C}^\infty_{X_\Sigma,p}$
is a flat
${\mathcal O}^\an_{X_\Sigma,p}$-module for each $p$ in $X_\Sigma(\C)$.
\end{lemma}

\proof
We consider for $n\geq 0$ the $\R$-algebra ${\mathcal O}_{\R^{2n},0}$
(resp. ${\mathcal E}_{\R^{2n},0}$) of germs of
real analytic (resp. real valued $\mathcal{C}^\infty$) functions
around $0$ in
$\R^{2n}$, and the $\C$-algebra ${\mathcal O}^\an_{\C^n,0}$ of germs of
holomorphic
functions around $0$ in $\C^n$.
The canonical map from ${\mathcal O}_{\R^{2n},0}$  to ${\mathcal
E}_{\R^{2n},0}$ is
flat by \cite[VI Cor. 1.3]{tougeron72}.
We have ${\mathcal
E}_{\R^{2n},0}\otimes_{\R}\C=\mathcal{C}_{\C^n,0}^\infty$
under the canonical identification of $\C^n$ with $\R^{2n}$.
Therefore $\kappa$ is flat if we show that
${\mathcal O}_{\R^{2n},0}\otimes_{\R}\C$ is
flat over  ${\mathcal O}^\an_{\C^n,0}$.
This can be checked on completions (which are faithfully flat).
We have
\[
\widehat{{\mathcal O}_{\R^{2n},0}\otimes_{\R}\C}=\C[[z_1,..,z_n,\ol{z}_1,..,\ol{z}_n]]
\]
and $\widehat{{\mathcal O}^\an_{\C^n,0}}=\C[[z_1,..,z_n]]$.
Our claim follows.
\qed

\subsubsection{}
Let $F_\infty$ denote the anti-holomorphic involution of the complex
manifold $X_\Sigma(\C)$ which maps $s:\Spec\,\C\to X$ to the
composition of complex conjugation in $\C$ with $s$.
We obtain an induced $\C$-antilinear involution
\[
F_\infty:A^k(X_\Sigma(\C),\C)\longrightarrow
A^k(X_\Sigma(\C),\C)\,\,,\,\,\alpha\mapsto
\overline{F_\infty^*(\alpha)}
\]
on the space of smooth complex valued $k$-forms on $X_\Sigma(\C)$.
One checks easily that this map is $\C$-anti-linear and
$\Gamma(X,{\mathcal O}_X)$-linear.
Furthermore it respects the $(p,q)$-type and commutes with $d$,
$\partial$, and $\overline{\partial}$.

\subsubsection{}\label{aforms}
For any ${\mathcal O}_X$-module $F$ on $X$, we consider the sheaf
\[
A^k(\ubar,F):=F_\C^\an \otimes_{\OanX}A^k(\ubar,\C) 
= F_\C \otimes_{\CinfX} A^k(\ubar,\C)
\]
on $X_\Sigma(\C)$.
It may be decomposed according to types:
\begin{equation}\label{types}
A^k(\ubar,F) = \bigoplus_{p+q=k} A^{p,q}(\ubar,F)
\end{equation}
where, for any two non-negative integers $p$ and $q$:
$$
A^{p,q}(\ubar,F):=F_\C^\an \otimes_{\OanX}A^{p,q}(\ubar,\C) 
= F_\C \otimes_{\CinfX} A^{p,q}(\ubar,\C).
$$

The space of sections $A^k(X_\Sigma(\C),F)$ is endowed with the
$\C$-antilinear involution  $F_\infty$ 
(which specializes to the one considered above when $F={\mathcal O}_{X}$),
defined by complex conjugation, acting both on 
$X_{\Sigma}(\C)$ and on the coefficients (namely, $k$-forms and fibers
of $F_{\C}$).

This involution is compatible with the decomposition into types
(\ref{types}) and with the Dolbeault operator.
We define
\[
A^k(\XR,F)=A^k(X_\Sigma(\C),F)^{F_\infty}
\mbox{ and }
A^{p,q}(\XR,F)=A^{p,q}(X_\Sigma(\C),F)^{F_\infty}
\]
and we obtain an induced Dolbeault operator
\[
\overline{\partial}_F:A^{p,q}(\XR,F)\longrightarrow  A^{p,q+1}(\XR,F).
\]
Its kernel will be denoted
$Z_{\overline{\partial}}^{p,q}(X_\R, F),$
and the $p$-th cohomology group 
$$Z_{\overline{\partial}}^{0,p}(X_\R, F)/ \overline{\partial}_F (A^{0, p-1}(\XR,F))$$ 
of the Dolbeault complex
$(A^{0,\cdot}(X_\R,F),\dolb_F),$ will be denoted $H^p_{\rm Dolb}(X_\R, F)$.
The Dolbeault isomorphism (see Appendix \ref{Dolbeault})
$${\rm Dolb}_{F_{\C}^\an}: H^p(X_{\Sigma}(\C),F_{\C}^\an) \longrightarrow 
\HD^p(X_{\Sigma}(\C),F_{\C}^\an):=H^p(X_{\Sigma}(\C),{\mathcal Dolb}(F_{\C}^\an))$$
yields an isomorphism
$${\rm Dolb}_{F_{\R}}: H^p(X_{\Sigma}(\C),F_{\C}^\an)^{F_{\infty}} 
\longrightarrow 
H^p_{\rm Dolb}(X_{\R},F)$$
between $H^p_{\rm Dolb}(X_{\R},F)$ and the real vector subspace in the cohomology
group $H^p(X_{\sigma}(\C),F_{\C}^\an)$ of elements invariant under
complex conjugation (acting both on $X_{\Sigma}(\C)$ and on
coefficients).

We shall also denote
$${\rm Dolb}_{F_{\R}}: H^p(X,F)
\longrightarrow
H^p_{\rm Dolb}(X_{\R},F)$$
the composition of the above isomorphism and of the canonical map
$j^\ast: H^p(X,F)
\rightarrow
H^p(X_{\sigma}(\C),F_{\C}^\an)^{F_{\infty}}$ defined by pulling back
through the morphism of ringed spaces $j:(X_\Sigma(\C),{\mathcal O}_{X_\Sigma}^{\,\an})
{\rightarrow } (X,{\mathcal O}_X).$ More generally, if $E$ and $F$ are ${\mathcal O}_X$-modules such that $E_\C^\an$ is a locally
free of finite rank (\ie, a holomorphic vector bundle), the base change by $j$ and the Dolbeault isomorphism
${\rm Dolb}_{E_\C,F_\C}$ (\cf \ref{Dolbeault}) define a map
$${\rm Dolb}_{E_\R,F_{\R}}: \Ext^p_{{\mathcal O}_X}(E,F)
\longrightarrow
H^p_{\rm Dolb}(X_{\R},E^\lor \otimes F).$$

Let $\Omega_{X/S}^r$ denote the sheaf of $r$-th relative Kähler differentials
of $X$ over $S$.
Then $(\Omega_{X/S}^{r})_\C^{\an}$ is the sheaf of holomorphic
$r$-forms
on $X_\Sigma(\C)$. This allows a natural identification
\[
\nu:A^{p,q}(X_\Sigma(\C),F)\stackrel{\sim}{\longrightarrow}
A^{0,q}(X_\Sigma(\C), F\otimes \Omega_{X/S}^p)\,\,,\,\,
f\otimes (\alpha\wedge\beta)\mapsto  (f\otimes \alpha)\otimes \beta
\]
for differential forms $\alpha$ and $\beta$ of type $(p,0)$
and $(0,q)$ respectively. Observe that the space $A^{p,p}(\XR,{\mathcal O}_{X})$
does \emph{not} coincide with the space $A^{p,p}(\XR)$, 
considered 
in \cite[3.2.1]{gilletsoule90}, of
real $(p,p)$-forms $\alpha$ on $X_\Sigma(\C)$
satisfying $F_\infty^*\alpha=(-1)^p\alpha$.
Instead, we have an embedding
\[
A^{p,p}(X_\R)\hookrightarrow A^{p,p}(\XR,{\mathcal O}_{X}) \simeq  A^{0,p}(X_\R,\Omega_{X/S}^p)\,\,
,\,\,\alpha\mapsto  (-2\pi i)^p \nu(\alpha),
\]
the image of which is the $(-1)^p$-eigenspace of the involution on
$A^{p,p}(\XR,{\mathcal O}_{X})$ defined by complex conjugation acting on
coefficients only.

\subsection{Hermitian coherent sheaves}
Let $X$ be an arithmetic scheme.
A vector bundle on $X$ is a locally free $\mathcal {O}_X$-module $E$
of
finite rank.
The dual vector bundle $\ihom(E,{\mathcal O}_X)$ is denoted by $E^\lor$.
Following \cite[Def. 25]{gilletsoule92} we define a {\it hermitian
coherent
sheaf} $\ol{E}$ on $X$ as a pair $(E,h)$ consisting of a
coherent ${\mathcal O}_X$-module $E$ whose restriction to the generic fiber
$X_F$ is locally free, together with a $F_\infty$-invariant
$\mathcal{C}^\infty$-hermitian metric $h$ on the holomorphic vector
bundle
$E_\C^\an$.
A {\it hermitian vector bundle} on $X$ is a hermitian coherent sheaf
whose
underlying coherent ${\mathcal O}_X$-module is locally free.
There are natural hermitian structures on tensor products, exterior
powers, and inverse
images of hermitian coherent sheaves, and on the dual of hermitian
vector bundles.

Observe also that, if $\ol{E}$ and $\ol{F}$ are two hermitian vector
bundles over $X$, then the canonical isomorphism $\ihom(E,F) 
\simeq E^\lor \otimes F$ allows us to equip $\ihom(E,F)$ with
a structure of hermitian vector bundle, which makes it canonically
isomorphic with $\ol{E}^\lor \otimes \ol{F}$. For any section $T$ of 
$\ihom(E,F)_{\C}$ and any $x \in X_{\Sigma}(\C),$ the
so-defined norm $\left\|T(x)\right\|_{\ol{E}^\lor \otimes \ol{F}}$ is 
the Hilbert-Schmidt norm of the $\C$-linear map between the hermitian
vector spaces $(E_{x}, \left\|.\right\|_{\ol{E}})$ and $(F_{x},
\left\|.\right\|_{\ol{F}})$. Occasionally  we shall also use the operator 
norm of such maps, and when confusion may arise, we shall denote
$\left\|T(x)\right\|_{\infty}$ or $\left\|T(x)\right\|^{\infty}$ the
latter, and $\left\|T(x)\right\|_{HS}$ the former.

\subsubsection{Direct image}\label{defdim}
Let $f:Y \to  X$ be a finite flat morphism of arithmetic varieties
such that $f_{F}:Y_{F} \to X_{F}$ is \'etale --- or
equivalently, such that $f_{\Sigma}:Y_{\Sigma} \to X_{\Sigma}$
is an \'etale covering --- then, for any hermitian coherent sheaf
$\ol{E}$ over $Y,$ we may consider its direct image
$f_{\ast}\ol{E},$ namely the hermitian coherent sheaf on $X$
defined by the coherent sheaf $f_{\ast}E$ equipped with the
hermitian structure which, for any $x \in X_{\Sigma}(\C),$ is given on
the fiber
\[
(f_{\ast}E)_{x}\simeq \bigoplus_{y\in f_\Sigma^{-1}(x)} E_{y}
\]
by the
direct sum of the hermitian structures on the $E_{y},$ $y\in
f_\Sigma^{-1}(x).$

\subsubsection{}\label{arithm.proj.formula}
For $f:Y\to X$ as above and hermitian coherent sheaves $\ol{E}$ on
$Y$ and
$\ol{F}$ on $X$, adjunction defines a natural morphism of
${\mathcal O}_Y$-modules
\[
f^*(f_*E\otimes_{{\mathcal O}_X} F)=(f^*f_*E)\otimes_{{\mathcal O}_Y} f^*F
\longrightarrow E\otimes_{{\mathcal O}_Y} f^*F.
\]
By adjunction it induces a canonical isomorphism of ${\mathcal O}_X$-modules
\[
f_*(\ol{E})\otimes_{{\mathcal O}_X}\ol{F}
\stackrel{\sim}{\longrightarrow }
f_*(\ol{E}\otimes_{{\mathcal O}_Y}f^*\ol{E})
\]
which is an isometry as a direct consequence of our definitions.

\subsubsection{}\label{arithm.base.change}
Let
\[
\begin{array}{ccc}
Y'& \stackrel{f'}{\longrightarrow } & X'\\
\downarrow \scriptstyle{h'}& & \downarrow \scriptstyle{h} \\
Y & \stackrel{f}{\longrightarrow } &X
\end{array}
\]
be a cartesian square of arithmetic varieties.
Let $\ol{E}$ be a hermitian coherent sheaf on $Y$ and assume
that $f$ is as in \ref{defdim}.
By adjunction we obtain a natural morphism of ${\mathcal O}_X$-modules
\[
h'^*f^*h_*E=f'^*h^*h_*E\longrightarrow f'^*E
\]
which induces a canonical isomorphism of ${\mathcal O}_X$-modules
\[
f^*h_*\ol{E}\stackrel{\sim}{\longrightarrow }h'_*f'^*\ol{E}.
\]
The latter is an isometry as a direct consequence of our definitions.

\subsection{Extensions}\label{algbaersum1}
We briefly recall some basic facts concerning $1$-extensions of
sheaves of modules. For more details and references, we refer
the reader to Appendix A.

Let $F$, $G$ denote ${\mathcal O}_X$-modules on a ringed space $(X,{\mathcal O}_X)$.
An {\it extension of $F$ by $G$} is a short exact sequence of
${\mathcal O}_X$-modules
\[
\mathcal{E}:0\longrightarrow  G\longrightarrow  E\longrightarrow
F\longrightarrow  0.
\]
A morphism of extensions
\begin{equation}\label{morph1}
(\alpha, \beta,\gamma):\mathcal{E}_1\longrightarrow  \mathcal{E}_2
\end{equation}
is given by a commutative diagram
\[
\begin{array}{ccccccccc}
{\mathcal E}_1\,:\,0 &\longrightarrow  &G_1& \longrightarrow
&E_1&\longrightarrow  &F_1& \longrightarrow  & 0\\
  & & \downarrow \scriptstyle{\alpha} & &\,\,\, \downarrow
{\scriptstyle \beta}
& & \downarrow\scriptstyle{\gamma}&\\
{\mathcal E}_2\,:\,0 &\longrightarrow  &G_2& \longrightarrow
&E_2&\longrightarrow  &F_2& \longrightarrow  & 0.
\end{array}
\]
Recall that two extensions $\mathcal{E}_1$ and $\mathcal{E}_2$ of $F$
by
$G$ are called {\it isomorphic} iff
there exists a morphism \refg{morph1} as above with $\alpha=\id_G$
and $\gamma=\id_F$.

Given an extension $\mathcal{E}$ as above,
we consider the boundary map
\[
\Hom_{{\mathcal O}_X}(F,F)\stackrel{\partial}{\longrightarrow
}\Ext^1_{{\mathcal O}_X}(F,G)
\]
where $\Ext^p_{{\mathcal O}_X}(F,\,.\,)$ denotes the $p$-th right
derived functor of $\Hom_{{\mathcal O}_X}(F,\,.\,)$.
It is well known (compare \cite[Th. 3.4.3]{weibel94},
\cite[Ex. III 6.1]{hartshorne77}, or Proposition \ref{extisext}) that
\begin{equation}
\genfrac{\{}{\}}{0pt}{}{\mbox{isomorphism classes of}}
{\mbox{extensions of $F$ by $G$}} \longrightarrow
\Ext^1_{{\mathcal O}_X}(F,G)\,\,,\,\,[\mathcal{E}]\mapsto  \partial(\id_F)
\end{equation}
defines a group isomorphism if we equip the left-hand side with
the group structure  induced by the Baer sum of extensions.
Recall that the Baer sum of two extensions
\begin{equation}\label{Baersum}
\mathcal{E}_j:0\longrightarrow  G\stackrel{i_j}{\longrightarrow }
E_j\stackrel{p_j}{\longrightarrow } F\longrightarrow  0\,\,\,(j=1,2)
\end{equation}
is the extension
\[
\mathcal{E}:\,0\longrightarrow  G\stackrel{i}{\longrightarrow }
E\stackrel{p}{\longrightarrow } F\longrightarrow  0
\]
where
\begin{equation}\label{algbaersum2alt}
E=\frac{\Ker(p_1-p_2:E_1\oplus E_2 { \longrightarrow  } F)}
{{\rm Im}\,((i_1,-i_2): G {\longrightarrow  } E_1\oplus E_2)},
\end{equation}
and $p$ and $i$ are given as $p(e_1,e_2)=p_1(e_1)=p_2(e_2)$ and
$i(g)=(i_1(g),0))=(0,i_2(g))$.

\section{The arithmetic extension group $\widehat{\rm Ext}_X^1(F,G)$}
This section is devoted to the basic definitions and properties of arithmetic 
extensions and of the corresponding groups of $1$-extensions on an arithmetic 
scheme $X$.
We associate a canonical differential form $\Psi(\mathcal{E},s)$ with
an
arithmetic extension $(\mathcal{E},s)$, namely its ``second fundamental
form" $\overline{\partial}s$.
The arithmetic extension group fits into two exact sequences which are
formally similar to corresponding sequences for arithmetic Chow
groups.
We discuss functorial properties of our extension groups and relate
arithmetic extensions to admissible extensions of hermitian coherent
sheaves and arithmetic torsors in the sense of Chambert-Loir and
Tschinkel. We also discuss an interpretation of the group of arithmetic extensions 
$\widehat{\rm Ext}_X^1(F,G)$
as suitable group of morphisms in the derived category of ${\mathcal O}_X$-modules.

In this section, $(R,\Sigma)$ denotes an arithmetic ring, $K$ the 
field of fractions of $R,$ and $X$ an
arithmetic scheme over $S:= \Spec R.$

\subsection{Basic definitions}\label{basicdef}

Let $F$ and $G$ be ${\mathcal O}_X$-modules.
An \emph{arithmetic extension $(\mathcal{E},s)$ of $F$ by $G$} is by
definition an
extension
\begin{equation}\label{ext}
\mathcal{E}:\,0\longrightarrow  G\stackrel{i}{\longrightarrow }
E\stackrel{p}{\longrightarrow } F\longrightarrow  0
\end{equation}
of ${\mathcal O}_X$-modules together with an $F_\infty$-invariant
$\mathcal{C}^\infty$-splitting
\begin{equation}\label{section}
s:F_{\C}\longrightarrow  E_\C
\end{equation}
of the associated extension of $\mathcal{C}^\infty_{X_\Sigma}$-modules
\[
\mathcal{E}_\C:\,0\longrightarrow
G_\C\stackrel{i_\C}{\longrightarrow }
E_\C\stackrel{p_\C}{\longrightarrow } F_\C\longrightarrow  0.
\]
There exists a unique map
such that the relation
\begin{equation}\label{st}
{\rm Id}_{E_\C}=s\circ p_\C+i_\C\circ t
\end{equation}
holds. It is $\mathcal{C}^\infty$ and
 $F_\infty$-invariant,
and the sections $s$ and $t$ determine each other uniquely.
We sometimes write $(\mathcal{E},s,t)$ to emphasize that $t$ is
defined
by \refg{st}.

A {\it morphism between arithmetic extensions}
$(\mathcal{E}_1,s_1,t_1)$ and
$(\mathcal{E}_2,s_2,t_2)$ is given by a morphism
$(\alpha,\beta,\gamma):\mathcal{E}_1\rightarrow \mathcal{E}_2$
of extensions such that
$\beta_\C\circ s_1=s_2\circ \gamma_\C$ holds.
Observe that this condition implies already
$t_2\circ \beta_\C =\alpha_\C\circ t_1$.

Two arithmetic extensions $(\mathcal{E}_1,s_1)$ and
$(\mathcal{E}_2,s_2)$ of $F$ by $G$ are called {\it isomorphic}
iff there exists a morphism $(\alpha,\beta,\gamma)$ from
$(\mathcal{E}_1,s_1)$ to $(\mathcal{E}_2,s_2)$
such that $\alpha={\rm id}_G$ and $\gamma={\rm id}_F$.
Any such morphism is automatically an isomorphism and defines an
isomorphism between
the arithmetic extensions $(\mathcal{E}_1,s_1)$ and
$(\mathcal{E}_2,s_2)$.
We denote
the set of isomorphism classes of arithmetic extensions of $F$ by
$G$ by $\widehat{\rm Ext}_X^1(F,G)$.

Let $(\mathcal{E}_1,s_1)$ and $(\mathcal{E}_2,s_2)$ be two arithmetic
extensions of $F$ by $G$ on $X$.
Let $\mathcal{E}$ denote the algebraic Baer sum
of $\mathcal{E}_1$ and $\mathcal{E}_2$.
The $\mathcal{C}^\infty$-splittings $s_1$ and $s_2$ induce an
$F_\infty$-invariant $\mathcal{C}^\infty$-splittings of $\mathcal{E}$
which maps a section $f$ of $F_\C$ to the class $s(f)$ of $(s_1(f), s_2(f))$
in
the quotient \refg{algbaersum2alt}.
We obtain an arithmetic extension $(\mathcal{E},s)$ which we
call the {\it arithmetic Baer sum} of $(\mathcal{E}_1,s_1)$ and
$(\mathcal{E}_2,s_2)$.

It is a straightforward exercise to check that the arithmetic Baer sum
defines on $\widehat{\rm Ext}_X^1(F,G)$ the structure of an abelian
group. Actually, the opposite of the class of $(\mathcal{E},s)$ is the
one of $(\tilde{\mathcal{E}},-s),$ if we let 
$$\tilde{\mathcal{E}}:\,0\longrightarrow  G\stackrel{i}{\longrightarrow }
E\stackrel{-p}{\longrightarrow } F\longrightarrow  0.$$ Moreover, with
the notation of \ref{sec:BaerSum}, if
$(\mathcal{E}_{1},s_{1}),\ldots,(\mathcal{E}_{k},s_{k})$ are
arithmetic extensions over $X,$ one defines a $F_{\infty}$-invariant 
$\mathcal{C}^\infty$-splitting
$s$ of $\mathcal{E}_{1}+\ldots+\mathcal{E}_{k}$ by sending a section
$f$ of $F_{\C}$ to the class of $(s_{1}(f),\ldots,s_{k}(f))$, and the 
class of 
$(\mathcal{E}_{1}+\ldots+\mathcal{E}_{k}, s)$ in $\widehat{\rm Ext}_X^1(F,G)$
is the sum of the classes
$[(\mathcal{E}_{1},s_{1})],$\ldots,$[(\mathcal{E}_{k},s_{k})]$.

\subsection{The first exact sequence}
Let $F$ and $G$ be  ${\mathcal O}_X$-modules on the arithmetic scheme $X$.
We consider the group ${\rm Hom}_{{\mathcal O}_X}(F,G)$
of homomorphisms of $\mathcal{O}_X$-modules from $F$ to $G$ and
the space
\[
{\rm Hom}_{\CinfX}(F_{\C},G_{\C})^{F_\infty}
=\{T=(T_\sigma)_\sigma\in \bigoplus\limits_{\sigma\in
   \Sigma}
{\rm Hom}_{\mathcal{C}^\infty_{X_\sigma(\C)}}
(F_{\sigma},G_{\sigma})\,|\,\overline{T}_\sigma=T_{\overline{\sigma}}\}.
\]
of $F_\infty$-invariant $\mathcal{C}^\infty$-homomorphisms
from $F_{\C}$ to $G_{\C}$.
There is a canonical map
\[
b:{\rm Hom}_{\CinfX}(F_\C,G_\C)^{F_\infty} \longrightarrow
\widehat{\rm Ext}^1_X(F,G)
\]
which sends any $T$ in $\Hom_{\CinfX}(F_\C,G_\C)^{F_\infty}$ to
the class of
the arithmetic
extension $(\mathcal{E},s)$ where $\mathcal{E}$ is the trivial
algebraic extension, defined by (\ref{ext}) with 
$E:=G\oplus F$ and $i$ and $p$ the obvious injection and projection
morphisms, and $s$ is given by $s(f)=(T(f),f))$.

\begin{theorem}\label{longexseq1}
The map $b$ is a group homomorphism. It fits into  an exact sequence
\begin{multline}\label{eq:longexseq1}
{\rm Hom}_{{\mathcal O}_X}(F,G)
\stackrel{\iota}{\longrightarrow }
{\rm Hom}_{\CinfX}(F_\C,G_\C)^{F_\infty}
\stackrel{b}{\longrightarrow }
\widehat{\rm Ext}_X^1(F,G)
\stackrel{\nu}{\longrightarrow } {\rm Ext}_{{\mathcal O}_X}^1(F,G)\\
\stackrel{\digamma}{\longrightarrow } {\rm
   Ext}^1_{\CinfX}(F_\C,G_\C)
\end{multline}
where $\iota$ is the canonical map, $\nu$ maps the class of an
arithmetic extension $(\mathcal{E},s)$ to the class of the underlying
 extension $\mathcal{E}$ of ${\mathcal O}_{X}$-modules,
and $\digamma$ maps the class of an extension of ${\mathcal O}_{X}$-modules to the class of 
the
associated extension of $\mathcal{C}^\infty_{X_{\Sigma}}$-modules.
Furthermore
\begin{equation}\label{vanish1}
{\rm Ext}^1_{\CinfX}(F_\C,G_\C)=0
\end{equation}
if $F_\C$ is a vector bundle.
\end{theorem}

Observe that, when this last assumption is satisfied, we may also
identify the real vector space ${\rm
Hom}_{\CinfX}(F_\C,G_\C)^{F_\infty}$ with
$A^0(X_{\R},F^\lor\otimes G).$

\proof
We first show that $b$ is a group homomorphism.
Consider $T_1,T_2$ in  $\Hom_{\CinfX}(F_\C,G_\C)^{F_\infty}$, and the 
associated arithmetic extension classes
$(\mathcal{E}_1,s_1)$ and
$(\mathcal{E}_2,s_2)$ that define $b(T_1)$ and $b(T_2)$ respectively.
Let $(\mathcal{E},s)$ be the arithmetic Baer sum of
$(\mathcal{E}_1,s_1)$ and $(\mathcal{E}_2,s_2)$
With self-explanatory notation, we get
\[
E_1\oplus E_2=G\oplus F\oplus G \oplus F.
\]
and obtain an isomorphism
\[
E=\frac{\Ker\,(p_1-p_2:E_1\oplus E_2\longrightarrow  F)}
{\Im\,((i_1,-i_2):G\longrightarrow  E_1\oplus E_2)}\longrightarrow
G\oplus F.
\]
which maps the class of $(g_1, f_1,g_2, f_2)$ to $(g_1+g_2,f_1)$. This
isomorphism defines an isomorphism of arithmetic extensions from
$(\mathcal{E},s)$ to the arithmetic extension of $F$ by $G$ that
defines $b(T_{1}+T_{2}).$
Consequently $b(T_{1})+b(T_{2})$ and $b(T_{1}+T_{2})$ coincide.

Let us check the exactness at $\widehat{\rm Ext}_X^1(F,G).$
We clearly have $\nu\circ b=0$. Conversely, 
let cl$(\mathcal{E},s)$ be an arithmetic extension class in the
kernel of
$\nu$.
One can chose an isomorphism $\tau:E\stackrel{\sim}{\to}
G\oplus F$ between $\mathcal{E}$ and the trivial algebraic
extension of $F$ by $G$.
Then the class of $(\mathcal{E},s)$ is precisely the image of
${\rm pr}_1\circ \tau_\C\circ s\in {\rm Hom}_{\CinfX}(F_\C,
G_\C)^{F_\infty}$
under $b$.

Let us show that $\ker\, b={\rm Im}\, \iota$.
Given $T\in\Hom_{{\mathcal O}_X}(F,G)$, the map
\[
G\oplus F\longrightarrow  G\oplus F\,\,\,,\,\,\,(g,f)\mapsto
(g+T(f),f)
\]
defines an isomorphism of arithmetic extensions which gives
$b(0)=b(T)$.
Hence we have $b\circ \iota=0$.
Conversely, let $T\in\Hom_{\CinfX}(F_\C,G_\C)^{F_\infty}$ be in the kernel of $b$.
One gets an isomorphism between the arithmetic extensions
representing $b(0)$
and $b(T)$ which is given by a matrix
\[
\left(\begin{array}{cc}1&\tilde{T}\\0&1\end{array}\right).
\]
for some $\tilde{T}\in\Hom_{{\mathcal O}_X}(F,G)$. It follows that
$T$ equals $\iota(\tilde{T}),$ and consequently
belongs to the image of $\iota$.

The image of $\nu$ is the kernel of $\digamma$ by the very definition of an
arithmetic extension class. This establishes the exactness at ${\rm Ext}^1_{{\mathcal O}_X}(F,G)$.
Finally
(\ref{vanish1}) follows from the existence of partitions of unity.
\qed

\begin{corollary}\label{cor.affine}
     When $X$ is an affine scheme,  $F$ is a vector bundle
\footnote{\ie, locally free of finite rank.}, and $G$ quasi-coherent, the map
     $b$ induces an isomorphism of abelian groups
     \[
     \Exthat^1_{X}\,(F,G) \simeq
     \frac{{\rm Hom}_{\CinfX}
     (F,G)^{F_\infty}}{\iota ({\rm Hom}_{{\mathcal O}_X}(F,G))}.
      \]
\end{corollary}

Indeed, under these assumptions on $X$ and $F$, the groups
$\Ext^1_{X}(F,G)$ and ${\rm Ext}^1_{\CinfX}(F_\C,G_\C)$ vanish.

\begin{example}\label{adelic}{\rm
Let $K$ be a number field.
We work over the arithmetic ring ${\mathcal O}_K$ (as defined in \ref{aring}).
Let $\A_K$ denote the ring of ad\`eles of $K$.
Let $E$ and $F$ be vector bundles on $S=\Spec\,{\mathcal O}_K$.
Consider the vector bundle $G={\mathcal Hom}_{{\mathcal O}_S}(E,F)$ as a
vector group scheme on $S$.
There is a canonical isomorphism
\begin{equation}
\widehat{\rm Ext}_S^1(E,F)=G(K)\backslash
G(\A_K)/\mathbb{K}
\end{equation}
where
\[
\mathbb{K}=\prod_{v\in {\scriptstyle \rm Spec}\,{\mathcal O}_K}
G(\hat{{\mathcal O}}_{K,v})× \prod_{v|\infty}\{0\}
\]
is the standard compact subgroup of the commutative group
$G(\A_K^{F_\infty})$.
In order to see this, we observe that  Corollary \ref{cor.affine} gives a
canonical isomorphism
\[
\widehat{\rm Ext}_S^1(E,F) \simeq
\frac{\bigl[ \prod_{\sigma}{\rm
Hom}_{\C}(E_\sigma,F_\sigma)\bigr]^{F_\infty}}
{{\rm Hom}_{{\mathcal O}_K}(E,F)}
=G({\mathcal O}_K)\backslash \Bigl[\prod_{\sigma}G_\sigma (\C)\Bigr]^{F_\infty}.
\]
Moreover $\Bigl[\prod_{\sigma}G_\sigma (\C)\Bigr]^{F_\infty}$ may be identified with  
$\prod_{v|\infty} G(K_v)$,  and the canonical map
\[
G({\mathcal O}_K)\backslash \Bigl[\prod_{\sigma}G_\sigma (\C)\Bigr]^{F_\infty}
\longrightarrow  G(K)\backslash G(\A_K^{F_\infty})/\mathbb{K}
\]
is an isomorphism by the strong approximation theorem (see for
example \cite[Chapter II.15]{Casselsetal67}).
}
\end{example}

\subsection{The second exact sequence}
In this paragraph, we consider ${\mathcal O}_X$-modules $F$ and $G$ over 
$X$ such that $F_\C^{\an}$ and $G_\C^{\an}$ are 
\emph{holomorphic vector bundles} over $X_\Sigma(\C)$. 
This is for instance the case when $F$ and $G$ are coherent ${\mathcal O}_X$-modules, 
and $F_K$ and $G_K$ are locally free. 

We use the notations introduced in \ref{aforms}, i.e. 
$Z_{\dolb}^{0,1}(X_\R, \ihob)$ denotes the subspace of forms in
$A^{0,1}(X_\R,\ihob)$ which are $\dolb_{\ihob}$-closed. 
For any
arithmetic extension $(\mathcal{E},s)$ of the ${\mathcal O}_X$-modules
$F$ by $G$ as in \ref{basicdef}, the Dolbeault operator
\[
\dolb_E:A^0(X_\R,E)\longrightarrow  A^{0,1}(X_\R,E)
\]
has a matrix representation
\begin{equation}\label{matrix}
\overline{\partial}_{E}=\left(\begin{array}{cc}\overline{\partial}_{G}
&\alpha\\0&\ol{\partial}_{F}\end{array}\right)
\end{equation}
 with respect to the
direct sum decomposition $E_\C\cong
G_\C\oplus F_\C$ induced by $s$, for some form $\alpha$ 
in
$Z_{\overline{\partial}}^{0,1}(X_\R,\ihob)$.
In other words, the form $\alpha$ is defined by the equality in 
$A^{0,1}(X_\R, \ihob)$: 
$$i_\C\circ\alpha= \dolb_{\ihob} s,
$$
or, equivalently, by the following identity, valid for any local 
$\mathcal{C}^\infty$ section $f$ of $F_\C$ over $X_\Sigma(\C)$:
$$\dolb_E(s\circ f)=s\circ \dolb_F f + i_\C\circ\alpha\circ f,$$
where, with a slight abuse of notation, we have denoted $s\circ .$ 
and $i_\C\circ.$ the action by composition of $s$ and $i_\C$ on 
sections of $F_\C$ and $E_\C$ with form coefficients. If, as before, $t$ is 
defined by \refg{st}, we also have:
$$\ol{\partial}_{E^\lor \otimes G}t=-\alpha \circ p_\C.$$

Following the terminology discussed in Appendix  \ref{secondfundamental},
we call
$\alpha$ the
{\it second fundamental form} of the arithmetic extension
$(\mathcal{E},s)$. We shall denote it $\Psi(\mathcal{E},s).$

Let us denote
$$p: Z_{\overline{\partial}}^{0,1}(X_\R,\ihob)
\longrightarrow 
H_{\rm Dolb}^1(X_\R, \ihob)
:= Z_{\overline{\partial}}^{0,1}(X_\R,\ihob)/ {\ol{\partial}_{\ihob}}(A^0(X_\R,\ihob))$$
the canonical quotient map.

\begin{proposition}\label{psioperator}
There is a well defined group homomorphism
\[
\Psi:\widehat{\rm Ext}_X^1(F,G)\longrightarrow
Z_{\overline{\partial}}^{0,1}(X_\R,\ihob)
\]
which associates the second fundamental form $\Psi(\mathcal{E},s)$
to the class of an
arithmetic extension $(\mathcal{E},s)$.
The map $\Psi$ fits into the commutative diagrams
\begin{equation}\label{dia1}
\begin{array}{ccc}
{\rm Hom}_{\CinfX}(F_\C,G_\C)^{F_\infty}
&\stackrel{b}{\longrightarrow }&\widehat{\rm Ext}_X^1(F, G)\\
\|&&\downarrow {\scriptstyle \Psi}\\
A^0(X_\R,\ihob)&\stackrel{\ol{\partial}_{\ihob}}{\longrightarrow }
&Z_{\ol{\partial}}^{0,1}(X_\R,\ihob).
\end{array}
\end{equation}
and
\begin{equation}\label{dia2}
\begin{array}{rcc}
\widehat{\rm Ext}_X^1(F,G)&\stackrel{\Psi}{\longrightarrow }
&Z_{\ol{\partial}}^{0,1}(X_\R,\ihob)\\
\downarrow{\scriptstyle \nu}&&\downarrow{\scriptstyle p}\\
{\rm Ext}_{{\mathcal O}_X}^1(F,G)&\stackrel{{\rm Dolb}_{F_\R, G_\R}}{\longrightarrow
}&H^{1}_{\rm Dolb}(X_\R, \ihob).
\end{array}
\end{equation}

\end{proposition}

\proof
Two isomorphic arithmetic extensions yield obviously the same form
$\alpha$.
Consider a morphism $T$ in $\Hom_{\CinfX}(F_\C,G_\C)^{F_\infty}$, and 
the
associated
extension $(\mathcal{E},s)$ that defines $b(T)$. The identity
$$\ol{\partial}_G(T\circ f)=\ol{\partial}_{\Hom (F,G)}(T)\circ f+
T\circ\dolb_F f,$$
valid for any local 
$\mathcal{C}^\infty$ section $f$ of $F_\C$ over $X_\Sigma(\C)$,
shows that
$\Psi(\mathcal{E},s)=\ol{\partial}(T)$.
This proves the commutativity of \refg{dia1}.
We can check locally on $X_\Sigma$ that $\Psi$ is a homomorphism.
Hence it suffices to consider elements in the image of $b$.
In this case the commutativity of \refg{dia1} implies that $\Psi$ is a
homomorphism.

The commutativity of \refg{dia1} follows from the classical fact that 
the ``second fundamental form" provides a representative in Dolbeault 
cohomology of the extension class of a $1$-extension of holomorphic vector 
bundles (\cf \cite{griffiths66} and Appendix A, Proposition \ref{prop:2f}).
\qed

As a consequence of
 Theorem  \ref{longexseq1} and Proposition \ref{psioperator}, we have:

\begin{theorem}\label{longexseq2} Let $F$ and $G$ be
two ${\mathcal O}_X$-modules such that
$X$ such that are $F_\C^{\an}$ and $G_\C^{\an}$ are holomorphic vector bundles.
We have an exact sequence
\begin{multline}\label{eq:longexseq2}
{\rm Hom}_{{\mathcal O}_X}(F,G)\stackrel{\iota}{\longrightarrow }
H^0(X_\R,\ihob)\stackrel{b}{\longrightarrow }
\Exthat^1_X(F,G)\longrightarrow  \\
\stackrel{{\genfrac{(}{)}{0pt}{}{\nu}{\Psi}}}{\longrightarrow }
\Ext_{{\mathcal O}_X}^1(F,G)\oplus Z^{0,1}_{\ol{\partial}}(X_\R,\ihob)
\stackrel{({{\rm Dolb}_{{F^\lor \otimes G}_\R}},-p)}{\longrightarrow }  
H^1_{\rm Dolb}(X_\R,\ihob)
\longrightarrow  0.
\end{multline}
\end{theorem}

\proof 
To establish the exactness at $H^0(X_\R,\ihob)$ and at $\Exthat^1_X(F,G),$ we 
use Theorem \ref{longexseq1}: according to the commutativity of \refg{dia1}, 
we may ``replace"
$\Hom_{\CinfX}(F_\C,G_\C)^{F_\infty}$ by
\begin{eqnarray*}
\Hom_{\CinfX}(F_\C,G_\C)^{F_\infty}
\cap \Ker \,\ol{\partial}_{\ihob}
&=&\Hom_{\OanX}(F_\C^\an,G_C^\an)^{F_\infty}\\
&=&H^0(X_\R,\ihob).
\end{eqnarray*}
The exactness at $H^1_{\rm Dolb}(X_\R,\ihob)$ is clear. It remains to show exactness at
$\Ext_{{\mathcal O}_X}^1(F,G)\oplus Z^{0,1}_{\ol{\partial}}(\XR,\ihob)$. The fact that the composition
 of ${\genfrac{(}{)}{0pt}{}{\nu}{\Psi}}$ and 
$({{\rm Dolb}_{{F^\lor \otimes G}_\R}},-p)$ vanishes is precisely the 
commutativity of \refg{dia2}. Conversely,
consider ${\rm cl}(\mathcal{E})$ in $\Ext_{{\mathcal O}_X}^1(F,G)$ and $\alpha$ in 
$Z^{0,1}_{\ol{\partial}}(\XR,\ihob)$ such that 
$${{\rm Dolb}_{{F^\lor \otimes G}_\R}}({\rm cl}(\mathcal{E})) -p(\alpha)=0.$$
Choose $s$ such that $(\mathcal{E},s)$ becomes an arithmetic
extension.
Then the commutativity of \refg{dia2} implies that
\[
\Psi(\mathcal{E},s)=\alpha+\ol{\partial}T
\]
for some $T\in \Hom_{\CinfX}(F_\C,G_\C)^{F_\infty}$, and the commutativity 
of \refg{dia1} that
$({\rm cl}(\mathcal{E}),\alpha)$ is the image of
cl$(\mathcal{E},s)-b(T)$.
\qed

\subsection{Pushout, pullback, and inverse image}
\label{subsec.funct}\label{compat}
Arithmetic extensions have the expected functorial properties and
behave well with respect
to the maps $\iota$, $b$, $\nu$, and $\Psi$ introduced above.

\subsubsection{Pushout} \label{funct2}
Let $g:G\to G'$ a morphism of ${\mathcal O}_X$-modules on an arithmetic scheme
$X$.
For an arithmetic extension $({\mathcal E},s)$ of an ${\mathcal O}_X$-module
$F$ by $G$,
we obtain a pushout diagram
\[
\begin{array}{ccccccccc}
{\mathcal E}\,:\,0 &\longrightarrow  &G& \longrightarrow  &E
&\longrightarrow  &F&
\longrightarrow  & 0\\
  & & \downarrow {\scriptstyle g}& &\,\,\, \downarrow {\scriptstyle
q} & & ||&\\
{\mathcal E}'\,:\,0 &\longrightarrow  &G'& \longrightarrow  &E'
&\longrightarrow  &F&
\longrightarrow  & 0.
\end{array}
\]
The section $s'=q_\C\circ s$ turns $({\mathcal E}',s')$ into an
arithmetic extension which we denote as $g \circ (\mathcal{E},s)$.
We obtain a morphism of arithmetic extensions
\[
(g,q,\id_F):(\mathcal{E},s)\longrightarrow  (\mathcal{E}',s')
\]
which is universal in the sense that every morphism of arithmetic
extensions
\[
(\alpha,\beta,\gamma):(\mathcal{E},s)\longrightarrow
(\mathcal{E}'',s'')
\]
with $\alpha=g$ factors uniquely through $(g,q,\id_F)$.
Our construction defines a canonical $\Z$-bilinear pairing
\[
\widehat{\rm Ext}_X^1(F,G)× \Hom_{{\mathcal O}_X}(G,G')\longrightarrow
\widehat{\rm Ext}_X^1(F,G')\,\,\,,
\,\,\,(\alpha={\rm cl}(\mathcal{E},s),g)\mapsto
g \circ \alpha={\rm cl}(\mathcal{E}',s').
\]
The formation of $g \circ \alpha$ is functorial in $g$.
It is compatible with the maps $\iota$, $b$, $\nu$, and $\Psi$
in the sense that we have the equalities
\begin{equation}\label{compat1}
\nu(g \circ \alpha)=g\circ \nu(\alpha)\,\,,\,\,
g \circ b(f)=b(\iota(g)\circ f)\,\,,\,\,
\Psi(g\circ \alpha)=g_\C\circ \Psi(\alpha)
\end{equation}
for all $\alpha$, $g$ as above and $f$ in 
$\Hom_{\CinfX}(F_\C,G_\C)^{F_\infty}$.
These are direct consequences of our definitions.

\subsubsection{Pullback}\label{funct3}
Let $h:F'\to F$ a morphism of ${\mathcal O}_X$-modules on an arithmetic scheme
$X$.
For an arithmetic extension $({\mathcal E},s,t)$ of $F$ by an
${\mathcal O}_X$-module $G$,
we obtain a pullback diagram
\[
\begin{array}{ccccccccc}
{\mathcal E}'\,:\,0 &\longrightarrow  &G& \longrightarrow  &E'
&\longrightarrow  &F'&
\longrightarrow  & 0\\
  & & || & & \,\,\,\downarrow {\scriptstyle p} & &
\downarrow {\scriptstyle h}&\\
{\mathcal E}\,:\,0 &\longrightarrow  &G& \longrightarrow  &E&
\longrightarrow  &F&
\longrightarrow  & 0.
\end{array}
\]
The section $t'=t \circ p_\C$ induces a section $s'$ of ${\mathcal
   E}_\C$ (as in (\ref{st})) which turns $({\mathcal E}',s',t')$ into
an
arithmetic extension which we denote as cl$(\mathcal{E},s)\circ h$.
We obtain a morphism of arithmetic extensions
\[
(\id_G,p,h):(\mathcal{E}',s')\longrightarrow  (\mathcal{E},s)
\]
which is universal in the sense that every morphism of arithmetic
extensions
\[
(\alpha,\beta,\gamma):(\mathcal{E}'',s'')\longrightarrow
(\mathcal{E},s)
\]
with $\gamma=h$ factors uniquely through $(\id_G,p,h)$.
Our construction defines a canonical $\Z$-bilinear pairing
\[
  \Hom_{{\mathcal O}_X}(F',F)×\widehat{\rm Ext}_X^1(F,G)\longrightarrow
\widehat{\rm Ext}_X^1(F',G)\,\,\,,
\,\,\,(h,\alpha={\rm cl}(\mathcal{E},s))\mapsto  \alpha\circ
h={\rm cl}(\mathcal{E}',s').
\]
The formation of $\alpha\circ h$ is functorial in $h$.
It is compatible with the maps $\iota$, $b$, $\nu$, and $\Psi$
in the sense that we have the equalities
\begin{equation}\label{compat2}
\nu(\alpha\circ h)= \nu(\alpha)\circ h\,\,,\,\,
b(f)\circ h=b(f \circ \iota(h))\,\,,\,\,
\Psi(\alpha\circ h)=\Psi(\alpha)\circ h_\C
\end{equation}
for all $\alpha$, $h$ as above and $f$ in 
$\Hom_{\CinfX}(F_\C,G_\C)^{F_\infty}$.
These are again a direct consequences of our definitions.

\subsubsection{Compatibilities}
Given an arithmetic extension $(\mathcal{E},s)$
and morphisms $g:G\rightarrow G'$ and $f:F'\rightarrow F$ as above,
there exists a natural
isomorphism
\begin{equation}\label{cpt1}
(g\circ (\mathcal{E},s)) \circ h\stackrel{\sim}{\longrightarrow }
g\circ ((\mathcal{E},s)\circ h).
\end{equation}
Given an arithmetic extension $(\mathcal{E},s)$ as in \refg{ext},
there exists a natural isomorphism
\begin{equation}\label{cpt2}
i\circ (\mathcal{E},s)\stackrel{\sim}{\longrightarrow }
(\mathcal{E},s)\circ p.
\end{equation}
Every morphism of arithmetic extensions
$(\alpha,\beta,\gamma):(\mathcal{E},s)\to (\mathcal{E}',s')$
determines a natural isomorphism
\begin{equation}\label{cpt3}
\alpha \circ (\mathcal{E},s)\stackrel{\sim}{\longrightarrow }
(\mathcal{E}',s')\circ \gamma.
\end{equation}
The isomorphisms \refg{cpt1}, \refg{cpt2}, and \refg{cpt3} are
constructed precisely as in the
algebraic case using the universal properties of pushout and pullback
(compare for example \cite[III.1]{maclane95}).

\subsubsection{Inverse image}\label{funct1}
Let $f:Y\to X$ be a morphism of arithmetic varieties.
Let $F$ and $G$ be ${\mathcal O}_X$-modules.
We assume either that $f$ is flat or that $F$ is a flat ${\mathcal O}_X$-module.
The pullback along $f$ defines a canonical homomorphism
\[
f^*:\widehat{\rm Ext}_X^1(F,G)\longrightarrow  \widehat{\rm
Ext}_Y^1(f^*F,f^*G)
\]
and the formation of $f^*$ is functorial in $f$ and compatible with
$\iota$, $b$, $v$ and $\Psi$ in the
expected way.

\subsubsection{An application}\label{appl}
Given an arithmetic extension $(\mathcal{E},s)$ of $F$ by $G$
and a vector bundle $E$ on the arithmetic scheme $X$, we obtain a
natural
arithmetic extension
\[
(\mathcal{E},s)\otimes E:=(\mathcal{E}\otimes E, s\otimes \id_{E_\C}).
\]
We obtain an induced map
\begin{equation}
  \widehat{\rm Ext}^1_X(F,G)\longrightarrow
\widehat{\rm Ext}^1_X(F\otimes E,G\otimes E)\,\,,\,\,{\rm
cl}(\mathcal{E},s)
\mapsto{\rm cl}\bigl((\mathcal{E},s)\otimes E\bigr)
\end{equation}
which is easily seen to be a group homomorphism.
The pushout and pullback constructions above allow
with the previous remark to construct the following canonical
isomorphism.

\begin{proposition}
For any ${\mathcal O}_X$-modules $F$, $G$ and any vector bundle ${E}$ on an
arithmetic
variety $X$, there is a canonical isomorphism
\[
\widehat{\rm
Ext}^1_X({F},{G}\otimes{E})\stackrel{\sim}{\longrightarrow }
\widehat{\rm Ext}^1_X({F}\otimes {E}^\lor, {G}).
\]
which maps the class of an arithmetic extension $(\mathcal{M},s_M)$
of $F$
by $ G\otimes E$ to the pushout of $\mathcal{M}\otimes E^\vee$
along the evaluation (or trace) map
${\rm id}_G\otimes {\rm ev}:G\otimes E\otimes E^\vee\to G$.
The inverse map sends the class of an arithmetic extension
$(\mathcal{N},s_N)$ to the pullback of $(\mathcal{N},s_N)\otimes E$
along
the canonical morphism
${\rm id}_F\otimes \Delta :F\to F\otimes E^\lor\otimes E$, defined by 
the ``identity" section $\Delta$ of $E^\lor\otimes E.$
\end{proposition}

\proof
Let $(\mathcal{M},s_M)$ be an arithmetic extension with underlying
algebraic extension
\[
\mathcal{M}:0\longrightarrow  G\otimes E\longrightarrow
M\stackrel{p}{\longrightarrow } F\longrightarrow  0.
\]
Consider the pushout
\[
(\mathcal{N},s_N)=\bigl((\mathcal{M},s_M)\otimes E^\vee\bigr) \circ
(\id_G\otimes {\rm ev})
\]
with associated morphism $f:M\otimes E^\vee\rightarrow N$ and section
$s_N=f_{\C}\circ (s_M\otimes \id_{E_\C^\vee})$.
Define $g$ as the composition
\[
g= (f \otimes \id_E)\circ (\id_M\otimes \Delta):
M{\longrightarrow } M\otimes E^\vee\otimes E\longrightarrow N\otimes
E.
\]
We claim that
\begin{equation}\label{hiso1}
(\id_{G\otimes E},g, \id_F\otimes\Delta):(\mathcal{M},s_M)
\longrightarrow   (\mathcal{N},s_N)\otimes E
\end{equation}
is a morphism of arithmetic extensions.
It is easily seen that \refg{hiso1} is a morphism
of the underlying algebraic extensions.
Hence $\mathcal{M}$ is the pullback of $\mathcal{N}\otimes E$
along $ \id_F\otimes\Delta$
and we obtain an induced section $\tilde{s}_M$ of $\mathcal{M}_\C$
as in \ref{funct3}.
We have to show that $s_M$ equals $\tilde{s}_M$.
We conclude from
\begin{eqnarray*}
g_\C\circ \tilde{s}_M\circ p_\C
&=&(s_N\otimes \id_{E_\C})\circ (\id_{F_\C}\otimes\Delta_\C)\circ
p_\C\\
&=&(f_\C\otimes\id_{E_\C})\circ(s_M\otimes\id_{E_\C^\vee\times E_\C})
\circ ( \id_{F_\C}\otimes\Delta_\C)\circ p_\C\\
&=&(f_\C\otimes \id_{E_\C})\circ( \id_{M_\C}\otimes\Delta_\C)
\circ s_M \circ p_\C\\
&=&g_\C\circ {s}_M\circ p_\C
\end{eqnarray*}
that $g_\C\circ (\tilde{s}_M-s_M)\circ p_\C=0$
which proves our claim.
Hence $(\mathcal{M},s_M)$ is the pullback of
$(\mathcal{N},s_N)\otimes E$
along $ \id_F\otimes\Delta$.
This proves one half of our proposition.
The other half follows by a dual argument which we leave to the
reader.
\qed

\subsection{Arithmetic extensions as homomorphisms in the derived
category}\label{derived}

In this paragraph, we present an interpretation of the arithmetic
extension group $\widehat{\rm Ext}_{X}^1(F,G)$ in terms of homological
algebra, in the spirit of the well-known identification of the
``classical" extension group ${\rm Ext}^1_{{\mathcal O}_{X}}(F,G)$, defined by classes
of 1-extensions of ${\mathcal O}_{X}$-modules
equipped with the Baer sum, with the ``cohomological"
extension group $\Hom_{\DOX)}(F, G[1])$ (see Appendix
\ref{extext}).

We follow the notation and conventions concerning derived categories  discussed in Appendix \ref{appa}.

To any ${\mathcal O}_X$-module $E$ is naturally attached the
adjunction map
\[
{\rm ad}_E\colon E\longrightarrow (\rho_* E_\C)^{F_\infty}
\]
with respect to the flat map
\[
\rho=j\circ \kappa
\colon (X_\Sigma(\C),\mathcal{C}_{X_\Sigma}^\infty)\longrightarrow
(X,{\mathcal O}_X)
\]
from \ref{arbach}. It maps any section $e$ of $E$ over an open subscheme 
$U$ of $X$ to the (analytic, hence $\mathcal{C}^\infty)$-section $s_\C$ of 
$E_\C$ over $U_\Sigma(\C),$ which is indeed a section over $U$ of 
$(\rho_* E_\C)^{F_\infty}$. Moreover, for any two 
${\mathcal O}_X$-modules $F$ and $G$, the map
\[
 {\rm
Hom}_{\,\mathcal{C}_{X_\Sigma}^\infty}(F_\C,G_\C)^{F_\infty}
\longrightarrow {\rm Hom}_{{\mathcal O}_X}\bigl(F,(\rho_* G_\C)^{F_\infty}\bigr)
\]
which maps an $F_\infty$-invariant $\mathcal{C}^\infty$-morphism 
$\varphi:F_\C\rightarrow G_\C$ to the composition
\[
\tilde{\varphi} \colon F\stackrel{{\rm ad}_F}{\longrightarrow} (\rho_*F_\C)^{F_\infty}
\stackrel{\rho_*(\varphi)}{\rightarrow} (\rho_*G_\C)^{F_\infty}
\]
is an isomorphism. (Its inverse 
 maps a morphism of ${\mathcal O}_X$-modules
$\psi:F\rightarrow (\rho_*G_\C)^{F_\infty}$ to the composition
\[
F_\C=\rho^*F\stackrel{\rho^*(\psi)}{\longrightarrow} \rho^*\rho_*G_\C\stackrel{{\rm
ad}^{G_\C}}{\rightarrow}G_\C.)
\]

Let $(\mathcal{E},s)$ be an arithmetic extension, 
\ie, an extension of ${\mathcal O}_X$-modules
\[
\mathcal{E}\colon 0\longrightarrow G\stackrel{i}{\longrightarrow}
E\stackrel{p}{\longrightarrow}F\rightarrow 0
\]
together with a $F_\infty$-invariant
$\mathcal{C}^\infty$-splitting $s:F_\C\rightarrow E_\C$.
Let $t:G_\C\rightarrow E_\C$ be the
$F_\infty$-invariant
$\mathcal{C}^\infty$-section of $i_\C$ such that such 
that ${\rm Id}_{E_\C}=s\circ p_\C+i_\C\circ t.$
Let us consider the following diagram of complexes (written horizontally, 
with morphisms written vertically) of ${\mathcal O}_X$-modules:
\begin{equation}\label{arext-qis}
\begin{array}{ccc}
& & F \\
& &\hspace*{3mm}\uparrow  {\scriptstyle p}  \\
G & \stackrel{i}{\longrightarrow} & E \\
\hspace*{6mm}\downarrow  {\scriptstyle{\rm Id}_G}
& &  \hspace*{3mm}\downarrow {\scriptstyle \tilde{t}}  \\
G & \stackrel{{\rm ad}_G}{\longrightarrow}
&  (\rho_* G_\C)^{F_\infty}
\end{array}
\end{equation}
where $F$, $E$, and  $(\rho_* G_\C)^{F_\infty}$ sit in degree zero, and $G$ in 
degree $-1$. The last two lines are precisely the cones $C(i)$ and 
$C({\rm ad}_G)$ of $i$ and ${\rm ad}_G$, and the  
map $p$ defines a quasi-isomorphism $\pp:C(i)\rightarrow F$.
Instead of (\ref{arext-qis}), we may write as well
\[
F\stackrel{\pp}{\longleftarrow} C(i)
\stackrel{({\rm Id}_G,\,\tilde{t})}{\longrightarrow} C({\rm ad}_G).
\]
This diagram defines a morphism
\[
\partial_{(\mathcal{E},s)}:= ({\rm Id}_G,\,
\tilde{t})\circ \pp^{-1} \colon F\longrightarrow C({\rm ad}_G)
\]
in the derived category $D^+({\mathcal O}_X\mbox{\rm 
-mod})$ of bounded below complexes of
${\mathcal O}_X$-modules. Clearly, if $(\mathcal{E},s)$ and $(\mathcal{E}',s')$ are 
isomorphic arithmetic extensions of $F$ by $G,$ the associated morphisms 
$\partial_{(\mathcal{E},s)}$ and $\partial_{(\mathcal{E}',s')}$ coincide.

\begin{proposition}\label{clariso}
For any two ${\mathcal O}_X$-modules $F$ and $G$, the mapping
\[
\clar_{F,G} \colon \widehat{\rm Ext}_X^1(F,G)\longrightarrow
{\rm Hom}_{D^+({\mathcal O}_X\mbox{\scriptsize\rm -mod})}\bigl(F,C({\rm
ad}_G)\bigr)\,,
\,\,\bigl[(\mathcal{E},s)\bigr]\mapsto \partial_{(\mathcal{E},s)}.
\]
is an isomorphism of abelian groups.
\end{proposition}

As hinted to above, the map $\clar_{F,G}$ is an avatar of the
classical isomorphism
\begin{equation}\label{canisoext}
\cl_{F,G} \colon {\rm Ext}_{{\mathcal O}_{X}}^1(F,G)\stackrel{\sim}{\longrightarrow}
{\rm Hom}_{D^+({\mathcal O}_X\mbox{\scriptsize\rm -mod})}\bigl(F,G[1]\bigr)\,,
\,\,[\mathcal{E}]\mapsto \partial_{\mathcal{E}}
\end{equation}
described in Appendix \ref{extext}. Recall that, to the class
$[\mathcal{E}]$ of a $1$-extension of ${\mathcal O}_X$-modules
\[
\mathcal{E}\colon 0\longrightarrow G\stackrel{i}{\longrightarrow}
E\stackrel{p}{\longrightarrow}F\longrightarrow 0,
\]
it associates its ``boundary operator"
$\partial_{\mathcal{E}}$ in 
${\rm Hom}_{D^+({\mathcal O}_X\mbox{\scriptsize\rm -mod})}\bigl(F,G[1]\bigr)$
that is defined as follows.
We may consider 
the diagram of complexes of ${\mathcal O}_{X}$-modules (written horizontally)
\begin{equation}\label{ext-qis}
\begin{array}{ccc}
& & F \\
& &\hspace*{3mm}\uparrow  {\scriptstyle p}  \\
G & \stackrel{i}{\longrightarrow} & E  \\
\hspace*{5mm}\downarrow  {\scriptstyle{\rm -Id}_G}  & &\\
G & &
\end{array}
\end{equation}
where $E$ and $F$ sit in degree 0, and $G$ in degree $-1.$ Its last
two lines are $C(i)$ and $G[1]$, the upper map defines a
quasi-isomorphism $\pp:C(i) \rightarrow F,$ and \refg{ext-qis} may be 
also written
$$F\stackrel{\pp}{\longleftarrow} C(i) \stackrel{(-{\rm
Id}_{G},0)}{\longrightarrow} G[-1].$$ Then, by definition, 
$\partial_\mathcal{E}$ is the morphism in 
$D^+({\mathcal O}_X\mbox{\rm -mod})$ defined by this last diagram:
\[
\partial_{\mathcal{E}}:= (-{\rm Id}_{G},0) \circ \pp^{-1}
\colon F\rightarrow G[1].
\]

Observe that, if we are given two  ${\mathcal O}_X$-modules $F$ and $G$, we may 
consider the distinguished triangle in $D^+({\mathcal O}_X\mbox{\rm -mod})$ attached 
to ${\rm ad}_G$:
\begin{equation}
G\stackrel{{\rm ad}_G}{\longrightarrow} (\rho_* G_\C)^{F_\infty}
\stackrel{\alpha}{\longrightarrow}
C({\rm ad}_G)\stackrel{\beta}{\longrightarrow} G[1]
\end{equation}
--- we use the sign conventions discussed in Appendix \ref{sign}; thus 
$\alpha$ (resp. $\beta$) is the morphism defined by 
${\rm Id}_{(\rho_* G_\C)^{F_\infty}}$ (resp. by ${\rm -Id}_G$) --- and apply 
the functor ${\rm Hom}_{D^+({\mathcal O}_X\mbox{\scriptsize\rm -mod})}\bigl(F,.)$. 
This yields the following long exact sequence:
\begin{multline}\label{longexa3}
{\rm Hom}_{D^+({\mathcal O}_X\mbox{\scriptsize\rm -mod})}(F,G)\stackrel{{\rm ad}_G\, \circ .}{\longrightarrow}
{\rm Hom}_{D^+({\mathcal O}_X\mbox{\scriptsize\rm 
-mod})}\bigl(F,(\rho_* G_\C)^{F_\infty}\bigr)
\stackrel{\alpha \,\circ .}{\longrightarrow}
{\rm Hom}_{D^+({\mathcal O}_X\mbox{\scriptsize\rm -mod})}\bigl(F,C({\rm
ad}_G)\bigr) \\
\stackrel{\beta\, \circ.}{\longrightarrow} {\rm 
Hom}_{D^+({\mathcal O}_X\mbox{\scriptsize\rm 
-mod})}\bigl(F,G[1]\bigr)
\stackrel{{\rm ad}_G[1]\, \circ.}{\longrightarrow}
{\rm Hom}_{D^+({\mathcal O}_X\mbox{\scriptsize\rm 
-mod})}\bigl(F,(\rho_* G_\C)^{F_\infty}[1]\bigr).
\end{multline}

We shall show simultaneously that $\clar_{F,G}$ is an isomorphism and that 
the exact sequence \refg{longexa3} is naturally isomorphic with the following 
variant of the long exact sequence in Theorem \ref{longexseq1}:
\begin{multline}\label{longexseq1bis}
{\rm Hom}_{{\mathcal O}_X}(F,G)
\stackrel{\iota}{\longrightarrow }
{\rm Hom}_{\CinfX}(F_\C,G_\C)^{F_\infty}
\stackrel{b}{\longrightarrow }
\widehat{\rm Ext}_X^1(F,G)
\stackrel{\nu}{\longrightarrow } {\rm Ext}_{{\mathcal O}_X}^1(F,G)\\
\stackrel{{\rm ad}_G\, \circ.}{\longrightarrow}  {\rm Ext}^1_{{\mathcal O}_X}(F_\C,(\rho_*G_\C)^{F_\infty}).
   \end{multline}
   
   The exactness of \refg{longexseq1bis} follows from Theorem \ref{longexseq1} 
combined with the following observation:
   
   \begin{lemma} For any extension of ${\mathcal O}_X$-modules
\[
\mathcal{E}\colon 0\longrightarrow G\stackrel{i}{\longrightarrow}
E\stackrel{p}{\longrightarrow}F\longrightarrow 0,
\]
an $F_\infty$-invariant $\mathcal{C}^\infty$ splitting of $\cE_\C$  may be 
equivalently described by:
\begin{itemize}
\item 
a morphism $t$ in ${\rm Hom}_{\CinfX}(E_\C,G_\C)^{F_\infty}$ such that  
$t\circ i_\C={\rm id}_{G_\C};$

\item a morphism   
$\tilde t$ in ${\rm Hom}_{{\mathcal O}_X}\bigl(E,(\rho_*G_\C)^{F_\infty}\bigr)$ such that 
$\tilde t\circ i={\rm ad}_G;$

\item a splitting of the extension of ${\mathcal O}_X$-modules defined as the pushout 
of $\cE$ by ${\rm ad}_G:$
\[
\begin{array}{ccccccccc}
0 & \longrightarrow & G&\stackrel{i}{\longrightarrow}&
E&\stackrel{p}{\longrightarrow}&F&\longrightarrow& 0\\
&& \hspace*{6mm}\downarrow  {\scriptstyle {\rm ad}_G} & & \downarrow
& & ||\\
0 & \longrightarrow & (\rho_*G_\C)^{F_\infty}&{\longrightarrow} &
E'&{\longrightarrow}&F&\longrightarrow& 0.
\end{array}
\]
\end{itemize}
Consequently, the homomorphisms $$\digamma:
{\rm Ext}_{{\mathcal O}_X}^1(F,G)
{\longrightarrow } {\rm
   Ext}^1_{\CinfX}(F_\C,G_\C)$$
and
$${\rm ad}_G \,\circ . :{\rm Ext}_{{\mathcal O}_X}^1(F,G)
{\longrightarrow}
 {\rm Ext}^1_{{\mathcal O}_X}(F_\C,(\rho_*G_\C)^{F_\infty})$$
   have the same kernel.
   \end{lemma}
   \endproof

Actually, the five lemma shows that Proposition \ref{clariso} follows from 
the long exact sequences \refg{eq:longexseq1} and \refg{longexseq1bis}, and 
Lemma \ref{clarhomo} and Proposition \ref{prep5lem} below.

\begin{lemma}\label{clarhomo}
The mapping $\clar_{F,G}$ in Proposition \ref{clariso} is a homomorphism of groups.
\end{lemma}

\proof
Consider arithmetic extensions $(\mathcal{E}_j,s_j)$, $j=1,2$ with
underlying extensions
\[
\mathcal{E}_j\colon 0\longrightarrow G\stackrel{i_j}{\longrightarrow}
E_j\stackrel{p_j}{\longrightarrow}F\longrightarrow 0,
\]
and $F_\infty$-invariant $\mathcal{C}^\infty$ splittings $s_j:F_\C\rightarrow
E_{j,\C}$.
Let $(\mathcal{E},s)$ denote the arithmetic Baer 
sum of $(\mathcal{E}_1,s_1)$ and
$(\mathcal{E}_2,s_2)$ as defined in \ref{basicdef}.
We have to show
\begin{equation}\label{eqindercat}
\partial_{(\mathcal{E},s)}=\partial_{(\mathcal{E}_1,s_1)}
+\partial_{(\mathcal{E}_2,s_2)}.
\end{equation}

Consider the following diagrams
\begin{equation}\label{arext-d1}
\begin{array}{ccc}
G\oplus G& \stackrel{w}{\longrightarrow}&
\ker(E_1\oplus E_2\stackrel{p_1-p_2}{\longrightarrow} F) \\
\hspace*{4mm}\downarrow  {\scriptstyle q_j}
& &\hspace*{3mm}\downarrow  {\scriptstyle q_j}  \\
G & \stackrel{i_j}{\longrightarrow} & E_j
\end{array}
\end{equation}
for $j=1,2$, and
\begin{equation}\label{arext-d2}
\begin{array}{ccc}
G\oplus G& \stackrel{w}{\rightarrow}&
\ker(E_1\oplus E_2\stackrel{p_1-p_2}{\longrightarrow} F) \\
\hspace*{4mm}\downarrow  {\scriptstyle +}
& &\hspace*{5mm}\downarrow  {\scriptstyle {\rm can}}  \\
G & \stackrel{i}{\longrightarrow} & E
\end{array}
\end{equation}
where $q_j$ denotes the $j$-th projection, $w$ the direct sum $i_1 \oplus i_2$, and `can' 
the
canonical projection.
These diagrams 
are commutative
and induce
quasi-isomorphisms $$\qq_j:C(w)\rightarrow C(i_j) \mbox{ and } a:C(w)\rightarrow C(i).$$

As before, let us introduce
 $t:G_\C\rightarrow E_{\C}$ (resp.
$t_j:G_\C\rightarrow E_{j,\C}$) such that  
${\rm Id}_{E_\C}=s\circ p_\C+i_\C\circ t$ (resp.
${\rm Id}_{E_{j,\C}}=s_j\circ p_{j,\C}+i_{j,\C}\circ t_j$).
According to the very definition of the arithmetic Baer sum,
\[
t\circ {\rm can}_\C = t_1\oplus t_2 .
\]
With obvious notation, we have:
\[
\partial_{(\mathcal{E},s)}\colon
F\stackrel{\pp}{\longleftarrow} C(i) \stackrel{a}{\longleftarrow} C(w)
\stackrel{a}{\longrightarrow} C(i)
\stackrel{({\rm Id}_G,\,\tilde{t})}{\longrightarrow} C({\rm ad}_G)
\]
and
\[
\partial_{(\mathcal{E}_j,s_j)}\colon
F\stackrel{\pp_j}{\longleftarrow} C(i_j) \stackrel{\qq_j}{\longleftarrow} C(w)
\stackrel{\qq_j}{\longrightarrow} C(i)
\stackrel{({\rm id}_G,\,\tilde{t_j})}{\longrightarrow} C({\rm ad}_G).
\]

We finally obtain (\ref{eqindercat}) from the straightforward equalities
\[
\pp_1\circ \qq_1=\pp\circ a=\pp_2\circ \qq_2   \colon C(w)\longrightarrow F
\]
and
\[
({\rm Id}_G,\tilde{t_1})\circ \qq_1+({\rm Id}_G,\tilde{t_2})\circ \qq_2
=({\rm Id}_G,\tilde{t})\circ a \colon C(w)\longrightarrow C({\rm ad}_G)
\]
of homomorphisms of complexes of ${\mathcal O}_X$-modules.
\qed

Recall that the category of ${\mathcal O}_X$-modules may be identified with a full 
subcategory of its derived category, and that consequently, for any two 
${\mathcal O}_X$-modules $F_1$ and $F_2$, we may identify
${\rm Hom}_{{\mathcal O}_X}(F_1,F_2)$ and ${\rm Hom}_{D^+({\mathcal O}_X\mbox{\scriptsize\rm -mod})}(F_1,F_2).$

Besides, as a special instance of the canonical isomorphism $\cl_{F,G}$ 
(with  $(\rho_*G_\C)^{F_\infty}$ instead of $G$), we may consider the
isomorphism
\[
\cl_{F,(\rho_*G_\C)^{F_\infty}}\colon{\rm Ext}_{{\mathcal O}_X}^1(F,
(\rho_*G_\C)^{F_\infty})\stackrel{\sim}{\longrightarrow}
{\rm Hom}_{D^+({\mathcal O}_X\mbox{\scriptsize\rm -mod})}
\bigl(F,(\rho_*G_\C)^{F_\infty}[1]\bigr).
\]

\begin{proposition}\label{prep5lem}
The following four diagrams are commutative:
\begin{equation}\label{prems}
\begin{array}{ccc}
{\rm Hom}_{{\mathcal O}_X}(F,G) & \stackrel{\iota}{\longrightarrow} &
{\rm Hom}_{\mathcal{C}_X^\infty}(F_\C,G_\C)^{F_\infty}\\
\| \hspace*{3mm}
& &  {\scriptstyle \sim}\downarrow {\scriptstyle T \mapsto \tilde{T}}  \\
{\rm Hom}_{D^+({\mathcal O}_X\mbox{\scriptsize\rm -mod})}(F,G) &\stackrel{{\rm ad}_G\,\circ.}{\longrightarrow} &
{\rm Hom}_{D^+({\mathcal O}_X\mbox{\scriptsize\rm 
-mod})}\bigl(F,(\rho_* G_\C)^{F_\infty}\bigr)
\end{array}
\end{equation}

\begin{equation}\label{deuz}
\begin{array}{ccc}
{\rm Hom}_{\mathcal{C}_X^\infty}(F_\C,G_\C)^{F_\infty} &
\stackrel{-b}{\longrightarrow} &
  \widehat{\rm Ext}_X^1(F,G)\\
{\scriptstyle \sim}\downarrow {\scriptstyle T \mapsto \tilde{T}}
& &  
\downarrow {\scriptstyle \clar_{F,G}}  \\
{\rm Hom}_{D^+({\mathcal O}_X\mbox{\scriptsize\rm -mod})}
\bigl(F,(\rho_* G_\C)^{F_\infty}\bigr) &\stackrel{\alpha\, \circ .}{\longrightarrow} &
{\rm Hom}_{D^+({\mathcal O}_X\mbox{\scriptsize\rm -mod})}\bigl(F,C({\rm
ad}_G)\bigr)\\
\end{array}
\end{equation}

\begin{equation}\label{troiz}
\begin{array}{ccc}
\widehat{\rm Ext}_X^1(F,G) & \stackrel{\nu}{\longrightarrow}
& {\rm Ext}^1_{{\mathcal O}_X}(F,G) \\
\downarrow  {\scriptstyle\clar_{F,G}}
& & {\scriptstyle \sim} 
\downarrow {\scriptstyle \cl_{F,G}}\\
{\rm Hom}_{D^+({\mathcal O}_X\mbox{\scriptsize\rm -mod})}\bigl(F,C({\rm
ad}_G)\bigr)
&\stackrel{\beta\, \circ .}{\longrightarrow}&
{\rm Hom}_{D^+({\mathcal O}_X\mbox{\scriptsize\rm -mod})}\bigl(F,G[1]\bigr)\\
\end{array}
\end{equation}

\begin{equation}\label{der}
\begin{array}{ccc}
{\rm Ext}^1_{{\mathcal O}_X}(F,G) & \stackrel{{\rm ad}_G\, \circ.}{\longrightarrow }& {\rm
Ext}^1_{{\mathcal O}_X}
\bigl(F,(\rho_* G_\C)^{F_\infty}\bigr) \\
{\scriptstyle \sim}\downarrow {\scriptstyle \cl_{F,G}}
& & {\scriptstyle \sim}\downarrow {\scriptstyle \cl_{F,(\rho_*G_\C)^{F_\infty}}} \\
{\rm Hom}_{D^+({\mathcal O}_X\mbox{\scriptsize\rm 
-mod})}\bigl(F,G[1]\bigr) & \stackrel{{\rm ad}_G[1]\, \circ .}{\longrightarrow} &
{\rm Hom}_{D^+({\mathcal O}_X\mbox{\scriptsize\rm 
-mod})}\bigl(F,(\rho_* G_\C)^{F_\infty}[1]\bigr).
\end{array}
\end{equation}

\end{proposition}

\proof
The commutativity of \refg{prems} is immediate. Let us check the one of \refg{deuz}.
Given
\[
T\in {\rm Hom}_{\mathcal{C}_{X_\Sigma}^\infty}(F_\C,G_\C)^{F_\infty},
\]
we equip the trivial extension
\[
\mathcal{E}\colon
0\longrightarrow G\stackrel{\genfrac{(}{)}{0pt}{}{{\rm
id}_G}{0}}{\longrightarrow}
G\oplus F \stackrel{(0\,{{\rm id}_F})}{\longrightarrow} F\longrightarrow 0
\]
with the $\mathcal{C}^\infty$-splitting defined by 
\[
s=\genfrac{(}{)}{0pt}{}{-T}{{\rm id}_{F_\C}} \mbox{ and } t=({{\rm
id}_{G_\C}}\; T).
\]
Then $\tilde{t}=({\rm ad}_G \,\, \tilde{T})$, the arithmetic extension
$(\mathcal{E},s)$
admits $-b(T)$ as class in $ \widehat{\rm Ext}_X^1(F,G)$, 
and $\partial_{(\mathcal{E},s)}$ is the morphism in 
${\rm Hom}_{D^+({\mathcal O}_X\mbox{\scriptsize\rm -mod})}\bigl(F,G[1]\bigr)$ 
defined by the diagram:
\[
\begin{array}{ccc}
&& F \\
& &\hspace*{10mm}\uparrow  {\scriptstyle (0\,{{\rm id}_F})}  \\
G & \stackrel{\genfrac{(}{)}{0pt}{}{{\rm id}_G}{0}}{\longrightarrow} &
G\oplus F \\
\hspace*{5mm}\downarrow  {\scriptstyle{\rm id}_G}
& &  \hspace*{4mm}\downarrow {\scriptstyle \tilde{t}}  \\
G & \stackrel{{\rm ad}_G}{\longrightarrow}
&  (\rho_* G_\C)^{F_\infty}.
\end{array}
\]
A quasi-inverse of $(0\,{{\rm id}_F})$ in this diagram is given by
\[
\begin{array}{ccc}
&& F \\
& &\hspace*{8mm}\downarrow  {\scriptstyle 
\genfrac{(}{)}{0pt}{}{0}{\rm id_F}}  \\
G & \stackrel{\genfrac{(}{)}{0pt}{}{{\rm id}_G}{0}}{\longrightarrow} &
G\oplus F.
\end{array}
\]
Consequently  $\partial_{(\mathcal{E},s)}$ is also defined by the map of complexes:
\[
\begin{array}{ccc}
& & F\\
&&\hspace*{2mm}\downarrow  {\scriptstyle \tilde{T}}\\
  G & \stackrel{{\rm ad}_G}{\longrightarrow}
&  [\rho_* G_\C]^{F_\infty} .
\end{array}
\]
This shows that $\clar_{F,G}(-b(T))= \alpha \,\circ \tilde{T}$, as required.

The commutativity of \refg{troiz}
is a direct consequence of our definitions and sign conventions.

Finally the diagram  \refg{der} commutes since
pushout on extension groups corresponds
under the isomorphism (\ref{canisoext}) to the natural functoriality
of the
${\rm Hom}_{D^+({\mathcal O}_X\mbox{\scriptsize\rm -mod})}$ functor in the
second
argument (compare Appendix \ref{extext}).
\qed

\subsubsection{}
We give a homological description of the operator $\Psi$ which
computes the
second fundamental form.

For any ${\mathcal O}_X$-module $G$ and any non-negative integer $j$, we shall denote
\[
A^{0,j}(-_\R,G)=\bigl(\rho_* A^{0,j}(-,G)\bigr)^{F_\infty}
\]
as in \ref{aforms}.
We have in particular $A^{0,0}(-_\R,G)=(\rho_* G_\C)^{F_\infty}$, and, for any 
coherent ${\mathcal O}_X$-module $F$, 
${\rm Hom}_{{\mathcal O}_X}(F, A^{0,j}(-_\R,G))$ may be identified with
$A^{0,j}(X_\R, \iho).$
The Dolbeault operator $\overline{\partial}$ induces
a complex of ${\mathcal O}_X$-modules
\[
(A^{0,.}(-_\R,G), \overline{\partial}_G) \colon
0\longrightarrow
A^{0,0}(-_\R,G)\stackrel{\overline{\partial}_G}{\longrightarrow}
A^{0,1}(-_\R,G)\stackrel{\overline{\partial}_G}{\longrightarrow}
A^{0,2}(-_\R,G)\stackrel{\overline{\partial}_G}{\longrightarrow}
\ldots
.
\]
The diagram
\[
\begin{array}{ccccccc}
G & \stackrel{{\rm ad}_G}{\longrightarrow} &A^{0,0}(-_\R,G)&&\\
& &  \hspace*{4mm}\downarrow {\scriptstyle \overline{\partial}_G} \\
& & A^{0,1}(-_\R,G) & \stackrel{-\overline{\partial}_G}{\longrightarrow} &
A^{0,2}(-_\R,G) & \stackrel{-\overline{\partial}_G}{\longrightarrow} &
\ldots
\end{array}
\]
defines a morphism of complexes
\[
\Psi'\colon C({\rm ad}_G)\longrightarrow 
(\sigma_{\ge 1}(A^{0,.}(-_\R,G),{\overline{\partial}_G}))[1]
\]
where $\sigma_{\ge 1}$ denotes the naive truncation.

Let $F$ and $G$ be two ${\mathcal O}_X$-modules, with $F$ coherent. Then the identification of 
$$A^{0,1}(X_\R, \iho)$$
and
$$
{\rm Hom}_{{\mathcal O}_X}(F, A^{0,1}(-_\R,G))$$
defines, by restriction, a canonical isomorphism
\begin{equation}\label{canhommap1}
Z^{0,1}_{\overline{\partial}}\bigl(X_\R,{\mathcal
Hom}_{{\mathcal O}_X}(F,G)\bigr)
={\rm Hom}_{{\mathcal O}_X}\Bigl(F,\,{\rm Ker}\,\bigl(\overline{\partial}_G:
A^{0,1}(-_\R,G)\rightarrow  A^{0,2}(-_\R,G)\bigr)\Bigr),
\end{equation}
hence a morphism of abelian groups
\begin{equation}\label{canhommap2}
z: Z^{0,1}_{\overline{\partial}}\bigl(X_\R,{\mathcal
Hom}_{{\mathcal O}_X}(F,G)\bigr)
\rightarrow
{\rm Hom}_{D^+({\mathcal O}_X-{\scriptstyle\rm mod})}\bigl(F,(\sigma_{\ge 1}(A^{0,.}(-_\R,G)
,{\overline{\partial}_G}))[1])\bigr).
\end{equation}

\begin{lemma}

i)
The sheaves $A^{0,j}(\ubar,G)$ over $X_\Sigma(\C)$ are fine. For any $k>0$, we have
\begin{equation}\label{finevanish1}
R^k\rho_*A^{0,j}(_\ubar,G)=0
\end{equation}
and
\begin{equation}\label{finevanish2}
R^k\rho_*\sigma_{\ge 1}(A^{0,.}(-_\R,G),{\overline{\partial}_G})=0.
\end{equation}
ii)
The morphism $z$ in  (\ref{canhommap2}) is an isomorphism
if $F$ is a vector bundle over $X$.
\end{lemma}

\proof
i) The sheaves $A^{0,j}(-, \C)$ on the paracompact Hausdorff space
$X_\Sigma(\C)$ are fine. It follows that the tensor products
$A^{0,j}(-,G)$ are fine.
We can compute $R^k\rho_*A^{0,j}(-,G)$
as the sheaf associated with the presheaf
\[
U\mapsto H^k\bigl((U_\Sigma(\C),A^{0,j}(\ubar,G)\bigr).
\]
The latter groups vanish as fine sheaves are 
$\Gamma\bigl(U_\Sigma(\C),-\bigr)$-acyclic. This proves \refg{finevanish1}, 
which in turn immediately implies \refg{finevanish2}.

ii) To prove that $z$ is an isomorphism when $F$ is a vector bundle, we may 
replace $F$ by ${\mathcal O}_X$ and $G$ by
$F^\lor \otimes G$. Hence we may assume that $F$ is ${\mathcal O}_X$.
To simplify notations, let us write 
$\sigma_{\geq 1}A^{0,\cdot}(X_\R,G)$
for $\sigma_{\ge 1}(A^{0,.}(-_\R,G),{\overline{\partial}_G}$.
The Leray spectral sequence
\[
E_2^{j,k}=H^{j}\bigl(X,R^k\rho_*
\sigma_{\geq 1}A^{0,\cdot}(-,G)\bigr)
\Rightarrow H^{j+k}\bigl(X_\Sigma(\C),\sigma_{\geq 1}A^{0,\cdot}(-,G)\bigr)
\]
degenerates at $E_2$ by i).
Hence
$$H^1\bigl(X,\sigma_{\geq 1} \rho_*A^{0,\cdot}(-,G)\bigr)
= H^1\bigl(X_\Sigma(\C),\sigma_{\geq 1}A^{0,\cdot}(-,G)\bigr),$$
and consequently,
\begin{eqnarray*}
{\rm Hom}_{D^+({\mathcal O}_X-{\scriptstyle\rm mod})}\bigl({\mathcal O}_X,(\sigma_{\ge 1}
A^{0,.}(-_\R,G))[1]\bigr)
&=& H^1\bigl(X,\sigma_{\geq 1} \rho_*A^{0,\cdot}(-,G)\bigr)^{F_\infty}\\
&=& H^1\bigl(X_\Sigma(\C),\sigma_{\geq 1}A^{0,\cdot}(-,G)\bigr)^{F_\infty}\\
&=&Z^{0,1}_{\overline{\partial}}(X_\Sigma(\C), G_\C)^{F_\infty}\\
&=& Z^{0,1}_{\overline{\partial}}\bigl(X_\R,G\bigr)
\end{eqnarray*}
as the sheaves $A^{0,\cdot}(-,G)$ are
$\Gamma\bigl(X_\Sigma(\C),-\bigr)$-acyclic.
\endproof

\begin{proposition} For any two vector bundles $F$ and
$G$ on $X$, the following diagram is commutative:
\[
\begin{array}{ccc}
\widehat{\rm Ext}_X^1(F,G) & \stackrel{\Psi}{\longrightarrow} &
Z^{0,1}_{\overline{\partial}}\bigl(X_\R,{\mathcal
Hom}_{{\mathcal O}_X}(F,G)\bigr) \\
\hspace*{3mm}\downarrow  {\scriptstyle\clar_{F,G}}
& & {\scriptstyle \sim}\downarrow  {-z}    \\
{\rm Hom}_{D^+({\mathcal O}_X\mbox{\scriptsize\rm -mod})}\bigl(F,C({\rm
ad}_G)\bigr)
&\stackrel{\Psi'\circ\, .}{\longrightarrow} &
{\rm Hom}_{D^+({\mathcal O}_X\mbox{\scriptsize\rm -mod})}\bigl(F,(\sigma_{\ge 1}(A^{0,.}(-_\R,G),{\overline{\partial}_G}))[1]\bigr)
\end{array}
\]
is commutative.
\end{proposition}

\proof  Let  $(\mathcal{E},s)$ be an arithmetic extension and, as before, define $t$
by (\ref{st}). 
The composite map $\Psi'\circ \clar_{F,G}$ sends the class of
$(\mathcal{E},s)$ in $\widehat{\rm Ext}_X^1(F,G)$ to 
the morphism 
$\Psi'\circ \partial_{(\mathcal{E},s)}$ in $D^+({\mathcal O}_X\mbox{\scriptsize\rm -mod})$, 
which is defined by the diagram
\[
\begin{array}{ccccccc}
& & F \\
& &\hspace*{3mm}\uparrow  {\scriptstyle p}  \\
G & \stackrel{i}{\longrightarrow} & E \\
\hspace*{6mm}\downarrow  {\scriptstyle\rm Id_G}
& &  \hspace*{2mm}\downarrow {\scriptstyle \tilde{t}}  \\
G & \stackrel{{\rm ad}_G}{\longrightarrow}
&  (\rho_* G_\C)^{F_\infty} \\
& &  \hspace*{4mm}\downarrow {\scriptstyle \overline{\partial}_G} \\
& & A^{0,1}(-_\R,G) &
\stackrel{-\overline{\partial}_G}{\longrightarrow} &
A^{0,2}(-_\R,G) & \stackrel{-\overline{\partial}_G}{\longrightarrow} &
\ldots
\end{array}
\]
It may also be written as
\begin{equation}\label{mordercat}
\Psi'\circ \partial_{(\mathcal{E},s)} : F\stackrel{\pp}{\longleftarrow} C(i)
\stackrel{({\rm Id}_G, \tilde{t})}{\longrightarrow}C({\rm ad}_G)
\stackrel{\overline{\partial}_G}{\longrightarrow}
(\sigma_{\ge 1}A^{0,\cdot}(-_\R,G))[1].
\end{equation}

Let us denote $
\alpha:=\Psi (\mathcal{E}, s).$ This is the element of 
$Z_{\overline{\partial}}^{0,1}\bigl(X_\R,\ihob\bigr)$ characterized by the relation
$$\overline{\partial}_E (s.f)= s.\overline{\partial}_F +i_\C (\alpha.f),$$
for any local  section $f$ of 
$A^{0,\cdot}(-,F_\C)$ over $X_\Sigma(\C).$ We also have
\begin{equation}\label{dt}
\ol{\partial}_{E^\lor \otimes G}t=-\alpha \circ p_\C.
\end{equation}

To any
 section $\beta$ of $A^{0,j}(U_\R,F)$ over some open $U$ in
$X$,
we may attach the element 
$$w(\beta):=(-\alpha(\beta),s(\beta))$$
of degree $j$
in the cone of the morphism $A^{0,\cdot}(-_\R,i)$ from 
$\bigl(A^{0,\cdot}(-_\R,G), \ol{\partial}_G\bigr)$ to
$\bigl(A^{0,\cdot}(-_\R,E), \ol{\partial}_E\bigr)$ defined by $i_\C$.
The differential of $w(\beta)$ in the complex $C\bigl(A^{0,\cdot}(-_\R,i)\bigr)$ is
\begin{eqnarray*}
d(-\alpha.\beta,s.\beta)
&=& \bigl(-\overline{\partial}_G(-\alpha.\beta), i_\C(-\alpha.\beta)
+\overline{\partial}_E (s.\beta)\bigr)\\
&=&\bigl(-\alpha . (\overline{\partial}_F\beta),s.
(\overline{\partial}_F.\beta)\bigr).
\end{eqnarray*}
Consequently $w$ defines a homomorphism of complexes
\[
w\colon A^{0,\cdot}(-_\R,F)\longrightarrow
C\bigl(A^{0,\cdot}(-_\R,i)\bigr),
\]
where, as before, we write $A^{0,\cdot}(-_\R,F)$ for $(A^{0,\cdot}(-_\R,F),\ol{\partial}_F).$
It is straightforward that it is a right inverse of the morphism of complexes 
$$ A^{0,\cdot}(-_\R,\pp)\colon
C\bigl(A^{0,\cdot}(-_\R,i)\bigr)
\longrightarrow
 A^{0,\cdot}(-_\R,F)$$
deduced from the quasi-isomorphism $\pp: C(i) \rightarrow F$ by considering 
the associated Dolbeault complexes.

Observe that --- since $\alpha$ is $\ol{\partial}$ closed --- one defines a 
morphism of complexes of ${\mathcal O}_X$-modules
$${A}:  
C\bigl(A^{0,\cdot}(-_\R,i)\bigr) \longrightarrow
 (\sigma_{\ge 1}A^{0,\cdot}(-_\R,G))[1]$$
 by mapping a sections $(\gamma, \delta)$ of 
 $$C\bigl(A^{0,\cdot}(-_\R,i)\bigr)^k:=
A^{0,k+1}(-_\R,G)
 \oplus
A^{0,k}(-_\R,E)$$
to the section $\alpha.p(\delta)$ of $A^{0,k+1}(-_\R,G)$.
Moreover the relation \refg{dt} shows that the right-hand square of the 
following diagram is commutative:
\[
\begin{array}{ccccc}
F& \stackrel{\pp}{\longleftarrow} &C(i)&
\stackrel{({\rm Id}_G,\tilde{t})}{\longrightarrow}&C({\rm ad}_G) \\
\hspace*{6mm}\downarrow  {\scriptstyle {\rm ad}_F} &&
\hspace*{16mm}\downarrow  {\scriptstyle ({\rm ad}_G, {\rm ad}_E)} &&
\hspace*{6mm}\downarrow  {\scriptstyle \overline{\partial}_G} \\
A^{0,\cdot}(-_\R,F)& \stackrel{ A^{0,\cdot}(-_\R,\pp)}{\longleftarrow}
&C\bigl(A^{0,\cdot}(-_\R,i)\bigr)&
\stackrel{-A}{\longrightarrow}
& (\sigma_{\ge 1}A^{0,\cdot}(-_\R,G))[1].
\end{array}
\]
The commutativity of the left-hand  square is straightforward.

Finally, together with 
(\ref{mordercat}) and the relation 
$\pp\circ w= {\rm Id}_{A^{0,\cdot}(-_\R,F)}$, the commutativity of the last 
diagram shows the equality  of morphisms in $D^+({\mathcal O}_X-{\rm mod})$:
\[
\Psi'\circ \partial_{(\mathcal{E},s)} =
-A \circ w\circ {\rm ad}_F.
\]
It is straightforward that this is precisely $-z(\alpha).$
\qed

\subsection{Admissible extensions}\label{admext}
Let $\overline{F}$ and $\overline{G}$ be hermitian coherent sheaves
over an arithmetic scheme $X$.

Given a hermitian coherent sheaf $\ol{E}$ on $X$ and an extension
\[
{\mathcal E}:\,0\longrightarrow G\longrightarrow
E\longrightarrow  F\longrightarrow  0,
\]
we write
\[
\overline{\mathcal E}:\,0\longrightarrow  \overline{G}\longrightarrow
\overline{E}\longrightarrow  \overline{F}\longrightarrow  0
\]
and call
$\ol{\mathcal{E}}$ an {\it admissible} extension of  hermitian
coherent sheaves
if $\overline{F}$ and $\overline{G}$ carry the induced hermitian
metrics
from $\overline{ E}$.
In this case, orthogonal projection on the orthogonal complement of 
$G_\C$ in $E_\C$ determines an $F_\infty$-invariant
$\mathcal{C}^\infty$-splitting $s^\perp:F_{\C}\to E_\C$
of the extension $\mathcal{E}_\C$ on $X_\Sigma(\C)$, which we shall 
call the \emph{orthogonal splitting} of $\cE$. Amongst the 
 splittings 
of $\mathcal{E}_\C$ over $X_\Sigma(\C)$, it is characterized as being 
fiberwise isometric.
In this way, each admissible extension $\ol{\mathcal{E}}$ determines an
arithmetic
extension $(\mathcal{E},s^\perp)$ of $F$ by $G$.

Conversely, if 
\[
{\mathcal E}:\,0\longrightarrow  G\stackrel{i}{\longrightarrow}
E\stackrel{p}{\longrightarrow}F\longrightarrow  0
\] 
is an extension of ${\mathcal O}_X$-modules, then any 
$F_\infty$-invariant
$\mathcal{C}^\infty$-splitting $s:F_{\C}\to E_\C$
of the extension $\mathcal{E}_\C$ determines an hermitian structure 
on $E.$ Namely, the map
$$\varphi:=(i_\C,s):G_\C \oplus F_\C \longmapsto E_\C$$
is a $F_\infty$-invariant
$\mathcal{C}^\infty$-isomorphism of vector bundles over 
$X_\Sigma(\C)$, and the hermitian metric on $G_\C \oplus F_\C$ defined as 
the orthogonal direct sum of the ones defining $\ol{G}$ and $\ol{F}$ 
may be transported by $\varphi$ to a $F_\infty$-invariant
$\mathcal{C}^\infty$-metric $\|.\|$ on $E$. If we let $\ol{E}:=(E, 
\|.\|),$ then
\[
\overline{\mathcal E}:\,0\longrightarrow  \overline{G}\stackrel{i}{\longrightarrow}
\overline{E}\stackrel{p}{\longrightarrow}\overline{F}\longrightarrow  0
\] 
is an admissible extension, the orthogonal splitting of which is 
precisely $s.$

An isomorphism of admissible extensions $\ol{\mathcal{E}}_1$ and
$\ol{\mathcal{E}}_2$ of $\ol{F}$ by $\ol{G}$ is,
by definition, an isomorphism from $\ol{E}_1$ to $\ol{E}_2$ which induces the
identity
on $F$ and $G$.
The constructions above induces a one to one correspondence
\[
\genfrac{\{}{\}}{0pt}{}{\mbox{isomorphism classes of admissible}}
{\mbox{ extensions of $\ol{F}$ by $\ol{G}$}}
\longleftrightarrow
\genfrac{\{}{\}}{0pt}{}{\mbox{isomorphim classes of arithmetic}}
{\mbox{ extensions of $F$ by $G$}},
\]
and we obtain a new interpretation of the group $\widehat{\rm
Ext}_X^1(F,G)$ as
group of isomorphism classes of admissible extensions of $\ol{F}$ by
$\ol{G}$.
The corresponding ``admissible Baer sum" of two admissible extensions
\[
\ol{\mathcal{E}}: 0\longrightarrow  \ol{G}
\stackrel{i_j}{\longrightarrow }
\ol{E}_j\stackrel{p_j}{\longrightarrow }\ol{F}
\longrightarrow  0\,\,,\,(j=1,2)
\]
has the following explicit description.
We define a hermitian metric on
\[
\ker(p_1-p_2:E_1\oplus E_2\longrightarrow  F)\subseteq E_1\oplus E_2
\]
by the formula
\begin{equation}\label{drole de somme}
\|(e_1,e_2)\|_\sigma^2
=2(\|e_1\|^2_{\ol{E}_1,\sigma}+\|e_2\|^2_{\ol{E}_2,\sigma})
-3\|p_1(e_1)\|^2_{\ol{F},\sigma}
\end{equation}
for $\sigma$ in  $\Sigma$ and $(e_1,e_2)$ a vector in a fiber of
 $E_{1,\sigma}\oplus E_{2,\sigma}$ such that $p_1(e_1)=p_2(e_2).$
We equip the algebraic Baer sum $E$ in (\ref{Baersum})($=$\refg{algbaersum2}) with the
quotient metric.
It is straightforward to check that the resulting extension
\[
\ol{\mathcal{E}_j}: 0\longrightarrow  \ol{G}\longrightarrow
\ol{E} \longrightarrow  \ol{F}\longrightarrow  0
\]
is admissible.
It corresponds to the arithmetic Baer sum.
Indeed we have
\[
\|[(s_1^\perp(f),s_2^\perp(f))]\|_{\ol{E},\sigma}
=\|f\|_{\ol{F},\sigma}
\]
for any vector $f$ in a fiber of $F_\sigma$, if $s_j^\perp$ denotes the orthogonal splitting of
$\ol{\mathcal{E}}_{j}$.

We leave the details to the interested reader, and just want to 
emphasize that this correspondence would \emph{not} hold with $\ol{E}$ defined as $E$ 
equipped with the hermitian structure it naively inherits as a 
subquotient of the orthogonal direct sum $\ol{E}_1 \oplus \ol{E}_2.$

\subsection{Arithmetic torsors}\label{arithmetic torsors}

Let $X$ be an arithmetic scheme over some arithmetic ring $R$.
Let $G$ be a connected flat linear group scheme over $S=$ Spec $R$ of
finite type.
For each $\sigma$ in $\Sigma$, we fix a maximal compact subgroup
$K_\sigma$ of $G_\sigma(\C)$ such that the family
$K_\infty=(K_\sigma)_{\sigma\in\Sigma}$ is invariant under $F_\infty$.
Chambert-Loir and Tschinkel define an arithmetic
$(G,K_\infty)$-torsor on $X$
as a pair $(\mathcal{T},s)$ given by a $G$-torsor $\mathcal{T}$ on $X$
together with a family of sections $s=(s_\sigma)_{\sigma\in \Sigma}$
where $s_\sigma$ is a $\mathcal{C}^\infty$-section of the
$K_\sigma\backslash G_\sigma(\C)$-fibration on $X_\sigma(\C)$ obtained
as the quotient of $\mathcal{T}_\sigma(\C)$ by the action of
$K_\sigma$.
The set of isomorphism classes of arithmetic  $(G,K_\infty)$-torsors
on
$X$ is denoted by $\widehat H^1\bigl(X,(E,K_\infty)\bigr)$
\cite[1.1]{chamberttschinkel01}.

Let $V$ be a vector bundle on $S$.
Let $(\mathcal{E}, s)$ be an arithmetic extension over $X$ with
underlying extension
\[
\mathcal{E}\colon\, 0\longrightarrow  f^*V\longrightarrow
E\stackrel{\pi}{\longrightarrow } {\mathcal O}_X\longrightarrow  0.
\]
where $f:X\to S$.
We consider $E$, $f^*V$ and ${\mathcal O}_X=\mathbb{G}_a$ as vector group
schemes and
denote by $1\colon\, X \to \mathbb{G}_a$ the section of $\mathbb{G}_a$
associated to the unit in ${\mathcal O}_X$.
The scheme $\mathcal{T}=\pi^{-1}(1)=E ×_{\mathbb{G}_a} X$ is the
$V$-torsor defined by the splittings of the extension $\mathcal{E}$.
The $\mathcal{C}^\infty$ splitting $s$ of $\mathcal{E}_\C$
induces a $\mathcal{C}^\infty$ section of the $V_\mathbb{C}$-torsor
$\mathcal{T}_\mathbb{C}$ on $X_\mathbb{C}$.
Let $K_\sigma = \{0\}$ be the maximal compact subgroup of
$G_\sigma(\C)$.
The pair $(\mathcal{T},s)$ is an arithmetic
$\bigl(G,K_\infty=(K_\sigma)_\sigma\bigr)$-torsor and
determines an element in $\widehat H^1\bigl(X,(G,K_\infty)\bigr)$.
One checks easily that the construction which associates $\bigl(\mathcal{T}=\pi^{-1}(1), s\bigr)$ to $(\mathcal{E},s)$ induces an isomorphism of groups:
\begin{equation} \widehat{\rm Ext}_X^1(\mathcal{O}_X,f^*V)
\longrightarrow
\widehat H^1(X,(G,K_\infty)).
\end{equation}

One recovers the first exact sequence \ref{longexseq1} when $F= {\mathcal O}_X$ and $G= f^\ast V$ as a special
case of
the exact sequence \cite[1.2.1 and 1.2.3]{chamberttschinkel01}, and
the adelic
description \ref{adelic} as a special case of 
\cite[1.2.6]{chamberttschinkel01}.

\section{Slopes of hermitian vector bundles and splittings of
extensions
over arithmetic curves}
Consider an extension
\[
\mathcal{E}\colon \, 0 \longrightarrow  G \longrightarrow
E \longrightarrow  F \longrightarrow  0
\]
of a line bundle $F$ by a line
bundle $G$ on a smooth projective geometrically connected curve $C$ of
genus $g$ over a field $k$. A straightforward application of Serre
duality
shows that
\emph{if $\mathcal{E}$ does not split, the following  inequality
holds:
}
\begin{equation} \deg G -\deg F \leq 2g-2. \label{extlines}
\end{equation}
Indeed, in this situation, the class of $\mathcal{E}$ provides a
non-zero element in
\[
\text{Ext}^1_C(F,G) \simeq H^1 (C,\,F^\lor \otimes G) \simeq
H^0 (C,\,\Omega_{C/k}^1 \otimes F \otimes G^\lor)^\lor.
\]
Hence $\Omega_{C/k}^1 \otimes F \otimes G^\lor$ has a non-trivial
regular section on $C$, and its degree is consequently non-negative;
this yields (\ref{extlines}).

In this section, our aim is primarily to establish an arithmetic
analog of the
inequality (\ref{extlines}), valid for an 
admissible extension $\overline{\mathcal
E}$  of hermitian vector bundles $\overline{F}$ and $\overline{G}$
of arbitrary ranks  over an ``arithmetic curve" $S:= \Spec {\mathcal O}_K$
defined by some number field $K$.

In order to formulate this analog, we need to introduce some
quantitative measure for the non-triviality of
an arithmetic
extension $\overline{\mathcal{E}}$: it will be given by its
\emph{size}
$\size(\overline{\mathcal{E}})$, defined as
the logarithm of the distance to zero of the 
corresponding point on the real torus
\[
\widehat{\rm 
Ext}^1_S(F,G)=\frac{\Hom_{{\mathcal O}_K}{(F,G)}\otimes_\Z\R}{\Hom_{{\mathcal O}_K}{(F,G)}}.
\]
Indeed, the hermitian structures of $\overline{F}$ and
$\overline{G}$ induce an euclidean norm on the real vector space
$\Hom_{{\mathcal O}_K}(F,G)\otimes_\Z\R$,
which may be seen as a flat
Riemannian metric on the torus $\widehat{\rm Ext}^1_S(F,G)$ and
defines a
distance on it (see \ref{subsec:size} \emph{infra} for the formal
definition of $\size(\overline{\mathcal{E}})$).

Our arithmetic analog of the  inequality (\ref{extlines}) will take
the form:
\begin{equation}
     \label{sizeineq}
\widehat{\mu}_{\rm min}(\ol{G})
- \widehat{\mu}_{\rm max}(\ol{F})
+\size(\overline{\mathcal{E}})\leq \frac{\log |\Delta_{K}|}{[K:\Q]}
+ \log\,\frac{\rk_K F_K\cdot\rk_K G_K}{2}.
\end{equation}
Here, $\widehat{\mu}_{\rm max}(\ol{F})$ (resp. $\widehat{\mu}_{\rm
min}(\ol{G})$) denotes the maximal (resp. minimal) normalized slope
of $\ol{F}$ (resp. of $\ol{G}$)
(\cf \cite{stuhler76}, \cite{grayson84}, and \ref{Adeg},
\emph{infra}), and  $\Delta_{K}$ the discriminant of the number field
$K$.
Our proof of (\ref{sizeineq}) will rely on (i) some upper bound on
the Arakelov degree of a sub-line bundle in
the tensor product of two hermitian vector bundles over $\Spec {\mathcal O}_K$
(Proposition \ref{prop.udegtens}), and (ii) on some ``transference"
result from the
geometry of numbers, relating minima of an euclidean lattice and of
its dual lattice, in the precise form
obtained by Banaszczyk (\cite{banaszczyk93}).

\subsection{Arithmetic degree and slopes}\label{Adeg}
We now discuss a few results concerning hermitian vector bundles on
``arithmetic curves" and their arithmetic degree and slopes. We refer
the reader to
\cite{szpiro85},\cite{lang88},\cite{stuhler76},\cite{grayson84}, and
\cite{neukirch99} for more extensive discussions of these notions.

Let $K$ be a number field, and $\Sigma$ the set of its fields
embeddings $\sigma: K \hookrightarrow \C.$
Let $\overline L=(L,(\|.\|_\sigma)_{\sigma \in \Sigma})$ be a
hermitian line bundle on the arithmetic curve
$S=\Spec {\mathcal O}_K$.
The expression
\[
\widehat{\rm deg}\ \overline{L} = \log \,\#\, \bigl(L/{\mathcal O}_K\cdot
s\bigr) - \sum\limits_{\sigma\in \Sigma}\log\,\|s\|_\sigma
\]
does not depend on the choice of a non-zero section $s\in\Gamma(S,L)$
and defines the arithmetic degree of  $\overline L$.

For an arbitrary  hermitian coherent sheaf $\overline E = (E,
(\|.\|_\sigma)_{\sigma \in \Sigma})$  over $S$
with ${\mathcal O}_K$-torsion subsheaf $E_{\rm tors}$, the quotient $E/E_{\rm
   tors}$ equipped with $(\|.\|_\sigma)_{\sigma \in \Sigma}$ defines
a hermitian
   vector bundle $\overline{ E/E_{\rm tors}}$ on $S$.
We define the {\it arithmetic degree} of the hermitian coherent sheaf
$\ol{E}$ by the formula \cite[2.4.1]{gilletsoule91}
\[
\widehat{\rm deg}\, \overline{E}
= \widehat{\rm deg}\,(\land^{\rm  max} \overline{ E/E_{\rm tors}})
+\log\,\#\,E_{\rm tors}.
\]
This is also the arithmetic degree $\dega {\rm det} \overline{E}$ of
the
determinant line of $\overline{E}$.

For every extension of hermitian coherent sheaves on $S$
\[
0\longrightarrow  \ol{G}\longrightarrow  \ol{E}\longrightarrow
\ol{F}\longrightarrow  0
\]
which is admissible in the sense of Section \ref{admext}, we have
\begin{equation}\label{formel}
\widehat{\rm deg}\,\ol{E}=\widehat{\rm deg}\,\ol{G}+\widehat{\rm
deg}\,\ol{F}
\end{equation}
(see for instance \cite[V 2.1]{lang88}).

It is also convenient to
introduce the \emph{normalized arithmetic degree} of an hermitian
coherent sheaf $\overline E$ over $S$, namely
\[\degan \ol{E}:= \frac{1}{[K:\Q]}\, \dega \ol{E}.\]
Then the (normalized)  {\it arithmetic slope} of an hermitian
coherent sheaf
$\ol{E}$ of positive rank ${\rm rk}\, E:={\rm rk}_K\, E_K$ is defined as
\begin{equation}
\label{eq:arslope}
\widehat{\mu}(\ol{E})= \frac{\degan(\ol{E})}{{\rm rk}\, E}
=\frac{1}{[K:\Q]}\cdot\frac{\widehat{\rm deg}(\ol{E})}{{\rm rk}\, E}.
\end{equation}

We define the {\it maximal slope}  $\widehat{\mu}_{\rm max}(\ol{E})$
of an hermitian vector bundle $\ol{E}$
of positive  rank
as the maximal arithmetic slope of a subbundle of positive rank in $E$
equipped with the induced metric,
  and its {\it minimal slope} as
  \begin{equation}
      \widehat{\mu}_{\rm
   min}(\ol{E}):=-\widehat{\mu}_{\rm max}(\ol{E}^\lor).
      \label{eq:maxminbis}
  \end{equation}
   Observe that $\muamax(\ol{E})$ (resp. $\muamin(\ol{E})$) coincides
with the maximal (resp. minimal) slope of a saturated subbundle (resp.
of a quotient coherent sheaf) of $E$, equipped with the hermitian
structure induced by restriction (resp. by quotient) from the one of
$\ol{E}$.

It is easily seen that, for any hermitian line bundle $\ol{L}$ on $S$,
we have
\[
\widehat{\mu}_{\rm max}(\ol{E}\otimes\ol{L})=\widehat{\mu}_{\rm
   max}(\ol{E})+\degan\ol{L}
\]
and
\[
\widehat{\mu}_{\rm min}(\ol{E}\otimes\ol{L})=\widehat{\mu}_{\rm
   min}(\ol{E})+\degan\ol{L},
\]
and that, if $\ol{F}$ is any hermitian vector bundle of positive rank
on $S$,
$$\mua (\ol{E}\otimes\ol{F})= \mua (\ol{E}) +
\mua (\ol{F}),$$
\begin{equation}
     \muamax (\ol{E}\otimes\ol{F}) \geq  \muamax (\ol{E}) +
\muamax (\ol{F}),
     \label{eq:mumaxtens}
\end{equation}
and
\begin{equation}
     \muamin (\ol{E}\otimes\ol{F}) \leq  \muamin (\ol{E}) +
\muamin (\ol{F}).
     \label{eq:mumintens}
\end{equation}

Finally, recall that the normalized degree is invariant under field
extension. Namely, if $K'$ is a number field containing $K$ and if
$$g: S':=\Spec {\mathcal O}_{K'} \longrightarrow S:=\Spec {\mathcal O}_{K}$$ is the
morphism of
schemes defined by the inclusion ${\mathcal O}_{K} \hookrightarrow {\mathcal O}_{K'},$
then for any hermitian vector bundle $\ol{E}$ over $S$, the hermitian
vector bundle $g^\ast\ol{E}$ over $S'$ satisfies
$$\degan  g^\ast\ol{E} = \degan \ol{E}.$$
Consequently, when $E$ has positive rank:
$$\mua  (g^\ast\ol{E}) = \mua (\ol{E}).
     \label{eq:muainv}$$
A simple descent argument, using the Harder-Narasimhan filtrations
of $\ol{E}$ and $g^\ast\ol{E}$, shows that this also holds for the
maximal and minimal slopes:
\begin{equation}
    \muamax  (g^\ast\ol{E}) = \muamax (\ol{E}),
     \label{eq:muamaxinv}
\end{equation}
\begin{equation}
    \muamin  (g^\ast\ol{E}) = \muamin (\ol{E}).
     \label{eq:muamininv}
\end{equation}

\begin{lemma}\label{elelem}
Let $\ol{E}$ be a hermitian vector bundle over $S$.
Let $\ol{F}$ be an hermitian coherent subsheaf of $E$, equipped with
the  metrics of $\ol{E}$, such
that the quotient $Q$ of $E$ by $F$ is torsion. 
Let $I$ be the annihilator ideal
of $Q$ in ${\mathcal O}_K$, and let $N(I):=\# {\mathcal O}_K/I$ denote its norm.
Then
\begin{eqnarray*}
\widehat{\mu}_{\rm min}(\ol{E})&\leq&
\muamin(\ol{F})+\frac{1}{[K:\Q]\cdot {\rm rk} E} \log\, \#Q \\
    &\leq &\widehat{\mu}_{\rm
     min}(\ol{F})+\frac{1}{[K:\Q]}\log\,\# N(I).
     \end{eqnarray*}
\end{lemma}

\proof
The various coherent sheaves on $S$ we shall consider will be
equipped with the hermitian structures deduced from the
one of $\ol{E}$ by quotient.

Let $F'$ be a quotient vector bundle of $F$ such that
\[
\widehat{\mu}(\ol{F'})=  \widehat{\mu}_{\rm min}(\ol{F}),
\]
and let us form the pushout diagram:
\[
\begin{array}{ccccccccc}
0 &\longrightarrow  &F& \longrightarrow  &E&\longrightarrow  &Q&
\longrightarrow  & 0\\
  & & \downarrow & &\,\,\, \downarrow & & ||&\\
0 &\longrightarrow  &F'& \longrightarrow  &E'&\longrightarrow  &Q&
\longrightarrow  & 0.
\end{array}
\]
Observe that $E'$ is a quotient coherent sheaf of $E$.

As a special case of (\ref{formel}) we get:
\begin{equation}\label{torsion}
\dega\ol{E'}=\dega\ol{F'}+\log\#Q.
\end{equation}
Besides, since the ${\mathcal O}_K$-module $Q$ is killed by $I$, and is a
quotient of the locally free ${\mathcal O}_K$-module $E$, we have:
\begin{equation}\label{rktorsion}
\log\#Q\leq \rk E \log N(I).
\end{equation}
From (\ref{torsion}), we get the first of the desired
inequalities:
\[
\widehat{\mu}_{\rm min}(\ol{E})\leq\widehat{\mu}(\ol{E'})
=\mua(\ol{F'})+\frac{1}{[K:\Q]\cdot {\rm rk} E} \log\, \#Q.
\]
The second one follows from (\ref{rktorsion}).
\qed

\begin{lemma}\label{lem.mumaxadm}
     For any admissible extension
     $$0\longrightarrow  \ol{N}\longrightarrow \ol{E} \longrightarrow
     \ol{Q} \longrightarrow 0$$
     of hermitian vector bundles of positive ranks over $S,$ the
     following inequality holds:
     \begin{equation}
         \muamax(\ol{E}) \geq \frac{\rk N}{1+\rk N}\, \mua(\ol{N})
         + \frac{1}{1+\rk N}\, \muamax(\ol{Q}).
         \label{eq:mumaxadm}
     \end{equation}
\end{lemma}

\proof let $F$ be any ${\mathcal O}_{K}$-submodule of positive rank in $Q$,
and $F'$ its inverse image in $E.$ Consider the hermitian vector
bundles $\ol{F}$ and $\ol{F'}$ defined by $F$ and $F'$ equipped with
the restrictions of the
hermitian metrics of $\ol{Q}$ and $\ol{E}$ respectively. We have:
$$\rk F'= \rk N + \rk F$$
and
$$\degan \ol{F'}= \degan \ol{N} + \degan \ol{F}.$$
Therefore
\begin{eqnarray*}
     \mua (\ol{F}) & = & \frac{1}{\rk F}\, \bigl(\degan \ol{F'} -
\degan
     \ol{N}\bigr)  \\
      &\leq &\frac{1}{\rk F}\,\bigl[\rk F'\cdot\muamax(\ol{E})
                   -\rk N\cdot\mua(\ol{N})\bigr]  \\
      & = & \muamax(\ol{E}) + \frac{\rk N}{\rk
F}\cdot\bigl(\muamax(\ol{E})
                   - \mua(\ol{N})\bigr).
\end{eqnarray*}
Since $\rk F \geq 1$ and $\muamax(\ol{E}) - \mua(\ol{N}) \geq 0,$
this shows:
$$\muamax(\ol{Q})\leq  \muamax(\ol{E}) + \rk
N\cdot\bigl(\muamax(\ol{E}) -
\mua(\ol{N})\bigr).$$
This inequality is equivalent to (\ref{eq:mumaxadm}).
\qed

\subsection{Euclidean lattices and direct images}\label{subsec.euclid}
By an \emph{euclidean lattice}, we shall mean a pair $(\Gamma, \|.\|)$
where $\Gamma$ is a free $\Z$-module of finite rank and $\|.\|$ a
euclidean norm on $\Gamma_{\R}:=\Gamma\otimes_{\Z}\R.$ such a norm
uniquely extends to an hermitian norm, invariant under complex
conjugation, on $\Gamma_{\C}:=\Gamma\otimes_{\Z}\C\simeq
\Gamma_{\R}\otimes_{\R}\C.$ Accordingly, euclidean lattices may be
identified with hermitian vector bundles over $\Spec \Z,$ and the
classical
invariants of the former, such as their successive minima, may be
interpreted as invariant of the latter.

For instance, if $\ol{V}:=(V, \|.\|)$ is any hermitian vector bundle
of positive rank over $\Spec \Z$, the first of its successive minima
is by definition:
$$\lambda_{1}(\ol{V}):=\min \left\{\|v\|, v \in V\setminus
\{0\}\right\}.$$
Recall also that, if $\mbox{covol}(\ol{V})$ denotes the covolume of
this euclidean lattice, then
$$\dega \ol{V}= -\log \mbox{covol}(\ol{V}).$$

For any positive integer $r$, let
$$v_{r}:=\frac{\pi^{r/2}}{\Gamma(\frac{r}{2}+1)}$$
be the volume of the unit euclidean ball in $\R^r$. Minkowski's First
Theorem on euclidean lattices may be reformulated as follows:

\begin{proposition}\label{prop.Minkowski}
     For any hermitian vector bundle $\ol{V}$ of positive rank $r$
over
     $\Spec \Z,$ we have
     \begin{equation}
         \log \lambda_{1}(\ol{V}) \leq - \mua(\ol{V}) -\frac{1}{r}
         \log v_{r} + \log 2.
         \label{eq:Minkowski}
     \end{equation}
\end{proposition}

Observe that
$$\psi(r):=-\frac{1}{r}\log v_{r} + \log 2 =
\frac{1}{r} \log \Gamma(\frac{r}{2}+1) -\frac{1}{2} \log
\frac{\pi}{4}$$
considered as a function of $r\in ]0, +\infty[,$
  is increasing (since $\log \Gamma$ is convex on $]0, +\infty[$ and
  vanishes at $1$).
Moreover, for any $x \in ]0, +\infty[$,  we have
$$\Gamma(x+1)=\sqrt{2\pi x} \left(\frac{x}{e}\right)^x e^{\theta(x)/12
x},$$
where $\theta(x)$ belongs to $]0,1[.$ This shows that, for positive
$r$,
$$\psi(r)=\frac{1}{2}\log r -\frac{1}{2}\log \frac{e\pi}{2} +
\varepsilon (r)$$
where
$$\varepsilon (r) := \frac{1}{2r}\log (\pi r) +\frac{1}{6r^2} \,
\theta
(r/2)$$
goes to $0$ when $r$ increases to infinity. Finally, since the unit
euclidean ball in $\R^r$ contains the ``cube" $[-1/\sqrt{r},
1/\sqrt{r}]^r$ for any positive integer $r$, we have
$$v_{r}\geq \left(\frac{2}{\sqrt{r}}\right)^r,$$
and consequently
$$\psi(r) \leq \frac{1}{2}\log r.$$
Moreover this inequality is strict if $r>1.$

If $K$ is a number field and $\ol{E}= (E,
(\|.\|_{\sigma})_{\sigma:K\hookrightarrow \C})$
an hermitian vector bundle over $\Spec {\mathcal O}_{K},$ we have defined its
direct image
$\pi_{\ast}\ol{E}$ --- where $\pi$ denotes the morphism
$\pi: \Spec {\mathcal O}_{K} \to \Spec \Z$ --- in \ref{defdim} as the euclidean
lattice $(\pi_{\ast}E, \|.\|),$ where $\pi_{\ast}E$ denotes $E$ seen
as a $\Z$-module, and $\|.\|$ the euclidean norm on
$E\otimes_{\Z}\R$ restriction of the hermitian scalar product
$\langle .,.
\rangle$
on $E\otimes_{\Z}\C\simeq \bigoplus_{\sigma:K\hookrightarrow \C}
E_{\sigma}$  defined by the direct sum of the hermitian scalar
products $\langle .,.
\rangle_{\sigma}$ on the $E_{\sigma}$'s attached to the norms
$\|.\|_{\sigma}.$ In other words, for any $v, w$ in $\pi_{\ast}E=E,$
$$\langle v,w
\rangle :=\sum_{\sigma  :K\hookrightarrow \C} \langle v,w
\rangle _{\sigma}.$$

If $\Delta_{K}$ denotes the discriminant of the number field $K,$
then the arithmetic degree of $\pi_{\ast}\ol{E}$ is given by:
$$\dega \pi_{*}\ol{E}= \dega \ol{E} - \frac{\rk_{{\mathcal O}_{K}} E}{2} \log
|\Delta_{K}|.$$
(see for instance \cite[III.7-8]{neukirch99} or
\cite[(2.1.13)]{bostetal94}).
Consequently, if $E$
has positive rank,
\begin{equation}
     \mua (\pi_{\ast}\ol{E})=\mua (\ol{E})-\frac{\log
|\Delta_{K}|}{2\, [K:\Q]}.
     \label{eq:mupi}
\end{equation}

Let $\omega_{{\mathcal O}_K}:=\Hom_\Z({\mathcal O}_K,\Z)$ denote the canonical module, or
inverse of the different, of the
number field $K$.
The formula
\begin{equation}\label{modul}
(af)(b)=f(ab)
\end{equation}
defines a ${\mathcal O}_K$-module structure on $\omega_{{\mathcal O}_K}$.
It is an invertible ${\mathcal O}_K$-module which is generated up to
torsion by the trace map ${\rm tr}_{K/\Q}: K \to \Q.$ It becomes
an hermitian line bundle
$\ol{\omega}_{{\mathcal O}_{K}}:=
(\omega_{{\mathcal O}_K},(\|.\|_{\sigma})_{\sigma:K\hookrightarrow \C})$ over
$\Spec {\mathcal O}_{K}$ if we equip it with the hermitian norms defined by
$\|{\rm tr}_{K/\Q}\|_{\sigma}=1.$
Then the arithmetic degree of $\ol{\omega}_{{\mathcal O}_K}$ is given by the
well-known formula
\begin{equation}
   \widehat{\rm deg}\,\ol{\omega}_{{\mathcal O}_K}=\log\,\#
\,\bigl(\omega_{{\mathcal O}_K}/{\mathcal O}_K\cdot{\rm tr}_{K/\Q}\bigr)
=\log|\Delta_K|.
     \label{eq:degom}
\end{equation}
Moreover, the hermitian line bundle $\ol{\omega}_{{\mathcal O}_K}$ satisfies the
following duality property:

\begin{proposition}\label{aradual}
     For any hermitian vector bundle $\ol{E}$ over $\Spec {\mathcal O}_{K}$,
     the $\Z$-linear map
     $$\begin{array}{crcl}
        I: & E^\lor \otimes_{{\mathcal O}_{K}} \omega_{{\mathcal O}_{K}} &
\longrightarrow &
         \Hom_{\Z}(E, \Z)  \\
         &\xi\otimes\lambda & \longmapsto & \lambda \circ \xi
     \end{array}
     $$
     defines an isometric isomorphism of hermitian vector bundles
     over $\Spec \Z$:
     \begin{equation}
     \pi_{\ast}(\ol{E}^\lor\otimes
     \ol{\omega}_{{\mathcal O}_K})\stackrel{\sim}{\longrightarrow}
     (\pi_{\ast}\ol{E})^\lor.
      \label{eq:dual}
     \end{equation}
     \end{proposition}

     \proof To check that $I$ is an isomorphism of $\Z$-modules, we
may
     work locally over $\Spec \Z.$ Then we may assume that $E$ is a
     trivial vector bundle, and the assertion is clear.

     Moreover, we have canonical isomorphisms:
     \begin{equation}
         \pi_{\ast}(E^\lor\otimes \omega_{{\mathcal O}_K})_{\C} \simeq
\bigoplus_{\sigma:K\hookrightarrow\C} (E_{\sigma}^\lor\otimes_{\C}
\omega_{{\mathcal O}_K,\sigma})
         \label{eq:dual1}
     \end{equation}
and
\begin{equation}
         (\pi_{\ast}E)^\lor_{\C} \simeq
\Hom_{\C}\left(\bigoplus_{\sigma:K\hookrightarrow\C} E_{\sigma},
\C\right) \simeq
\bigoplus_{\sigma:K\hookrightarrow\C} E_{\sigma}^\lor.
         \label{eq:dual2}
     \end{equation}
Since $\omega_{{\mathcal O}_{K}}\otimes_{{\mathcal O}_{K}}K= K\cdot{\rm tr}_{K/\Q},$ for
any
embedding $\sigma \hookrightarrow \C,$
the $\C$-vector space $\omega_{{\mathcal O}_{K},\sigma}$ is one dimensional
with basis the image ${\rm tr}_{K/\Q,\sigma}$ of ${\rm tr}_{K/\Q},$
and a vector in $E_{\sigma}^\lor\otimes_{\C}
\omega_{{\mathcal O}_K,\sigma}$ may be written uniquely as $\xi \otimes {\rm
tr}_{K/\Q,\sigma}$ for some $\xi$ in $E_{\sigma}^\lor.$

For any $\xi$ in $E^\lor,$ the image of
the element $\xi \otimes {\rm tr}_{K/\Q}$
of $\pi_{\ast}(E^\lor\otimes \omega_{{\mathcal O}_K})_{\C}$ by (\ref{eq:dual})
is $(\xi_{\sigma}\otimes {\rm tr}_{K/\Q,\sigma})_{\sigma},$ and the
image of
the element ${\rm tr}_{K/\Q}\circ \xi$ of $\pi_{\ast}E^\lor$ by
(\ref{eq:dual2}) is $(\xi_{\sigma})_{\sigma}.$ Indeed, the $\C$-vector
space $(\pi_{\ast}E)_{\C}:=E\otimes_{\Z}\C$ is generated by the image
of the inclusion $$E\hookrightarrow E\otimes_{\Z}\C
\simeq \bigoplus_{\sigma:K\hookrightarrow\C} E_{\sigma},$$
and
for any $v$ in
$E,$
\begin{eqnarray*}
     ({\rm tr}_{K/\Q}\circ \xi) (v) & = &  {\rm tr}_{K/\Q}(\xi ((v))
\\
      & = & \sum_{\sigma:K\hookrightarrow\C} \sigma (\xi(v))  \\
      & = & \sum_{\sigma:K\hookrightarrow\C} \xi_{\sigma}(v_{\sigma}).
\end{eqnarray*}

Since $\pi_{\ast}(E^\lor\otimes \omega_{{\mathcal O}_K})_{\C}$ is generated, as
a $\C$-vector space, by its ${\mathcal O}_{K}$-submodule $E^\lor\otimes {\rm
tr}_{K/\Q},$
this shows that, using the identifications (\ref{eq:dual1}) and
\ref{eq:dual2}), the $\C$-linear map
$$I_{\C}: \pi_{\ast}(E^\lor\otimes \omega_{{\mathcal O}_K})_{\C}  \longmapsto
(\pi_{\ast}E)^\lor_{\C}$$
maps an arbitrary element $(\xi(\sigma)\otimes {\rm
tr}_{K/\Q,\sigma})_{\sigma}$ of $\pi_{\ast}(E^\lor\otimes
\omega_{{\mathcal O}_K})_{\C}$
to $(\xi(\sigma))_{\sigma}.$ This makes clear that $I_{\C}$ is an
isometry
with respect to the hermitian structures defining
$\pi_{\ast}(\ol{E}^\lor\otimes \ol{\omega}_{{\mathcal O}_K})$
and $(\pi_{\ast}\ol{E})^\lor.$
\qed

\subsection{First minimum, upper degree, and maximum slope}
If $K$ is a number field and
$\ol{E}= (E, (\|.\|_{\ol{E}, \sigma})_{\sigma:K\hookrightarrow \C})$
an hermitian vector bundle of  positive rank over
$\Spec {\mathcal O}_{K}$, we may define its (normalized) \emph{upper arithmetic
degree} $\udega \ol{E}$ as the maximum of the normalized degree
$\degan \ol{L}$ of a
sub-line bundle $L$ of $E$ equipped with the hermitian structure
induced by the one of $\ol{E}$. Clearly, in this
definition, we may restrict to saturated sub-line bundles. Moreover,
if for any prime ideal ${\bf p} \neq (0)$ in ${\mathcal O}_{K},$ we denote
$\|.\|_{E, {\bf p}}$ the ${\bf p}$-adic norm on $E_{K_{\bf
p}}:=E\otimes_{{\mathcal O}_{K}}K_{\bf p}$ attached to its ${\mathcal O}_{\bf p}$-lattice
$E\otimes_{{\mathcal O}_{k}}{\mathcal O}_{\bf p},$ we have:
\begin{equation}
     \udega \ol{E}=-\min_{v\in E_{K}\setminus \{0\}}\frac{1}{[K:\Q]}
     \left( \sum_{{\bf p}\in \Spec {\mathcal O}_{K}\setminus \{(0)\}} \log
     \|v\|_{E,{\bf p}} +   \sum_{\sigma: K\hookrightarrow \C} \log
     \|v\|_{\ol{E},\sigma} \right).
     \label{eq:defudega}
\end{equation}
Observe also that, for any hermitian line bundle $\ol{L}$ over $\Spec
{\mathcal O}_{K}$,
\begin{equation}
     \udega \ol{E}\otimes\ol{L}= \udega \ol{E} +\degan\ol{L}.
     \label{eq:udegatens}
\end{equation}

In this paragraph, we show that, for any hermitian vector bundle
$\ol{E}$ over $\Spec {\mathcal O}_{K}$ as above, the differences between
$\widehat{\mu}_{\rm max}(\ol{E}),$
$\udega(\ol{E}),$ and $-\log \lambda_{1}(\pi_{\ast}\ol{E})$ may be
bounded in terms of $K$ and $\rk_{K}E_{K}$ only. We refer the reader to 
\cite{soule97}, 1.1-2, and \cite{borek05} for related results 
concerning the successive minima of hermitian vector bundles over $\Spec {\mathcal O}_{K}$.

\begin{proposition}\label{nelem}
Let $K$ be a number field and $\pi$ the morphism from $\Spec {\mathcal O}_{K}$
to $\Spec \Z$.
For any hermitian vector bundle $\ol{E}$ of positive rank $r$ over
$\Spec {\mathcal O}_{K},$ the following inequalities hold:
\begin{equation}
     \udega \ol{E} \leq \widehat{\mu}_{\rm max}(\ol{E}),
     \label{eq:triv}
\end{equation}

\begin{equation}
     -\log \lambda_{1}(\pi_{\ast}\ol{E}) \leq \udega \ol{E} -
     \frac{1}{2}\, \log [K:\Q],
     \label{eq:geoquadr}
\end{equation}
and
\begin{equation}
     \widehat{\mu}_{\rm max}(\ol{E}) \leq -\log
     \lambda_{1}(\pi_{\ast}\ol{E}) + \psi([K:\Q]\cdot \rk_{{\mathcal O}_{K}}E)
                              + \frac{\log |\Delta_{K}|}{2\, [K:\Q]}.
     \label{eq:maxmin}
\end{equation}
\end{proposition}

Observe that the compatibility of these three estimates for hermitian
vector bundles of rank one imply the following lower bound on
$|\Delta_{K}|,$ \emph{à la} Hermite-Minkowski:
\begin{equation}
     \frac{\log |\Delta_{K}|}{[K:\Q]} \geq \log [K:\Q] - 2
\psi([K:\Q]).
     \label{eq:HM}
\end{equation}
The right-hand side of (\ref{eq:HM}) is positive if
$[K:\Q]>1$ and has the positive limit $\log \frac{e\pi}{2}$ when
$[K:\Q]$ goes to infinity.

In the sequel, we shall use the following consequence of
(\ref{eq:maxmin}) and (\ref{eq:geoquadr}):
\begin{equation}
     \widehat{\mu}_{\rm max}(\ol{E}) \leq
     \udega \ol{E} -
     \frac{1}{2}\, \log [K:\Q] + \psi([K:\Q]\cdot \rk_{{\mathcal O}_{K}}E)
                                 + \frac{\log |\Delta_{K}|}{2\,
[K:\Q]}.
     \label{mumaxudeg1}
\end{equation}
Using the inequality
$$\log ([K:\Q]\cdot \rk_{{\mathcal O}_{K}}E) -2 \psi([K:\Q]\cdot \rk E) \geq
0,$$
we also obtain a slightly weaker version of (\ref{mumaxudeg1}):
\begin{equation}
     \widehat{\mu}_{\rm max}(\ol{E}) \leq
     \udega \ol{E} +
     \frac{1}{2}\, \log \rk E + \frac{\log |\Delta_{K}|}{2\, [K:\Q]}.
     \label{mumaxudeg2}
\end{equation}

\proof The estimate (\ref{eq:triv}) is a trivial consequence of the
definitions of $\udega \ol{E}$ and $\widehat{\mu}_{\rm max}(\ol{E})$.

To prove (\ref{eq:geoquadr}), consider a non-zero element $v$ of
$E (=\pi_{*}E)$ such that
$$\|v\|_{\pi_{*}\ol{E}}^2= \sum_{\sigma:K\hookrightarrow
\C}\|v\|_{\sigma}^2$$
is minimal, and the rank one sub-bundle $L:={\mathcal O}_{K}\cdot v$ of $E$ it
generates. We have:
\begin{equation}
     - \log \lambda_{1}(\pi_{*}\ol{E})=\frac{1}{2}\log
     \|v\|_{\pi_{*}\ol{E}}^2
     =-\frac{1}{2}\log \left(\frac{1}{[K:\Q]}
\sum_{\sigma:K\hookrightarrow
\C}\|v\|_{\sigma}^2\right) -\frac{1}{2} \log [K:\Q],
     \label{eq:quadr}
\end{equation}
and
\begin{equation}
    -\frac{1}{2[K:\Q]} \sum_{\sigma:K\hookrightarrow
\C}\log \|v\|_{\sigma}^2 = \degan \ol{L} \leq \udega \ol{E}.
     \label{eq:geo}
\end{equation}
The inequality (\ref{eq:geoquadr}) follows from (\ref{eq:quadr}),
(\ref{eq:geo}), and the convexity of the function $-\log$.

To prove (\ref{eq:maxmin}), consider an arbitrary sub-vector bundle
$F$ of positive rank in $E$. Clearly,
\begin{equation}
     \lambda_{1}(\pi_{\ast}\ol{E}) \leq \lambda_{1}(\pi_{\ast}\ol{F}).
     \label{eq:EF}
\end{equation}
Moreover, Minkowski's theorem (\ref{eq:Minkowski}) shows that
\begin{equation}
   \log \lambda_{1}(\pi_{\ast}\ol{F}) \leq - \mua 
(\pi_{\ast}\ol{F}) + \psi([K:\Q]\cdot 
\rk_{{\mathcal O}_{K}}F).
\end{equation}
Using (\ref{eq:mupi}) and the fact that $\psi$ is increasing, this
shows:
\begin{equation}
     \log \lambda_{1}(\pi_{\ast}\ol{F})
     \leq
     - \widehat{\mu}(\ol{F})
    + \frac{\log |\Delta_{K}|}{2\, [K:\Q]}  + \psi([K:\Q]\cdot
\rk_{{\mathcal O}_{K}}E).
     \label{eq:maxminter}
\end{equation}
The estimate (\ref{eq:maxmin}) follows from (\ref{eq:EF})
and (\ref{eq:maxminter}).
\qed

Actually, the proof above
establishes a stronger form of
(\ref{eq:geoquadr}). Namely, it shows that, if $v$ is an element of
$E\setminus\{0\}$ such that $\|v\|_{\pi_{*}\ol{E}}$ is minimal, then
\begin{equation}
   \begin{split}
     -\log \lambda_{1}(\pi_{\ast}\ol{E}) = -\log \|v\|_{\pi_{*}\ol{E}}
     &\leq \udega \ol{{\mathcal O}_{K}\cdot v} - \frac{1}{2}\, \log [K:\Q]\\
     &\leq \udega \ol{E} - \frac{1}{2}\, \log [K:\Q].
   \end{split}
     \label{eq:geoquadrbis}
\end{equation}

When $\ol{E}$ has rank one, (\ref{eq:maxmin}) may be written:
\begin{equation}
     \degan(\ol{E}) \leq -\log
     \|v\|_{\pi_{*}\ol{E}} + \psi([K:\Q]) + 
\frac{\log |\Delta_{K}|}{2\, [K:\Q]},
     \label{eq:maxminencore}
\end{equation}
and the conjunction of (\ref{eq:geoquadrbis}) and
(\ref{eq:maxminencore}) may be seen as a ``quantitative version" of
the
fundamental theorems of Dirichlet about the rings of
algebraic integers --- the finiteness of the ideal class group, and
the ``unit theorem".

{\small
Indeed, these theorems easily lead to (\ref{eq:maxminencore}), with
$$\psi([K:\Q]) + \frac{\log |\Delta_{K}|}{2\, [K:\Q]}$$
replaced by a constant depending on $K$ only: if $E$ is the trivial
line bundle over ${\mathcal O}_{K},$ this follows from the unit theorem; the
general case follows using the finiteness of the ideal class group.

Conversely, (\ref{eq:geoquadrbis}) and
(\ref{eq:maxminencore}) easily imply these two theorems (compare
\cite{szpiro85}). To show this, let us define
$$A(K):= - \frac{1}{2}\log [K:\Q] +\psi([K:\Q]) +
\frac{\log |\Delta_{K}|}{2\, [K:\Q]}$$
and
$$B(K)=[K:\Q]\cdot A(K).$$
(One might observe that $A(K) \leq (\log |\Delta_{K}|)/(2 [K:\Q])$ and
$B(K)\leq  (\log |\Delta_{K}|)/2.)$
If, for every $\bp \in \Spec {\mathcal O}_{K}\setminus\{(0)\},$ we denote
$n_{\bp}$
the $\bp$-adic valuation of $v$ as a section of $E$, we get
from (\ref{eq:geoquadrbis}) and
(\ref{eq:maxminbis}):
\begin{eqnarray*}
     \sum_{\bp}n_{\bp} \log N_{\bp}  =  \dega \ol{E} -\dega
     \ol{{\mathcal O}_{K}\cdot v}
      & \leq  & [K:\Q] ( \dega \ol{E} + \log \|v\|_{\pi_{*}\ol{E}}
                                       - \frac{1}{2}\, \log [K:\Q])
\\
      & \leq & B(K).
\end{eqnarray*}
In other words, the divisor
$${\rm div}\, v:= \sum_{\bp}n_{\bp} \bp$$ of $v$
belongs to the finite set of effective
divisors on $\Spec {\mathcal O}_{K}$:
$${\mathcal D}:=\left\{\sum_{\bp}n_{\bp} \bp
       \mid n_{\bp}\geq 0\wedge \sum_{\bp}n_{\bp} \log
N_{\bp} \leq B(K)\right\}.$$

This already shows that any element of the ideal class group
$Cl_{K}:= {\rm Pic}(\Spec{\mathcal O}_{K})$ of $K$ is the class of some divisor
in ${\mathcal D}$ --- in particular, $Cl_{K}$ is finite.

Moreover, if $\ol{L}$ is any hermitian line bundle over $\Spec {\mathcal O}_{K}$
such that $L={\mathcal O}_{K}$ and $\dega \ol{L} =0,$ then
we may consider the hermitian line bundles $\ol{L}^{\otimes n}$,
$n\in \N,$
and  choose non-zero elements $v_{n}$ in $L^{\otimes n} \simeq {\mathcal O}_{K}$
such that $\|v_{n}\|_{\pi_{*}\ol{L}^{\otimes n}}$ is minimal. Since
the divisors ${\rm div}\, v_{n}$ lie in the finite set ${\mathcal
D}$, there
exists an increasing sequence $(n_{i})$ in $\N$ such that the
divisors ${\rm div}\, v_{n_{i}}$ coincide. Then the elements
$u_{i}:=v_{n_{0}}^{-1}\cdot v_{n_{i}}$ of $L^{\otimes (n_{i}-n_{0})}
\simeq
{\mathcal O}_{K}$ are units of ${\mathcal O}_{K},$ and their norms
$\|u_{i}\|_{\pi_{*}\ol{L}^{\otimes (n_{i}-n_{0})}}$ stay bounded when
$i$ goes to infinity.
If we let $t_\sigma:= \log \|1\|_{\ol{L}, \sigma}^{-1}$ for any
embedding $\sigma : K \hookrightarrow \C,$ and $m_i:=n_i-n_0,$ this
shows
that the sequences $(\log |\sigma(u_i)|-m_i.t_\sigma)_i$ are bounded
from above. Since $\sum_{\sigma : K 
\hookrightarrow \C} \log |\sigma(u_i)| $ and 
$\sum_{\sigma: K
\hookrightarrow \C}t_\sigma$ vanish, these sequences are also bounded
from below, and consequently, for any embedding $\sigma : K
\hookrightarrow
\C,$
$$\lim_{i \rightarrow + \infty}\frac{1}{m_i} \log |\sigma(u_i)| =
t_\sigma.$$
As any family $(t_\sigma)_{\sigma:K \hookrightarrow \C}$ of $[K:\Q]$
real numbers satisfying $t_{\overline{\sigma}}=t_\sigma$ and
$\sum_\sigma t_\sigma =0$ arises in this construction, this implies
Dirichlet's unit theorem.
}

\subsection{The upper degree of a tensor
product}\label{subsec.udegtens}

\begin{proposition}\label{prop.udegtens}
     For any two hermitian vector bundles $\ol{F}$ and $\ol{G}$ of
positive
     rank over $\Spec {\mathcal O}_{K}$, we have:
     \begin{equation}
         \udega (\ol{F}\otimes\ol{G}^\lor) \leq \widehat{\mu}_{\rm
         max}(\ol{F})-\widehat{\mu}_{\rm min}(\ol{G})
         \label{eq:udegtens1}
     \end{equation}
     and
     \begin{equation}
         \udega (\ol{F}\otimes\ol{G}) \leq \widehat{\mu}_{\rm
         max}(\ol{F})+\widehat{\mu}_{\rm max}(\ol{G}).
         \label{eq:udegtens2}
     \end{equation}
\end{proposition}

Observe that, from this proposition and the comparison of the upper
degree and the maximal slope in Proposition \ref{nelem}, we derive
estimates on the maximum and minimum slopes of tensor products which
complement (\ref{eq:mumaxtens}) and (\ref{eq:mumintens}), namely:
\begin{equation}
     \muamax (\ol{E}\otimes\ol{F}) \leq  \muamax (\ol{E})+
            \muamax (\ol{F}) +\frac{1}{2}\log(\rk E\cdot \rk F)
         + \frac{\log
         |\Delta_{K}|}{2[K:\Q]},
     \label{eq:mumaxtensbis}
\end{equation}
and
\begin{equation}
     \muamin (\ol{E}\otimes\ol{F})
       \geq \muamin (\ol{E}) + \muamin (\ol{F})
           - \frac{1}{2}\log(\rk E\cdot \rk F)
            - \frac{\log |\Delta_{K}|}{2[K:\Q]}.
     \label{eq:mumintensbis}
\end{equation}
Actually, it is possible to
establish similar estimates where the term $$\frac{1}{2}\log(\rk
E\cdot \rk F)
         + \frac{\log
         |\Delta_{K}|}{2[K:\Q]}$$ is replaced by a constant depending
         on $\rk E$ and $\rk F$ only (see \cite{bost96a}, Propositions
         A.4 and A.5, and also \cite{graftieaux01}, appendix), which
         however is larger than this term when $K=\Q.$

\proof According to the expression (\ref{eq:maxmin}) of $\muamin$ in
term of $\muamax$, its is enough to establish (\ref{eq:udegtens1}).
The expression (\ref{eq:defudega}) for the upper degree shows that, to
achieve this, we need to prove that any non-zero element $\phi$ of
$(F\otimes G^\lor)_{K}\simeq \Hom_{K}(G_{K},F_{K})$ satisfies
\begin{equation}
     -\frac{1}{[K:\Q]}
     \left( \sum_{{\bf p}\in \Spec {\mathcal O}_{K}\setminus \{(0)\}} \log
     \|\phi\|_{F\otimes G^\lor,{\bf p}} + \sum_{\sigma:
K\hookrightarrow \C}\log
     \|\phi\|_{\ol{F}\otimes \ol{G}^\lor,\sigma} \right) \leq
      \muamax(\ol{F}) - \muamin(\ol{G}) .
     \label{eq:sumphi}
     \end{equation}

     For any such $\phi,$ let $I\subset F$ (resp. $J\subset G$) be the
     saturated subbundle of $F$ (resp. of $G$) such that $I_{K}=\Im
     (\phi)$ (resp. $J_{K}=\Ker (\phi)$), and consider the canonical
     factorization of $\phi$:
     $$\begin{array}{ccc}
     G_{K} & \stackrel{\phi}{\longrightarrow} & F_{K}  \\
      \downarrow &  & \uparrow  \\
      (G/J)_{K} & \stackrel{\phit}{\longrightarrow} & I_{K}.
  \end{array}$$
By construction, $\phit$ is an isomorphism of $K$-vector spaces. We
shall denote $r$ the rank of $\phi$ and $\phit$; it is also the rank
of $I$ and of $G/J$.

The determinant $\det \phit$ of $\phit$ may be seen as a non-zero
meromorphic section of $(\Lambda^rG/J)^\lor\otimes\Lambda^rI$ over
$\Spec {\mathcal O}_{K}$. As such, it has a well defined valuation $v_{\bf
p}(\det\phit)$ at every ${\bf p}$ in $\Spec {\mathcal O}_{K}\setminus
\{(0)\},$ which vanishes for almost every $\bf p$. We may equip $G/J$
and $I$ with the hermitian structures defined by quotient and
restriction from the ones of $\ol{G}$ and $\ol{F}$. Then we have:
\begin{equation}
     \begin{split}
     -\dega \ol{G/J} + \dega \ol{I} &= \dega
     \Lambda^r\ol{G/J}^\lor\otimes\Lambda^r\ol{I}\\
     &= \sum_{{\bf p}}v_{\bf p}(\det \phit) \log N({\bf p})
     -\sum_{\sigma} \log \|\det
     \phit\|_{\Lambda^r\ol{G/J}^\lor\otimes\Lambda^r\ol{I},\sigma}.
     \end{split}
     \label{eq:degdet}
\end{equation}
Observe that, if $v_{\bf p}(\phi)$ denotes the valuation at $\bf p$
of $\phi$ seen as a meromorphic section of $F\otimes G^\lor$ over
$\Spec {\mathcal O}_{K}$, we have:
$$v_{\bf p}(\det \phit) \geq r v_{\bf p}(\phi),$$
or, equivalently:
\begin{equation}
     -\log \|\det \phit \|_{\Lambda^rG/J^\lor\otimes\Lambda^rI,\bf p}
     \geq -r\log \| \phi \|_{F\otimes G^\lor, \bf p}.
     \label{eq:divelem}
\end{equation}
This is a straightforward consequence of the theory of elementary
divisors
applied to $\phi_{K_{\bf p}}: G_{K_{\bf p}} \to F_{K_{\bf
p}}.$

Similarly, for any embedding $\sigma: K \hookrightarrow \C,$  by
considering
the polar decomposition of $\phi_{\sigma}$ we see that
$$\|\det \phit
\|_{\Lambda^r\ol{G/J}^\lor\otimes\Lambda^r\ol{I},\sigma}
\leq \| \phi \|_{\ol{F}\otimes\ol{G}^\lor, \sigma}^r,$$
or, equivalently,
\begin{equation}
     - \log \|\det \phit
     \|_{\Lambda^r\ol{G/J}^\lor\otimes\Lambda^r\ol{I},\sigma} \geq
     - r \log \| \phi \|_{\ol{F}\otimes\ol{G}^\lor, \sigma}.
     \label{eq:polar}
\end{equation}

From, (\ref{eq:degdet}), (\ref{eq:divelem}) and (\ref{eq:polar}), we
deduce:
$$- \sum_{{\bf p}\in \Spec {\mathcal O}_{K}\setminus \{(0)\}} \log
     \|\phi\|_{F\otimes G^\lor,{\bf p}}-\sum_{\sigma:
K\hookrightarrow \C}\log
     \|\phi\|_{\ol{F}\otimes \ol{G}^\lor,\sigma} \leq -
     \mua(\ol{G/J})+\mua(\ol{I}).$$
Since $\mua(\ol{G/J})\geq \muamin(\ol{G})$ and $\mua(\ol{I})\leq
\muamax(\ol{F}),$ this proves (\ref{eq:sumphi}).

   \qed

\subsection{The size of an arithmetic extension}\label{subsec:size}
\subsubsection{Definition of the size}\label{subsubsec.sizede}
Consider a number field $K$ with ring of integers $\mathcal{O}_K$ and
denote the spectrum of ${\mathcal O}_K$ by $S$.
Let $F$ and $G$ two vector bundles on the arithmetic curve $S$, defined
by
projective ${\mathcal O}_K$-modules we shall also $F$ and $G$.
The morphism
\[
b:\pi_{\ast}(F^\lor\otimes G)_{\R}=\Hom_{{\mathcal O}_K}(F,G)\otimes_\Z\R
\longrightarrow  \Exthat_{S}^1\,(F,G)
\]
induces an isomorphism of abelian groups
(compare Corollary
\ref{cor.affine}):
\begin{equation}
\Exthat_{S}^1\,(F,G)\simeq\frac{\Hom_{{\mathcal O}_K}(F,G)\otimes_\Z \R}{\Hom_{{\mathcal O}_K}(F,G)}
     \label{eq:extcurve}
\end{equation}
Moreover, as an abelian group, the right-hand side of
(\ref{eq:extcurve})
may be identified with the real torus
$$\frac{\pi_{\ast}(F^\lor\otimes G)_{\R}}{\pi_{\ast}(F^\lor\otimes
G)}.$$

From the hermitian structures on $\ol{F}$ and $\ol{G}$, we deduce an
hermitian structure on $\ol{F}^\lor \otimes \ol{G}$, that is, for
every field embedding $\sigma:K\hookrightarrow\C$, an hermitian
structure on the $\C$-vector space
$$(F^\lor\otimes G)_{\sigma}\simeq \Hom_{\C}(F_{\sigma},G_{\sigma}).$$
It is given by the ``Hilbert-Schmidt" hermitian scalar product
$\langle .,. \rangle _{\sigma}$, defined by
$$\langle T_{1},T_{2}
\rangle_{\sigma}:= {\rm tr} (T_{2}^{\ast}T_{1})$$
for any $T_{1}, T_{2} \in \Hom_{\C}(F_{\sigma},G_{\sigma}),$
where the adjoint $T_{2}^{\ast}$ is taken with respect to the
hermitian norms $\|.\|_{\ol{F}, \sigma}$ and $\|.\|_{\ol{G}, \sigma}$
on $F_{\sigma}$ and $G_{\sigma}$.

We may use these metrics to define the size of an arithmetic
extension $(\mathcal{E},s)$ over $S$ with underlying
algebraic extension
\[
{\mathcal E}:0\longrightarrow  {G} \longrightarrow  {E}
\longrightarrow
{F} \longrightarrow  0.
\]
Namely, we define the \emph{size of} (the class ${\rm cl}\,({\mathcal
E},s)$
in $\Exthat_{S}^1\,(F,G)$ of) \emph{the arithmetic extension}
$({\mathcal E},s)$ as
\begin{equation}\label{eq:defsize}
\size({\mathcal E},s):=  \log\,{\rm inf}\,
\Biggl\{
   \sqrt{\frac{1}{[K:\Q]}\sum\limits_{\sigma:K\hookrightarrow \C}
   \| h_{\sigma}\|_{\ol{F}^\lor \otimes \ol{G},\sigma}^2}\,
   \Bigg|
   \genfrac{}{}{0pt}{}{(h_{\sigma})\in
\bigl[\bigoplus_{\sigma} \Hom_{\C}(F_\sigma,G_\sigma)\bigr]^{F_\infty}
}{\wedge\,\,b\bigl((h_\sigma)\bigr) ={\rm cl}({\mathcal E},s)}
\Biggr\}.
\end{equation}
It takes values in $\R\cup\{-\infty\}$ and the equality
$\size({\mathcal E},s)=-\infty$ holds if and only
${\rm  cl}\,({\mathcal E},s)$ vanishes.
Clearly, the size of $({\mathcal E},s)$ is the logarithm of the
distance from $0$ to ${\rm cl}\,({\mathcal E},s)$ in the real
torus
$$\Exthat_{S}^1\,(F,G) = \frac{\pi_{\ast}(F^\lor\otimes
G)_{\R}}{\pi_{\ast}(F^\lor\otimes G)}$$
equipped with the translation invariant Riemannian metric defined by
the euclidean norm $[K:\Q]^{-1/2} \|.\|_{\pi_{\ast}(\ol{F}^\lor
\otimes
\ol{G})}$ on $\pi_{\ast}(F^\lor\otimes
G)_{\R}$. Moreover, the infimum in the right-hand side of
(\ref{eq:defsize}) is actually a minimum.

Let us emphasize that, to define the size $\size({\mathcal E},s)$, we
need to have chosen hermitian structures on $F$ and $G$.
It will  be
therefore
more appropriate to call (\ref{eq:defsize}) the \emph{size of
$({\mathcal
E},s)$ with respect to $\ol{F}$ and $\ol{G}$} and to denote it
$$\size_{\ol{F},\ol{G}}({\mathcal E},s)$$
when some ambiguity may arise.

\subsubsection{Size of admissible
extensions}\label{subsubsec.sizeadmissible}
We define the \emph{size $\size({\ol{\mathcal E}})$ of an admissible
extension
}
\begin{equation}
     {\ol{\mathcal E}}:0\longrightarrow  \ol{G} \longrightarrow
\ol{E}
\longrightarrow \ol{F} \longrightarrow  0
     \label{eq:admissibleencore}
\end{equation}
of hermitian vector bundles over $S$ as the size with respect
to $\ol{F}$ and $\ol{G}$ of its class in $\Exthat_{S}^1\,(F,G)$.
If
\[
s^\perp \in \Hom_{{\mathcal O}_K}(F,E)\otimes_\Z\R
\simeq
[ \bigoplus_{\sigma: K \hookrightarrow \C} \Hom(F_\sigma, 
E_\sigma)]^{F_\infty}
\]
denotes the orthogonal splitting of (\ref{eq:admissibleencore}), we
have:
\[
\size(\ol{\mathcal E})=\log {\min} \Bigl\{ [K:\Q]^{-1/2}\cdot
\|s-s^\perp \|_{\pi_{\ast}(\ol{F}^\lor \otimes \ol{G})} \,\Big|
\,\mbox{$s:F\to E$ an ${\mathcal O}_K$-linear splitting of $\mathcal{E}$}
\Bigr\}.
\]

Let us consider the set $\Triv_\OK (\cE)$ of 
trivializations
of $\cE$ over $\Spec \OK$, namely of $\OK$-modules isomorphisms
$$\varphi: E \stackrel{\sim}{\longrightarrow} G \oplus F$$
such that ${\rm pr}_2 \circ \phi =p$ and $\phi \circ i = (Id_G, 
0_F).$ The map which sends $\varphi \in \Triv_\OK (\cE)$ to 
$\varphi^{-1}\circ (0, Id_F) \in \Hom_\OK(F,G)$ defines a bijection
from 
$\Triv_\OK (\cE)$ onto the set of splittings of the extension $\cE$ 
of $\OK$-modules. In particular, $\Triv_\OK (\cE)$ is non-empty, and
becomes a torsor under the abelian group 
$\Hom_\OK (F,G)$ thanks to the action defined by letting, for any 
$\psi$ in $\Hom_\OK (F,G)$ and any $\varphi$ in $\Triv_\OK (\cE)$:
$$\psi + \varphi := i \circ \psi \circ p + \varphi.$$

For any $\varphi$ in $\Triv_\OK (\cE)$, the difference 
$$s^\perp - ((\varphi^{-1}\circ (0, Id_F))_\sigma) \in
[ \bigoplus_{\sigma: K \hookrightarrow \C} \Hom(F_\sigma, 
E_\sigma)]^{F_\infty}$$
factorizes through the morphism $(i_\sigma)
\in
[ \bigoplus_{\sigma: K \hookrightarrow \C} \Hom(G_\sigma, 
E_\sigma)]^{F_\infty},$ and consequently may be considered as 
an element of 
$[ \bigoplus_{\sigma: K \hookrightarrow \C} \Hom(F_\sigma, 
G_\sigma) ]^{F_\infty}.$
Moreover its image under the map
$$b:  [ \bigoplus_{\sigma: K \hookrightarrow \C} \Hom(F_\sigma, 
G_\sigma) ]^{F_\infty} \rightarrow   \Exthat^1_S (F,G)$$
is precisely the class $[\ol{\cE}]$ of the admissible extension
$\ol{\cE}$.

In the sequel, for any embedding $\sigma \hookrightarrow \C,$ the
$\C$-vector spaces $E_\sigma,$ $F_\sigma,$ and  $G_\sigma$ are
equipped with the hermitian structures defining $\ol{E},$ $\ol{F}$
and $\ol{G}$, and $F_\sigma\oplus G_\sigma$ with the  direct sum of
the hermitian structures of $F_\sigma$ and $G_\sigma.$
We shall also denote $\|.\|_{\infty,\sigma}$ (resp. 
$\|.\|_{HS,\sigma}$) the operator norm (resp.  the Hilbert-Schmidt 
norm) on $\C$-linear maps between 
some of the hermitian spaces $E_\sigma,$ $F_\sigma,$ $G_\sigma$, or
$F_\sigma\oplus G_\sigma.$ Finally we define  
$$\theta_\sigma: E_\sigma \simeq F_\sigma \oplus G_\sigma$$
as the orthogonal trivialization over $\C$ of the extension
$\cE_\sigma$ of $\C$-vector spaces. By definition, it is unitary when
$E_\sigma$ and $F_\sigma \oplus G_\sigma$ are equipped with the above
hermitian structures, and we have:
$$s^\perp =(\theta_\sigma^{-1}\circ (0, Id_{F_\sigma})).$$

\begin{proposition}\label{proposition.splittingsize}
1) The map 
$$
\begin{array}{crcl}
\beta : & \Triv_\OK (\cE)  & \longrightarrow & b^{-1} ( [\ol{\cE}]) \\
 & \varphi & \longmapsto &  s^\perp - ((\varphi^{-1}\circ (0,
Id_F))_\sigma)
\end{array}
$$
is a bijection.

2) For any $\varphi$ in $\Triv_\OK (\cE),$ the norms of 
$(\varphi_\sigma)$ and of its image $(h_\sigma):=\beta(\varphi)$ 
satisfy the relations:
$$\|\varphi_\sigma\|_{HS,\sigma}^2=  
\|\varphi_\sigma^{-1}\|_{HS,\sigma}^2 = \rk E +   
\|h_\sigma\|_{HS,\sigma}^2,$$
$$\|\varphi_\sigma -\theta_{\sigma}\|_{HS,\sigma}^2=  
\|\varphi_\sigma^{-1}-\theta_{\sigma}^{-1}\|_{HS,\sigma}^2 =   
\|h_\sigma\|_{HS,\sigma}^2,$$
and
$$\|\varphi_\sigma -\theta_{\sigma}\|_{\infty,\sigma}=  
\|\varphi_\sigma^{-1}-\theta_{\sigma}^{-1}\|_{\infty,\sigma} \leq  
\|h_\sigma\|_{\infty,\sigma} \leq \|h_\sigma\|_{HS,\sigma}.$$

\end{proposition}

\proof This is a straightforward consequence of the definitions
and of the following elementary lemma:

\begin{lemma} 

Let $p$ and $q$ be non-negative integers, and $n:= p+q.$ 
For any matrix $A$ in $M_{p,q}(\C),$ the matrix
$$\tilde{A}:=
\begin{pmatrix} I_p & A \\ 0 & I_q \end{pmatrix}$$
belongs to ${\rm GL}_n(\C),$ its inverse is 
$$\tilde{A}^{-1}:=
\begin{pmatrix} I_p & -A \\ 0 & I_q \end{pmatrix},$$ and
 the Hilbert-Schmidt and operator norms of $A,$ $\tilde{A},$ 
and
$\tilde{A}^{-1}$ satisfy:
$$\|\tilde{A}\|_{HS}^2 = \|\tilde{A}^{-1}\|_{HS}^2= n +
\|A\|_{HS}^2,$$
$$\|\tilde{A}-I_n\|_{HS}^2 = \|\tilde{A}^{-1}-I_n\|_{HS}^2= 
\|A\|_{HS}^2,$$
and
$$\|\tilde{A}- I_n\|_{\infty} = \|\tilde{A}^{-1}-I_n\|_{\infty}\leq
\|A\|_{\infty} \leq
\|A\|_{HS}.$$

\end{lemma}

\qed

\begin{corollary}\label{cor.splittingsize} 
For any $\varphi$ in $\Triv_\OK (\cE)$, we have:
\begin{equation}\label{b1}
\begin{split}
[K:\Q]^{-1} \sum_{\sigma: k \hookrightarrow \C} \|\varphi_\sigma
-\theta_{\sigma}\|_{HS,\sigma}^2 & =  
[K:\Q]^{-1} \sum_{\sigma: k \hookrightarrow \C}
\|\varphi_\sigma^{-1}-\theta_{\sigma}^{-1}\|_{HS,\sigma}^2
\\ &  =   
[K:\Q]^{-1} \sum_{\sigma: k \hookrightarrow \C}
\|h_\sigma\|_{HS,\sigma}^2.
\end{split}
\end{equation}
When $\varphi$ runs over $\Triv_\OK (\cE)$, this quantity achieves
$\exp (2  \size(\ol{\mathcal E}))$ as 
minimal value.  Moreover any   $\varphi \in \Triv_\OK (\cE)$ which
achieves it also satisfies:
\begin{equation}\label{b2}
\begin{split}
[K:\Q]^{-1} \sum_{\sigma: k \hookrightarrow \C} \|\varphi_\sigma
-\theta_{\sigma}\|_{\infty,\sigma}^2 & =  
[K:\Q]^{-1} \sum_{\sigma: k \hookrightarrow \C}
\|\varphi_\sigma^{-1}-\theta_{\sigma}^{-1}\|_{\infty,\sigma}^2
\\ &  \leq   
\exp (2  \size(\ol{\mathcal E})).
\end{split}
\end{equation}\end{corollary} 

Observe that, from any admissible extension over $\Spec \OK$ as above
$$ {\ol{\mathcal E}}:0\longrightarrow  \ol{G}
\stackrel{i}{\longrightarrow}  \ol{E}
\stackrel{p}{\longrightarrow} \ol{F} \longrightarrow  0,$$
we derive an other one by duality:
$$ {\ol{\mathcal E}}^\lor : 0\longrightarrow  \ol{F}^\lor
\stackrel{p^t}{\longrightarrow} \ol{E}^\lor
\stackrel{i^t}{\longrightarrow} \ol{G}^\lor \longrightarrow  0.$$
The sets of trivializations of $\cE$ and $\cE^\lor$ are in bijection
\emph{via} the map
$$
\begin{array}{ccc}
\Triv_\OK (\cE)  & \longrightarrow   & \Triv_\OK (\cE^\lor)   \\
  \varphi & \longmapsto   & (\varphi^t)^{-1}.   
\end{array}
$$
As a consequence of Proposition \ref{proposition.splittingsize}, we
derive:

\begin{corollary}\label{sizedual} For any admissible extension
$\ol{\cE}$ over $\Spec  \OK,$ we have:
$$\size(\ol{\cE}^\lor) = \size( \ol{\cE}).$$
\end{corollary}

\begin{proposition}\label{sizedirectsummand} For any admissible
extension 
${\ol{\mathcal E}}:0\longrightarrow  \ol{G} \longrightarrow  \ol{E}
\longrightarrow \ol{F} \longrightarrow  0$ over $\Spec \OK$ 
such that $m:= \min (\rk F, \rk G)$ is positive, there exists a 
submodule $F'$ of the $\OK$-module $E$ satisfying $E = F' \oplus G$
and $$\degan \ol{F}' \geq \degan \ol{F} - \frac{m}{2} \log (1 +
e^{2\size (\ol{\cE})/m}).$$
\end{proposition}

As usual, $\ol{F}'$ denotes the hermitian vector bundle over $\Spec 
\OK$ defined by $F'$ equipped with the restriction of the hermitian 
structure of $\ol{E}.$

\proof Actually the image $F' := s(F)$ of any splitting $s$ of $\cE$ 
over $\OK$ such that
$$\size(\ol{\mathcal E})=\log ( [K:\Q]^{-1/2}\cdot
\|s-s^\perp \|_{\pi_{\ast}(\ol{F}^\lor \otimes \ol{G})})$$
satisfies the above conditions.  This follows from the concavity of 
the logarithm, combined with the following simple consequence of the 
polar decomposition of linear maps between finite dimensional 
hermitian vector spaces:

\begin{lemma}
Let $p$ and $q$ be positive integers, $n:= p+q,$ and $m:= \min
(p,q).$ 
For any matrix $S$ in $M_{p, q}(\C),$ the  matrix
$$\tilde{S}:=
\begin{pmatrix} S \\ I_q \end{pmatrix}$$
in $M_{n,q}(\C)$ defines a $\C$-linear map from $\C^q$ to $\C^n$ 
whose (operator) norm of the $q$-th exterior power satisfy the 
following upper bound:
$$ \log \|\wedge^q \tilde{S}\|^2 \leq m
 \log (1 + 
\|S\|_{HS}^2/m),$$
when $\C^q$, $\C^n$, and there $q$-th exterior powers are equipped with their 
standard hermitian structures.
\end{lemma}
\qed

\subsection{Size and operations on
extensions}\label{subsubsec.sizeoper}
As observed in \ref{compat}, the description (\ref{eq:extcurve})
of the group $\Exthat_{S}^1\,(F,G)$ is compatible with various
functorial operations on arithmetic extensions. This allows us to
study the behavior of sizes of extensions under these operations.

\subsubsection{Pushout and Pullback}
   For instance, if $\ol{F'}$ and $\ol{G'}$ are hermitian vector
   bundles over $\Spec {\mathcal O}_{K}$, and if
   $$\alpha: F' \longrightarrow F \;\; \mbox{ and } \;\;
   \beta: G \longrightarrow G'$$
    are morphisms of ${\mathcal O}_{K}$-modules, we easily get, for any
extension
    class $e$ in $\Exthat_{S}^1\,(F,G)$:
    \begin{equation}
      \size_{\ol{F'},\ol{G}}(e\circ\alpha) \leq
      \size_{\ol{F},\ol{G}}(e) + \log \max_{\sigma:K\hookrightarrow\C}
      \|   \alpha\|_{\sigma}^\infty,
        \label{eq:sizecomp1}
    \end{equation}
    and
    \begin{equation}
      \size_{\ol{F},\ol{G'}}(\beta\circ e) \leq
      \size_{\ol{F},\ol{G}}(e) + \log \max_{\sigma:K\hookrightarrow\C}
      \|   \beta\|_{\sigma}^\infty,
        \label{eq:sizecomp2}
    \end{equation}
    where $\|.\|_{\sigma}^\infty$ denotes the operator norm on
    $\Hom_{\C}(F'_{\sigma},F_{\sigma})$ (resp. on
    $\Hom_{\C}(G_{\sigma},G'_{\sigma})$) deduced from the hermitian
    norms $\|.\|_{\ol{F'}, \sigma}$ and $\|.\|_{\ol{F}, \sigma}$
(resp.
    from
    $\|.\|_{\ol{G}, \sigma}$ and $\|.\|_{\ol{G'}, \sigma}$).
    Applied to $F'=F$ and $G'=G$, (\ref{eq:sizecomp1}) and
    (\ref{eq:sizecomp2})
give a control on the variation of $\size_{\ol{F},\ol{G}}(e)$ when
the hermitian structures on $\ol{F}$ and $\ol{G}$ are modified.

\subsubsection{Inverse image}
   Consider now $K'$ a number field containing $K$ and
   \[
   g:S'=\Spec\,{\mathcal O}_{K'}\longrightarrow  S=\Spec\,{\mathcal O}_K
   \]
   the associated morphism of arithmetic curves. Define
   $\ol{F'}:=g^\ast\ol{F}$ and $\ol{G'}:= g^\ast\ol{G}$. Then the
pullback
   morphism
   \begin{equation}\label{map1}
   g^*:\widehat{\rm Ext}_S^1(F,G)\longrightarrow \widehat{\rm
Ext}_{S'}^1(F',G')
   \end{equation}
gets identified with the morphism of real tori
\[
\frac{\left(\bigoplus_{\sigma:K\hookrightarrow\C}
\Hom_{\C}(F_{\sigma},G_{\sigma})\right)^{F_{\infty}}}{\Hom_{{\mathcal O}_{K}}(F,G)}
\longrightarrow
\frac{\left(\bigoplus_{\sigma':K'\hookrightarrow\C}
\Hom_{\C}(F'_{\sigma'},G'_{\sigma'})\right)^{F_{\infty}}}
{\Hom_{{\mathcal O}_{K'}}(F',G')}
\]
deduced from the  diagram
\begin{equation}
    \begin{CD}
     \Hom_{{\mathcal O}_{K}}(F,G)@>>>\Hom_{{\mathcal O}_{K'}}(F',G')\\
     @VV{\iota}V                    @VV{\iota}V \\
    \left(\bigoplus_{\sigma:K\hookrightarrow\C}
\Hom_{\C}(F_{\sigma},G_{\sigma})\right)^{F_{\infty}} 
@>>> \left(\bigoplus_{\sigma':K'\hookrightarrow\C}
\Hom_{\C}(F'_{\sigma'},G'_{\sigma'})\right)^{F_{\infty}},
\end{CD}
     \label{eq:bigcart}
\end{equation}
where the upper horizontal arrow is the ``extension of scalars" from
${\mathcal O}_{K}$ to ${\mathcal O}_{K'}$, and the lower one, the ``diagonal" $\R$-linear
map
$$(T_{\sigma})_{\sigma:K\hookrightarrow\C} \longmapsto
(T_{\sigma'_{|K}})_{\sigma':K'\hookrightarrow\C}.$$

Observe that an ${\mathcal O}_{K'}$-linear map $T':F'\to G'$ descends
to an ${\mathcal O}_{K}$-linear map $T:F \to G$ iff the $K'$-linear map
$T'_{K'}:F'_{K'}\to G'_{K'}$
descends to a $K$-linear map $T_{K}:F_{K} \to G_{K}.$ Together
with the separability of $K'$ over $K$, this shows that the diagram
(\ref{eq:bigcart}) is cartesian.  Moreover, its upper horizontal
arrow  is an
isometry, in the sense that,
for any $(T_{\sigma})_{\sigma:K\hookrightarrow\C}$ in
$\bigoplus_{\sigma:K\hookrightarrow\C}
\Hom_{\C}(F_{\sigma},G_{\sigma})$,
$$[K':\Q]^{-1}\sum\limits_{\sigma':K'\hookrightarrow
   \C}\|T_{\sigma'_{|K}}\|_{\ol{F'}^\lor \otimes \ol{G'},\sigma'}^2
   =
   [K:\Q]^{-1}\sum\limits_{\sigma:K\hookrightarrow
   \C}\|T_{\sigma}\|_{\ol{F}^\lor \otimes \ol{G},\sigma}^2.$$

From these remarks, we get:

\begin{proposition}\label{prop.sizedecr}
     With the above notation, the base change morphism (\ref{map1}) is
injective, and
for any $e$ in $\widehat{\rm Ext}_S^1(F,G)$, we have:
\begin{equation}
     \size_{\ol{F'},\ol{G'}}(g^\ast e) \leq
     \size_{\ol{F},\ol{G}}(e).
     \label{eq:sizedecr}
\end{equation}
\end{proposition}

By considering a geometric analogue of the above notion of size,
concerning extensions of vector bundles over a smooth projective
curve over a field, one is led to wonder whether the inequality
(\ref{eq:sizedecr}) would not actually be an equality (see
  \ref{subsec.geoarg} \emph{infra}). We shall
investigate this issue more closely in Section
\ref{sect:inv} below.

\subsection{Covering radius and size of admissible extensions}

Recall that, if $\Gamma$ is a lattice and $B$ a convex, symmetric body
in some finite dimensional real vector space $V$, the corresponding
\emph{inhomogeneous minimum}, or \emph{covering radius}, is defined
as\footnote{In \cite{banaszczyk93}, it is denoted $\mu (\Gamma, B)$.
We depart
from this notation to avoid confusion with slopes.}
$$\rho (\Gamma, B):= \min \{ r \in \R^\ast_{+} \mid \Gamma + rB =
V\}.$$
In particular, if $\ol{E}=(E, \|.\|)$ is an hermitian vector bundle
over $\Spec \Z$ --- or equivalently an euclidean lattice --- we may
consider the covering radius attached to the lattice $E$ and to the
unit ball $B_{\ol{E}}:=\{v\in E_{\R}\mid \|v\|\leq 1 \}$ in $E_{\R}$.
We shall
denote it $\rho(\ol{E})$.

This invariant of euclidean lattices is relevant for estimating the
sizes of arithmetic extensions over arithmetic curves. Indeed, by the
very definitions of these sizes and of the covering radius of an
euclidean
lattice, we have:

\begin{proposition}\label{prop.obvious}
     For any number field $K$ and any two hermitian vector bundles
     $\ol{F}$ and $\ol{G}$ over $S:=\Spec {\mathcal O}_{K},$ we have\footnote{As
     before, we denote $\pi$ the morphism from $\Spec {\mathcal O}_{K}$ to
     $\Spec \Z$.}:
     \begin{equation}
         \max \Bigl\{\size (\ol{\mathcal E}) \,\Big|\,
\genfrac{}{}{0pt}{}{\, \ol{\mathcal E} \, \mbox{\rm admissible
extension }}
{\mbox{\rm of } \ol{F}\mbox{ \rm by }\ol{G}} \Bigr\}
         = \log \left([K:\Q]^{-1/2}\cdot\rho(\pi_{\ast}
(\ol{F}^\lor\otimes
         \ol{G}))\right).
         \label{eq:obvious}
     \end{equation}
\end{proposition}

The so-called ``transference theorems" of the geometry of numbers
relate minima of various kinds attached to some lattice to minima of
the dual lattice. To get upper bounds on sizes of arithmetic
extensions, we shall use the following transference result of
Banaszczyk
(\cite[Theorem 2.2]{banaszczyk93}):
\begin{theorem}\label{thm.ban}
For any hermitian vector bundle $\ol{E}$ of rank $n$ over $\Spec \Z$,
we have:
   \begin{equation}
       \rho(\ol{E})\cdot\lambda_{1}(\ol{E}^\lor) \leq \frac{n}{2}.
       \label{eq:Ban}
   \end{equation}
\end{theorem}

  The covering radius of an hermitian vector bundle $\ol{E}$ of
  positive rank $n$ over $\Spec \Z$ satisfies also the following
  elementary lower bounds:
  \begin{equation}
      \rho(\ol{E})\cdot \lambda_{1}(\ol{E}^\lor) \geq \frac{1}{2}
      \label{eq:low1}
  \end{equation}
  and
   \begin{equation}
      v_{n} \rho(\ol{E})^n \geq {\rm covol}(\ol{E}).
      \label{eq:low2}
  \end{equation}
Indeed, (\ref{eq:low1}) follows from the surjectivity of the map
$$
\begin{array}{rcl}
     E_{\R}/E & \longrightarrow & \R / \Z  \\
{[v]} & \longmapsto & \xi(v) \mbox{ mod }1 \,\,,
\end{array}
$$
where $\xi \in E^\lor \setminus\{0\}$ is an element of minimal norm,
and (\ref{eq:low2}) from the surjectivity of the composition
$$B_{\ol{E}} \hookrightarrow E_{\R} \longrightarrow  E_{\R}/E.$$

Let us us emphasize that, as pointed out in \cite{banaszczyk93}, p.
633, Banaszczyk's upper bound (\ref{eq:Ban}) is basically optimal.
Indeed, as explained in \cite{Milnor73} (Chapter II, Theorem 9.5), by
using
Siegel's formula on
integral quadratic forms in a given genus,
Conway and Thompson have shown that one may construct a sequence
$\ol{CT}_{n}$ of rank $n$ euclidean lattices, $n\in \N$, such that
there exist (symmetric) isometric isomorphisms
\begin{equation}
     \ol{CT}_{n} \simeq \ol{CT}_{n}^\lor,
     \label{eq:CT1}
\end{equation}
and their first minima satisfy
\begin{equation}
     \lambda_{1}(\ol{CT}_{n}) \geq \sqrt{\frac{n}{2\pi e}} (1+o(1))
     \;\;\mbox{ as $n$ goes to infinity.}
     \label{eq:CT2}
\end{equation}
From (\ref{eq:CT1}), it follows that
$$\mbox{covol}(\ol{CT}_{n})=1,$$
and, by (\ref{eq:low2}), that
\begin{equation}
     \rho(\ol{CT}_{n})\geq v_{n}^{-\frac{1}{n}}\sim \sqrt\frac{n}{2
\pi e}
                        (1+o(1)) \;\;\mbox{ as $n$ goes to infinity.}
     \label{eq:CT3}
\end{equation}
Finally,
\begin{equation}
     \rho(\ol{CT}_{n})\lambda_{1}(\ol{CT}_{n}^\lor)
          =\rho(\ol{CT}_{n})\lambda_{1}(\ol{CT}_{n})
     \sim \frac{n}{2
     \pi e} (1+o(1))
     \;\;\mbox{ as $n$ goes to infinity,}
     \label{eq:CT4}
\end{equation}
and for the euclidean lattices $\ol{CT}_{n}$, Banaszczyk's
transference upper bound (\ref{eq:Ban}) is an equality up to a
bounded multiplicative factor.

With the notation of Proposition (\ref{prop.obvious}), when $\ol{F}$
and $\ol{G}$ have positive rank,we obtain from  Banaszczyk's upper
bound (\ref{eq:Ban}):
\begin{equation}
\log \rho(\pi_{\ast}(\ol{F}^\lor\otimes\ol{G}))\leq -\log
\lambda_{1}(\pi_{\ast}(\ol{F}^\lor\otimes\ol{G})^\lor) + \log
(\frac{1}{2}\cdot[K:\Q]\cdot\rk F \cdot \rk G).
  \label{eq:up1}
\end{equation}

Moreover, using successively  Propositions \ref{aradual} and
\ref{nelem} and (\ref{eq:udegatens}), we derive:

\begin{equation}
  \begin{split}
        -\log
        \lambda_{1}(\pi_{\ast}(\ol{F}^\lor\otimes\ol{G})^\lor) & =
        -\log
\lambda_{1}(\pi_{\ast}(\ol{F}\otimes\ol{G}^\lor\otimes\ol{\omega}_{K}))\\
        & \leq \udega (\ol{F}\otimes\ol{G}^\lor\otimes\ol{\omega}_{K})
        -\frac{1}{2}\log [K:\Q]\\
        &=\udega (\ol{F}\otimes\ol{G}^\lor) + \degan \ol{\omega}_{K}
                                              -\frac{1}{2}\log [K:\Q].
      \end{split}
    \label{eq:up2}
\end{equation}

From (\ref{eq:obvious}), (\ref{eq:up1}), and (\ref{eq:up2}), using
(\ref{eq:degom}) and Proposition \ref{prop.udegtens},
we finally obtain the upper bound on sizes of arithmetic extensions
announced in (\ref{sizeineq}):

\begin{theorem}\label{prop.sizeup}
     For any number field $K$, for any two hermitian vector bundles
     $\ol{F}$ and $\ol{G}$ of positive rank over $S:=\Spec {\mathcal O}_{K},$
     and for any admissible extension $\ol{\mathcal E}$ of $\ol{F}$
     by $\ol{G},$   we have:
     \begin{equation}
         \begin{split}
\size (\ol{\mathcal E})
         & \leq \udega (\ol{F}\otimes\ol{G}^\lor) + \frac{\log
         |\Delta_{K}|}{[K:\Q]}+ \log \frac{\rk F\cdot\rk G}{2}\\
          & \leq  \widehat{\mu}_{\rm 
max}(\ol{F})-\widehat{\mu}_{\rm min}(\ol{G})
          + \frac{\log
         |\Delta_{K}|}{[K:\Q]}+ \log \frac{\rk F\cdot\rk G}{2} .
              \end{split}
         \label{eq:sizeup}
     \end{equation}
     \end{theorem}

     Observe that, when $S$ is $\Spec \Z$ and $\ol{F}$ the trivial
hermitian
     line bundle, the first upper bound in
     (\ref{eq:sizeup}) is equivalent (by taking the logarithm) to
     Banaszczyk's bound (\ref{eq:Ban}). In particular, when moreover
     $\ol{G}=\ol{CT}_{n},$  it becomes an equality, up to a bounded
additive
     error term. Observe also that the right hand side of the second
     upper bound is invariant by unramified extension of the number
     field $K$.

     Similarly, the lower bounds (\ref{eq:low1}) and (\ref{eq:low2})
     lead to lower bounds on the maximal size of admissible extensions
     of $\ol{F}$ by $\ol{G}$. We leave this to the interested reader.

\subsection{The geometric case I}\label{geomext}
The results in the previous sections admit (simpler) analogs in the
geometric case where the number field $K$ is replaced by the function
field $k(C)$ defined by
a smooth projective geometrically connected curve $C$ over a field
$k$, and hermitian vector bundles over $\Spec {\mathcal O}_{K}$ by vector
bundles over $C$.

Recall that the slope of a such vector bundle $E$ of positive rank is
defined as the quotient
\begin{equation}
\label{eq:geoslope} \mu(E):=\frac{\deg E}{\rk E},
\end{equation}
its maximal slope $\mu_{\rm max}(E)$ as the maximum of the slope of a
subvector bundle of positive rank in $E$, and its minimal slope as
$$\mu_{\rm min}(E):=-\mu_{\rm max}(E^\lor).$$

We may also introduce the upper degree of $E$:
$${\rm udeg} (E):= \max\{\deg L\,|\, L \mbox{ sub-line bundle of } E
\}.$$
For any line bundle $M$ on $C$, we have:
\begin{equation}
     {\rm udeg} (E\otimes M)={\rm udeg} (E) + \deg M.
     \label{eq:udegtens}
\end{equation}

Observe that we are considering ``un-normalized" slopes: their
behavior
under some base extension involves the degree of the latter. Namely,
if $C'$ is  another smooth projective geometrically connected curve
$C$
over $k$, and $f:C' \to C$ a dominant $k$-morphism, then
$$\deg f^{\ast} E =\deg f \cdot \deg E,$$
and consequently,
$$\mu(f^{\ast}E) =\deg f \cdot \mu(E).$$
The maximal and minimal slopes satisfy similar formulae when $f$ is a
separable morphism (in particular when $\mbox{char}(k)=0$).

A (simplified) variant of the proof of Proposition \ref{prop.udegtens}
establishes:

\begin{proposition}\label{prop.udegtensgeo}
      For any two  vector bundles $F$ and $G$ of positive
     rank over $C$, we have:
     \begin{equation}
         {\rm udeg} (F\otimes G^\lor) \leq \mu_{\rm max}(F) -\mu_{\rm
         min}(G).
         \label{eq:udegtensgeo1}
     \end{equation}
     and
     \begin{equation}
         {\rm udeg} (F\otimes G) \leq \mu_{\rm max}(F) +\mu_{\rm
         max}(G).
         \label{eq:udegtensgeo2}
     \end{equation}
     \end{proposition}
In characteristic zero, (\ref{eq:udegtensgeo2}) and consequently
(\ref{eq:udegtensgeo1}) follows form the trivial upper bound
$${\rm udeg} (F\otimes G) \leq \mu_{\rm max}(F\otimes G)$$
combined with the equality
$$\mu_{\rm max}(F\otimes G)=\mu_{\rm max}(F) +\mu_{\rm
         max}(G),$$
which is
nothing but a reformulation of the classical fact that the tensor
product of two semi-stable vector bundles is semi-stable. However, in
positive characteristic, this equality does not hold in general.

Using Proposition \ref{prop.udegtensgeo}, we may establish a
generalization of the upper bound (\ref{extlines}) concerning
non-trivial extensions of vector bundles of positive ranks over $C$:

\begin{proposition}\label{prop.extvectgeo} Let $g$ denote the genus
of $C$.
For any extension $$
\mathcal{E}\colon \, 0\longrightarrow  G\longrightarrow
E\longrightarrow  F \longrightarrow  0
$$
of vector bundles of positive ranks over $C$ which does not split,
we have:
$$\label{geomest}
{\mu}_{\rm min}(G)
- {\mu}_{\rm max}(F) \leq 2g-2.$$
\end{proposition}

\proof The class of $\mathcal{E}$ provides a non-zero element in $$
\text{Ext}^1_C(F,G) \simeq H^1 (C,\,F^\lor \otimes G) \simeq
H^0 (C,\,\Omega_{C/k}^1 \otimes F \otimes G^\lor)^\lor.
$$
Hence $\Omega_{C/k}^1 \otimes F \otimes G^\lor$ has a non-trivial
regular section on $C$, and its upper degree is consequently
non-negative.
As (\ref{eq:udegtens}) and Proposition (\ref{prop.udegtensgeo}) show
that
$${\rm udeg} (\Omega_{C/k}^1 \otimes F \otimes G^\lor)=
  {\rm udeg} (F \otimes
G^\lor)+ 2g-2 \leq \mu_{\rm max}(F) -\mu_{\rm
         min}(G) + 2 g-2,  $$
this proves (\ref{geomest}).
\qed

\subsection{The geometric case II}\label{geomextII}
It turns out that the geometric analogues of our arithmetic results
discussed in the previous sections possess some ``refined versions"
that
more closely
parallel our investigation of the arithmetic case, where
\emph{archimedean places} of  number fields play a special role.

To formulate it, let  $C$ be as above a smooth 
geometrically connected projective curve over
some field $k$, and let
$$D=\sum_{P \in |D|} n_P.P$$ be an effective
divisor on $C$ with non-empty support $|D|$, and let
$\Co$ be the affine curve $C\setminus |D|.$

The function field $K:=k(C)=k(\Co)$ is a geometric analogue of a
number field --- in this analogy, the curve $\Co$ 
plays the role of the affine scheme
$\Spec {\mathcal O}_{K}$ where $K$ is a number field, and the points of $|D|$
of its
archimedean places. Moreover a vector bundle over $C$ (resp., its
restriction $\Eo$ to $\Co$) is the
counterpart of an hermitian vector bundle
$(E,(\|.\|_\sigma)_{\sigma:K\rightarrow C})$ (resp. of the vector
bundle $E$) over $\Spec {\mathcal O}_{K}$. In this section, we want to extend
this dictionary by describing the analogues, in 
this geometric setting, of our arithmetic
extension groups and of the size function on them.

For any closed point $P$ of $C$, let  $\Oc_{C,P}$ the completion of
the local ring ${\mathcal O}_{C,P}$ of $C$ at $P,$  $\m_P$ its maximal ideal,
and $\Kc_P$ its field of fractions.

It is natural to define an analogue of the 
arithmetic extension group over an arithmetic 
curve
by mimicking its description in Proposition \ref{cor.affine}.
Namely,
if$F$ and $G$ are two vector bundles over $C$, and $\Fo$ and $\Go$
their restrictions to $\Co,$
we let
\begin{equation}\label{eq:defgeo}
\Exthat^1_{\Co}(\Fo,\Go):=\frac{\bigoplus_{P\in |D|} \Hom_\Co (\Fo,
\Go)\otimes_{{\mathcal O}(\Co)}\Kc_P}{\iota (\Hom _\Co (\Fo,  \Go))},
\end{equation}
where $$\iota: \Hom _\Co (\Fo,  \Go) 
\longrightarrow \bigoplus_{P\in |D|} \Hom_\Co 
(\Fo,
\Go)\otimes_{{\mathcal O}(\Co)}\Kc_P$$
denotes the diagonal embedding, defined by $$\iota(x):= (x\otimes
1_{\Kc_{P}})_{P\in
|D|}.$$

Observe that its image is discrete in
$\bigoplus_{P\in |D|} \Hom_\Co (\Fo,
\Go)\otimes_{{\mathcal O}(\Co)}\Kc_P$ equipped with the product topology deduced
from
the $\m_P$-adic topologies on the finite dimensional $\Kc_P$-vector
spaces  $$\Hom_\Co (\Fo,
\Go)\otimes_{{\mathcal O}(\Co)}\Kc_P \simeq (\check{F}\otimes G)_K\otimes_K
\Kc_P\simeq (\check{F}\otimes G)_{{\mathcal O}_{X,P}}\otimes_{{\mathcal O}_{X,P}}
\Kc_P.$$
Consequently the  quotient topology on the 
abelian group $\Exthat^1_{\Co}(\Fo,\Go)$  is
separated and complete.

Besides, a neighbourhood basis of zero in 
$\Exthat^1_{\Co}(\Fo,\Go)$ equipped with this 
topology
is formed by the images  in $\Exthat^1_{\Co}(\Fo,\Go)$ of the
subgroups
\begin{equation}\label{eq:Nthnbd}
\bigoplus_{P \in |D|} (\check{F}\otimes
G)_{{\mathcal O}_{X,P}}\otimes_{{\mathcal O}_{X,P}} \m_P^{N.n_P}, \,\,\, N \in \N
\end{equation}
of
$$\bigoplus_{P \in |D|} (\check{F}\otimes
G)_{{\mathcal O}_{X,P}}\otimes_{{\mathcal O}_{X,P}}
\Kc_P,$$
and the quotient of  $\Exthat^1_{\Co}(\Fo,\Go)$ by the image  of
(\ref{eq:Nthnbd}) may be
identified with
$$\frac{\bigoplus_{P\in |D|} (\check{F}\otimes
G)_{{\mathcal O}_{X,P}}\otimes_{{\mathcal O}_{X,P}}
\Kc_P}{\iota (\Hom _\Co (\Fo,  \Go)) + \bigoplus_{P\in |D|}
(\check{F}\otimes
G)_{{\mathcal O}_{X,P}}\otimes_{{\mathcal O}_{X,P}} \m_P^{N.n_P}}.$$
In turn, this space is canonically isomorphic to
$$H^1(C, \check{F}\otimes G (-N.D)) \simeq \Ext^1_C(F, G(-N.D)).$$
This follows from the long exact sequence of cohomology groups
associated to the ``localization" short exact sequence of sheaves of
${\mathcal O}_{C}$-modules:
\begin{multline}
0\longrightarrow \check{F}\otimes 
G(-N.D)\longrightarrow 
j_{\ast}j^\ast(\check{F}\otimes G) \oplus 
\bigoplus_{P\in |D|}
i_{P \ast}((\check{F}\otimes
G)_{{\mathcal O}_{X,P}}\otimes_{{\mathcal O}_{X,P}} \m_P^{N.n_P})\\ \longrightarrow
\bigoplus_{P \in |D|}i_{P \ast}((\check{F}\otimes
G)_{{\mathcal O}_{X,P}}\otimes_{{\mathcal O}_{X,P}}
\Kc_P )\longrightarrow 0,
\end{multline}
where $j$ (resp. $i_{P}$) denotes the open (resp. closed) immersion
$\Co \hookrightarrow C$ (resp. $\{P\} \hookrightarrow C$).

Finally, the topological group $\Exthat^1_{\Co}(\Fo,\Go)$ is
canonically isomorphic to the projective limit
$$\lim_{\stackrel{\longleftarrow}{N \in \Z}} \Ext^1_C(F,G(-N.D))$$
of the finite dimensional $k$-vector spaces  $\Ext^1_C(F,G(-N.D))$
equipped
with the discrete topology.  In particular, it is 
a linearly compact topological
$k$-vector space. (This is
similar to the compactness of $\Exthat^1_{\Spec {\mathcal O}_K} (F,G)$ in the
number field case.)

For any integer $N \in \Z,$ let
$$\pi_N: \Exthat^1_{\Co}(\Fo,\Go) \simeq 
\lim_{\stackrel{\longleftarrow}{M \in \Z}} 
\Ext^1_C(F,G(-M.D))
\longrightarrow  \Ext^1_C(F,G(-N.D))$$
be the projection on the $N$-th component. To define a geometric
counterpart of
the size of arithmetic extensions, we let, for every $e \in
\Exthat^1_{\Co}(\Fo,\Go)$:
$$\size_{F,G}(e):= -\inf \{N\in \Z \mid \pi_N(e) \neq 0\}.$$
It is straightforward that it is an element of $\Z \cup \{-\infty\},$
finite for any $e\neq 0,$ and that Proposition
\ref{prop.extvectgeo} applied to the vector 
bundles $F$ and $G(-N.D)$ yields the
following geometric analogue of Theorem \ref{prop.sizeup}:
\begin{equation}  \label{eq:sizeupgeo}
    \deg D. \size_{F,G}(e)\leq \mu_{\max} (F) -\mu_{\min}(G) +2g-2,
\end{equation}
when $F$ and $G$ have positive rank\footnote{The occurrence of $\deg
D$
in the right hand side of (\ref{eq:sizeupgeo}) is related to the use
of ``non-normalized" slopes (\ref{eq:geoslope}) instead of normalized
slopes (\ref{eq:arslope}) in the arithmetic case.}.

\section{Sizes of admissible extensions: explicit computations
  and an application to reduction theory}

In this section, we want to show how  evaluating the size of
admissible
  extensions is related to basic questions in the geometry of
lattices.
  Firstly we compute it explicitly in some elementary examples ---
  notably in the most basic case of extension of hermitian line
  bundles over $\Spec \Z$. Then we consider
  the size of the restriction of the universal admissible
extension over $\P^1_\Z$ at rational points in $\P^(\Q)$ ($\simeq
\P^1(\Z)$) and relate it to the
usual logarithmic height of these points, and to the geometry of
the Ford
circles\footnote{Namely, the horocycles image
of $\{ \Im z =1 \}$ under the action of $SL^2(\Z)$.} in the upper
half-plane.
Finally, using the upper bound (\ref{sizeineq}) on the sizes of
admissible
extensions over $\Spec {\mathcal O}_{K}$, we derive the existence of some
``almost-splitting" for any hermitian vector bundle $\ol{E}$ over
$\Spec {\mathcal O}_{K}$, namely the existence of $n:= \rk E$ hermitian lines
bundles
$\ol{L}_{1},$ \ldots, $\ol{L}_{n},$  and of an isomorphism
of ${\mathcal O}_{K}$-modules
\[
\phi: E \stackrel{\sim}{\longrightarrow}
     \bigoplus_{i=1}^{n}L_{i}
\]
     such that the archimedean norms of $\phi$ and $\phi^{-1}$
     (defined by means of the hermitian metrics on $\ol{E}$ and
     $\bigoplus_{i=1}^{n}\ol{L}_{i})$ are bounded in terms of $K$ and
     $n$ only. When $K=\Q,$ this is basically the main result of the
     classical reduction theory of positive quadratic forms. For
     general number fields $K$, our method yields explicit bounds on
     the norms of $\phi$ and $\phi^{-1}$, which improve on the
     qualitative results which may be derived from
     the general reduction theory for reductive algebraic groups
     over number fields.

\subsection{Some explicit computations of
size}
We discuss various examples of admissible extensions
over arithmetic curves, the sizes of which can be ``explicitly"
computed.

   \begin{example}\label{exAn}\rm
   For any positive integer $n$, the morphism of ${\mathcal O}_{K}$-modules
   $$
     \begin{array}{cccc}
         p: & {\mathcal O}_{K}^{n+1} & \longrightarrow & {\mathcal O}_{K}  \\
          & (x_{i})_{0\leq i\leq n} & \longmapsto &
\sum_{i=0}^{n}x_{i}
     \end{array}
     $$
     defines an admissible extension:
     $$\overline{\mathcal A}_{n,K}:
     \,0\longrightarrow \overline{\ker p}\longrightarrow
     \overline{{\mathcal O}_{S}}^{\oplus (n+1)}\stackrel{p}{\longrightarrow}
     \overline{L}_{n,S}\longrightarrow  0,$$
     where $\overline{{\mathcal O}_{S}}^{\oplus (n+1)}$
     denotes the trivial hermitian vector bundle of rank $n+1$
     over $S$, and  $\overline{L}_{n,S}$  the hermitian line bundle
     $({\mathcal O}_{K}, (\|.\|_{n,\sigma})_{\sigma:K\hookrightarrow\C}))$
     defined by the hermitian metrics:
     $$\|1\|_{n,\sigma}:=\frac{1}{\sqrt{n+1}}.$$
     This admissible extension is the base change to $S$ of the
     admissible extension over $\Spec \Z$ defining the root lattice
     $A_{n}$ (see \ref{subsec.rootlattice} below).

     The orthogonal splitting $s^\perp$ of $\overline{\mathcal
     A}_{n,K}$ satisfies

$$s^\perp_{\sigma}(1)=(\frac{1}{n+1},\ldots,\frac{1}{n+1})\,\,\,\mbox{
     for every embedding }\sigma:L\hookrightarrow\C.$$
     An algebraic splitting $s$ of ${\mathcal A}_{n,K}$ over ${\mathcal O}_{K}$
     is given by $$s(1):=(1,0,\ldots,0).$$
     Since the euclidean norm of
    $$s^\perp_{\sigma}(1)-s(1)
       =(-\frac{n}{n+1},\frac{1}{n+1},\ldots,\frac{1}{n+1})$$
    is $\sqrt{n/n+1}$, and $\|1\|_{n,\sigma}=1/\sqrt{n+1},$ this
    shows that
    $$\size(\overline{\mathcal A}_{n,K}) \leq \frac{1}{2}\log n.$$

    It is easy to check that this inequality is an equality when
    $K=\Q$. Our results on the invariance of size by base change in
    section \ref{subsubsec.bcheuclidean} below will show that, for any
    number field $K$,
    $$\size(\overline{\mathcal A}_{n,K}) =\frac{1}{2}\log n.$$
    \end{example}

\begin{example}\rm
    Let $S={\rm Spec}\,\Z$, and let $\ol{F}$ and $\ol{G}$ be hermitian
    line bundles. We may choose generators $f$ and $g$ for $F$ and $G$
    over $\Z$; they are well defined up-to-sign, and determine a group
    isomorphism:
    $$
       \begin{array}{rcl}
        \Exthat^{1}_{S}\, (F,G) & \longrightarrow  & \Exthat^1_{\Spec
        \Z}(\Z, \Z) \simeq \R/\Z  \\
        e & \longmapsto & \tilde{e}:=g^{-1}\circ e \circ f^{-1}.
    \end{array}$$
    Then, if we let
    $$\begin{array}{rcl}
        \delta: \R/\Z & \longrightarrow & [0, 1/2] \\
        {[x]} &\longmapsto  & \min_{k \in \Z} |x-k| \;(=|x| \mbox{ if
}
        |x| \leq 1/2),
    \end{array}$$
     the size of an extension class $e \in \Exthat^{1}_{S}\, (F,G)$
     may be expressed as
     \begin{equation}
         \begin{split}
     \size_{\ol{F},\ol{G}}(e) & = \log \delta (\tilde{e}) + \log
     \|g\| - \log \|f\|\\
     &= \log \delta (\tilde{e}) + \dega \ol{F} -\dega \ol{G}.
     \end{split}
     \end{equation}
     This clearly implies (\ref{sizeineq}) in this special situation.

     \end{example}

   \begin{example}\label{ext11}\rm
Let $\ol{E}$ be an hermitian vector bundle over $S:= \Spec \Z.$
Consider the attached projective space over $S$:
$$\pi:\P(E):= \mathbf{Proj}\, \mathrm{Sym}(E) \longrightarrow S,$$
and the tautological exact sequence of vector bundles over $\P(E)$:
$$ 0 \longrightarrow V \longrightarrow \pi^\ast E \longrightarrow
{\mathcal O}_{E}(1)\longrightarrow 0.$$
We may equip $V$ (resp. ${\mathcal O}_{E}(1)$) with the induced (resp. quotient)
hermitian structure deduced from $\pi^\ast \ol{E}$. In such a way, we
define the \emph{tautological admissible extension} over $\P(E)$:
\begin{equation}
     \overline{\mathcal E}:\
     0 \longrightarrow \ol{V} \longrightarrow \pi^\ast \ol{E}
\longrightarrow
\ol{{\mathcal O}}_{E}(1)\longrightarrow 0.
     \label{eq:tautadmis}
\end{equation}

The function
$$
\begin{array}{rcl}
     \P(E)(\Q) \simeq \P(E)(\Z) & \longrightarrow  & [-\infty,
+\infty[  \\
     P & \longmapsto & \size (P^\ast \overline{\mathcal E})
\end{array}
$$
may be described as follows.

For any $P \in \P(E)(\Z)$, the line bundle $P^\ast \ol{{\mathcal O}}_{E}(1)$
over $\Spec \Z$ is trivial, and $P^\ast \overline{\mathcal E}$ is
isomorphic to an admissible extension of the form
\begin{equation}
     0 \longrightarrow \ol{\ker a} \longrightarrow \ol{E}
    \stackrel{a}{\longrightarrow} (\Z, \|.\|_{P}) \longrightarrow 0,
     \label{eq:adma}
\end{equation}
where $a$ denotes a surjective morphism in $E^\lor:=\Hom_{\Z}(E,\Z)$,
and $\|.\|_{P}$ the quotient norm on $\Z_{\R}\simeq \R$. Let
$\tilde{a}$ denote the  image of $a$ by the isomorphism
$E^\lor_{\R}\simeq E_{\R}$ determined by the euclidean structure
$\|.\|_{\ol{E}}$ on $E_{\R}.$ The orthogonal splitting $s^\perp: \R
\to E_{\R}$ of (\ref{eq:adma}) satisfies
$$s^\perp(1) = \frac{\at}{\|\at\|_{\ol{E}}^2},$$
and therefore
$$ \|1\|_{P}=
\|s^\perp(1)\|_{\ol{E}}=
\|\at\|_{\ol{E}}^{-1}.$$
Consequently,
\begin{equation}
   \size (P^\ast \overline{\mathcal E}) =\log \min \left\{
\|m-\|\at\|_{\ol
   E}^{-2} \at \|_{\ol{E}} \cdot
        \|\at\|_{\ol{E}}\,\,, \, m \in E \mbox{ such that }
   a(m)=1
   \right\}.
     \label{eq:sizePast1}
\end{equation}

Observe that, for any $v \in E_{\R},$
$$ \|v\|_{\ol{E}}^2 \cdot \|\at\|^2_{\ol{E}}=\langle v,\at
\rangle _{\ol{E}}^2 +
\|v\land \at\|_{\Lambda^2\ol{E}}^2.$$
For any $m$ in $E$ such that $a(m)=1,$ the vector $m-\|\at\|_{\ol
E}^{-2} \at$
is orthogonal to $\at$, and the identity above shows that
$$\|m-\|\at\|_{\ol E}^{-2} \at \|_{\ol{E}}\cdot \|\at\|_{\ol{E}}=
   \|m \land \at \|_{\Lambda^2 \ol{E}}.$$
   This shows that the expression (\ref{eq:sizePast1}) for the size
   of $P^\ast \overline{\mathcal E}$ may also be written:
\begin{equation}
   \size (P^\ast \overline{\mathcal E}) =\log \min \left\{ \|m \land
\at
   \|_{\Lambda^2 \ol{E}}\,\,,\, m \in E \mbox{ such that }
   a(m)=1
   \right\}.
     \label{eq:sizePast2}
\end{equation}
   \end{example}

   \subsection{Universal extensions and heights 
over $\P^1_\Z$}\label{subsec.universal}

Let us specialize the previous example \ref{ext11} to the
    situation where $\ol{E}$ is the trivial hermitian vector bundle of
    rank 2 over $\Spec \Z$. In other words, $\ol{E}=(E, \|.\|)$ where
    $\|.\|$ denotes the ``standard" hermitian metric on $\Z^2\otimes
    \C= \C^2,$ and $\P(E)$ is the ``projective line" $\P_{\Z}^1.$

    For any $P$ in $\P(E)(\Q)=\P^1(\Q),$ we may choose homogeneous
    coordinates $(a_{0}:a_{1})$ such that $a_{0}$ and $a_{1}$ are
    coprime integers, and $m_{0}$ and $m_{1}$ be coprime integers
such that
    $a_{0}m_{0}+a_{1}m_{1}=1.$ Then the class in $\R/\Z$ of
    $$\frac{m_{0}a_{1}-m_{1}a_{0}}{a_{0}^2+a_{1}^2}$$ depends only
    on $P$, and the discussion above applied to
    $a=(a_{0},a_{1})$ and $\tilde{a}=\tbinom{a_{0}}{a_{1}}$ shows
    that
  $$\size (P^\ast \overline{\mathcal E}) =\log \delta
        \left(\frac{m_{0}a_{1}-m_{1}a_{0}}{a_{0}^2+a_{1}^2}\right) +
\log
        ({a_{0}^2+a_{1}^2}).$$
    Observe that the second term in the right-hand side is the
    \emph{height of} $P$ with respect to $\ol{{\mathcal O}}_{E}(1)$:
    $$\log \sqrt{a_{0}^2+a_{1}^2} = h_{\ol{{\mathcal O}}_{E}(1)}(P)
    :=\dega P^\ast \ol{{\mathcal O}}_{E}(1)
    =- \dega P^\ast \ol{V}.$$
    Therefore, if we let
$$s(P):=\left[\frac{m_{0}a_{1}-m_{1}a_{0}}{a_{0}^2+a_{1}^2}\right]
\,\,(\in \R/\Z),$$
we have
    \begin{equation}
        \size (P^\ast \overline{\mathcal E}) =\log \delta (s(P))
        + 2 h_{\ol{{\mathcal O}}_{E}(1)}(P).
        \label{eq:sizePast3}
    \end{equation}

The matrix
\[
\gamma:=\begin{pmatrix}
a_1 & -m_0 \\
a_0 & m_1
\end{pmatrix}
\]
belongs to $SL_{2}(\Z)$and we have
\footnote{We let
       $SL_{2}(\Z)$ act on $\P^1(\Q)=\Q \cup \{\infty\}$ and
     the upper half-plane by the usual formula $\left(
     \begin{smallmatrix} a & b\\ c& d \end{smallmatrix}
     \right)\cdot z:=(az+b)/(cz+d).$}
     \[
     {\gamma}\cdot\infty=\frac{a_{1}}{a_{0}}.
     \]
    The inverse matrix of ${\gamma}$,
    $${\gamma}^{-1}:=
    \begin{pmatrix} m_1 & m_0 \\ -a_0 & a_1 \end{pmatrix},
        $$
        satisfies:
     $$\mathrm{Re}({\gamma}^{-1}\cdot i)
         =\frac{m_{0}a_{1}-m_{1}a_{0}}{a_{0}^2+a_{1}^2}.
     $$
     Consequently, the map $s:\P^1(\Q)\to \R/\Z$ admits the
     following description in terms of the actions of
     $\Gamma:=SL_{2}(\Z)$ on $\P^1(\Q)$ and the upper half-plane:
     if
       \begin{displaymath}
     \Gamma_{\infty}
     := \left\{\gamma \in \Gamma \mid \gamma\cdot\infty
     =\infty \right\} 
     =\left\{\left(
     \begin{smallmatrix} \epsilon & n\\ 0& \epsilon \end{smallmatrix}
     \right), \epsilon=± 1, n \in \Z \right\},
    \end{displaymath}
     then the map $s$ is characterized by the commutativity of the
     following diagram:
   \begin{equation}
    \begin{array}{cccc}
              & {[} \gamma ]  & \longmapsto & \gamma
              \cdot\infty\\
            {[} \gamma ] & \Gamma/\Gamma_{\infty} &
             \stackrel{\sim}{\longrightarrow} & \P^1(\Q)  \\
            \Big\downarrow & \Big\downarrow  &  &
           \Big\downarrow\vcenter{%
\rlap{$s$ }}
             \\
        {[} \Re (\gamma^{-1}\cdot i)] & \R/\Z & = & \R/\Z.
         \end{array}
         \label{sgroup}
  \end{equation}

     \begin{proposition}\label{prop.sizeProj}
         With the above notation, for any $P \in \P^1(\Z),$ the
         extension
         class of the admissible extension $P^\ast \overline{\mathcal
E}$
         vanishes iff $P$ is $0:=(1:0)$ or $\infty:=(0:1)$.

         Moreover, for any $P$ in $\P^1(\Z)\setminus\{0,\infty\}$, we
have:
         \begin{equation}
             -\frac{1}{2}\log 2 + h_{\ol{{\mathcal O}}_{E}(1)}(P) \leq
\size (P^\ast \overline{\mathcal E}) \leq -\log 2 + 2
h_{\ol{{\mathcal O}}_{E}(1)}(P).
             \label{eq:sizeProj}
         \end{equation}
     \end{proposition}

     \proof As above, let us choose $a_{0}$ and $a_{1}$ prime
     together such that $P=(a_{0}:a_{1}).$

     To establish the first assertion, observe that the
     following conditions are successively equivalent:

     - the admissible extension $P^\ast \overline{\mathcal E}$ is
split;

     - there exists $(m_{0}, m_{1})$ in $\Z^2$ such that
     $a_{0}m_{0}+a_{1}m_{1}=1$ and $m_{0}a_{1}-m_{1}a_{0}=0;$

     - there exists $k \in \Z$ such that
     $(m_{0},m_{1})=k(a_{0},a_{1})$ and $a_{0}m_{0}+a_{1}m_{1}=1;$

     - $a_{0}^2+a_{1}^2=1;$

     - $(a_{0},a_{1})$ belongs to $\{(1,0), (-1,0), (0,1), (0,-1) \}.$

     The second estimate in (\ref{eq:sizeProj}) follows from
     (\ref{eq:sizePast3}), since the values of $\delta$ lie in $[0,
1/2]$.

     To
     derive the first one, without loss of generality, we may assume
     $a_{0}$ and $a_{1}$ positive. Let us choose $b_{0}$ and $b_{1}$
two
     integers satisfying
     $a_{0}b_{1}-a_{1}b_{0}=1$ such that $b_{0}$ is non-negative and
     minimal, or equivalently, in $\{0, \ldots, a_{0}-1 \}$. Then
     $a_{0} b_{1}=1+a_{1} b_{0}$ is positive and $\leq 1+
     a_{1}(a_{0}-1) \leq a_{0} a_{1}.$ Consequently $b_{1}$ belongs
     to $\{1, \ldots, a_{1}\},$ and the following inequality holds:
     $$0 \leq \frac{b_{0}}{a_{0}} \leq
\frac{a_{0}b_{0}+a_{1}b_{1}}{a_{0}^2+a_{1}^2}<\frac{b_{1}}{a_{1}}\leq
     1.$$

     If $b_{0}\neq 0$ and $b_{1}\neq a_{1},$ this implies
     $$\frac{1}{a_{0}} \leq
\frac{a_{0}b_{0}+a_{1}b_{1}}{a_{0}^2+a_{1}^2}
     \leq 1- \frac{1}{a_{1}},$$
     and consequently,

$$\delta\left(\frac{a_{0}b_{0}+a_{1}b_{1}}{a_{0}^2+a_{1}^2}\right)
     \geq \min (\frac{1}{a_{0}},\frac{1}{a_{1}}) \geq
     \frac{1}{\sqrt{a_{0}^2+a_{1}^2}}.$$

     If $b_{0}=0,$ then $a_{0}b_{1}=1,$ and necessarily
     $a_{0}=b_{1}=1,$
     hence

$$\frac{a_{0}b_{0}+a_{1}b_{1}}{a_{0}^2+a_{1}^2}=\frac{a_{1}}{a_{1}^2+1}.$$
     This is a number in $[0,1/2],$ and

$$\delta\left(\frac{a_{0}b_{0}+a_{1}b_{1}}{a_{0}^2+a_{1}^2}\right)
     = \frac{a_{1}}{a_{1}^2+1}
     \geq \frac{1}{\sqrt{2} \sqrt{a_{1}^2+1}} = \frac{1}{\sqrt{2}
     \sqrt{a_{0}^2+a_{1}^2}}.$$

     If $b_{1}=a_{1}$, then $a_{1}(b_{0}-a_{0})=-1,$ and necessarily
     $a_{1}=1$  and $b_{0}-a_{0}=-1.$ Then
     $$\frac{a_{0}b_{0}+a_{1}b_{1}}{a_{0}^2+a_{1}^2}=
     \frac{a_{0}(a_{0}- 1) +1}{a_{0}^2+1}= 1-
     \frac{a_{0}}{a_{0}^2+1}.$$
     Since $a_{0}/(a_{0}^2+1)$ belongs to $[0, 1/2],$ we obtain:

$$\delta\left(\frac{a_{0}b_{0}+a_{1}b_{1}}{a_{0}^2+a_{1}^2}\right)
     = \frac{a_{0}}{a_{0}^2+1}
     \geq  \frac{1}{\sqrt{2}
     \sqrt{a_{0}^2+a_{1}^2}}.$$

     We have shown that the lower bound

$$\delta\left(\frac{a_{0}b_{0}+a_{1}b_{1}}{a_{0}^2+a_{1}^2}\right)
     \geq  \frac{1}{\sqrt{2}
     \sqrt{a_{0}^2+a_{1}^2}}$$
     always holds. This may also be written
     $$\log s(P)\geq -h_{\ol{{\mathcal O}}_{E}(1)}(P) - \frac{1}{2} \log 2,$$
and,
     according to (\ref{eq:sizePast3}), is equivalent to the first
inequality
     in (\ref{eq:sizeProj}).
     \qed

     Observe also that the estimates (\ref{eq:sizeProj}) are basically
     optimal.

     Indeed, if for any positive integer $n$, we let $P_{n}:=(1:n),$
     then
     $$h_{\ol{{\mathcal O}}_{E}(1)}(P_{n})=\log \sqrt{n^2+1}=\log n +
O(1/n^2),$$
     and, since $\left(
     \begin{smallmatrix} 1 & 1\\ n & n+1 \end{smallmatrix}\right)$
     belongs to $SL_{2}(\Z),$ we have
     $$s(P_{n})=\frac{n^2+n+1}{n^{2}+1}=\frac{n}{n^2+1} \mod{\Z}$$
     and
     $$\size (P^\ast \overline{\mathcal E}) - 2
     h_{\ol{{\mathcal O}}_{E}(1)}(P_{n})= \log \frac{n}{n^2+1} = -\log n +
     O(1/n^2).$$
     In particular, if $\mathcal F$ denotes the filter of complements
     in $\P^1(\Q)$ of finite subsets, we have:
     $$\liminf_{\mathcal F} \left(\size (P^\ast \overline{\mathcal
E}) -
     h_{\ol{{\mathcal O}}_{E}(1)}(P)\right)
     \leq 0.$$

     Moreover, the interpretation of $s$ given by the diagram
     (\ref{sgroup}) shows that $s(\P^1(\Q))$ is dense in $\R/\Z$.
     Indeed, the closure of $\Gamma\cdot i$ in $\P^1(\C)$ contains the
     limit set $\P^1(\R)$ for the action of $\Gamma$ in $\P^1(\C)$;
     therefore, $\Re (\Gamma\cdot i)$ is dense in $\R.$ This implies
that
     the set of values of
     $$\size (P^\ast \overline{\mathcal E}) - 2 h_{\ol{{\mathcal O}}_{E}(1)}(P)
=\log
     \delta (s(P))$$
     when $P$ runs over $\P^1(\Q)$ is dense in $\log \delta (\R/\Z) =
     [-\infty, -\log 2].$ In particular, with $\mathcal F$ as above,
      $$\limsup_{\mathcal F} \left(\size (P^\ast) \overline{\mathcal
E}
                - 2 h_{\ol{{\mathcal O}}_{E}(1)}(P)
      \right) \geq -\log 2.$$

\subsubsection{A geometric interpretation by means of Ford circles}
    It may be worth noting that the function $s:\P^1(\Q) \mapsto
    \R/\Z$ in terms of which we have expressed $\size (P^\ast
    \overline{\mathcal E})$ admits a geometric interpretation by means
    of the so-called \emph{Ford circles}.\footnote{These circles also
    appear in the literature under the names of Farey or Speiser
    circles. We refer the reader to \cite{ford38} and
\cite{rademacher64},
    chapter 6, for more information and references about their
properties.}

    Recall these are circles $C(q)$ in the upper half-plane $\{ \Re z
    \geq 0 \}$ in $\C$ associated to rational numbers $q\in \Q$: if
    $q=b/a$, where $a$ and $b$ are integers which are prime together,
    $C(q)$ is defined as the circle of center $q+i/2a^2$ and radius
    $1/2a^2.$ It meets tangentially  the line $\R$ at $q$.  It is
    convenient to define also $C(\infty)$ as the subset $(i+\R) \cup
    \{\infty\}$ in $\P^1(\C)=\C\cup\{\infty\}.$ Then, for any $q\in
    \P^1(\Q)$, $C(q)\setminus\{q\}$ is an horocycle in the Poincaré
    upper half-plane, and a straightforward computation shows that,
for
    any $\gamma$ in $SL_{2}(\Z),$
    $$\gamma\cdot C(q)=C(\gamma\cdot q).$$

\begin{center}                                           %
\begin{figure}[h]                                            %
\begin{center}                                               %
   
\begin{picture}(0,0)%
\includegraphics{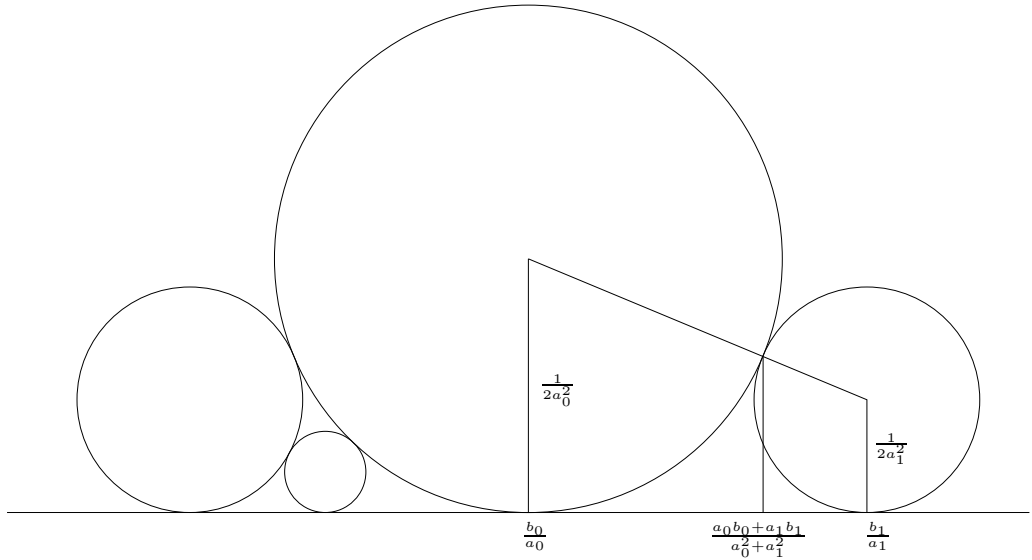}%
\end{picture}%
\setlength{\unitlength}{1243sp}%
\begingroup\makeatletter\ifx\SetFigFont\undefined%
\gdef\SetFigFont#1#2#3#4#5{%
  \reset@font\fontsize{#1}{#2pt}%
  \fontfamily{#3}\fontseries{#4}\fontshape{#5}%
  \selectfont}%
\fi\endgroup%
\begin{picture}(20409,10749)(-11,-10264)
\put(13996,-10186){\makebox(0,0)[lb]{\smash{{\SetFigFont{5}{6.0}{\rmdefault}
{\mddefault}{\updefault}{\color[rgb]{0,0,0}$\frac{a_0b_0+a_1b_1}{a_0^2+a_1^2}$}}}}}
\put(10621,-7171){\makebox(0,0)[lb]{\smash{{\SetFigFont{5}{6.0}{\rmdefault}
{\mddefault}{\updefault}{\color[rgb]{0,0,0}$\frac{1}{2a_0^2}$}}}}}
\put(17101,-10186){\makebox(0,0)[lb]{\smash{{\SetFigFont{5}{6.0}{\rmdefault}
{\mddefault}{\updefault}{\color[rgb]{0,0,0}$\frac{b_1}{a_1}$}}}}}
\put(10261,-10186){\makebox(0,0)[lb]{\smash{{\SetFigFont{5}{6.0}{\rmdefault}
{\mddefault}{\updefault}{\color[rgb]{0,0,0}$\frac{b_0}{a_0}$}}}}}
\put(17281,-8341){\makebox(0,0)[lb]{\smash{{\SetFigFont{5}{6.0}{\rmdefault}
{\mddefault}{\updefault}{\color[rgb]{0,0,0}$\frac{1}{2a_1^2}$}}}}}
\end{picture}%
\caption{The tangency point of two adjacent Ford circles}   %
\label{fig:Ford}                                             %
\end{center}                                                 %
\end{figure}                                                 %
\end{center}

    This easily implies that for any two distinct points
    $q_{0}=b_{0}/a_{0}$ and $q_{1}=b_{1}/a_{1}$ in $\P^1(\Q)$ --- with
    $a_{0}$ and $b_{0}$ (resp. $a_{1}$ and $b_{1}$) prime together ---
    the ``circles" $C(q_{0})$ and $C(q_{1})$ are disjoint if
    $|a_{0}b_{1}-a_{1}b_{0}|\neq 1,$ and are externally tangent if
    $|a_{0}b_{1}-a_{1}b_{0}|=1.$ Moreover, when the latter possibility
    arises, the abscissa of their tangency point is
    $$\frac{a_{0}b_{0}+a_{1}b_{1}}{a_{0}^2+a_{1}^2}$$
    (see Figure \ref{fig:Ford}).

    Consequently, for any $P$ in $\P^1(\Q)\setminus\{0, \infty\},$ its
    image $s(P)$ in $\R/\Z$ may be constructed as follows. Let $a_{0}$
    and $a_{1}$ be two integers prime together such that
    $P=(a_{0}:a_{1})$, and let $C(q_{0})$ and $C(q_{1})$ be two
    tangent Ford circles such that $q_{0}<q_{1},$ attached to rational
    points $q_{0}=b_{0}/a_{0}$ and $q_{1}=b_{1}/a_{1}$ with $b_{0}$
    (resp. $b_{1}$) an integer prime to $a_{0}$ (resp. $a_{1}$).
Denote by $\varepsilon$ the sign of $a_{0}a_{1}.$
    Then
    $\varepsilon.s(P)$ is the abscissa of the tangency point of $C(q_{0})$ and
    $C(q_{1})$ in $\R/\Z.$ Observe that this construction makes
    geometrically obvious
    the density of $s(\P^1(\Q))$ in $\R/\Z.$

\subsection{An application to reduction theory}\label{subsec.reduc}
In this section, we want to discuss the relation between
the classical reduction theory and
our results concerning
  extensions of hermitian vector bundles over arithmetic curves and
  their sizes (see for instance \cite{waerden56} and
  \cite{borel69} for classical expositions and references, and
\cite{lagariasetal90} for more recent applications to lattice
geometry ).

In one direction, observe that, prior to Banaszczyk's
contribution --- which is based on properties of ``Gaussian-like"
measures on lattices (\cite{banaszczyk95}, Section1) ---
transference results similar to Theorem \ref{thm.ban}
were classically established by means of reduction
theory (see for instance \cite{lagariasetal90}, Sections 3 and 5).

Conversely, using our basic upper bound on sizes in Theorem
\ref{prop.sizeup}, we are going to establish the following

\begin{theorem}\label{thm.reducArak}
     For any number field $K$ and any positive integer $n$, there
     exists a constant $r(n,K)$ satisfying the following property:
     for every hermitian vector bundle $\ol{E}$ of rank $n$ over
     $\Spec {\mathcal O}_{K},$ there exist hermitian line bundles $\ol{L}_{1},$
     \ldots, $\ol{L}_{n}$, and an isomorphism of ${\mathcal O}_{K}$-modules
     $$\phi: E \stackrel{\sim}{\longrightarrow}
     \bigoplus_{i=1}^{n}L_{i}$$
     such that, for any embedding $\sigma: K\hookrightarrow \C,$
     \begin{equation}
        \log \|\phi_{\sigma}\|_{\ol{E}, \bigoplus_{i=1}^{n}\ol{L}_{i},
         \sigma}^\infty \leq r(n,K)\,\,\,
         \mbox{ and }\,\,\,\log \|\phi_{\sigma}^{-1}\|_{
        \bigoplus_{i=1}^{n}\ol{L}_{i}, \ol{E}, \sigma}^\infty\leq
r(n,K).
         \label{eq:reducArak}
     \end{equation}
\end{theorem}

As above $\|.\|_{\ol{E}, \bigoplus_{i=1}^{n}\ol{L}_{i},
\sigma}^\infty$
        (resp. $\|.\|_{\bigoplus_{i=1}^{n}\ol{L}_{i}, \ol{E}, \sigma
        }^\infty$)
        denotes the operator norm on the space
    $\Hom_{\C}(E_{\sigma},\bigoplus_{i=1}^{n}L_{i,\sigma})$ (resp. on
    $\Hom_{\C}(\bigoplus_{i=1}^{n}L_{i,\sigma},(E_{\sigma})$) deduced
    from the hermitian norms $\|.\|_{\ol{E}, \sigma}$ and
    $\|.\|_{\bigoplus_{i=1}^{n}\ol{L}_{i}, \sigma}$.

When $K=\Q,$ this follows from the classical reduction theory of
quadratic forms of Lagrange, Gau{ß} ($n=2$), Hermite, and
Korkin-Zolotarev ($n$ arbitrary). For arbitrary number fields,
Proposition (\ref{eq:reducArak}) follows from the general theory of
fundamental domains for arithmetic groups (for instance, from
\cite{borel69}, Th\' eor\`eme 13.1, applied to the Weil restriction
from $K$ to $\Q$ of $GL_{n,K}$).

An effective control
of the constants $r(n,K)$ does not seem to follow simply from
results in the literature, except when $K=\Q$ (in which case the best 
available estimates would follow from the results in
\cite{lagariasetal90} concerning bases of euclidean lattices which
are ``Korkin-Zolotarev reduced"). An
interesting feature of the proof below is the explicit values of
$r(n,K)$ it provides for arbitrary number fields.

Actually we will establish a significantly more precise version of
Theorem \ref{thm.reducArak}. 
Before we state it, we need to introduce some preliminary definitions.

Consider an hermitian vector bundle $\ol{E}$ of rank $n$ over $S:=
\Spec \OK$ as above, and a \emph{complete flag}
$$E_\bullet: E_{0}=\{0\} \subset E_{1}\subset \cdots \subset
E_{n-1}\subset
E_{n}=E,$$
namely a filtration by saturated $\OK$-submodules of ranks $\rk
E_{i}=i,$ $0\leq i \leq n.$ To these data, we may attach the vector
bundle of rank $n$ over $S$
\newcommand{\Gr}{{\rm{Gr}}}
$$\Gr E_{\bullet}:=\bigoplus_{1\leq i \leq n} E_{i}/E_{i-1}.$$
We define a
\emph{splitting over $\OK$} of the flag $E_{\bullet}$ as an
isomorphism 
of $\OK$-modules
$$\varphi: E \stackrel{\sim}{\longrightarrow} \Gr
E_{\bullet}$$
such that, for any $k \in \{0, \ldots, n\},$
$$\varphi (E_{k})=\bigoplus_{1\leq i \leq k} (E_{i}/E_{i-1}).$$
For any such splitting $\varphi$  and any $i \in \{1, 
\ldots,n\},$ the inverse image
\begin{equation}\label{ailes}
L_i := \varphi^{-1}((E_i/E_{i-1}))
\end{equation}
is a direct summand of rank $1$ in $E$ such that 
$$E_i= E_{i-1} \oplus L_i.$$
Conversely any family $(L_i)_{1 \leq i \leq n}$ of $\OK$-submodules
of 
$E$ that satisfies these conditions is deduced by the above 
construction from a unique splitting $\varphi$ of $E_{\bullet}$ over 
$\OK.$ Observe also that to any complete flag $E_{\bullet}$ in $E$ is attached by duality a complete flag 
$E_{\bullet}^{\perp}:=((E/E_{n-i})^\lor)_{0\leq i \leq n}$ in $E^\lor$, and that this construction establishes a one-to-one correspondence between complete flags in $E$ and $E^\lor.$

For any hermitian vector bundle $\ol{E}$ as above, we may perform the 
following inductive construction of a complete flag $E_{\bullet}$: we 
let $E_{0}:=\{0\},$ and for any integer $i \in \{0,\ldots,n-1\},$ we
choose $E_{i+1}$ as a (necessarily saturated) $\OK$-submodule of rank
$i+1$ 
in $E$ containing $E_{i}$ that satisfies
$$\degan \ol{E_{i+1}/E_{i}}= \udega \ol{E/E_{i}}.$$
We shall call a complete flag $E_{\bullet}$ obtained by the above
construction a \emph{reduced complete flag} associated to $\ol{E}.$
This terminology is meant to emphasize the analogy with classical
reduction theory \emph{\`a la} Hermite-Korkin-Zolotarev. Actually,
when $S= \Spec \Z,$ if $(b_{1},\ldots,b_{n})$ is a basis of the
euclidean lattice $\ol{E}$ which is ``Korkin-Zolotarev reduced" (see
for instance
\cite{lagariasetal90}, Section 2), then the complete flag
$(\bigoplus_{1\leq i \leq k}\Z{b_{i}})_{0\leq k \leq n}$ of $E$ is
reduced
in our sense.

We may now state a refined version of Theorem \ref{thm.reducArak}.

\begin{theorem}\label{thm.reducArakbis}
    Let $K$ be a number field, and let 
    $(c(n,K))_{n>0}$     be a sequences in $\R_{+}$ that satisfy, for
any $n>0,$
    \begin{equation}\label{definesc}
   c(n,K) \geq \frac{1}{2n} \sum_{1 \leq i \leq n-1} \log (1 +
|\Delta_{K}|^{4/[K:\Q]} (n-1)^2)
    \end{equation}

    For any hermitian vector bundle $\ol{E}$ of positive rank $n$ over
    $\Spec \OK$ and any complete flag $E_{\bullet}$ in $E$ such that the
     dual flag $E_\bullet^\perp$ is a reduced flag of  $\ol{E}^\lor,$  there exists 
  a splitting over $\OK$ of  $E_{\bullet}$ 
such that the associated $\OK$-submodules $(L_i)_{1\leq i \leq n}$ 
satisfy
\begin{equation} \label{delta}
\delta (\ol{E}; L_1, \ldots, L_n) := \mua (\ol{E}) - \frac{1}{n} 
\sum_{i=1}^n \degan \ol{L}_i \leq  c(n,K).
\end{equation}

Moreover, the archimedean operator norms of the tautological ``sum
map"
$$\Sigma : \bigoplus_{1 \leq i \leq n} L_i
\stackrel{\sim}{\longrightarrow} 
E$$
and of its inverse,
computed with respect to the hermitian structures of $\bigoplus_{1
\leq 
i \leq n} \ol{L}_i$ and $\ol{E},$
satisfy the following bounds:
\begin{equation}\label{obvious}
 \|\Sigma\|_\sigma \leq \sqrt{n} \quad \mbox{ for any embedding }\sigma:K\hookrightarrow \C,
 \end{equation}
and
\begin{equation}\label{lessobvious}
[K:\Q]^{-1} \sum_{\sigma: k \hookrightarrow \C}\log  
\|\Sigma^{-1}\|_\sigma  \leq \frac{n-1}{2} \log n + n. c(n,K).
\end{equation}

\end{theorem}

Condition (\ref{definesc}) is satisfied 
by 
$$ c(n,K):=  2 \frac{\log
     |\Delta_{K}|}{[K:\Q]} + \frac{1}{n}\log n!,$$
 and consequently also by
 \begin{equation}\label{simplevalue}
 c(n,K):=  2 \frac{\log
     |\Delta_{K}|}{[K:\Q]} + \log n.
     \end{equation}

Observe that, for any family $(L_i)_{1\leq i \leq n}$ of 
$\OK$-submodules of rank $1$ in $E$ such that $E_K$ is the direct 
sum $\bigoplus_{1 \leq i \leq n} L_{i,K},$ we may attach the real 
number
$$\delta (\ol{E}; L_1, \ldots, L_n) := \mua (\ol{E}) - \frac{1}{n} 
\sum_{i=1}^n \degan \ol{L}_i.$$
It is easily checked to be non-negative and to vanish iff the 
tautological ``sum map" defines an isomorphism of hermitian vector
bundles 
from $\bigoplus_{1 \leq i \leq n} \ol{L}_i$ to $\ol{E}.$

When $K=\Q$ and $(b_1, \ldots, b_n)$ is a $\Z$-basis of the euclidean 
lattice $\ol{E},$ then $\delta (\ol{E}; \Z b_1, \ldots, \Z b_n)$ is 
the logarithm of the ``orthogonality defect" $\prod_{1\leq i \leq n}
\left\| 
b_i \right\|/{\rm{covol}}(\ol{E})$ of the basis $(b_1, \ldots, b_n)$ 
(compare \cite{lagariasetal90}, p.336). 
Lagarias and his coworkers establish in \emph{loc. cit.} that the
submodules $L_{i}:= \Z.b_{i}$ generated by the vectors of a
Korkin-Zolotarev reduced base $(b_{1},\ldots,b_{n})$ of $\ol{E}$ satisfy the
upper bound (\ref{delta}) with, in place of $c(n,\Q)$, a function of $n$
which, like $c(n,\Q)$,   grows 
like $\log n + O(1)$ when $n$ goes to infinity.   

The existence of a splitting of $E_\bullet$ such that the associated 
family $(L_i)_{1 \leq i \leq n}$ satisfies (\ref{delta})  follows by
induction on
$n$ from Corollary \ref{sizedual}, Proposition \ref{sizedirectsummand}, and the following lemma 
applied to the dual of $\ol{E}$ and the submodule $(E/E_{n-1})^\lor$ of rank $1$ in $E^\lor$:


\begin{lemma}\label{lem.reducArak} Let $\ol{E}$
     be an hermitian vector bundle  of rank $n>1,$ and $L$ a saturated
     ${\mathcal O}_{K}$-submodule of rank 1 in $E$ such that
     $\degan \ol{L}$ $(=\mua (\ol{L}))$ equals $\udega \ol{E}.$

     We have:
     \begin{equation}
     \udega \ol{E/L} - \degan \ol{L} \leq \log 2 + \frac{\log
	  |\Delta_{K}|}{[K:\Q]}.
	  \label{eq:reducupper}
     \end{equation}
     Moreover the size of the admissible extension
     \begin{equation}
     {\ol{\mathcal E}}:0\longrightarrow  \ol{L} \longrightarrow
     \ol{E} \longrightarrow \ol{E/L} \longrightarrow  0
     \label{eq:admissibleL}
\end{equation}
satisfies
\begin{equation}
     \size({\ol{\mathcal E}}) \leq 2 \frac{\log
     |\Delta_{K}|}{[K:\Q]} + \log(n-1).
     \label{eq:sizeL}
\end{equation}

\end{lemma}

\proof  
Consider a saturated $\OK$-submodule $M$ of rank $1$ in $E/L$ such
that
\begin{equation}\label{MEL}
\degan \ol{M} =\udegan \ol{E/L},
\end{equation}
and let $p$ be its inverse image in $E$. It a saturated
$\OK$-submodule 
of rank 2 in $E$, which contains $L$ and defines an admissible 
extension over $\Spec \OK$:
$$0\longrightarrow  \ol{L} \longrightarrow
     \ol{P} \longrightarrow \ol{M} \longrightarrow  0.$$
In particular,
\begin{equation}\label{PLM}
\degan \ol{P} =\degan \ol{L} + \degan \ol{M}.
\end{equation}
Moreover, by the very definition of $L,$ the inequalities 
$$\degan \ol{L} \leq \udegan \ol{P} \leq \udegan \ol{E}$$
are equalities. Together with ((\ref{mumaxudeg2}),  this shows that:
\begin{equation}\label{LP}
\degan \ol{L} =\udegan \ol{P} \geq \widehat{\mu}(\ol{P}) - 
\frac{1}{2} \log 2  - \frac{\log |\Delta_{K}|}{2[K:\Q]}.
\end{equation}
Inequality  (\ref{eq:reducupper}) follows from 
(\ref{MEL}), (\ref{PLM}), and  (\ref{LP}).

According to the upper bound (\ref{eq:sizeup}) in
Theorem \ref{prop.sizeup} on the size of ${\ol{\mathcal E}},$ we
have:
\begin{equation}
  \size({\ol{\mathcal E}})  \leq   \udegan(\ol{E/L}\otimes
\ol{L}^\lor) 
  + \frac{\log |\Delta_{K}|}{[K:\Q]} + \log \bigl((n-1)/2\bigr).
      \label{eq:L01}
\end{equation}
Inequality (\ref{eq:admissibleL}) follows from (\ref{eq:L01}), 
((\ref{eq:reducupper}), and the equality
$$\udegan(\ol{E/L}\otimes \ol{L}^\lor) =
\udegan(\ol{E/L}) -\degan \ol{L},$$
which is a straightforward consequence of (\ref{eq:udegatens}).

\qed

The upper bound (\ref{obvious}) on the archimedean norms of $\Sigma$ 
is obvious. 
The upper bounds  (\ref{lessobvious}) on its inverse map 
then follow from the next observation applied to $\Sigma.$

\begin{lemma}
For any two  hermitian vector bundles $\ol{E}_1$ and $\ol{E}_2$ of 
positive rank $n$ over $\Spec \OK$ and any isomorphism of 
$\OK$-modules $\psi: E_1 
\stackrel{\sim}{\rightarrow} E_2,$ we have:
\begin{multline}
 [K:\Q]^{-1} \sum_{\sigma: K \hookrightarrow \C} \log \|
\psi^{-1}\|_{\ol{E}_2, 
 \ol{E}_1, \sigma} \\ \qquad \qquad \qquad
 =  [K:\Q]^{-1} \sum_{\sigma: K \hookrightarrow 
 \C} (\log \|  \wedge^{n-1}\psi\|_{\wedge^{n-1}\ol{E}_1, 
 \wedge^{n-1}\ol{E}_2, \sigma} 
 - \log \|  \wedge^{n}\psi\|_{\wedge^{n}\ol{E}_1, 
 \wedge^{n}\ol{E}_2, \sigma} )
 \\
 \leq 
(n-1) [K:\Q]^{-1} \sum_{\sigma: K \hookrightarrow \C} \log \|
\psi\|_{\ol{E}_1, 
 \ol{E}_2, \sigma}
 + n (\mua (\ol{E}_2) - \mua(\ol{E}_1)). \qquad \qquad \enskip \;
\end{multline}
\end{lemma}

Theorem \ref{thm.reducArakbis} --- or rather a variant, involving
different constants --- may also be deduced by induction from
Lemma \ref{lem.reducArak} by means of the bounds (\ref{b1}) and
(\ref{b2}) on the norms of trivializations of admissible extensions.
The above proof, which emphasizes the role of the invariant
$\delta(\ol{E}; L_{1},\ldots,L_{n})$ attached to a linearly
independent $n$-tuple of rank one subbundles, is computationally
simpler.

  \section{Invariance of the size under base change
and Voronoi cells of euclidean lattices}\label{sect:inv}

     An intriguing issue is the invariance property of the size under
     base changes $\Spec {\mathcal O}_{K'} \to \Spec {\mathcal O}_{K}$ associated
     to extensions of number fields $K\hookrightarrow K'.$ It seems
     plausible that, defined with the precise normalization we
     introduce in paragraph \ref{subsec:size} below, the size is
     invariant by any such base change, at least when $K=\Q$.
      In this section, we
     establish various results which support this expectation.

\subsection{A geometric consideration}\label{subsec.geoarg}

Let $C$ and $C'$ be two projective curves, smooth and geometrically
connected over some field $k,$ and let $f:C'\rightarrow C$ be a finite
(or equivalently dominant) $k$-morphism.

\begin{proposition}\label{prop.nonsplitgeo}
When the degree $\deg f$ of $f$ is prime to the characteristic
exponent
of $k$ \emph{(for instance, when $k$ is a field of characteristic
zero)},
then, for any vector bundle $E$ over $C,$ the $k$-linear map
\begin{equation}
     f^\ast:H^1(C,E) \longrightarrow H^1(C',f^\ast E)
     \label{eq:FF1}
\end{equation}
is injective.
\end{proposition}

\proof Since the  morphism $f$ is affine and $E$ is
locally free, we have a
canonical isomorphism
$$ H^1(C',f^\ast E) \simeq H^1(C,E \otimes f_\ast{\mathcal O}_{C'})$$
and the map (\ref{eq:FF1}) may be identified with the one deduced from
the canonical morphism of ${\mathcal O}_{C}$-modules
\begin{equation}
     {\mathcal O}_{C} \longrightarrow f_{\ast}{\mathcal O}_{C'}
     \label{eq:FF2}
\end{equation}
by applying
the functors $E\otimes .$ and $H^1(E,.).$ Under the above hypothesis,
the morphism (\ref{eq:FF2})
is split by the morphism of ${\mathcal O}_{C}$-modules:
\begin{equation}
   \frac{1}{\deg f}{\rm Tr}_{f}:  f_{\ast}{\mathcal O}_{C'}
\longrightarrow {\mathcal O}_{C},
     \label{eq:FF3}
     \end{equation}
    where  ${\rm Tr}_{f}$ denotes the ``naive" 
trace morphism\footnote{For any $U$ open in $C$,
the map
    ${\rm Tr}_{f\mid U}:f_{\ast}{\mathcal O}_{C'}
(U)={{\mathcal O}}_{C'}(f^{-1}(U)) \rightarrow {\mathcal O}_{C}(U)$
     is defined as the trace on  the ${\mathcal O}_{C}(U)$-algebra
${{\mathcal O}}_{C'}(f^{-1}(U)),$ which  is  a finitely generated projective
${\mathcal O}_{C}(U)$-module since $f$ is finite and flat.}  from
$f_\ast{\mathcal O}_{C'}$ to ${\mathcal O}_{C}.$

  Applied  to (\ref{eq:FF3}), the functors $E\otimes .$ and $H^1(C,.)$
produce
  a splitting of (\ref{eq:FF1}).
  \qed

  {\small However, when the characteristic $p$ of $k$ is positive and
  divides $\deg f$, the pullback map (\ref{eq:FF1}) may not be
  injective. Consider for instance any smooth projective,
geometrically
  connected curve $C$ over $k$ and the (relative) Frobenius morphism
  $$F:C \longrightarrow C^{(p)}.$$
  Then
  $$ F^\ast:H^1(C^{(p)}, {\mathcal O}_{C^{(p)}}) \longrightarrow
  H^1(C,{\mathcal O}_{C})$$
  is injective iff the Hasse-Witt matrix of $C$ is invertible. (This
  does not hold for instance when $C$ is a supersingular elliptic
  curve.)

  Actually, when $p$ divides $\deg f$, (\ref{eq:FF1}) may not be
  injective even when $f$ is separable, as demonstrated by the
  following observation, applied to an ordinary elliptic curve and its
  quotient by an \'etale subgroup of order $p$:

  \begin{lemma}
      If $E$ and $E'$ are two elliptic curves over a field $k$ of
      positive characteristic $p$ and if $f:E' \rightarrow E$ is a
      separable $k$-isogeny of degree divisible by $p$, then the
      pullback map:
      $$f^\ast:H^1(E, {\mathcal O}_{E}) \longrightarrow
H^1(E',{\mathcal
      O}_{E'})$$
      vanishes.
  \end{lemma}

\proof Consider the ``trace -map"
$$t_f:H^1(E',{\mathcal O}_{E'}) \simeq H^1(E,f_{\ast}{\mathcal
O}_{E'}) \longrightarrow H^1(E, {\mathcal
O}_{E})$$
induced by the morphism of ${\mathcal O}_E$-modules
$$
  {\rm Tr}_{f}:  f_{\ast}{\mathcal O}_{E'} \longrightarrow {\mathcal
O}_{E}.
$$
Using  Serre duality, it may be identified with the
transpose of $f$  acting by pullback on  regular $1$-forms:
$$f^\ast_{\Omega^1}: \Omega^1(E) \longrightarrow \Omega^1(E').$$
Since $f$ is separable, $f^\ast_{\Omega^1}$ and  consequently $t_f$
are
isomorphisms. Besides, we have:
$$t_f\circ f^\ast= \deg f. Id_{H^1(E,{\mathcal O}_E)}.$$
Since $p$ divides $\deg f$, this vanishes, and consequently
$f^\ast$ vanishes too.
\qed
  }

Consider now $D$ a non-zero effective divisor on $C,$ $D':=f^\ast(D)$
its inverse image in $C',$
and the affine curves $\Co:=C\setminus |D|$ and
${\Co}':=C'\setminus |D'|.$ Let also $F$ and $G$ be two vector bundles
on $C$, $\Fo$ and $\Go$ their restrictions to $\Co$,  $F'$ and $G'$
their pullback to $C'$, and ${\Fo}'$ and ${\Go}'$ the restrictions
of the
latter to ${\Co}',$ or equivalently, the pullback of
$\Fo$ and $\Go$ by $f_{\mid {\Co}'}$. With the
notation of Section \ref{geomextII}, pulling back extensions by $f$ defines a
natural $k$-linear map:
$$f^\ast : \Exthat^1_{\Co}(\Fo,\Go) \longrightarrow
\Exthat^1_{{\Co}'}({\Fo}',{\Go}').$$
By applying Proposition \ref{prop.nonsplitgeo} to the vector bundles
$E:= \check{F}\otimes G(-N.D),$ $N\in \Z,$ we obtain:

\begin{corollary}\label{sizeinvgeo}
     When the degree of $f$ is prime to the characteristic exponent of
     $k,$ then, for any $e\in \Exthat^1_{\Co}(\Fo,\Go),$
     $$\size_{F',G'}(f^\ast(e))=\size_{F,G}(e).$$
\end{corollary}

\subsection{Size and base change}\label{subsec.sizechange}
We give two equivalent formulation of the problem of the invariance of
size under base change.
\subsubsection{The condition
\condp{K}{K'}{\ol{E}}}\label{subsubsec.PK}
As before, consider a number field $K$,  a number field $K'$
containing $K$, and
   \[
   g:S'=\Spec\,{\mathcal O}_{K'}\longrightarrow  S=\Spec\,{\mathcal O}_K
   \]
   the associated morphism of arithmetic curves. Let us also denote
   by $\pi$ (resp. $\pi'$) the morphism
   from $\Spec {\mathcal O}_{K}$ (resp. $\Spec {\mathcal O}_{K'}$) to $\Spec \Z$.

   \noindent \textbf{Problem 1}. \emph{Let $\ol{F}$ and $\ol{G}$ be
   two hermitian vector bundles over $S$. Is it true that, for any
   extension class $e$ in $\Exthat^1_{S}\, (F,G),$ the inequality
   $$ \size_{g^{\ast}\ol{F}, g^{\ast}\ol{G}} \, g^{\ast}(e) \leq
   \size_{\ol{F}, \ol{G}}\, (e)$$
   is indeed an equality ?}

   Let $\ol{E}$ be an hermitian vector bundle over $S$. The
map of ``extension of scalars" from
${\mathcal O}_{K}$ to ${\mathcal O}_{K'}$
$$
\begin{array}{rcl}
     E & \longrightarrow & E\otimes_{{\mathcal O}_{K}}{\mathcal O}_{K'}  \\
     v & \longmapsto & v\otimes 1,
\end{array}$$
seen as a $\Z$-linear map, defines a morphism of $S$-vector bundles
$$\Delta:\pi_{\ast}E \longrightarrow \pi'_{\ast}g^\ast E.$$
This is also the morphism deduced from the natural morphism
$$\delta: E \longrightarrow g_{\ast}g^{\ast}E$$
by taking its direct image by $\pi$.

The linear maps
$$\delta_{\sigma}: E_{\sigma} \longrightarrow
(g_{\ast}g^{\ast}E)_{\sigma}$$
(where $\sigma$ denote an embedding of $K$ in $\C$) and
$$\Delta=\Delta_{\R}:(\pi_{\ast}E)_{\R}\simeq E\otimes_{\Z}\R
\longrightarrow (\pi'_{\ast}g^\ast  E)_{\R}
\simeq E\otimes_{{\mathcal O}_{K}}{\mathcal O}_{K'}\otimes_{\Z}\R$$
are isometric
up to a factor $[K':K]^{1/2}$ when we equip these vector spaces with
the norms defining respectively the hermitian structures of $\ol{E},$
$g_{\ast}g^{\ast}\ol{E},$ $\pi_{\ast}\ol{E},$ and
$\pi'_{\ast}g^\ast\ol{E}.$ Namely, for any $v$ (resp. $w$) in
$E_{\sigma}$ (resp. $(\pi_{\ast}E)_{\R}$), we have:
\begin{equation}
     \|\delta_{\sigma}(v)\|_{g_{\ast}g^{\ast}\ol{E},\sigma}=
[K':K]^{1/2}
     \|v\|_{\ol{E}, \sigma} \;\;\;(\mbox{resp. }\;\;
     \|\Delta(w)\|_{\pi'_{\ast}g^\ast\ol{E}}= [K':K]^{1/2}
     \|w\|_{\pi_{\ast}\ol{E}}.)
     \label{eq:isom}
\end{equation}

\noindent \textbf{Problem 2}. \emph{With the notation above, is it
true that, for any $w\in (\pi_{\ast}E)_{\R},$ the inequality
\begin{equation}
\min_{\beta \in g^\ast E }\|\Delta(w)-\beta
\|_{\pi'_{\ast}g^\ast\ol{E}} \leq
\min_{\alpha \in E }
\|\Delta(w)-\Delta(\alpha)\|_{\pi'_{\ast}g^\ast\ol{E}}
= [K':K]^{1/2}
\min_{\alpha \in E}       \|w-\alpha\|_{\pi_{\ast}\ol{E}}
     \label{eq:P2}
\end{equation}
is actually an equality ?}

These two problems are equivalent. More precisely, Problem 1 for some
pair $(\ol{F}, \ol{G})$ of hermitian vector bundles is equivalent
to Problem 2 for $\ol{E}=\ol{F}^\lor \otimes \ol{G}.$

We shall say that \emph{Condition \condp{K}{K'}{\ol{E}} holds}
when Problem 2 has an affirmative answer.

Observe that, as a straightforward consequence of definitions, this
condition satisfies the following compatibility with the operation of
direct sum:

\begin{lemma}\label{lem.Psum}
     For any finite family $(\ol{E}_{i})_{i \in I}$ of hermitian
     vector bundles over $\Spec {\mathcal O}_{K}$, the condition
     \condp{K}{K'}{\bigoplus_{i\in I}\ol{E}_{i}} holds
    iff  the condition \condp{K}{K'}{\ol{E}_{i}} holds for any $i\in
I$.
     \end{lemma}

     The following property also is immediate:

     \begin{lemma}\label{lem.Phomothety}
         Let
$\ol{E}:=(E,(\|.\|_{\sigma})_{\sigma:K\hookrightarrow\C}$ be an
         hermitian vector bundle over $\Spec {\mathcal O}_{K}$, and $\lambda$ a
positive
         real number, and let us consider
$\ol{E}':=(E,(\lambda\|.\|_{\sigma})_{\sigma:K\hookrightarrow\C}.$
         Then  \condp{K}{K'}{\ol{E}} holds  iff
\condp{K}{K'}{\ol{E}'} holds.
     \end{lemma}

\begin{lemma}\label{primediscr}
Let $L/K$ and $K'/K$ be two extensions of number 
fields which are disjoint (i.e.
$L'=L\otimes_K K'$ is again a number field) and 
have coprime discriminant ideals.
Let $h$ denote the natural map from Spec ${\mathcal O}_L$ to Spec ${\mathcal O}_K$.
Let $\ol{E}$ be an hermitian vector bundle over Spec ${\mathcal O}_L$.
Then \condp{K}{K'}{h_*\ol{E}} implies \condp{L}{L'}{\ol{E}}.
\end{lemma}

\proof
The condition on discriminants assures that the middle square in
the commutative diagram
\[
\begin{array}{ccccccc}
\Spec \Z & \stackrel{\tau'}{\leftarrow } & \Spec{\mathcal O}_{L'}
& \stackrel{g'}{\longrightarrow } & \Spec {\mathcal O}_L
& \stackrel{\tau}{\longrightarrow }& \Spec \Z \\
\| & &\downarrow {\scriptstyle h'} &
&\downarrow {\scriptstyle h} & &\|\\
\Spec \Z & \stackrel{\pi'}{\leftarrow } & \Spec{\mathcal O}_{K'}
& \stackrel{g}{\longrightarrow } & \Spec {\mathcal O}_K
& \stackrel{\pi}{\longrightarrow }& \Spec \Z \\
\end{array}
\]
is cartesian (see for example \cite[I 2.11]{neukirch99}).
The base change isomorphism
$g^*h_*\ol{E}\stackrel{\sim}{\to}h'_*g'^*\ol{E}$ is an isometry
by \ref{arithm.base.change}.
For each element $w\in (\tau_*\ol{E})_\R=(\pi_*h_*\ol{E})_\R$, we get
\[
\min_{\beta \in g'^* E } \|\Delta(w) -\beta \|_{\tau'_*g'^*\ol{E}}
=\min_{\beta \in g^*h_* E } \|\Delta(w) -\beta
\|_{\pi'_*g^*h_*\ol{E}}.
\]
The condition \condp{K}{K'}{h_*\ol{E}} implies that the last term
equals
\[
[K':K]^{1/2}
\min_{\alpha \in E}       \|w-\alpha\|_{\pi_{\ast}h_*\ol{E}}
= [L':L]^{1/2}
\min_{\alpha \in E}       \|w-\alpha\|_{\tau_{\ast}\ol{E}}.
\]
Hence \condp{L}{L'}{\ol{E}} holds.
\qed

\subsection{The condition \condp{K}{K'}{\ol{E}} and Voronoi cells}
Let us now turn to a reformulation of condition \condp{K}{K'}{\ol{E}}
  involving the geometry of  Voronoi
cells of  euclidean lattices.

     Recall that the \emph{Voronoi cell} ${\mathcal V}(\ol{F})$ of an
     euclidean lattice $\ol{F}:=(F, \|.\|)$ consists of those points
     of $F_{\R}$ that are at least as close to the origin as to  any
     element $e$ in $F$:
  $${\mathcal V}(\ol{F}):= \{x \in F_{\R} \mid \|x\|   \leq
  \|x-e\| \mbox{ for all } e \in F\}.$$
  If $\langle .,.
\rangle $ denotes the euclidean scalar product on $F_{\R}$
  associated to $\|.\|$, the condition
  \begin{equation}
    \|x\|   \leq
  \|x-e\|
      \label{eq:balls}
  \end{equation}
  may also be expressed as
  \begin{equation}
      \langle x-\frac{e}{2},e
\rangle  \leq 0.
      \label{eq:halfspace}
  \end{equation}
  This shows that the Voronoi cell ${\mathcal V}(\ol{F})$ is a
compact,
  convex, symmetric neighborhood of the origin. Actually, it is
  defined by a finite number of the conditions (\ref{eq:halfspace}).
  In other words, it is a polytope, and consequently posses a finite
  number of extremal points, its \emph{vertices}. Moreover, we have:
  \begin{equation}
      F_{\R}=\bigcup_{e \in F}(e + {\mathcal V}(\ol{F})).
      \label{eq:Vtiling}
  \end{equation}

  \begin{lemma}\label{lem.CNSVor}
      With the above notation, condition \condp{K}{K'}{\ol{E}} holds
      iff, for any vertex $P$ of ${\mathcal V}(\pi_{\ast}\ol{E})$ and
      any $\beta \in g^\ast E$, the following equivalent conditions
      holds:
           \begin{equation}
           \|\Delta(P) -\beta \|_{\pi'_{\ast}g^\ast\ol{E}} \geq
           \|\Delta(P)\|_{\pi'_{\ast}g^\ast\ol{E}} \,\,(=
           [K':K]^{1/2} \|P\|_{\pi_{\ast}\ol{E}});
          \label{eq:CNSVor1}
      \end{equation}
\begin{equation}
           2\langle \Delta(P), \beta
\rangle _{\pi'_{\ast}g^\ast\ol{E}}
          \,\,\leq\,  \|\beta \|^2_{\pi'_{\ast}g^\ast\ol{E}}.
          \label{eq:CNSVor2}
      \end{equation}
\end{lemma}

\proof For any $P$ in $(\pi_{\ast}E)_{\R}\, (=E\otimes \R)$, the
inequalities
(\ref{eq:CNSVor1}) and (\ref{eq:CNSVor2}) are equivalent (this is the
special case of the equivalence of (\ref{eq:balls}) and
(\ref{eq:halfspace}) when $\|.\| = \|.\|_{\pi'_{\ast}g^\ast\ol{E}},$
$v=\Delta(P),$ and $e=\beta$).

Since $(\pi_{\ast}E)_{\R}$ is the union of the translates of
${\mathcal
V}(\pi_{\ast}\ol{E})$ by elements of $\pi_{\ast}E\, (=E),$ the
validity of (\ref{eq:P2}) for all $v \in (\pi_{\ast}E)_{\R}$ is
equivalent to its validity when $v$ is
a point $P$ of ${\mathcal V}(\pi_{\ast}\ol{E})$.
Consequently, \condp{K}{K'}{\ol{E}} holds iff
(\ref{eq:CNSVor1}), or (\ref{eq:CNSVor2}), holds
for any $P$ in ${\mathcal V}(\pi_{\ast}\ol{E})$ and any $\beta$ in
$g^\ast E$. Observe finally that, for any given $\beta$,
(\ref{eq:CNSVor2}) holds for all $P$ in the
polytope ${\mathcal V}(\pi_{\ast}\ol{E})$
iff it holds for its extremal points.
\qed

Observe that (\ref{eq:CNSVor2}) is fulfilled as soon as
$$
\|\beta \|_{\pi'_{\ast}g^\ast\ol{E}} \geq
           2\|\Delta(P)\|_{\pi'_{\ast}g^\ast\ol{E}} \,\,(=
           2[K':K]^{1/2} \|P\|_{\pi_{\ast}\ol{E}}).$$
Consequently, Lemma \ref{lem.CNSVor} reduces the proof of
\condp{K}{K'}{\ol{E}} to checking the \emph{finite set} of
inequalities
(\ref{eq:CNSVor2}) for $P$ a vertex of ${\mathcal
V}(\pi_{\ast}\ol{E})$ and $\beta$ a point in $g^\ast E$ such that
$\|\beta
\|_{\pi'_{\ast}g^\ast\ol{E}}$ is at most
$$2[K':K]^{1/2} \max_{P\in{\mathcal V}(\pi_{\ast}\ol{E})}
\|P\|_{\pi_{\ast}\ol{E}}=
2[K':K]^{1/2} \rho (\pi_{\ast}\ol{E}).$$
This shows that there are algorithms that would theoretically enable
one
to check the validity of \condp{K}{K'}{\ol{E}} for any
``explicitly given" $K$, $K',$ and $\ol{E}$. However, one should be
aware that determining the vertices of the Voronoi cell of a euclidean
lattice is a ``hard" problem.

\subsection{Base change from euclidean
lattices}\label{subsubsec.bcheuclidean}
  Let $L$ be a number field, and $\ol{E}$ an hermitian vector bundle
  over $\Spec \Z$.  Lemma \ref{lem.CNSVor} may be used to
  reformulate the  condition \condp{K}{K'}{\ol{E}}
  for
  $K=\Q$ and $K'=L$.

  Indeed, for any embedding $\sigma:L\hookrightarrow\C$, we may define
  $$
   \begin{array}{cccc}
       \sigma_{E}:=id_{E}\otimes \sigma: & g^\ast
E=E\otimes_{\Z}{\mathcal O}_{L}  &
       \longrightarrow & E_{\C}=E\otimes_{\Z}\C   \\
        & e\otimes\lambda & \longmapsto & e\otimes \sigma(\lambda).
   \end{array}
   $$
   Then the morphism of $\Z$-modules
   \[
   (\sigma_{E})_{\sigma:L\hookrightarrow\C }: g^\ast E\longmapsto
   E_{\C}^{(\Spec L)(\C)}=\oplus_{\sigma:L\hookrightarrow \C}E_\C
   \]
   extends to an isomorphism of $\C$-vector spaces:
   $$g^\ast E \stackrel{\sim}{\longrightarrow} E_{\C}^{(\Spec
L)(\C)},$$
and
condition
(\ref{eq:CNSVor2}) may be written
\begin{equation}
     2 \sum_{\sigma:L\hookrightarrow\C } {\rm
     Re}\langle P,\sigma_{E}(\beta)
\rangle _{\ol{E}}\, \leq\,
     \sum_{\sigma:L\hookrightarrow\C }
     \|\sigma_{E}(\beta)\|^2_{\ol{E}}.
     \label{eq:CNSVorQ}
\end{equation}

We shall establish that this holds for any vertex $P$ of ${\mathcal
V}(\ol{E})$ and any $\beta$ in $E\otimes_{\Z}{\mathcal O}_{L}$ --- in other
words that condition \condp{\Q}{L}{\ol{E}} holds --- for various
special lattices.

Our proof will rely on the following observation:
\begin{lemma}\label{lem.sumnorm}
     For any number field $L$ and any element $\alpha$ of its ring of
     integers ${\mathcal O}_{L}$, we have:
     $$  \sum_{\sigma:L\hookrightarrow\C } |\sigma(\alpha)|^2
     - \sum_{\sigma:L\hookrightarrow\C } {\rm Re} \, \sigma(\alpha)
\geq
     \sum_{\sigma:L\hookrightarrow\C } |\sigma(\alpha)|^2
     - \sum_{\sigma:L\hookrightarrow\C } |\sigma(\alpha)| \geq 0.$$
     \end{lemma}

     \proof We may assume that $\alpha \neq 0.$ Then
     $$\prod_{\sigma:L\hookrightarrow\C } |\sigma(\alpha)| =
     |N_{L/\Q}(\alpha)| \geq 1,$$
     and therefore
     \begin{equation}
       \frac{1}{[L:\Q]} \sum_{\sigma:L\hookrightarrow\C }
       |\sigma(\alpha)| \geq 1.
         \label{eq:sumnorm1}
     \end{equation}
Moreover, by Cauchy-Schwarz inequality,
     \begin{equation}
     \left(\sum_{\sigma:L\hookrightarrow\C} 
|\sigma(\alpha)|\right)^2 \leq [K:\Q]
     \sum_{\sigma:L\hookrightarrow\C} |\sigma(\alpha)|^2.
         \label{eq:sumnorm2}
     \end{equation}
     The second inequality in Lemma \ref{lem.sumnorm} follows from
     (\ref{eq:sumnorm1}) and
     (\ref{eq:sumnorm2}). The first one is obvious.

     \qed

     \begin{proposition}\label{prop.dirsum}
         For any hermitian  vector bundle $\ol{E}$  over
         $\Spec \Z$ which is a direct sum of
         hermitian line bundles, condition
         \condp{\Q}{L}{\ol{E}} holds for any number field $L$.
     \end{proposition}
     \proof  Proposition \ref{lem.Psum} shows that we may suppose that
     $\ol{E}$ has rank one.

     Let $e_{0}$ denote a generator of the rank one free
     $\Z$-module $E$. Then the Voronoi cell of $\ol{E}$ is the
     "interval"
     $${\mathcal
V}(\ol{E})= [-1/2,1/2]\cdot e_{0},$$
and its vertices $P$ are the points $e_{0}/2$ and $-e_{0}/2$. The
elements $\beta$ of $E\otimes{\mathcal O}_{L}$ may be written $e_{0}\otimes
\alpha$
with $\alpha \in {\mathcal O}_{L}$,
and condition (\ref{eq:CNSVorQ}) reduces to the inequality in Lemma
\ref{lem.sumnorm} applied to $\alpha$ and $-\alpha.$
\qed

The last proposition leads to ask the following:

\begin{question}\label{invQ?} Is it true that
condition \condp{\Q}{L}{\ol{E}} holds for
any hermitian vector bundle
$\ol{E}$ over $\Spec \Z$ and any number field $L$ ?
     \end{question}

     A positive answer  would imply that the size of any
     admissible extension over $\Spec \Z$ is invariant by the base
     change from $\Z$ to ${\mathcal O}_{L}$ for any number field $L$.

    In the next sections, we shall prove
the following theorem which, together with Proposition
     \ref{prop.dirsum}, points towards a positive answer to question
     \ref{invQ?}.

\begin{theorem}\label{theorem.bc}
For any hermitian vector bundle $\ol{E}$ over Spec $\Z$ and
any number field $L$, the condition $P(L/\Q,\ol{E})$ holds
in the following cases:
\begin{enumerate}
\item[i)]
If $L/\Q$ is an abelian extension.
\item[ii)]
If $\ol{E}$ is a root lattice.
\item[iii)]
If $\ol{E}$ is a lattice of Voronoi's first kind.
\end{enumerate}
\end{theorem}

\proof
Case i) is treated in Section \ref{subsec.cyclo}.
Recall that an integral euclidean lattice 
$\ol{E}$ over Spec $\Z$ is called a {\it root 
lattice} if $E$ is generated by its subset 
$\{e\in E\,|\, \|e\|^2\in \{1,2\} \}$.
The proof of ii) is given in Section \ref{subsec.rootlattice}.
The definition of a lattice of Voronoi's first 
kind is recalled in the appendix.
The proof of iii) can be found in Section \ref{proof.iii}.
\qed

\subsubsection{An application of reduction
theory}\label{subsec.applred}

     Observe that, by combining Proposition \ref{prop.dirsum} and
     our ``reduction"  Theorem \ref{thm.reducArak}, we obtain that
     Question \ref{invQ?} has a positive answer up to an additive
error
     term bounded in terms of the
     ranks of the considered hermitian vector bundles:

     \begin{proposition}\label{prop:almostinvariant}
    There exist non-negative real numbers $s(n),$
       $n\in \N,$ which satisfy the following property.

       Let
       $\ol{E}$  be an hermitian vector bundle over $\Spec \Z$, and
       $n$ its rank. Let $L $ be a number field, and $g$ the unique
morphism
       $\Spec {\mathcal O}_L \to \Spec \Z,.$

       For any $w\in E_{\R},$
       the following inequalities hold
       \begin{equation}
e^{-s(n)}       [L:\Q]^{1/2}
\min_{\alpha \in E}       \|w-\alpha\|_{\ol{E}}  \leq
  \min_{\beta \in g^\ast E } \|\Delta(w) -\beta
\|_{g_{\ast}g^\ast\ol{E}} \leq
  [L:\Q]^{1/2}
\min_{\alpha \in E}       \|w-\alpha\|_{\ol{E}}
     \label{eq:invapprox}
\end{equation}
  \end{proposition}

\proof The second inequality in (\ref{eq:invapprox}) is clear, and,
according to Proposition \ref{prop.dirsum}, is an equality when
$\ol{E}$ is an orthogonal direct sum of hermitian line bundles over
$\Spec \Z.$

Let $\ol{E}_{0}$ be such a direct sum and $\varphi: E \rightarrow E_{0}$
an isomorphism of $\Z$-modules, and denote
$\Delta_{0}:E_{0}\rightarrow g^\ast E_{0}:=E_{0}\otimes_{\Z}{\mathcal O}_{L}$
the map $.\otimes 1,$ and $\left\|\varphi\right\|_{\infty}$
(resp. $\left\|\varphi^{-1}\right\|_{\infty}$) the operator norm of
$\varphi_{\C}$, or equivalently of $\varphi_{\R}$
(resp. of its inverse). Then, for any $w$ in
$E_{0,\R},$ we have
$$
\min_{\alpha \in E}
\left\|w - \alpha \right\|_{\ol{E}}
\leq 
\left\|\varphi^{-1}\right\|_{\infty}
\min_{\alpha_{0}\in E_{0}}
\left\|\varphi_{\R}(w)-\alpha_{0}\right\|_{\ol{E}_{0}}
$$
and
$$\left\|\varphi\right\|_{\infty}^{-1}
\min_{\beta_{0}\in g^\ast E_{0}}
\left\|\Delta_{0}(\varphi_{\R}(w))-\beta_{0}\right\|_{g_{\ast}g^\ast\ol{E}_{0}}
\leq
\min_{\beta \in g^\ast E}
\left\|\Delta(w) - \beta \right\|_{g_{\ast}g^\ast\ol{E}}.
$$
Besides, since $\ol{E}_{0}$ is an orthogonal direct sum of hermitian line
bundles, we have:
$$\min_{\beta_{0}\in g^\ast E_{0}}
\left\|\Delta_{0}(\varphi_{\R}(w))-\beta_{0}\right\|_{g_{\ast}g^\ast\ol{E}_{0}}
=
[L:\Q]^{1/2} \min_{\alpha_{0}\in E_{0}}
\left\|\varphi_{\R}(w)-\alpha_{0}\right\|_{\ol{E}_{0}}.$$
Consequently,
$$\left\|\varphi\right\|_{\infty}^{-1}
\left\|\varphi^{-1}\right\|_{\infty}^{-1}
\min_{\alpha \in E}
\left\|w - \alpha \right\|_{\ol{E}}
\leq
\min_{\beta \in g^\ast E}
\left\|\Delta(w) - \beta \right\|_{g^\ast\ol{E}}.
$$

To complete the proof, just recall that, as shown in Theorem
\ref{thm.reducArak}, we may find $\ol{E}_{0}$ and $\varphi$ as above
with $\left\|\varphi\right\|_{\infty}$ 
and $\left\|\varphi^{-1}\right\|_{\infty}$ bounded in terms of $n$
alone. 

\qed

The above proof establishes that, 
with the notation of Theorem \ref{thm.reducArak} and 
Theorem \ref{thm.reducArakbis}, Proposition \ref{prop:almostinvariant} 
holds with
$$s(n) =2 r(n,\Q),$$
or
$$s(n)= n (\frac{1}{2} \log n + c(n, \Q)).$$
Using for instance (\ref{simplevalue}), we may take:
$$s(n) =\frac{3}{2} n \log n.$$

Using the numbers $s(n)$, $n\in\N$, from the previous proposition, we
get:

       \begin{corollary} Let $L$ be any number field, and let $g$ be
       the unique morphism
       $\Spec {\mathcal O}_L \to \Spec \Z.$
           For any two hermitian vector bundles
             $\ol{F}$ and $\ol{G}$ over $S:=\Spec \Z,$ and any
             extension class $e$ in $\Exthat^1_{S}\, (F,G),$  we have:
\[
\size_{\ol{F}, \ol{G}}\, (e) - s(\rk F\cdot \rk G)
\leq \size_{g^{\ast}\ol{F}, g^{\ast}\ol{G}} \, (g^{\ast}(e))
\leq \size_{\ol{F}, \ol{G}}\,(e).
\]
             \end{corollary}

\subsection{Cyclotomic base change}\label{subsec.cyclo}

Let $\ol{V}$ be an hermitian vector bundle over $\Spec \Z$ and
$\epsilon$ a non-zero-vector in $V$. We shall say that $\epsilon$ is
\emph{\Vs} in $\ol{V}$ if, for every hermitian vector bundle $\ol{E}$
over $\Spec \Z$, the map
$$
\begin{array}{cccc}
     \Delta: & E_{\R} & \longrightarrow & E_{\R} \otimes_{\R} V_{\R}
\\
      & v & \longmapsto & v \otimes \epsilon
\end{array}$$
satisfies the following compatibility conditions with the Voronoi
cells of $\ol{E}$ and $\ol{E}\otimes \ol{V}$:
\begin{equation}
     \Delta(\cV (\ol{E})) \subset \cV (\ol{E}\otimes \ol{V}).
     \label{eq:Vs}
\end{equation}

Observe that, for any $v \in E_{\R},$
$$\|\Delta (v) \|_{\ol{E}\otimes \ol{V}}=
\|v \otimes \epsilon \|_{\ol{E}\otimes \ol{V}}
=\|\epsilon\|_{\ol{V}} \cdot \|v\|_{\ol{E}}$$
and that $\Delta (E) \subset E \otimes V.$ Therefore
\begin{equation}
     \min_{\beta \in E \otimes V} \|\Delta (v) -
\beta\|_{\ol{E}\otimes \ol{V}}
     \leq
     \|\epsilon\|_{\ol{V}}\cdot\min_{\alpha \in E } \|v -
\alpha\|_{\ol{E}},
     \label{eq:Vsbis}
\end{equation}
and condition (\ref{eq:Vs}) holds iff this inequality is actually an
equality. Observe also that (\ref{eq:Vsbis}) precisely means that
\begin{equation}
     \Delta^{-1}(\cV (\ol{E} \otimes \ol{V})) \subset \cV (\ol{E}),
\end{equation}
and that consequently (\ref{eq:Vs}) is equivalent to

\begin{equation}
     \Delta^{-1}(\cV (\ol{E} \otimes \ol{V})) = \cV (\ol{E}).
     \label{eq:Vster}
\end{equation}

\begin{lemma}\label{lem.Vssorite}
     1) For any euclidean norm $ \|.\|$ on $\R,$  the element $1$ is
     \Vs $(\Z, \|.\|)$.

     2) If $\epsilon$ is \Vs in $\ol{V}$, then, for any hermitian
     vector bundle $\ol{W}$, $\epsilon \oplus 0$ is \Vs in
     $\ol{V}\oplus \ol{W}$.

     3) 
     Let $V$ be a free $\Z$-module of finite rank, $\epsilon$ a
     non-zero element in $V$, and $\|.\|$ and $\|.\|'$ two euclidean
     norms on $V_{\R}$, and let us assume that
     $$\|.\|\leq \|.\|'\,\,\,\mbox{ and } \|\epsilon\| =
\|\epsilon\|'.$$
     Then, if $\epsilon$ is \Vs in $\ol{V}:=(V,\|.\|),$
     it is \Vs in $\ol{V}':=(V,\|.\|').$

     4) If $\epsilon$ (resp. $\epsilon'$) is \Vs in $\ol{V}$ (resp.
     $\ol{V}'$), then
     $\epsilon \otimes \epsilon'$ is \Vs in $\ol{V} \otimes \ol{V}'.$
\end{lemma}

\proof Assertions 1), 2), and 4) are straightforward consequences of
the definition of \Vs vectors. To prove 3), observe that, since
$\|.\|\leq \|.\|',$ for any hermitian vector bundle $\ol{E}$ over
$\Spec \Z,$ we have $\|.\|_{\ol{E}\otimes \ol{V}}
\leq \|.\|_{\ol{E}\otimes \ol{V'}}.$ This is clear by decomposing
$V_{\R}$ as a direct sum of rank-one $\R$-vector spaces which is
orthogonal with respect to both $\|.\|$ and $\|.\|'.$ \qed

Let $L$ be a number field, $\pi$ the morphism from $\Spec {\mathcal O}_{L}$
to $\Spec \Z,$ and $\pi_{\ast}\ol{{\mathcal O}_{L}}$ the direct image of the
trivial hermitian line bundle over $\Spec {\mathcal O}_{L}.$ For any hermitian
vector bundle $\ol{E}$ over $\Spec \Z,$ we have a canonical
isomorphism of hermitian vector bundles over $\Spec \Z$
$$\ol{E} \otimes \pi_{\ast}\ol{{\mathcal O}_{L}} \simeq \pi_{\ast}\pi^{\ast}
\ol{E},$$
under which the canonical map
$$E \hookrightarrow \pi_{\ast}\pi^{\ast}E \simeq E \otimes_{\Z}
{\mathcal O}_{L}$$
gets identified to $v\mapsto v \otimes 1$ by
\ref{arithm.proj.formula}.

Consequently,\emph{ Condition \condp{\Q}{L}{\ol{E}} holds for
any hermitian vector bundle
$\ol{E}$ over $\Spec \Z$ iff $1$ is \Vs in $\pi_{\ast}\ol{{\mathcal O}_{L}}.$}

For any positive integer $n$, let $\zeta_{n}:=e^{2\pi i/n}$ and
$K_{n}:=\Q(\zeta_{n})$ the number field generated by the $n$-th roots
of unity. Its ring of integers ${\mathcal O}_{K_{n}}$ is
$$\Z[\zeta_{n}]= \bigoplus_{i=0}^{\varphi(n)-1}\Z\cdot \zeta_{n}^i.$$
We shall denote $\pi_{n}$ the morphism from $\Spec {\mathcal O}_{K_{n}}$ to
$\Spec \Z$. Observe that, for any two elements $x$ and $y$ in
${\mathcal O}_{K_n}$, their scalar product with respect to the euclidean
structure of $\pi_{n\ast}\ol{{\mathcal O}_{K_{n}}}$ is
\begin{equation}
     \langle x,y
\rangle _{\pi_{n\ast}\ol{{\mathcal O}_{K_{n}}}}=
\sum_{\sigma:K_{n}\hookrightarrow\C}
\sigma(x) \ol{\sigma(y)} = {\rm Tr}_{K_{n}/\Q} x\ol{y}.
\label{eq:scalarcyclo}
\end{equation}
In particular, for any pair of integers $(i,j),$
\begin{equation}
     \langle \zeta_{n}^i,\zeta_{n}^j
\rangle _{\pi_{n\ast}\ol{{\mathcal O}_{K_{n}}}}=
     \sum_{d\in (\Z/n\Z)^\ast} \zeta_{n}^{d(i-j)}.
\label{eq:scalarroots}
\end{equation}

\begin{proposition}\label{prop.cyclo}
     For any positive integer $n$, $1$ is a \Vs vector in
     $\pi_{n\ast}\ol{{\mathcal O}_{K_{n}}}.$
\end{proposition}

Thanks to the theorem of Kronecker-Weber, this establishes that
\emph{ Condition \condp{\Q}{L}{\ol{E}} holds
for any \emph{abelian} extension $L$ of $\Q$ and any hermitian vector
bundle
$\ol{E}$ over $\Spec \Z$.} Using Lemma \ref{primediscr}, we also
obtain that, \emph{
for any number field $L$ and any positive integer 
$n$ prime to the absolute discriminant of $L$, the cyclotomic
extension $L_{n}:=L(\zeta_{n})$ satisfies \condp{L}{L_{n}}{\ol{E}}
for any hermitian vector bundle $\ol{E}$ over $\Spec L.$
} 

The proof will rely on some auxiliary results if 
Kitaoka (\cite{kitaoka93}, Section 7.1 and
Theorem 7.1.3), which he established when investigating minimal
vectors
in tensor products of euclidean lattices.

When $n=p^k$ is a prime power, Proposition (\ref{prop.cyclo}) follows
from the following two lemma:

\begin{lemma}\label{cyclop}
     For any prime number $p$, $1$ is  \Vs  in
$\pi_{\ast}\ol{{\mathcal O}_{K_{p}}}.$
\end{lemma}

\begin{lemma}\label{lem.Epn} \emph{(\cite{kitaoka93}, p. 196, Lemma
7.1.7)}
     For any prime number $p$ and any positive integer $k$, consider
     the $\Z$-submodule
     $$E_{p^k}:=\bigoplus_{i=0}^{p^{k-1}-1} \Z\cdot \zeta_{p^k}^i$$
     in $\Z[\zeta_{p^k}]={\mathcal O}_{K_{p^k}},$ and let $\ol{E}_{p^k}$ be the
     hermitian vector bundle over $\Spec \Z$ defined by $E_{p^k}$
     equipped with the restriction of the euclidean structure  of
     $\pi_{\ast}\ol{{\mathcal O}_{K_{p^k}}}.$

   The vectors in the basis $(\zeta_{p^k}^i)_{0\leq i \leq
     p^{k-1}-1}$ of $E_{p^k,\R}$  are pairwise orthogonal
     with respect to this euclidean structure.

     Moreover
     the morphism of $\Z$-modules
     $$ \psi_{p,k}: E_{p^k} \otimes \Z[\zeta_{p}]
         \longrightarrow
         \Z[\zeta_{p^k}] $$
         which maps $a \otimes b$ to $ab$
     defines an isomorphism of hermitian vector bundles over $\Spec
     \Z$:
     \begin{equation}
         \ol{E}_{p^k} \otimes 
\pi_{p\ast}\ol{{\mathcal O}_{K_{p}}} 
\stackrel{\sim}{\longrightarrow}
         \pi_{p^k\ast}\ol{{\mathcal O}_{K_{p^k}}}.
         \label{eq:Epn}
     \end{equation}
\end{lemma}

Indeed, according to Lemma \ref{lem.Vssorite}, 1) and 2), the first
assertion in Lemma \ref{lem.Epn} shows that $1$ is a \Vs vector in
$\ol{E}_{p^k}.$ Using Lemma \ref{cyclop}
and  Lemma \ref{lem.Vssorite}, 4),
  we deduce that $1\otimes 1$ is \Vs in $\ol{E}_{p^k} \otimes
  \pi_{p\ast}\ol{{\mathcal O}_{K_{p}}}.$ Since $\psi_{p,k}$ maps $1\otimes 1$
  to $1,$ the second
assertion in Lemma (\ref{lem.Epn}) finally establishes that $1$ is
\Vs in
  $\pi_{p^k\ast}\ol{{\mathcal O}_{K_{p^k}}}$.

  Taking the above two lemma for granted,  Proposition
\ref{prop.cyclo}
  then follows --- by writing $n$ as a product of prime powers ---
from  Lemma
  \ref{lem.Vssorite}, 4), and from

  \begin{lemma}\label{lem.n1n2}\emph{(\cite{kitaoka93}, p. 197, Proof
of Theorem
  7.1.3)}
      Let $n_{1}$ and $n_{2}$ be two positive integers that are prime
      together, and let $n:= n_{1}\cdot n_{2}.$

      The morphism of rings
      \begin{equation}
          \begin{array}{ccc}
              \Z[\zeta_{n_{1}}] \otimes
              \Z[\zeta_{n_{2}}]&\longrightarrow & \Z[\zeta_{n}]  \\
             a \otimes b  & \longmapsto   & ab
          \end{array}
          \label{eq:n1n2}
          \end{equation}
          is an isomorphism and defines an isomorphism of hermitian
vector
          bundles over $\Spec \Z$:
          $$\pi_{n_{1}\ast}\ol{{\mathcal O}_{K_{n_{1}}}} 
\otimes \pi_{n_{2}\ast}\ol{{\mathcal O}_{K_{n_{2}}}}
          \simeq \pi_{n\ast}\ol{{\mathcal O}_{K_{n}}}.$$
  \end{lemma}

\noindent\emph{Proof of the Lemma.}
That the morphism $\psi_{p,k}$ in Lemma \ref{lem.Epn} is an
isomorphism is straightforward by considering the bases
  $(\zeta_{p^k}^i)_{0\leq i \leq
     p^{k-1}-1}$ of $E_{p^k}$, $(\zeta_{p}^i)_{0\leq i \leq p-2}$ of
     $\Z[\zeta_{p}]$, and $(\zeta_{p^k}^i)_{0\leq i \leq
     (p-1)p^{k-1}-1})$ of $\Z[\zeta_{p^k}].$ Moreover, from
     (\ref{eq:scalarroots}), one obtain that, for any pair of
     integers $(i,j)$,
     \begin{equation}
         \begin{array}{ccll}
        \langle \zeta_{p^k}^i,\zeta_{p^k}^j
\rangle _{\pi_{p^k\ast}\ol{{\mathcal O}_{K_{p^k}}}} & = &
       (p-1) p^{k-1} & \,\,\mbox{ if }\,\, i=j \mod{p^k}
        \\
          & = & -p^{k-1} &
          \,\,\mbox{ if }\,\, i=j \mod{p^{k-1}}
          \,\,\mbox{ and }\,\, i\neq j \mod{p^k}
          \\
          & = & 0 &
          \,\,\mbox{ if }\,\, i=j \mod{p^{k-1}}.
          \end{array}
         \label{eq:scalarrootspk}
     \end{equation}
     This establishes in particular the first assertion in Lemma
     \ref{lem.Epn}. The isometry property (\ref{eq:Epn}) of
     $\psi_{p,k}$ also is easily derived from (\ref{eq:scalarrootspk})
     and its special case where $k=0$, which shows that
     \begin{equation}
         \begin{array}{ccll}
        \langle \zeta_{p}^i,\zeta_{p}^j
\rangle _{\pi_{p\ast}\ol{{\mathcal O}_{K_{p}}}} & = &
       (p-1)  & \,\,\mbox{ if }\,\, i=j \mod{p}
        \\
          & = & -1&
          \,\,\mbox{ if }\,\, i\neq j \mod{p}.
          \end{array}
         \label{eq:scalarrootsp}
     \end{equation}

     Similarly, by considering the bases $(\zeta_{n_{1}}^i)_{0\leq i
     <\phi(n_{1})}$ of $\Z[\zeta_{n_{1}}],$
     $(\zeta_{n_{2}}^i)_{0\leq i
     <\phi(n_{2})}$ of $\Z[\zeta_{n_{2}}],$
     and $(\zeta_{n}^i)_{0\leq i
     <\phi(n)}$ of $\Z[\zeta_{n}],$ one proves that the morphism
     (\ref{eq:n1n2}) is an isomorphism. Its compatibility with the
     euclidean structures on $\pi_{n_{1}\ast}\ol{{\mathcal O}_{K_{n_{1}}}},$
     $\pi_{n_{2}\ast}\ol{{\mathcal O}_{K_{n_{2}}}},$ and
$\pi_{n\ast}\ol{{\mathcal O}_{K_{n}}}$
     directly follows from the expressions (\ref{eq:scalarcyclo}) for
     the associated scalar products.

     Let us finally establish Lemma \ref{cyclop}. Let $\ol{E}$ be any
     hermitian vector bundle over $\Spec \Z$ and $v$ an element of
     $\cV (\ol{E})$.  To prove that $v\otimes 1$ belongs to $\cV
     (\ol{E}\otimes \pi_{p\ast}\ol{{\mathcal O}_{K_{p}}})$, we need to show
     that, for any $\beta$ in $E\otimes {\mathcal O}_{K_{p}},$
     \begin{equation}
       \|v\otimes 1 - \beta\|_{\ol{E}\otimes
\pi_{p\ast}\ol{{\mathcal O}_{K_{p}}}}
       \geq \|v \otimes 1\|_{\ol{E}\otimes
\pi_{p\ast}\ol{{\mathcal O}_{K_{p}}}}.
         \label{eq:vbeta}
     \end{equation}

     To achieve this, observe that a vector $u$ in
     $E_{\R}\otimes_{\R}{\mathcal O}_{K_{p}}$ may be uniquely written as
     $$u= \sum_{i=0}^{p-2} u_{i}\otimes \zeta_{p}^i,$$
     where $u_{0},$\ldots,$u_{p-2}$ are some vectors in $E_{\R}.$
     Moreover, we have (compare \cite{kitaoka93}, p. 196):
     $$
         \begin{array}{ccll}
         \|u\|_{\ol{E}\otimes \pi_{p\ast}\ol{{\mathcal O}_{K_{p}}}}^2
         & = &  \sum\limits_{0\leq i,j \leq p-2}
         \langle
\zeta_{p}^i,\zeta_{p}^j\rangle_{\pi_{p\ast}\ol{{\mathcal O}_{K_{p}}}}
          \cdot\langle u_{i},u_{j}\rangle _{\ol{E}} & \\
          & = & (p-1)\cdot \sum\limits_{0\leq i\leq p-2}
\|u_{i}\|_{\ol{E}}^2
          - \sum\limits_{0\leq i\neq j \leq p-2}
          \langle u_{i},u_{j} \rangle _{\ol{E}}&
\bigl(\mbox{according to
          (\ref{eq:scalarrootsp})}\bigr)   \\
          & = & \sum\limits_{0\leq i\leq p-2}
          \|u_{i}\|_{\ol{E}}^2 + \sum\limits_{0\leq i< j \leq p-2}
          \|u_{i}-u_{j}\|_{\ol{E}}^2 &
     \end{array}$$
by an elementary computation.

     In particular, we get:
     $$ \|v \otimes 1\|_{\ol{E}\otimes \pi_{p\ast}\ol{{\mathcal O}_{K_{p}}}}^2=
     (p-1) \|v \|_{\ol{E}}^2, $$
     and, if $\beta_{0},$ \ldots, $\beta_{p-2}$ are the elements of
$E$
     such that
     $$\beta= \sum_{i=0}^{p-2} \beta_{i}\otimes \zeta_{p}^i,$$
     we also have:
   $$\|v\otimes 1 - \beta\|_{\ol{E}\otimes
\pi_{p\ast}\ol{{\mathcal O}_{K_{p}}}}^2
   = \|v-\beta_{0}\|_{\ol{E}}^2 + \sum_{1\leq i\leq p-2}
          \|\beta_{i}\|_{\ol{E}}^2 +
          \sum_{0< j \leq p-2}
          \|v-\beta_{0}+\beta_{j}\|_{\ol{E}}^2
          + \sum_{1\leq i\neq j \leq p-2}
          \|\beta_{i}-\beta_{j}\|_{\ol{E}}^2.$$
   Since $v$ belongs to the Voronoi cell $\cV (\ol{E})$,  the $p-1$
terms
\[
\|v-\beta_{0}\|_{\ol{E}}^2,\,
\|v-\beta_{0}+\beta_{1}\|_{\ol{E}}^2,\ldots,
\,\|v-\beta_{0}+\beta_{p-2}\|_{\ol{E}}^2
\]
are all $\geq \|v  \|_{\ol{E}}^2$. This establishes (\ref{eq:vbeta}).
\qed

\subsection{Base change from root lattices}\label{subsec.rootlattice}
The aim of this section is to establish the following:
   \begin{proposition}\label{prop.root}
         If $\ol{E}$ is the hermitian vector bundle over $\Spec \Z$
         defined by any one of the root lattices $A_{n}$ $(n\geq 1)$,
$D_{n}$
         $(n\geq 4)$,
         $E_{6},$ $E_{7},$ or $E_{8},$ then
         condition
         \condp{\Q}{L}{\ol{E}} holds for any number field $L$.
     \end{proposition}

     The definitions of these root lattices are recalled in the proof
     below. We refer the reader to \cite{martinet03}, chapter 4, for
     more informations concerning them. Let us only recall that,
     according to a theorem of Witt (\cite{witt41}), the integral
     lattices\footnote{\emph{i.e.}, the euclidean lattices $(\Gamma,
\|.\|)$
     such that the scalar product of any two vectors in $\Gamma$ is
     an integer. } which are generated by vectors of square norms 1 or
     2 are precisely the orthogonal sums of lattices isometric to one
     of the lattices $(\Z, |.|),$ $A_{n}$
     $(n\geq 1)$, $D_{n}$
         $(n\geq 4)$,
         $E_{6},$ $E_{7},$ or $E_{8}.$
Hence Proposition \ref{prop.root}, Proposition \ref{prop.dirsum},
and Lemma \ref{lem.Psum} yield a proof of Theorem \ref{theorem.bc}
ii).

     {\small

     Our proof of Proposition \ref{prop.root}  relies on the
description of the
     vertices of the Voronoi cell of such  euclidean lattices $\ol{E}$
     appearing in the work of Conway and Sloane (\cite{conwaysloane99}
     Chapter 21, and \cite{conwaysloane91}): these vertices are the
     images under the Weyl group associated to the 
root lattice of a ``small" set
     of vertices ${\mathcal F}(\ol{E})$, whose coordinates are
     explicitly given in \emph{loc. cit.}. Observe that, to prove that
     \condp{\Q}{L}{\ol{E}} holds, it is enough to establish the
     validity of (\ref{eq:CNSVor1}) (or equivalently of
(\ref{eq:CNSVor2})
     or (\ref{eq:CNSVorQ})) for any $P$ in ${\mathcal
     F}(\ol{E})$ and any $\beta$ in $E\otimes_{\Z}{\mathcal O}_{L}.$ Indeed the
     action of the Weyl group extends, by base change, to an isometric
     action on $g^\ast \ol{E}$.

     To check these conditions, we shall use Lemma \ref{lem.sumnorm}
     and the following related inequalities, valid for any number
     field $L$:

     (i) \emph{For any non-zero $\alpha$ in} ${\mathcal O}_{L}$,
     \begin{equation}
         \frac{1}{[L:\Q]}\sum_{\sigma:L\hookrightarrow\C}
      |\sigma(\alpha)|^2 \geq 1.
         \label{eq:sumnorm3}
     \end{equation}
This follows for instance from (\ref{eq:sumnorm1}) and
(\ref{eq:sumnorm2}).

(ii) \emph{For any $\alpha \neq 1/2$ in $\frac{1}{2}{\mathcal O}_{L}$,
     \begin{equation}
         \sum_{\sigma:L\hookrightarrow\C}
         (|\sigma(\alpha)|^2 -\Re \sigma(\alpha)) \geq 0,
         \label{eq:sumnorm4}
     \end{equation}
     and}
     \begin{equation}
         \frac{1}{[L:\Q]}\sum_{\sigma:L\hookrightarrow\C}
      |\sigma(\alpha)-1/3|^2 \geq 1/9.
         \label{eq:sumnorm5}
     \end{equation}
     The lower bound (\ref{eq:sumnorm4}) follows from
     (\ref{eq:sumnorm3}) applied to $2\alpha -1$. To prove
     (\ref{eq:sumnorm5}), observe that, for any embedding
     $\sigma:L\hookrightarrow\C$,
     \begin{eqnarray*}
         4 |\sigma(\alpha) - \frac{1}{3}|^2 & = & 
|\sigma(2\alpha -1) + \frac{1}{3}|^2  \\
          & = & |\sigma(2\alpha -1)|^2 +\frac{2}{3} \Re
\sigma(2\alpha -1) +
         \frac{1}{9}   \\
          & \geq & |\sigma(2\alpha -1)|^2 - \frac{2}{3}
|\sigma(2\alpha -1)| +
          \frac{1}{9}    \\
          & = & |\sigma(2\alpha -1)|^2 - |\sigma(2\alpha -1)|+
          \frac{1}{3} |\sigma(2\alpha -1)| +
          \frac{1}{9},
      \end{eqnarray*}
      and apply Lemma \ref{lem.sumnorm}
and (\ref{eq:sumnorm1}) to $2\alpha -1.$

     For any two elements $a$ and $b$ of some set $\mathcal E$ and any
     two positive integers $i$ and $j$, we shall denote
     $$(a^{× i}, b^{× j}):=(\underset{i\mbox{
     times}}{\underbrace{a,\ldots,a}}, \underset{j\mbox{
     times}}{\underbrace{b,\ldots,b}})\,\, (\in {\mathcal E}^{i+j}).$$

     \noindent \emph{Proof for} $A_{n}.$
     The  lattice $A_{n}$ is the lattice $\Z^{n+1}\cap A_{n,\R}$
     in the hyperplane
     $$A_{n,\R}:=\{(x_k)_{0\leq k \leq n} \in 
\R^{n+1}\mid \sum_{k=0}^{n}x_{k}=0\}$$
     of $\R^{n+1}$ equipped with the standard euclidean norm. The
     corresponding hermitian vector bundle over $\Spec \Z$
     coincides with the hermitian vector bundle $\ol{\ker p}$
     considered in Example \ref{exAn} above when $S=\Spec \Z$.

     According to \cite{conwaysloane99}
     Chapter 21,3.B, or \cite{conwaysloane91}, section 4, the
     vertices of the Voronoi cell of the root lattice $A_{n}$ 
are the images under its
     Weyl group of the following $n$ points in $A_{n,\R}$:
     $$[i]:=\left(\left(\frac{j}{n+1}\right)^{× i},
     \left(-\frac{i}{n+1}\right)^{× j}\right), \,\,\,\, 1\leq
     i\leq n, \,\,\, j:= n+1-i.$$
     By considering the conditions (\ref{eq:CNSVorQ}) with $P=[i],$
     $1\leq i\leq n,$ we are reduced to proving the following

  \begin{lemma}\label{lem.An}
      Let $L$ be a number field, $i$ an integer in $\{1,\ldots,n\},$
      and $j:= n+1-i.$

      For any $(\beta_{k})_{0\leq k \leq n}$ in ${\mathcal O}_{L}^{n+1}$ such
      that $\sum_{k=0}^{n}\beta_{k}=0,$ we have:
      \begin{equation}
      2  \sum_{\sigma:L\hookrightarrow\C} \left[\sum_{0\leq k \leq
i-1}
      \frac{j}{n+1}\,{\rm Re}\,\sigma(\beta_{k})
      -\sum_{i\leq k \leq n}
      \frac{i}{n+1}\,{\rm Re}\,\sigma(\beta_{k}) \right]
      \leq  \sum_{\sigma:L\hookrightarrow\C} \sum_{0\leq k \leq n}
      |\sigma(\beta_{k}|^2.
         \label{eq:An}
     \end{equation}
  \end{lemma}

  To prove this lemma, observe that, for any embedding
  $\sigma:L\hookrightarrow\C$,
  $$\sum_{\sigma:L\hookrightarrow\C}\, {\rm Re}\,\sigma(\beta_{k})
=0.$$
Consequently the left-hand side of (\ref{eq:An}) is equal to:
\begin{eqnarray*}
      2  \sum_{\sigma:L\hookrightarrow\C}\left(\frac{j}{n+1} +
  \frac{i}{n+1}\right) \sum_{0\leq k \leq i-1}
      \,{\rm Re}\,\sigma(\beta_{k})
       & = &  2  \sum_{\sigma:L\hookrightarrow\C}  \sum_{0\leq k \leq
i-1}
      \,{\rm Re}\,\sigma(\beta_{k})
     \\
      & = &     \sum_{\sigma:L\hookrightarrow\C}\left( \sum_{0\leq k
\leq i-1}
      \,{\rm Re}\,\sigma(\beta_{k}) -    \sum_{i\leq k \leq n}
      \,{\rm Re}\,\sigma(\beta_{k})\right).
\end{eqnarray*}
The inequality (\ref{eq:An}) now follows from Lemma
(\ref{lem.sumnorm}) applied to
$\alpha$ in
$$\{\beta_{0},\ldots,\beta_{i-1},-\beta_{i},\ldots,-\beta_{n}\}.$$
\qed

     \noindent \emph{Proof for} $D_{n}.$
     The lattice $D_{n}$ is the sublattice of index 2 in the
     lattice $\Z^{n}$
     of $\R^{n}$ equipped with the standard euclidean norm consisting
     of vectors $(x_{1},\ldots,x_{n})$ for which the $x_{i}$ are
     integers with an even sum.

     According to \cite{conwaysloane99}
     Chapter 21, 3.C, and \cite{conwaysloane91}, section 6, the
     vertices of the Voronoi cell of $D_{n}$ are the images under its
     Weyl group of the points $((1/2)^{× n})$ and
     $(1, 0^{× (n-1)})$ in $\R^n$.
     By considering the condition (\ref{eq:CNSVorQ}) (resp.
     (\ref{eq:CNSVor1}) with $P=((1/2)^{× n})$ (resp. with
     $P=(1, 0^{× (n-1)})$),
     we are reduced to proving:

  \begin{lemma}\label{lem.Dn}
      Let $L$ be a number field. For any $(\beta_{k})_{1\leq k\leq n}$
      in ${\mathcal O}_{L}^n$ such that
      $\sum_{k=1}^{n} \beta_{k}$ belongs to $2{\mathcal O}_{L},$ the following
      inequalities hold:
      \begin{equation}
      2  \sum_{\sigma:L\hookrightarrow\C} \sum_{1\leq k \leq n}
\frac{1}{2}
     \,{\rm Re}\,\sigma(\beta_{k})
      \leq  \sum_{\sigma:L\hookrightarrow\C} \sum_{1\leq k \leq n}
      |\sigma(\beta_{k})|^2,
         \label{eq:Dn1}
     \end{equation}
     and
       \begin{equation}
           \frac{1}{[L:\Q]}\sum_{\sigma:L\hookrightarrow\C}
      |\sigma(\beta_{1})-1|^2 +
      \frac{1}{[L:\Q]}\sum_{\sigma:L\hookrightarrow\C}  \sum_{2\leq k
\leq n}
      |\sigma(\beta_{k})|^2 \geq 1.
         \label{eq:Dn2}
     \end{equation}
  \end{lemma}

  The first inequality (\ref{eq:Dn1}) directly follows from
  Lemma
(\ref{lem.sumnorm}). To prove (\ref{eq:Dn2}), observe that
that at least one of the algebraic integers $\beta_{1}-1,$
$\beta_{2},$ $\dots,$ $\beta_{n}$ is not zero, and use
(\ref{eq:sumnorm3}).
\qed

  \noindent \emph{Proof for} $E_{8}.$
     The lattice $E_{8}$ is the lattice in
$\R^{n}$ equipped with the standard euclidean norm consisting
     of the vectors $(x_{1},\ldots,x_{8})$ for which the $x_{i}$ are
all
     in $\Z$ or all in $1/2+\Z$, and have an even sum.

     According to \cite{conwaysloane91}, section 8, the
     vertices of the Voronoi cell of $E_{8}$ are the images under its
     Weyl group of the points $((0)^{× 7},1)$ and
     $((1/3)^{× 7},-1/3)$ in $\R^8$.
     By considering the condition
     (\ref{eq:CNSVor1}) with $P=((0)^{× 7},1)$ and with
     $P=((1/3)^{× 7},-1/3)$,
     we are reduced to proving:

  \begin{lemma}\label{lem.E8}
      Let $L$ be a number field.
      If $\beta_{1},\ldots,\beta_{8}$ in $\frac{1}{2}{\mathcal O}_{L}$ have the
      same class $\gamma$ in $(\frac{1}{2}{\mathcal O}_{L})/{\mathcal O}_{L}$ and if their
      sum $\beta_{1}+\ldots+\beta_{8}$ belongs to $2{\mathcal O}_{L},$ then the
      following inequalities hold:
    \begin{equation}
        \frac{1}{[L:\Q]}\sum_{\sigma:L\hookrightarrow\C}
      \left( \sum_{1\leq k \leq 7} |\sigma(\beta_{k})|^2
      +
      |\sigma(\beta_{8})-1|^2 \right)
     \geq 1,
         \label{eq:E81}
     \end{equation}
     and
       \begin{equation}
\frac{1}{[L:\Q]}\sum_{\sigma:L\hookrightarrow\C} \left( \sum_{1\leq k
\leq
      7}
      |\sigma(\beta_{k})-1/3|^2
      +
      |\sigma(\beta_{8})+1/3|^2 \right)
     \geq 8/9.
         \label{eq:E82}
     \end{equation}
  \end{lemma}

  When $\gamma \neq [0],$ each of the expressions
  $$\frac{1}{[L:\Q]}\sum_{\sigma:L\hookrightarrow\C}
      |\sigma(\beta_{k})|^2, \,\,\,\, 1\leq k\leq 7 $$
      and
      $$\frac{1}{[L:\Q]}\sum_{\sigma:L\hookrightarrow\C}
      |\sigma(\beta_{8})-1|^2$$
      is $\geq 1/4$, by (\ref{eq:sumnorm3}) applied to $\alpha$ in
      $\{2\beta_{1},\ldots,2\beta_{7},2\beta_{8}-2\}$, and
      (\ref{eq:E81}) follows. When $\gamma =[0]$, 
these expressions are $\geq 1,$
      unless $\beta_{k}=0$ (when $1\leq k \leq 7$), or $\beta_{8}=1.$
      Since $\beta_{1}+\ldots+\beta_{8}\neq 1,$ the vector
      $(\beta_{1},\ldots,\beta_{8})$ does not coincide with $((0)^{×
      7},1)$, and at least one of them is $\geq 1$, which
      implies (\ref{eq:E81}).

      If $\beta_{1},\ldots,\beta_{7}$ and $-\beta_{8}$ are all
distinct
      from $1/2,$ then (\ref{eq:E82}) directly follows from
      (\ref{eq:sumnorm5}). Otherwise, $\gamma$ is $[1/2],$ and the
      vector
      $(\beta_{k})_{1\leq k \leq 8} -(1/2)^{× 8}$ belongs to
      $D_{8}\otimes {\mathcal O}_{L}.$ The validity of \condp{\Q}{L}{D_{8}},
      applied to $v=((1/3)^{× 7},-1/3) - (1/2)^{× 8}$, shows that
      the left-hand side of (\ref{eq:E82}) is greater or equal to the
      minimum of the expression
      \begin{equation}
           \sum_{1\leq k \leq 7}
      (\alpha_{k}- 1/3)^2 + (\alpha_{8}+ 1/3)^2
    \label{eq:alpha8}
      \end{equation}
      when $(\alpha_{k})_{1\leq k \leq 8}$ varies in $D_{8}+(1/2)^{×
      8}$. This minimum is clearly greater or equal than the one of
      (\ref{eq:alpha8}) when $(\alpha_{k})_{1\leq k \leq 8}$ varies in
      $E_{8}$ --- which is $\|((1/3)^{× 7},-1/3)\|^2=8/9$ since
      $((1/3)^{× 7},-1/3)$ belongs to the Voronoi cell of $E_{8}$.
\qed

\noindent \emph{Proof for} $E_{7}.$
     The lattice $E_{7}$ may be realized as the 
lattice $E_{8}\cap E_{7,\R}$ in the
     hyperplane
     $$E_{7,\R}:=\{(x_k)_{1\leq k \leq 8} \in 
\R^{8}\mid \sum_{k=1}^{8}x_{k}=0\}$$
     of $\R^{8}$ equipped with the standard euclidean.

     According to \cite{conwaysloane91}, section 9, the
     vertices of the Voronoi cell of $E_{7}$ are the images under its
     Weyl group of the points $(7/8,(-1/8)^{× 7})$ and
     $((3/4)^{× 2},(-1/4)^{× 6})$ in $E_{7,\R}.$
     By considering the condition
     (\ref{eq:CNSVorQ}) with $P=(7/8,(-1/8)^{× 7})$ and
     $P=((3/4)^{× 2},(-1/4)^{× 6})$,
     we are reduced to proving:

  \begin{lemma}\label{lem.E7}
      Let $L$ be a number field.
      If $\beta_{1},\ldots,\beta_{8}$ in $\frac{1}{2}{\mathcal O}_{L}$ have the
      same class $\gamma$ in $(\frac{1}{2}{\mathcal O}_{L})/{\mathcal O}_{L}$ and if
      $\beta_{1}+\ldots+\beta_{8}=0,$ then the
      following inequalities hold:
    \begin{equation}
\sum_{\sigma:L\hookrightarrow\C}
      \left(|\sigma(\beta_{1})|^2 -\frac{7}{4} \Re \sigma(\beta_{1}))
      +
            \sum_{2\leq k \leq 8}
     (|\sigma(\beta_{k})|^2 +\frac{1}{4} \Re \sigma(\beta_{k}))\right)
      \geq 0,
         \label{eq:E71}
     \end{equation}
     and
       \begin{equation}
  \sum_{\sigma:L\hookrightarrow\C} \left(\sum_{1\leq k \leq
      2}
      (|\sigma(\beta_{k})|^2 -\frac{3}{2} \Re \sigma(\beta_{k})
      +
       \sum_{3\leq k \leq
      8}
      (|\sigma(\beta_{k})|^2 +\frac{1}{2} \Re \sigma(\beta_{k}))
\right)
     \geq 0.
         \label{eq:E72}
     \end{equation}
  \end{lemma}

  Observe that, with the notation of the proof for $A_{n}$ with $n=8$,
  we have
  $$(7/8,(-1/8)^{× 7})=[1],$$
  and
  $$((3/4)^{× 2},(-1/4)^{× 6})=[2].$$
  Therefore the computation in the proof of Lemma \ref{lem.An},
  with $n=8,$ shows that (\ref{eq:E71}) and (\ref{eq:E72}) may also be
  written
  \begin{equation}
\sum_{\sigma:L\hookrightarrow\C}
      (|\sigma(\beta_{1})|^2 - \Re \sigma(\beta_{1}))
      +
            \sum_{2\leq k \leq 8} \sum_{\sigma:L\hookrightarrow\C}
     (|\sigma(\beta_{k})|^2 + \Re \sigma(\beta_{k}))
      \geq 0,
         \label{eq:E71bis}
     \end{equation}
     and
       \begin{equation}
  \sum_{1\leq k \leq
      2} \sum_{\sigma:L\hookrightarrow\C}
      (|\sigma(\beta_{k})|^2 - \Re \sigma(\beta_{k}))
      +
       \sum_{3\leq k \leq
      8}\sum_{\sigma:L\hookrightarrow\C}
      (|\sigma(\beta_{k})|^2 + \Re \sigma(\beta_{k}))
     \geq 0.
         \label{eq:E72bis}
     \end{equation}
  These inequalities directly follow from (\ref{eq:sumnorm4}) when
  $\gamma\neq [1/2].$

  If $\gamma =[1/2],$ we may introduce
  \begin{equation}
      (\tilde{\beta}_{k})_{1\leq k  \leq 8}:=(\beta_{k})_{1\leq k \leq
  8}-
  ((1/2)^{× 4},(-1/2)^{× 4})).
      \label{eq:betatilde}
  \end{equation}
  Then $\tilde{\beta}_{1},\ldots,\tilde{\beta}_{8}$ belong to
${\mathcal O}_{L}$,
  and
  their sum vanishes. Moreover, expressed in terms of
  $\tilde{\beta}_{1},\ldots,\tilde{\beta}_{8},$
  (\ref{eq:E71bis}) takes the form:
  \begin{equation}
\sum_{1\leq k \leq 8} \frac{1}{[L:\Q]}\sum_{\sigma:L\hookrightarrow\C}
     (|\sigma(\tilde{\beta}_{k})|^2 - \epsilon_{k} \Re
     \sigma(\tilde{\beta}_{k})) +1
      \geq 0,
         \label{eq:E71ter}
     \end{equation}
     where $\epsilon_{k}:= 1$ if $k\in \{1,5,6,7,8\}$ and
     $\epsilon_{k}:=-1$ if $k\in \{2,3,4\}.$ Similarly,
     (\ref{eq:E72bis}) is equivalent to:
       \begin{equation}
  \sum_{1\leq k \leq 8} \sum_{\sigma:L\hookrightarrow\C}
     (|\sigma(\tilde{\beta}_{k})|^2 - \epsilon'_{k} \Re
     \sigma(\tilde{\beta}_{k}))
      \geq 0,
         \label{eq:E72ter}
     \end{equation}
     where $\epsilon'_{k}:= 1$ if $k\in \{1,2,5,6,7,8\}$ and
     $\epsilon'_{k}:=-1$ if $k\in \{3,4\}.$

     The inequalities (\ref{eq:E71ter}) and (\ref{eq:E72ter})
     directly follow from Lemma \ref{lem.sumnorm}.
\qed

\noindent \emph{Proof for} $E_{6}.$
     The lattice $E_{6}$ may be realized as the 
lattice $E_{8}\cap E_{6,\R}$ in the
     codimension 2 subspace
     $$E_{6,\R}:=\{(x_k)_{1\leq k \leq 8} \in \R^{8}\mid
     x_{1}+x_{8}=x_{2}+\ldots+x_{7}=0\}$$
     of $\R^{8}$ equipped with the standard euclidean norm.

     According to \cite{conwaysloane91}, section 10, the
     vertices of the Voronoi cell of $E_{6}$ are the images under its
     Weyl group of the point $(0,(2/3)^{× 2},(-1/3)^{× 4},0)$
     in $E_{6,\R}.$
     By considering the condition
     (\ref{eq:CNSVorQ}) with $P=(0,(2/3)^{× 2},(-1/3)^{× 4},0)$,
     we are reduced to proving:

  \begin{lemma}\label{lem.E6}
      Let $L$ be a number field.
      If $\beta_{1},\ldots,\beta_{8}$ in $\frac{1}{2}{\mathcal O}_{L}$ have the
      same class $\gamma$ in $(\frac{1}{2}{\mathcal O}_{L})/{\mathcal O}_{L}$ and if
      $\beta_{1}+\beta_{8}=0,$ and $\beta_{2}+\ldots+\beta_{6}=0,$
then the
      following inequality holds:
    \begin{equation}
\sum_{\sigma:L\hookrightarrow\C}
      \left(|\sigma(\beta_{1})|^2)
      +
            \sum_{2\leq k \leq 3}
     (|\sigma(\beta_{k})|^2 -\frac{4}{3} \Re \sigma(\beta_{k})) +
     \sum_{4\leq k \leq 7}
     (|\sigma(\beta_{k})|^2 +\frac{2}{3} \Re \sigma(\beta_{k}))
     +
      |\sigma(\beta_{8})|^2\right)
      \geq 0.
         \label{eq:E6}
     \end{equation}
  \end{lemma}

Again the computation in the proof of Lemma \ref{lem.An},
now with $n=8,$ $i=2$ and $j=4,$ allows us to replace
(\ref{eq:E6}) by the equivalent condition:
  \begin{equation}
            \sum_{1\leq k \leq 8}  \sum_{\sigma:L\hookrightarrow\C}
     (|\sigma(\beta_{k})|^2 - \eta_{k} \Re \sigma(\beta_{k}))
      \geq 0,
         \label{eq:E6bis}
     \end{equation}
  where $\eta_{1}=\eta_{8}:=0,$ $\eta_{2}=\eta_{3}:=1,$ and
  $\eta_{4}=\eta_{5}=\eta_{6}=\eta_{7}=-1.$
According to (\ref{eq:sumnorm4}), this holds if $\beta_{k}\neq
\eta_{k}/2$ for $2\leq k \leq 7,$ in particular when
  $\gamma\neq [1/2].$

  When $\gamma =[1/2],$ we  introduce again
$\tilde{\beta}_{1},\ldots,\tilde{\beta}_{8}$
  defined by (\ref{eq:betatilde}).
  They belong to ${\mathcal O}_{L}$,
  the sum $\tilde{\beta}_{2}+\ldots+\tilde{\beta}_{7}$ vanishes, and
  (\ref{eq:E6bis}) takes the form:
  \begin{equation}
\sum_{1\leq k \leq 8}  \sum_{\sigma:L\hookrightarrow\C}
     (|\sigma(\beta_{k})|^2 - \eta'_{k} \Re \sigma(\beta_{k}))
      \geq 0,
         \label{eq:E6ter}
     \end{equation}
     where $\eta'_{k}:= 1$ if $k\in \{2,3,5,6,7,8\}$ and
     $\eta'_{k}:=-1$ if $k\in \{1,4\}.$ Finally,
     (\ref{eq:E6ter}) follows from Lemma \ref{lem.sumnorm} applied to
     $\alpha=\eta'_{k} \tilde{\beta}_{k},$ $1\leq k \leq 8.$
\qed
}

\subsection{Base change from lattices of Voronoi's first
kind}\label{proof.iii}
From the description of the Voronoi cell of
$\overline V(p), \ p\in \bigl(\R^*_+\bigr)^{\frac{n(n+1)}{2}}$ in the
appendix,
we may easily derive:

\begin{proposition}\label{prop4.1}
If $\overline E$ is the hermitian vector bundle over
$\Spec\ \Z$ defined by an euclidean lattice of Voronoi's first kind,
then the
condition \condp{\Q}{L}{\ol{E}} holds for any
number field $L$.
\end{proposition}

\proof
A straightforward perturbation argument shows that we may restrict to
the case
where $ \overline E$ possesses a strictly obtuse superbase, that is,
to the case
where $\overline E = \overline V(p)$ for some $p \in
\bigl(\R^*_+\bigr)^{\frac{n(n+1)}{2}}$,
where $n:= {\rm rk}\, E$.
According to Proposition \ref{prop3.4} every vector of
$\mathcal{V}\bigl(\overline V(p)\bigr)$ is of the form $s_A$ for some
$A\in \mathcal{S}(n)$.
Moreover, for any such
\[A=\bigl\{\{i_1,\ldots,\, i_n\},\ \{i_2,\ldots,\,
i_n\},\ldots,\,\{i_n\}\bigr\}\]
and any $(x_0,\ldots,\, x_n) \in \R^{n+1}$ if we let
\[x:= \sum^n_{k=0} x_k\, v_{i_k},\]
then Lemma \ref{abs3.3} shows that
\begin{equation}
\label{gl4.1}
\begin{array}{rcl}
\Vert x\Vert^2_{\overline V(p)} - 2\, \langle s_A,
x\rangle_{\overline V(p)} & = & \Vert x-s_A \Vert^2_{\overline V(p)}
- \Vert s_A\Vert^2_{\overline V(p)} \\
& = & \sum\limits_{0\leq k< \ell\leq n} \,p_{i_ki_\ell}
\bigl[(x_\ell-x_k)^2 - (x_\ell
-x_k)\bigr].
\end{array}
\end{equation}

For any complex embedding $\sigma\colon L\to \C$ of some number field
$L$, consider
\[\begin{array}{rcl}
\sigma_{V(p)} \colon V(p)\otimes L & \longrightarrow & V(p)_\C =
V(p)\otimes \C \\
e\otimes \lambda & \longmapsto & e\otimes \sigma(\lambda).
\end{array}\]
As observed in Lemma \ref{lem.CNSVor}, to complete the proof of the
Proposition it is enough to check that,
for any $\beta \in V(p) \otimes_\Z \,\mathcal{O}_L$,
\begin{equation}
\label{gl4.2}
\sum_{\sigma:L\hookrightarrow\C} \,\Bigl[ \big\Vert \sigma_{V(p)}\,
(\beta)\Vert^2_{\overline V(p)} - 2 \mbox{ Re } \langle s_A,
\,\sigma_{V(p)}\,(\beta) \rangle_{\overline V(p)}\Bigr] \geq 0.
\end{equation}
Any such $\beta$ may be written
\[
\beta = \sum^n_{k=0} v_{i_k} \otimes \beta_k
\]
for some $(\beta_0,\ldots ,\,\beta_n) \in \mathcal{O}^{n+1}_L$, and
then (\ref{gl4.1})
shows that the left-hand side of (\ref{gl4.1}) is equal to
\[
\sum_{\sigma:L\hookrightarrow\C} \ \sum_{0\leq k< \ell\leq n} \,p_{i_k
i_\ell} \,\Bigl[\big\vert \sigma (\beta_k-\beta_\ell)\big\vert^2 -
\mbox{ Re } \sigma (\beta_\ell-\beta_k)\Bigr].
\]
This is indeed non-negative by Lemma \ref{lem.sumnorm}.
\qed

\begin{appendix}

\section{Extension groups}\label{appa}

In this Appendix, we discuss various
definitions of extension groups of sheaves of modules (notably, of
groups of $1$-extensions) which are used in this paper.
We pay some special
attention to the sign issues which arise in defining canonical
isomorphisms between various constructions of extension groups.

\subsection{Notation and sign conventions}\label{sign}
When dealing with categories of complexes in abelian categories and
their derived categories,  we follow the same conventions, notably
concerning  signs, as in  \cite{berthelotbreenmessing82},
\cite{illusie96}, and
\cite{conrad00}.  (The reader should be aware that these conventions
do \emph{not} agree with the ones in some classical
references such that \cite{hartshorne66}, 
\cite{deligne73}, or \cite{weibel94}.)

We
refer to \cite{conrad00}, sections 1.2-3, for a thorough discussion of
  signs issues, and we  just specify some basic definitions where some
  ambiguity may arise.

If $X= (X^k)_{k \in \Z}$ is a complex (with cohomological
grading) in some abelian category $\cC$, then, 
for any integer $i,$ $X[i]$ is the complex 
defined by
$$X[i]^k:=X^{k+i}\,\,\mbox{ and }\,\, d_{X[i]}=(-1)^i d_X.$$
For any map of complexes $f^\bullet: X^\bullet \rightarrow
Y^\bullet,$ the map of complexes $f^\bullet[i]: X[i]^\bullet \rightarrow
Y[i]^\bullet$ is defined to be $f^{k+i}$ in degree $i.$ The
isomorphisms
\begin{equation}\label{sanssigne}
H^k(X[i]) \simeq H^{k+i}(X)
\end{equation}
are defined without the intervention of signs.

A \emph{triangle} in some category of complexes (such as the usual
category
of complexes $C(\cC)$ in some  abelian category $\cC,$ or the
  category $K(\cC)$ of such complexes with morphisms up to homotopy,
  or its derived category $D(\cC)$) is a sequence of morphisms of
complexes of the form
$$X \rightarrow Y \rightarrow Z \rightarrow X[1].$$

The \emph{cone} $C(f)$ of a morphism of complexes $f: X
\rightarrow Y$ is defined by
$$C(f)^k:= (X[1] \oplus Y)^k =X^{k+1} \oplus Y^k \,\,\mbox{ and }
d_{C(f)}(x,y):=(-d_X(x), f(x) + d_Y(y)).$$
The \emph{standard triangle} associated to $f$ is:
$$X\stackrel{f}{\longrightarrow} Y \stackrel{i}{\longrightarrow}
C(f) \stackrel{p}{\longrightarrow} X[1],$$
where $i:Y\rightarrow C(f)$ is the natural injection, and
$p:C(f)\rightarrow X[1]$ the \emph{opposite} of the natural
projection.

A \emph{distinguished triangle} is a triangle isomorphic to such a
standard
triangle. For instance, a short exact sequence of complexes
$$\cE:\, 0 \longrightarrow
X\stackrel{f}{\longrightarrow} Y \stackrel{g}{\longrightarrow}
Z \longrightarrow 0$$
determines a distinguished triangle
$$X\stackrel{f}{\longrightarrow} Y \stackrel{g}{\longrightarrow} Z
\stackrel{\partial_{\cE}}{\longrightarrow} X[1]$$
in $D(\cC),$  where $\partial_{\cE}$ is defined 
as the composite of $p:C(f)\rightarrow
X[1]$ with the inverse of
the quasi-isomorphism $q: C(f) \rightarrow Z$ given by $g$ on $Y$ and
$0$ on $X[1]$:
$$\partial_{\cE}: Z \stackrel{q}{\longleftarrow} C(f)
\stackrel{p}{\longrightarrow} X[1].$$

Observe that, using these sign conventions, for 
any integer $i$, the map between cohomology 
groups determined by $\partial_{\cE}$:
$$H^i(\partial_{\cE}): H^i(Z) \longrightarrow H^i(X[1]) \simeq
H^{i+1}(X)$$
coincides with the usual\footnote{namely, defined without intervention
of signs, as in \cite{maclane95}, II.4, or 
\cite{weibel94}, 1.3 : the graph of 
$H^i(\partial_{\cE})$ is the image in
$H^i(Z) \times H^{i+1}(X)$ of $$\{(z,x)\in Z^i 
\times X^{i+1}\mid d_Z(z)=0 \mbox{ and } \exists 
y \in Y^i, g^i(y)= z \mbox{ and } d_Y(y) = 
f(x)\}.$$
}
boundary map in the long exact sequence of 
cohomology groups associated to $\cE.$

For any two complexes $X$ and $Y$ in $\cC,$ and any integer $i,$ we
let:
$$\Ext_\cC^i(X,Y):= \Hom_{D(\cC)}(X,Y[i]).$$
If $Z$ is again a complex in $\cC,$ and $j$ an integer, the
composition of two elements
$$f \in \Ext_\cC^i(X,Y)\,\,\mbox{ and }\,\,g \in \Ext_\cC^j(Y,Z)$$
is the element
$$g\circ f \in \Ext_\cC^{i+j}(X,Z)$$
defined by composing the arrows
$$f:X\longrightarrow Y[i] \,\,\mbox{ and }\,\, g[i]:Y[i]
\longrightarrow Z[i+j]$$
in $D(\cC).$

\subsection{Extension groups of sheaves of modules}

Let $X$ be any topological space, equip\-ped with a sheaf $\cA$ of
commutative rings. Let $\cAmod$ be the abelian category of sheaves of
$\cA$-modules over $X$ and $D^+(\cAmod ),$
the derived category of bounded below complexes of 
$\cA$-modules. It is a full triangulated subcategory of  
$D(\cAmod ),$
the derived category of  complexes of 
$\cA$-modules. Moreover, since the abelian category 
$\cAmod$ admits
enough injectives, the derived category $D^+(\cAmod )$ is equivalent
to
the category of bounded below complexes of 
injectives sheaves of $\cA$-modules over
$X$, with morphisms the homotopy classes of morphisms of complexes.

Any object in $\cAmod$ may be seen as a complex concentrated in degree
0. Accordingly, for any two objects $E$ and $F$ in $\cAmod$, the
associated
\emph{extension groups} are defined as
$$\Ext_{\cA}^i (E,F):=\Hom_{D^+(\cAmod )}(E, F[i]), \,\, \mbox{ for
any } i\in \N.$$
For any resolution
\begin{equation}
  0\rightarrow F \rightarrow F^0 \rightarrow F^1 \rightarrow \ldots
\rightarrow F^i \rightarrow F^{i+1} \rightarrow \ldots
     \label{eq:res}
\end{equation}
of $F$ by injective sheaves of $\cA$-modules, these are the
cohomology group of the complex of abelian groups deduced from
$$0\rightarrow F^0 \rightarrow F^1 \rightarrow \ldots
\rightarrow F^i \rightarrow F^{i+1} \rightarrow \ldots$$
by applying the functor $\Hom_{\cA}(E, .).$

When $\cA$ is the constant sheaf $\Z_{X},$ $\cAmod$ is the category of
sheaves of abelian groups on $X,$ and for any such sheaf $F$ and any
integer $i\in \N,$  the extension group $\Ext^i_{\Z_{X}}(\Z_{X}, F)$
coincides with the $i$-th cohomology group $H^i(X, F),$ by the very
definition of the latter.
More generally, for any sheaf of ring $\cA$ over $X,$ the injective
objects in $\cAmod$ are flabby sheaves of abelian groups, and, for any
$F$ in $\cAmod,$ the extension group $\Ext^i_{\cA}(\cA,F)$ coincides
with the cohomology group $H^i(X,F)$ of $F$ considered as a sheaf of
abelian groups.

If the sheaf of
$\cA$-modules  $E$ is locally free of finite rank,
and if $$\check{E}:={\mathcal Hom}_{\cA}(E,\cA)$$
denotes the dual sheaf, the functorial isomorphism of sheaf of
$\cA$-modules
\begin{equation}
   \check{E}\otimes_{\cA}G  \stackrel{\sim}{\longrightarrow}
{\mathcal Hom}_{\cA}(E,G),
     \label{eq:dualiso}
\end{equation}
where $G$ denotes any sheaf in $\cAmod$, leads to isomorphisms of
extension groups
\begin{equation}
     \Ext_{\cA}^i(\cA, \check{E}\otimes_{\cA}F) 
\stackrel{\sim}{\longrightarrow}
     \Ext_{\cA}^i(E,F),
     \label{eq:extisoA}
\end{equation}
defined by applying (\ref{eq:dualiso}) to  any resolution
(\ref{eq:res}), and taking  global sections and cohomology groups.
Finally, for any two sheaves of $\cA$-modules $E$ and $F,$ with $E$
locally free of finite rank, we get canonical isomorphisms
\begin{equation}
   \Ext_{\cA}^i(E,F)
   \stackrel{\sim}{\longrightarrow} H^i(X,\check{E}\otimes_{\cA}F).
     \label{eq:extisoB}
\end{equation}

\subsection{Extension groups of quasi-coherent sheaves of modules
over schemes}

When $\cX:=(X,\cA)$ is a scheme, we may consider the full subcategory
$\cXqc$
of the category $\cXmod$ of sheaves of ${\mathcal O}_{\cX}:=\cA$-modules over
$X$
  defined by the quasicoherent sheaves.
  When $\cX$ is the affine scheme $\Spec A$ defined by some
commutative
  ring $A,$ the functor $\Gamma(X, .)$ realizes an equivalence of
  category between $\cXqc$ and the category of $A$-modules.

  When the scheme $\cX$ is quasi-compact and quasi-separated (\eg,
when
  $\cX$ is noetherian), the category $\cXqc$ has enough injectives. If
  indeed $\cX$ is noetherian, the inclusion functor
  from $\cXqc$
to $\cXmod$ preserves injectives, and its extension, as a functor of
triangulated category from $D^+(\cXqc)$ to $D^+(\cXmod)$ is fully
faithful and its essential image is the subcategory of $D^+(\cXmod)$
defined by bounded below complexes in $\cXmod$ with quasi-coherent
cohomology (see for instance \cite{hartshorne66}, II.7,
\cite{illusie71a}, section 3, and
\cite{thomason90}, Appendix B).

We shall frequently use the following consequence of the above full
faithfulness:
\emph{for any two quasicoherent sheaves $E$ and 
$F$ on a noetherian scheme $\cX,$
the morphism of abelian groups
$$\Hom_{D^+(\cXqc )}(E, F[i])
\longrightarrow
\Hom_{D^+(\cXmod )}(E, F[i])$$
is an isomorphism for any integer $i.$} In other words, the extension
groups
$$\Ext_{\cXqc}^i(E,F):= \Hom_{D^+(\cXqc )}(E, F[i])$$
computed in the abelian category of quasicoherent sheaves of
${\mathcal O}_{\cX}$-modules over $\cX$ coincide
with the extension groups $\Ext^i_{{\mathcal O}_{\cX}}(E,F)$ computed in the
category
of all sheaves of ${\mathcal O}_{\cX}$-modules.

In particular, when $\cX$ is the
affine scheme $\Spec A$ defined by some noetherian commutative ring
$A,$
the extension group $\Ext^i_{{\mathcal O}_{\cX}}(E,F)$ may
be identified with the $A$-module
$\Ext^i_{A}(\Gamma(X,E),\Gamma(X,F)),$ computed in the abelian
category of $A$-modules.

\subsection{Groups of 1-extensions}

In this section, $(X, \cA)$ denotes as above a topological space
equipped with a sheaf $\cA$ of commutative ring, and $\cAmod$ the
abelian category of sheaves of $\cA$-modules over $X.$ Actually, the
constructions we shall now describe would still make sense with the
category $\cAmod$ replaced by a general abelian category, satisfying
suitable smallness assumptions.

Let us indicate that the approach to $\Ext^1$ groups in terms of
$1$-extensions described below
--- which directly inspired our definition of the arithmetic extension
group in terms of admissible extension of hermitian vector bundles ---
originates in the papers of Baer \cite{baer34} and Eilenberg--MacLane
\cite{eilenbergmaclane42}. See  \cite{maclane95}, p.103, notes of
Chapter III, for additional historical references.
\subsubsection{1-extensions}

Let $F$ and $G$ be two objects in $\cAmod.$ A 
\emph{1-extension} in $\cAmod$ (or
shortly, when no confusion may arise, an \emph{extension}) $\cE$ of
$F$ by $G$  is a short exact sequence in $\cAmod$ of the form
\begin{equation}
     \label{eq:ExtE}
     \cE: \,\, 0 \longrightarrow G \stackrel{i}{\longrightarrow} E
     \stackrel{p}{\longrightarrow} F \longrightarrow 0.
\end{equation}
A \emph{morphism of extensions} from the extension $\cE$ to a second
one
\begin{equation}
     \cE': \,\, 0 \longrightarrow G \stackrel{i'}{\longrightarrow} E'
     \stackrel{p'}{\longrightarrow} F \longrightarrow 0
     \end{equation}
     is a morphism of $\cA$-modules
     $$f:E \longrightarrow E'$$
     such that the diagram
     $$
     \begin{array}{ccccccccc}
0 &\longrightarrow  &G& \stackrel{i}{\longrightarrow}
&E&\stackrel{p}{\longrightarrow}  &F& \longrightarrow  & 0\\
  & & || & & \,\,\,\downarrow {\scriptstyle f} & &
||&\\
0 &\longrightarrow  &G& \stackrel{i'}{\longrightarrow}
&E'&\stackrel{p'}{\longrightarrow}  &F& \longrightarrow  & 0.
\end{array}
$$
is commutative. Such a morphism of extensions, if it exists, is an
isomorphism.
Moreover the isomorphism classes of extensions of $F$ by
     $G$ constitute a set, which we shall denote ${\mathbb E}{\rm
     xt}^1_{\cA}(F,G).$  The previous observation 
shows that any two extensions of $F$ by $G$
     related by a morphism of extensions define the same element in
     ${\mathbb E}{\rm
     xt}^1_{\cA}(F,G).$

     An extension (\ref{eq:ExtE}) of $F$ by $G$ is said to be
     \emph{split} or \emph{trivial} if it is isomorphic to the
extension
     $$0 \longrightarrow  G\stackrel{\binom{
{\rm id}_{G}}{0}}{\longrightarrow }
G \oplus F \stackrel{(0,{\rm id}_{F})}{\longrightarrow }
F\longrightarrow
0,$$
     or, equivalently, if $p$ admits a left 
inverse, or if $i$ admits a right inverse in 
$\cAmod.$

     \subsubsection{Pullback and pushout}\label{pbpo}

     Let $\cE$ an extension of $F$ by $G$ in $\cAmod$ as above
(\ref{eq:ExtE}) and $u:F'\longrightarrow F$
a morphism in $\cAmod.$ From $\cE$ and $u$, we may
define the $\cA$-module
$$E':=\Ker (E \oplus F' \stackrel{(p,-u)}{\longrightarrow } F)$$
and construct
a cartesian  diagram
\begin{equation}
    \begin{array}{cccc}
              & (e,f')  & \longmapsto & f'\\
           (e,f') & E' &
             \stackrel{\tilde{p}}{\longrightarrow}& F'  \\
            \Big\downarrow & \Big\downarrow\vcenter{%
\llap{ $\tilde{u}\,\,$ }}  &  &
           \Big\downarrow\vcenter{%
\rlap{$u$ }}
             \\
       e & E & \stackrel{p}{\longrightarrow}& F,
         \end{array}
         \label{pullback1}
  \end{equation}
  a morphism $\tilde{i}:G \rightarrow E'$ in $\cAmod$ defined by
  $\tilde{i}(g):=(i(g),0)$,
  and an extension $\cE \circ u$ of $F'$ by $G$:
  $$\cE \circ u:
  \,\, 0 \longrightarrow G \stackrel{\tilde{i}}{\longrightarrow}
  E'\stackrel{\tilde{p}}{\longrightarrow} F' \longrightarrow 0$$
By construction, the diagram
\begin{equation} \label{pullbackdiagram}
\begin{array}{ccccccccc}
0 &\longrightarrow  &G&\stackrel{\tilde{i}}{\longrightarrow} &E'
&\stackrel{\tilde{p}}{\longrightarrow}&F'&
\longrightarrow  & 0\\
  & & || & & \,\,\,\downarrow {\scriptstyle \tilde{u}} & &
\downarrow {\scriptstyle u}&\\
0 &\longrightarrow 
&G&\stackrel{i}{\longrightarrow}  &E& 
\stackrel{p}{\longrightarrow}
   &F&
\longrightarrow  & 0
\end{array}
\end{equation}
is commutative, and its right-hand square is cartesian. The pair
formed by the extension $\cE \circ u$ and the morphism $\tilde{u}$ are
characterized by these properties, and defines the \emph{pullback} of
$\cE$ by $u.$

For instance, the diagram
\begin{equation} 
\begin{array}{ccccccccc}
0 &\longrightarrow  &G&\stackrel{\binom{{\rm
id}_{G}}{0}}{\longrightarrow}
&G \oplus E \,
&\stackrel{\tilde{p}}{\longrightarrow}&F&
\longrightarrow  & 0\\
  & & || & & \;\;\;\downarrow {\scriptstyle (i, {\rm id}_{E})} & &
\downarrow {\scriptstyle p}&\\
0 &\longrightarrow 
&G&\stackrel{i}{\longrightarrow}  &E& 
\stackrel{p}{\longrightarrow}
   &F&
\longrightarrow  & 0
\end{array}
\end{equation}
establishes that the pullback $\cE \circ p$ is a split extension.

Symmetrically, if $v:G\longrightarrow
G'$ is a morphism in $\cAmod,$  we may
define the $\cA$-module
$$E':=\Coker (G \stackrel{\binom{i}{-v}}{\longrightarrow} E \oplus
G')$$
and construct
a cocartesian  diagram:
\begin{equation}
    \begin{array}{cccc}
   G &\stackrel{i}{\longrightarrow}& E & e  \\
            \Big\downarrow\vcenter{\llap{ $v\,\,$ }}& &
\Big\downarrow\vcenter{%
\rlap{ $\tilde{v}\,\,$ }}    &
           \Big\downarrow \\
       G' & \stackrel{\tilde{i}}{\longrightarrow}& E' &{[}(e,0)] \\
       g' & \longmapsto &{[}(0,g')],&
         \end{array}
         \label{pushout1}
  \end{equation}
  a morphism
  $\tilde{p}:E'\rightarrow F$ in $\cAmod$ defined by
   $\tilde{p}([(e,g')]):=p(e)$,
  and an extension $v \circ \cE$ of $F$ by $G'$:
  $$ v\circ\cE:
  \,\, 0 \longrightarrow G' \stackrel{\tilde{i}}{\longrightarrow}
  E'\stackrel{\tilde{p}}{\longrightarrow} F \longrightarrow 0.$$
By construction, the diagram
  $$
\begin{array}{ccccccccc}
0 &\longrightarrow 
&G&\stackrel{i}{\longrightarrow}  &E& 
\stackrel{p}{\longrightarrow}
   &F&
\longrightarrow  & 0
\\
  & &\downarrow {\scriptstyle v} & & \,\,\,\downarrow {\scriptstyle
  \tilde{v}} & &
|| &\\
0 &\longrightarrow  &G'&\stackrel{\tilde{i}}{\longrightarrow} &E'
&\stackrel{\tilde{p}}{\longrightarrow}&F&
\longrightarrow  & 0
\end{array}
$$
is commutative, and its left-hand square is cocartesian.
The pair
formed by the extension $v \circ \cE$ and the morphism $\tilde{v}$ are
characterized by these properties, and defines the \emph{pushout} of
$\cE$ by $v.$

The constructions of pullback and pushout of extensions satisfy
associativity properties. Namely, with the above notation,
  the extensions $v\circ (\cE \circ u)$ and $(v \circ \cE) \circ
u$ of $F'$ by $G'$ are canonically isomorphic. Moreover,
  if $u':F''\rightarrow F'$ (resp.  $v': G' \rightarrow G''$) is
another morphism in $\cAmod,$ the extensions $(\cE \circ u)\circ u'$
and $\cE \circ (u \circ u')$ (resp. $v'\circ (v \circ \cE)$ and
$(v'\circ
v )\circ \cE$) are canonically isomorphic.
This follows from the above characterization of pullback and
pushforward extensions, and legitimates the notation $v \circ \cE
\circ u,$ $\cE \circ u \circ u',$ and  $v'\circ v \circ \cE $ for
these extensions.

\subsubsection{Baer sum}\label{sec:BaerSum}

     The \emph{Baer sum} of two extensions of $F$ by $G$ in
     $\cAmod$
     $$\mathcal{E}_1:0\longrightarrow
G\stackrel{i_1}{\longrightarrow }
     E_1\stackrel{p_1}{\longrightarrow } F\longrightarrow  0\,\,\mbox{
     and }\,\,\mathcal{E}_2:0\longrightarrow
G\stackrel{i_2}{\longrightarrow }
E_2\stackrel{p_2}{\longrightarrow } F\longrightarrow  0$$
is the extension
$$\cE_{1}+\cE_{2}:\,0\longrightarrow  G\stackrel{i}{\longrightarrow }
E\stackrel{p}{\longrightarrow } F\longrightarrow  0 $$
where
\begin{equation}\label{algbaersum2}
E=\frac{\Ker((p_1,-p_2):E_1\oplus E_2 { \longrightarrow  } F)}
{{\rm Im}\,(\binom{i_1}{-i_2}: G {\longrightarrow  } E_1\oplus E_2)},
\end{equation}
and $p$ and $i$ are given by $$p[(e_1,e_2)]=p_1(e_1)=p_2(e_2)
\,\,\mbox{ and
}\,\,
i(g)=[(i_1(g),0))]=[(0,i_2(g))].$$

This construction defines a composition  law $+$ on ${\mathbb E}{\rm
     xt}^1_{\cA}(F,G)$. One easily checks that
     $({\mathbb E}{\rm
     xt}^1_{\cA}(F,G),+)$ is a commutative group in which the opposite
     of the class of an extension
     $$\mathcal{E} : 0\longrightarrow  G\stackrel{i}{\longrightarrow}
	E\stackrel{p}{\longrightarrow } F\longrightarrow  0$$
	is the class of
	$$\tilde{\mathcal{E}} : 0\longrightarrow 
G\stackrel{i}{\longrightarrow}
		E\stackrel{-p}{\longrightarrow } F\longrightarrow  0.$$
Actually, for any integer $k\geq 1$ and any $k$  extensions
$\cE_1,$ \ldots $\cE_k$ of $F$ by $G$ in $\cAmod,$
$$\mathcal{E}_j:0\longrightarrow  G\stackrel{i_j}{\longrightarrow }
     E_j\stackrel{p_j}{\longrightarrow } F\longrightarrow  0, \;\;\,\,
     1\leq j \leq k,$$the sum of their
classes in ${\mathbb E}{\rm
     xt}^1_{\cA}(F,G)$
     is the class of their Baer sum, which is defined as
\begin{equation}\label{algbaersumk}
     \cE_1 + \dots  +\cE_k: = \Sigma^k_G \circ (\cE_1 \oplus \dots
  \oplus \cE_k) \circ \Delta^k_F,
  \end{equation}
  where
  $$\Sigma^k_G: G^{\oplus k} \longrightarrow G$$
  denotes the ``sum" morphism,
and where $$\cE_1 \oplus \dots
  \oplus \cE_k : 0 \longrightarrow  G^{\oplus
  k}\stackrel{i_1\oplus\ldots\oplus i_k}{\longrightarrow }
     E_1\oplus \ldots \oplus 
E_k\stackrel{p_1\oplus \ldots\oplus 
p_k}{\longrightarrow }
     F^{\oplus k}\longrightarrow  0$$
     is the direct sum of the extensions $\cE_k,$ \ldots $\cE_k,$
     and
     $$\Delta^k_F: F \longrightarrow F^{\oplus k}$$
     the diagonal embedding of $F$ into $F^{\oplus k}.$ This
     definition is easily seen to coincide with the previous one when
     $k=2$. In general, for any integer $l$ such that $1< l < k,$ one
     sees that the extensions $(\cE_1 + \dots  +\cE_l) + (\cE_l +
\dots
     +\cE_k)$ and $\cE_1 + \dots  +\cE_k$ are isomorphic. This
establishes in particular the associativity of the Baer sum.

Using this expression for the Baer sum, it is straightforward to
derive  that it is compatible with composition. Namely, if $u: F'
\rightarrow F$ is a morphism in $\cAmod,$ the extensions $(\cE_1 +
\ldots
+\cE_k)\circ u$ and $\cE_1\circ u+ \ldots  +\cE_k\circ u$ are
canonically isomorphic. Similarly, if $v: G \rightarrow G'$ s a
morphism
in $\cAmod,$ the extensions $v \circ (\cE_1 + \ldots
+\cE_k)$ and $v \circ \cE_1 + \ldots
+v \circ \cE_k$ are canonically isomorphic.

\subsubsection{$\Ext^1$ and ${\mathbb E}{\rm xt}^1$}\label{extext}

Consider a 1-extension
$$  \cE: \,\, 0 \longrightarrow G \stackrel{i}{\longrightarrow} E
     \stackrel{p}{\longrightarrow} F \longrightarrow 0
     $$
     in $\cAmod$ as above.
     It may be seen as an exact sequence of complexes in $\cAmod$ (by
     identifying $E,$ $F,$ and $G$ with complexes concentrated in
     degree $0$), and defines a distinguished triangle
     $$G \stackrel{i}{\longrightarrow} E
\stackrel{p}{\longrightarrow}F\stackrel{\partial_{\cE}}{\longrightarrow}G[1]$$
     in $D(\cAmod)$ by the construction recalled in \ref{sign}.

     The arrow $\partial_{\cE}$ in
     $$\Hom_{D(\cAmod)}(F, G[1]) =: \Ext^1_{\cA}(F,G)$$
     will be denoted $\cl(\cE)$.

     {\footnotesize  Its definition boils down to the following.
     From $\cE,$ we  construct the cone of $i$ --- which simply is
      the complex in $\cAmod$
     $$C(i):= [ \ldots \longrightarrow 0 
\longrightarrow G \stackrel{i}{\longrightarrow} E
      \longrightarrow 0 \longrightarrow \ldots ],$$
      where $G$ (resp. $E$) is placed in degree $-1$ (resp. $0$) ---
and
      two morphisms of complexes of $\cA$-modules,
      $$\mathbf{p}: C(i) \longrightarrow F,$$
      defined by $p$ from $E$ to $F$ both placed in degree $0,$ and
      $$q: C(i) \longrightarrow G[1],$$ defined
      by the identity morphism of $G,$ placed in degree $-1.$ By the
      exactness of $\cE,$ $\mathbf{p}$ is a quasi-isomorphism and
      defines an isomorphism
      $$\mathbf{p}: C(i) \stackrel{\sim}{\longrightarrow} E \,\,\mbox{
      in } D(\cAmod).$$
      Then we have
$$\cl(\cE)= \partial_{\cE}:= - q \circ 
\mathbf{p}^{-1} \in \Ext^1_{\cA}(F,G):= 
\Hom_{D(\cAmod)}(F,G[1]).$$
}

In this way, one defines a map
$$\cl : {\mathbb E}{\rm xt}^1_{\cA}(F,G) \longrightarrow
\Ext^1_{\cA}(F,G).$$

Observe that, if $u:F' \longrightarrow F$ is a morphism in $\cAmod,$
then, from the diagram (\ref{pullbackdiagram}) defining $\cE \circ u,$
we derive a map of complexes
\[\scriptstyle
\tilde{u}\colon C(\tilde{i}):=[\,\ldots 
\,\longrightarrow\, 0\,\longrightarrow \,G\,
\stackrel{\tilde{i}}{\longrightarrow}\, E'\, \longrightarrow\, 0\,
\longrightarrow \,\ldots\,] \,\longrightarrow\, C(i)\,:=\,
[\,\ldots\, \longrightarrow\, 0\,\longrightarrow\, G\,
\stackrel{i}{\longrightarrow}\, E\, \longrightarrow\, 0\,
\longrightarrow\, \ldots\,]
\]
defined by ${\rm id}_{G}$ (resp. $\tilde{u}$) in degree $-1$ (resp.
$0$),
and a commutative diagram of complexes
$$
\begin{array}{ccccc}
     F' & \stackrel{\mathbf{\tilde{p}}}{\longleftarrow} &
C(\tilde{i}) &
     \stackrel{\tilde{q}}{\longrightarrow} & G  \\
     \hspace*{3mm}\downarrow {\scriptstyle u} &  & 
\hspace*{3mm}\downarrow {\scriptstyle
     \mathbf{u}} &  & ||  \\
     F &\stackrel{\mathbf{p}}{\longleftarrow} & C(i) &
     \stackrel{q}{\longrightarrow}  & G.
\end{array}$$
This shows that $\partial_{\cE \circ u} = \partial_{\cE}\circ u,$ or
in other terms, the equality $\cl(\cE \circ u) =\cl(\cE) \circ u$ in
$\Ext^1_{\cA}(F,G).$

Similarly, one shows that, for any morphism $v: G
\longrightarrow G'$ in $\cAmod,$ we have $\cl(v \circ \cE) =v \circ
\cl (\cE).$

Using the definition (\ref{algbaersumk}) of the Baer sum, it follows
that the map $\cl$ is indeed a morphism of abelian groups.

Observe
also that the class in ${\mathbb E}{\rm xt}^1_{\cA}(F,G)$ of the
extension $\cE$ lies in the kernel of $\cl$ precisely when
$\partial_{\cE}$ vanishes, that is when $\cE$ splits. (Indeed, the
vanishing of $\partial_{\cE}$ implies that the map
$p_{\ast}: \Hom_{\cA}(F,E) \rightarrow \Hom_{\cA}(F,F)$ is onto. Any
element in $p_{\ast}^{-1}({\rm id}_{F})$ then defines a splitting of
$\cE$.)
Consequently, the group morphism $\cl$ is injective.

Conversely, let $F$ and $G$ be any two objects in $\cAmod,$ and let
$$\iota: G \hookrightarrow I$$ be
a monomorphism from $G$ to an injective object in $\cAmod.$ Together
with the quotient map
$$\pi: I \longrightarrow \Coker \iota,$$
this monomorphism defines an extension
$$\cI:\,\,  0 \longrightarrow G \stackrel{\iota}{\longrightarrow} I
     \stackrel{\pi}{\longrightarrow} \Coker \iota \longrightarrow 0
     $$
     in ${\mathbb E}{\rm xt}^1_{\cA}(\Coker \iota, G).$
     The long exact sequence of $\Ext^._{\cA}(F, .)$'s deduced from
$\cI$
     starts as follows:
     \[\scriptstyle
     0\,\longrightarrow\, \Hom_{\cA}(F,G)\,
     \stackrel{\iota_{\ast}}{\longrightarrow} \,\Hom_{\cA}(F,I)\,
     \stackrel{\pi_{\ast}}{\longrightarrow}\, \Hom_{\cA}(F,\Coker
\iota)\,
     \stackrel{\partial}{\longrightarrow}\, \Ext_{\cA}^1(F,G)\,
     \longrightarrow\, \Ext_{\cA}^1(F,I)\,=\,\{0\}.
     \]
In particular, the boundary map
$$\partial: \Hom_{\cA}(F,\Oker \iota)\longrightarrow
\Ext_{\cA}^1(F,G)$$
induces an isomorphism of abelian groups
$$j_{\iota}: \Oker \pi_{\ast} 
\stackrel{\sim}{\longrightarrow} 
\Ext_{\cA}^1(F,G).$$
Since $\partial$  coincides with the composition
with  $\cl(\cI) \in \Ext^1_{\cA}(\Coker
\iota, G)$, the isomorphism $j_{\iota}$ maps the class in
$\Oker \pi_{\ast}$ of an element $u \in \Hom_{\cA}(F,I)$ to
\begin{equation}\label{comp}
\cl(\cI)\circ u = \cl(\cI \circ u).
\end{equation}
Besides, as the $1$-extension $\cI \circ \pi$ is split (\cf \ref{pbpo}
above), one may define a  morphism of abelian groups
$$j'_{\iota}: \Oker \pi_{\ast} \stackrel{\sim}{\longrightarrow}
{\mathbb E}{\rm
xt}_{\cA}^1(F,G)$$ by sending
the class in
$\Oker \pi_{\ast}$ of an element $u \in \Hom_{\cA}(F,I)$ to the class
of $\cI \circ u,$ and (\ref{comp}) asserts that this morphism
satisfies
$$\cl \circ j'_{\iota} = j_{\iota}.$$

This shows that the map $\cl$ is onto --- 
actually this constructs an inverse of $\cl$ --- 
and concludes the proof of the
following:

\begin{proposition}\label{extisext}
     The map
$$\cl : {\mathbb E}{\rm xt}^1_{\cA}(F,G) \longrightarrow
\Ext^1_{\cA}(F,G)$$
is an isomorphism from the group ${\mathbb E}{\rm xt}^1_{\cA}(F,G)$
equipped with the Baer sum onto the ``cohomological" extension group
$\Ext^1_{\cA}(F,G).$ 

\end{proposition}

\subsection{Extension groups of holomorphic vector bundles}

Let $X$ be a paracompact complex analytic manifold, and ${\mathcal O}_{X}$ the
sheaf of
$\C$-analytic functions on $X.$

\subsubsection{Dolbeault isomorphisms}\label{Dolbeault}

For any ${\mathcal O}_X$-module $F$, we obtain a quasi-isomorphism of
complexes of ${\mathcal O}_{X}$-modules
\begin{equation}
{\mathcal Dolb}_{F}: F \longrightarrow  {\mathcal 
Dolb}(F):=(F \otimes _{{\mathcal O}_X} A^{0,.}_X,
\overline{\partial}_F)
\label{DolbSheaf}
\end{equation}
--- the \emph{Dolbeault resolution of} $F$ ---
by applying the functor $F\otimes _{{\mathcal O}_X}.$ to the Dolbeault
resolution
$$0 \rightarrow {\mathcal O}_X \hookrightarrow A^{0,0}_X 
\stackrel{ \overline{\partial}}{\rightarrow}
A^{0,1}_X \rightarrow \cdots \rightarrow A^{0,i}_X
\stackrel{ \overline{\partial}}{\rightarrow} A^{0,i+1}_X \rightarrow
\cdots $$
of the sheaf ${\mathcal O}_{X}$ of holomorphic functions on $X.$
This construction is functorial and natural: to any morphism of ${\mathcal O}_X$-modules 
$\phi: F_1 \rightarrow F_2$, we may attach the morphism of Dolbeault complexes:
$${\mathcal Dolb}(\phi):= \phi \otimes_{{\mathcal O}_{X}}
{\rm Id}_{A_X^{0,.}}:= 
{\mathcal Dolb}(F_{1})
\longrightarrow {\mathcal Dolb}(F_{2}),$$
which fits into a commutative diagram:
$$
\begin{array}{ccc}
     F_{1} & \stackrel{{\mathcal Dolb}_{F_{1}}}{\longrightarrow} &
{\mathcal Dolb}(F_{1})  \\
     \hspace*{3mm}\downarrow {\scriptstyle \phi} &  & 
\hspace*{3mm}\downarrow {\scriptstyle
     {\mathcal Dolb}(\phi)} \\
     F_{2} &\stackrel{{\mathcal Dolb}_{F_{2}}}{\longrightarrow} &{\mathcal
     Dolb}(F_{2}).
\end{array}$$
Actually, the functor ${\mathcal Dolb}$ may be defined on bounded
below complexes of ${\mathcal O}_{X}$-modules, and is exact (by the flatness of 
$A^{0,.}_{X}$ over ${\mathcal O}_{X}$; \cf \ref{arbach}).

The sheaves
$F \otimes _{{\mathcal O}_X} A^{0,i}_X$ are sheaves of $A^{0,0}_X$-modules,
hence acyclic. Consequently the cohomology groups $H^i(X,
{\mathcal Dolb}(F))$ may be identified with the Dolbeault cohomology
groups $\HD^i(X,F),$ defined as the cohomology groups of the complex
of $\C$-vector spaces:
\[
\scriptstyle
0\, \rightarrow\, \Gamma(X,    F \otimes_{{\mathcal O}_X} A^{0,0}_X)\,
\stackrel{ \overline{\partial}_F}{\rightarrow} \,
\Gamma (X, F \otimes_{{\mathcal O}_X} A^{0,1}_X)\, \rightarrow\, \ldots\,
\rightarrow \,
\Gamma(X, F \otimes_{{\mathcal O}_X} A^{0,i}_X)\,
\stackrel{ \overline{\partial}_F}{\rightarrow}\, \Gamma(X,  F
\otimes_{{\mathcal O}_X}
A^{0,i+1}_X)\,\rightarrow\,
\ldots\, .
\]
Thanks to this identification, the isomorphisms of $\C$-vector spaces
$$H^i(X,F) \longrightarrow  H^i(X,{\mathcal Dolb}(F))$$
deduced from the quasi-isomorphism (\ref{DolbSheaf}) define the
\emph{Dolbeault isomorphisms}:
$${\rm Dolb}_F: H^i(X,F) \longrightarrow \HD^i(X,F).$$

More generally, let $E$ be a locally free ${\mathcal O}_X$-module of finite
rank  (or,
equivalently, the sheaf of $\C$-analytic sections of some
$\C$-analytic vector bundle over $X$). Then, for any ${\mathcal O}_X$-module
$F,$ the composition of the isomorphism 
(\ref{eq:extisoB}) (where $\cA={\mathcal O}_X$) and
${\rm Dolb}_{\check{E}\otimes F}$ defines an isomorphism of
$\C$-vector
spaces
$${\rm Dolb}_{E,F}: \Ext^i_{{\mathcal O}_X} (E,F)
\stackrel{\sim}{\longrightarrow} \HD^i(X, \check{E}\otimes F),$$
which coincides with ${\rm Dolb}_F$ when $E={\mathcal O}_{X}$.

\subsubsection{Second fundamental form}\label{secondfundamental}

Consider a short exact sequence of $\C$-analytic vector bundles over
$X$, or equivalently, of locally free ${\mathcal O}_X$-modules of finite rank:
$$\cE: \,\,\, 0 \longrightarrow G \stackrel{i}{\longrightarrow} E
\stackrel{p}{\longrightarrow}F \longrightarrow 0.$$
Let $s$ be a $C^\infty$-splitting of $\cE$, namely a
$C^\infty$-section
of $p.$
We may see it as an element of
$$\Gamma(X,   \check{F} \otimes_{{\mathcal O}_X} E\otimes_{{\mathcal O}_X}A^{0,0}_X)
\simeq \Hom_{A^{0,0}_X}(F \otimes_{{\mathcal O}_X} A^{0,0}_X, E\otimes_{{\mathcal O}_X}
A^{0,0}_X),$$ and consider its image
  by the Dolbeault operator
  $$\overline{\partial}_{\check{F}\otimes E}\, s \in
  \Gamma(X,\check{F} \otimes_{{\mathcal O}_X} E \otimes_{{\mathcal O}_X} A^{0,1}_X)
\simeq \Hom_{A^{0,0}_X}(F\otimes_{{\mathcal O}_X} A^{0,0}_X, E\otimes_{{\mathcal O}_X}
A^{0,1}_X).$$
The relation $p \circ s ={\rm id}_E$ implies that $(p
\otimes_{{\mathcal O}_X} {\rm id}_{A^{0,1}})
\overline{\partial}_{\check{F}\otimes E} s $
vanishes. Consequently there exists a unique $\alpha$ in
$$\Gamma(X,  \check{F} \otimes_{{\mathcal O}_X} G \otimes_{{\mathcal O}_X} A^{0,0}_X )
\simeq \Hom_{A^{0,0}_X}( F \otimes_{{\mathcal O}_X} A^{0,0}_X, G \otimes_{{\mathcal O}_X}
A^{0,0}_X)$$
such that
$$ \overline{\partial}_{\check{F}\otimes E}\, s
=(i \otimes_{{\mathcal O}_X} {\rm id}_{A^{0,1}_X})\,\alpha.$$
Moreover $s$, and therefore $\alpha$, is
$\overline{\partial}$-closed, namely
$\overline{\partial}_{\check{F}\otimes G}\,\alpha$ vanishes in
$\Gamma(X,\check{F} \otimes_{{\mathcal O}_X} G\otimes_{{\mathcal O}_X}A^{0,0}_X).$

When $s$ is the orthogonal splitting of $\cE$ deduced from some
$C^\infty$-hermitian metric $\|.\|$ on $E,$ the $(0,1)$-form $\alpha$
with coefficients in $\check{F}\otimes G$ is the adjoint of the
so-called
second
fundamental form associated to the extension of hermitian
vector bundles on $X$ defined by $\cE$ and $\|.\|$ (\cf
\cite{griffiths66}, VI.3; see also \cite{griffiths69}, 2.d-e, and
\cite{griffithsharris78}, pp. 72-73). In our more general context, we
shall  call the $(0,1)$-form $\alpha$ itself the \emph{second
fundamental form}
attached to the extension $\cE$ and its $C^\infty$-splitting
$s$.

Recall that the extension $\cE$ of analytic vector bundles over $X$
defines a class $\cl(\cE)$ in $\Ext^1_{{\mathcal O}_{X}}(F,G),$ and that the
$\overline{\partial}_{\check{F}\otimes G}$-closed $(0,1)$-form
$\alpha$
in $A^{0,1}(X, \check{F}\otimes G):= \Gamma(X, \check{F}\otimes G
\otimes_{{\mathcal O}_{X}}A^{0,1}_{X})$ defines a class
$[\alpha]$ in $\HD^1(X, \check{F}\otimes
G)$.

The following proposition is classical (see for instance
\cite{griffiths66}, Proposition p.422), up to the precise
determination of the sign.
\begin{proposition}\label{prop:2f}
     With the above notation, the following equality holds in
     the Dolbeault cohomology group $\HD^1(X,
     \check{F}\otimes G)$:
     $${\rm Dolb}_{F,G}(\cl(\cE))=[\alpha].$$
\end{proposition}

In other words, the image of $\cl(\cE)$ under the composition of
canonical isomorphisms
$$
\begin{array}{ccl}
     \Ext^1_{{\mathcal O}_{X}}(F,G) & = & \Hom_{D({\mathcal O}_{X}-\mathbf{mod})}(F, G[1])  \\
      & \simeq & \Hom_{D({\mathcal O}_{X}-\mathbf{mod})}(F, {\mathcal Dolb}(G)[1])   \\
      & \simeq & \Hom_{D({\mathcal O}_{X}-\mathbf{mod})}({\mathcal O}_{X}, 
\check{F}\otimes_{{\mathcal O}_{X}}{\mathcal Dolb}(G)[1]) 
\\
    & \simeq & \Hom_{D({\mathcal O}_{X}-\mathbf{mod})}({\mathcal O}_{X}, {\mathcal 
Dolb}(\check{F}\otimes_{{\mathcal O}_{X}}G)[1]) \\
      & \simeq & \HD^1(X,\check{F}\otimes_{{\mathcal O}_{X}}G)
\end{array}
$$
which defines ${\rm Dolb}_{F,G}$ coincides with the class of $\alpha.$
This follows from the following two straightforward lemmas:

\begin{lemma}\label{CDolb}
By applying the functor ${\mathcal Dolb}$ to ${\mathcal E}$, we obtain
an exact sequence of complexes of sheaves of ${\mathcal O}_{X}$-modules:
$${\mathcal Dolb}(\cE):
\,\,\, 0 \longrightarrow {\mathcal Dolb}(G) \stackrel{{\mathcal
Dolb}(i)}{\longrightarrow} {\mathcal Dolb}(E)
\stackrel{{\mathcal Dolb}(p)}{\longrightarrow}{\mathcal Dolb}(F) \longrightarrow 0.$$
Moreover the cone $C({\mathcal Dolb}(p))$ may be identified with ${\mathcal
Dolb}(C(p))$, and the following diagram is commutative:
$$
\begin{array}{rccccc}
 \partial_{\mathcal E}:&F & \stackrel{\mathbf{\pp}}{\longleftarrow} &
C(i) &
     \stackrel{(-{\rm Id}_{G},0)}{\longrightarrow} & G [1]  \\ &
     \hspace*{6mm}\downarrow {\scriptstyle {\mathcal Dolb}_{F}} &  & 
\hspace*{6mm}\downarrow {\scriptstyle
     {\mathcal Dolb}_{C(i)}} &  & \hspace*{8mm}\downarrow {\scriptstyle
     {\mathcal Dolb}_{G}[1]}  \\
\partial_{{\mathcal Dolb}({\mathcal E})}: &   {\mathcal Dolb}(F) &\stackrel{{\mathcal
Dolb}(\pp)}{\longleftarrow} & C({\mathcal Dolb}(i)) &
     \stackrel{(-{\rm Id}_{{\mathcal Dolb}(G)},0)}{\longrightarrow}  & 
{\mathcal Dolb}(G)[1].
\end{array}$$
\end{lemma}

\begin{lemma}\label{w}
    One defines a morphism of complexes of ${\mathcal O}_{X}$-modules
    $$w: {\mathcal Dolb}(F) \rightarrow C({\mathcal Dolb}(i))$$
    that is a right inverse of ${\mathcal Dolb}(\pp)$ by letting, for
    any local section $\beta$ of ${\mathcal Dolb}(F)^i:= F
    \otimes_{{\mathcal O}_{X}}A^{0,i}_{X}$,
    $$w(\beta):=(-\alpha.\beta, s.\beta).$$
\end{lemma}

Indeed Lemma 
\ref{CDolb} and \ref{w} show that  the image 
in $\Hom_{D({\mathcal O}_{X}-\mathbf{mod})}(F, {\mathcal Dolb}(G)[1])$ of $\cl(\cE)$ is defined by the
composition:
$$F \stackrel{{\mathcal Dolb}_{F} }{\longrightarrow} {\mathcal Dolb}(F)
\stackrel{w}{\longrightarrow} C({\mathcal Dolb}(i)) 
\stackrel{(-{\rm Id}_{{\mathcal Dolb}(G)},0)}{\longrightarrow}  
{\mathcal Dolb}(G)[1].$$
It is immediate that it corresponds in $\HD^1(X,\check{F}\otimes_{{\mathcal O}_{X}}G)$ to the class 
 of
$\alpha$.

\section{Lattices of Voronoi's first kind}

\subsection{Selling parameters}
We recall the definition of lattices of Voronoi's 
first kind and their associated
Selling parameters.

\subsubsection{Definitions}
Let $\overline E = \bigl(E,\Vert \cdot\Vert\bigr)$ be an euclidean
lattice, of positive rank $n$.
It is said to be a {\it lattice of Voronoi's first kind} if it
possesses
what Conway and Sloane (\cite[§2]{conwaysloane92}) call an {\it
obtuse
superbase}, namely if there exists a $n+1$-tuple
  $(v_0,\ldots,\,v_n)$ of vectors in $E$ such that
\begin{itemize}
\item[i)]
$(v_1,\ldots,\,v_n)$ is a $\Z$-basis of $E$, and
\[v_0+\cdots + v_n = 0\]
(this is a {\it superbase});
\item[ii)]
for any $(i,j) \in \{0,\ldots ,\,n\}^2, \ i\not= j$,
\begin{equation}
\label{gl1.1}
p_{ij} := -v_i\cdot v_j \geq 0\end{equation}
(the angle of the vectors $v_i$ and $v_j$ is {\it obtuse}).
\end{itemize}

When the inequalities (\ref{gl1.1}) are strict, the superbase
$(v_0,\ldots,\,v_n)$ is called
{\it strictly obtuse}.

Observe that, for any superbase $(v_0,\ldots,\, v_n)$ of $E$, the
$n(n+1)/2$-coefficients $$p_{ij} = p_{ji},\ 0\leq i< j\leq n,$$
defined by (\ref{gl1.1}) uniquely determine the euclidean structure
of $\overline E$.
Indeed, for any $i\in \{0,\ldots,\,n\}$ we have
\[
\Vert v_i\Vert^2
= - v_i\cdot \sum_{\genfrac{}{}{0pt}{}{0\leq k\leq n}{ k\not= i}} v_k
= \sum_{\genfrac{}{}{0pt}{}{0\leq k\leq n} {k\not= i}} p_{ik},
\]
and consequently, for any $(x_0,\ldots,\, x_n) \in \R^{n+1}$,
\begin{equation}
\label{gl1.2}
\Big\Vert \sum^n_{i=0} x_i v_i\big\Vert^2 = \sum_{0\leq i < j\leq n}
p_{ij}(x_i-x_j)^2.
\end{equation}
This formula goes back to Selling \cite{selling74}, at least when
$n=3$, and the coefficients
$$(p_{ij})_{\ 0\leq i < j\leq n}$$
will be called the {\it Selling parameters} attached to the superbase
$(v_0,\ldots,\,v_n)$.

Selling's formula (\ref{gl1.2}) shows in particular that a superbase
$(v_0,\ldots,\,v_n)$ of an
euclidean lattice $\overline E = (\overline E,\Vert \cdot \Vert)$ is
obtuse iff the quadratic
form $\Vert \cdot \Vert^2$ on $E_\R$ expressed in the base
$(v_1,\ldots,\,v_n)$, takes the form
\[\sum^n_{k=1} \lambda_k X^2_k + \sum_{1\leq i<j\leq n}
\lambda_{ij}(X_i-X_j)^2\]
for some $\lambda_i,\ \lambda_{ij}$ in $\R_+$ or equivalently, iff
the matrix
$(a_{ij})_{1\leq i,j\leq n} := (v_i\cdot v_j)_{1\leq i,j\leq n}$ of
this
quadratic form satisfies
\[\sum^n_{j=1} a_{ij} \geq 0 \ \mbox{ for any }\ i\in
\{1,\ldots,\,n\}\]
and
\[
a_{ij} = a_{ji} \leq 0 \ \mbox{ if }\ 1\leq i < j\leq n.
\]
(Indeed, with the above notation,
\[
p_{ij} = \lambda_{ij} = - a_{ij} \ \mbox{ if }\ 1\leq i < j\leq n
\]
and
\[
p_{0i} = \lambda_i = \sum^n_{j=1} a_{ij}.)
\]
These quadratic forms are precisely the ones in the domain associated
by Voronoi
to the ``forme parfaite principale''
\begin{equation}
\label{gl1.3}
\varphi := \sum_{1\leq i\leq j\leq n} X_i X_j
\end{equation}
(\cf \cite[29]{voronoi08a}).

Observe also that an euclidean lattice $\overline E$ of rank $n$ is
of Voronoi's first kind iff
there exists a basis $(\xi_1,\ldots,\,\xi_n)$ of the $\Z$-module
$E^\lor$ and a family
\[
(p_{ij})_{0\leq i< j \leq n} \in  \R^{\frac{1}{2}\cdot n(n+1)}_+
\]
such that, for any $x\in E_\R$,
\begin{equation}
\label{gl1.4}
\Vert x\Vert^2 = \sum_{0\leq i<j\leq n} p_{ij} \bigl(\xi_i(x) -
\xi_j(x)\bigr)^2
\end{equation}
where $\xi_0:= 0$. Indeed this identity is equivalent to the fact
that the superbase
$(v_0,v_1,\ldots,v_n)$ of $E$, defined by the dual basis
$(v_1,\ldots,\,v_n)$ of
$(\xi_1,\ldots,\,\xi_n)$ and $v_0:= - v_1-\cdots - v_n$, satisfies
(\ref{gl1.2}).
Any $\Z$-basis $(\xi_1,\ldots,\,\xi_n)$ of $E^\lor$ satisfying
(\ref{gl1.4}) will be said to
be {\it adapted} to the euclidean lattice $\overline E$ of Voronoi's
first kind
and the (unique) family $(p_{ij})_{\ 0\leq i < j\leq n}$ the {\it
associated Selling parameters.}

For any subset $S$ of $\bigl\{(i,j), \ 0\leq i < j\leq n\bigr\}$, let
$\gamma(S)$ be the (non-oriented)
graph whose set of vertices is $\{0,\ldots,\,n\}$, with an edge
between any $i$ and
$j$ such that $\,0\leq i < j\leq n$ iff $(i,j) \in S$. The following
lemma is left as an easy exercise for the reader:

\begin{lemma}\label{abs1.1}
Let $\overline E$ be an euclidean lattice of Voronoi's first kind,
$(\xi_1,\ldots,\,\xi_n)$ an adapted $\Z$-basis of $E^\lor$ and
$(p_{ij})_{0\leq i< j\leq n}$
the corresponding Selling parameters, defined by $(\ref{gl1.4})$
where $\xi_0:= 0$.

1) If the euclidean lattice $\overline E $ is indecomposable ---
\emph{namely, if it is not $($isomorphic to$)$
the direct sum of two euclidean lattices of positive rank} --- then
the graph $\gamma(S)$ attached to
\[
S:= \bigl\{(i,j)\big\vert p_{ij} \not= 0\bigr\}
\]
is connected.

2) For any subset $S'$ of $\bigl\{(i,j),\ 0\leq i < j\leq n\bigr\}$ such
that $\gamma(S')$ is connected, we have
\[
E^\lor = \sum_{(i,j)\in S'} \,\Z(\xi_i-\xi_j).
\]
\end{lemma}

\subsection{Examples}\label{abs2.1}
We now describe remarkable classes of euclidean lattices which are
and which are not of Voronoi's first kind.
\subsubsection{}
Any euclidean lattice of rank 2 is of Voronoi's first kind. Indeed,
if $(v_1,v_2)$ is
a base of such a lattice which is obtuse (i.e., such that $v_1\cdot
v_2\leq 0)$ and
reduced in the sense of Lagrange (i.e., which satisfies
$\Vert v_1\Vert \leq \Vert v_2 \Vert \leq \Vert v_2± v_1\Vert\,)$,
then the superbase $(-v_1-v_2, v_1, v_2)$ is obtuse.

\subsubsection{}\label{abs2.2}
Any euclidean lattice of rank 3 is of Voronoi's fist kind. This is a
classical result of
Selling  \cite{selling74}, which has been reproved by Voronoi
\cite[33]{voronoi08a} as a
consequence of his theory of ``Voronoi's reduction'' associated to
perfect forms.
A direct argument appears in \cite{conwaysloane92}, and may be
concisely
reformulated as follows.

To any superbasis $(v_0, v_1, v_2, v_3)$ of an euclidean lattice of
rank 3, we may
attach

\begin{equation}
\label{gl2.1}
\begin{array}{rrl}
  N(v_0, v_1, v_2, v_3) & := &\sum\limits_{S\subset \{0,1,2,3\}}
\,\Vert \sum\limits_{k\in S}\,v_k\Vert^2 \\
  & = & 2\, \bigl[ \Vert v_1\Vert^2 + \Vert v_2\Vert^2+ \Vert
v_3\Vert^2 + \Vert v_1+v_2\Vert^2 + \Vert v_2+v_3\Vert^2 \\
  & & + \Vert v_3 + v_1\Vert^2 + \Vert v_1 + v_2 + v_3\Vert^2\bigr].
\end{array}
\end{equation}
If we let, as above, $p_{ij} := -v_i \cdot v_j$ if $0\leq i < j\leq
3$,
a straightforward
computation shows that:
\[
N(v_0, v_1, v_2, v_3) = 8 \sum_{0\leq i< j\leq 3}\, p_{ij}.
\]
There exists a superbasis $(v_0, v_1, v_2, v_3)$ such that
$N(v_0, v_1, v_2, v_3)$ is minimal --- this follows from its very
definition
{\rm (\ref{gl2.1})} --- and any such ``minimal'' superbasis is
obtuse.
Indeed, if some $p_{ij}$, say $p_{01}$, were negative, then the
superbasis
\[(v'_0, v'_1, v'_2, v'_3) := (-v_0, v_1, v_0 + v_2, v_0 + v_3)\]
would satisfy

\begin{eqnarray*}
N(v'_0, v'_1, v'_2, v'_3) - N(v_0, v_1, v_2, v_3) & = & 2\,
\bigl(\Vert -v_0 + v_1\Vert^2 - \Vert v_0 + v_1\Vert^2\bigr) \\
& = & 8 \, p_{01} < 0.
\end{eqnarray*}

A similar argument shows that a superbase $(v_0, v_1, v_2)$ of an
euclidean lattice
of rank 2 is obtuse if it minimizes $\Vert v_0\Vert^2 + \Vert
v_1\Vert^2 + \Vert v_2\Vert^2$.

\subsubsection{}\label{abs2.3}
For any positive integer $n$, the euclidean lattice $A_n$,
defined in section \ref{subsec.rootlattice}, 
admits an obtuse superbase defined by
\[v_i := e_i-e_{i+1},\ 0\leq i\leq n,\]
where $(e_0, \ldots,\, e_n)$ is the standard basis of $\Z^{n+1}$ and
$e_{n+1} := e_0$.

The dual euclidean lattice $A^*_n$ may be identified with the lattice
in the real vector space
\[A^*_{n,\R} := A_{n,\R} = \bigl\{(x_0,\ldots,\, x_n) \in \R^{n+1}
\big\vert \sum^n_{k=0}
\,x_k=0\bigr\},\]
equipped with the restriction of the standard euclidean scalar
product on
$\R^{n+1},$ that is defined by

\begin{eqnarray*}
A^*_n & := & \bigl\{x\in A_{n,\R} \,\vert\, \forall\, y\in A_n,\
x\cdot y\in \Z\bigr\} \\
& \, = & \bigl\{ (x_0,\ldots , \, x_n) \in \R^{n+1} \big\vert
\sum^n_{k=0}\, x_k=0 \mbox{ and} \\
& & \qquad\qquad\qquad \qquad\qquad \forall\,(i,j) \in
\{0,\ldots,\,n\}^2,\ x_j - x_i \in \Z\bigr\}.
\end{eqnarray*}
One easily checks that the vectors
\[
v_i:= e_i - \frac{1}{n+1}\,\sum^{n+1}_{k=0} e_k\,,\,\,\,\ (0\leq i\leq
n),
\]
constitute a superbase of $A^*_n$.
It is strictly obtuse, since $v_i\cdot v_j = -1/(n+1)$ if $i\not= j$.

This shows that $A_n$ and $A^*_n$ are euclidean lattices of Voronoi's
first kind.
This property of $A^*_n$ also follows from Voronoi's theory
\cite{voronoi08a}.
Indeed, as observed above, an euclidean lattice $\overline E$ is of
Voronoi's first kind
iff it is isometric to a lattice $(\Z^n, \Vert \cdot \Vert_\Psi)$
where
$\Vert \cdot \Vert_\Psi := \Psi^{1/2}$ is the euclidean norm on
$\R^n$
defined by a quadratic form $\Psi$ on $\R^n$ in the polyhedral domain
associated by Voronoi's ``reduction of the first kind"
(\cite[15]{voronoi08a}) to the perfect form
$\varphi$ defined in {\rm (\ref{gl1.3})}. The euclidean lattice
$(\Z^n, \,\Vert \cdot \Vert_\varphi)$, up to a scaling, is isomorphic
to $A_n$.
(Indeed the isomorphism
\[
\begin{array}{ccccc}
\pr\colon A_n& \hookrightarrow & \Z^{n+1} & \longrightarrow & \Z^n \\
&&(x_i)_{0\leq i\leq n} & \longmapsto & (x_i)_{1\leq i\leq n}
\end{array}\]
satisfies $\Vert x\Vert^2_{A_n} = 2 \,\Vert \pr_\R (x)
\Vert^2_\varphi$ for any $x \in A_{n,\R}$.)
As the perfect lattice $A_n$ --- or equivalently, the perfect form
$\varphi$ --- is extreme
\footnote{Recall from \cite[3.4.6]{martinet03} that
a lattice is extreme iff it is perfect and eutactic and that the
irreducible root lattices $A_n$, $D_n$, $E_6$, $E_7$ and $E_8$ are
extreme
\cite[4.7.2]{martinet03}.}, it is eutactic,
\ie the adjoint
form\footnote{Namely up to a scaling, the quadratic form on $\R^n$
whose matrix is the inverse of the
one of $\varphi$.}
$\widetilde\varphi$ belongs to the interior of the polyhedral domain
attached to $\varphi$ (\cf {\cite[17]{voronoi08a} and
\cite{coxeter51}).
Since, up to a scaling, $(\Z^n,\,\Vert
\cdot\Vert_{\widetilde\varphi})$ is isometric
to $A^*_n$, this establishes again that $A^*_n$ is of Voronoi's first
kind.
Actually, the fact the $\widetilde\varphi$ lies in the {\it interior}
of the
domain attached to $\varphi$ shows that $A^*_n$ admits a {\it
strictly} obtuse superbase.

Since the Voronoi domains attached to two non-proportional perfect
forms
on $\R^n$ meet only along a common face (\cite[20]{voronoi08a}),
we see similarly that the only extreme form $\Psi$ on $\R^n$ whose
adjoint
form $\widetilde\Psi$ belongs to the domain of $\varphi$ is $\varphi$
itself, or a
multiple of $\varphi$.
  In other words,  {\it if an euclidean lattice
$\overline E$ of rank $n$ is extreme and if its dual $E^\lor$ is of
Voronoi's first kind, then, up to a scaling, $\overline E$ is
isometric to $A_n$. }

In particular, the dual root lattices $D^*_n (n\geq 4),\ E^*_6,\,
E^*_7$ and
$E^*_8$ are not of Voronoi's first kind.
Consequently the euclidean lattice $E_8$ (resp. $D_n$) which is
isometric
to $E_8^*$ (resp. to $D^*_n$, up to a scaling) is {\it not} of
Voronoi's first kind.

\subsubsection{}\label{abs2.4}
From Voronoi's theory, one may also derive some ``upper bound'' on
automorphism groups
of euclidean lattices of Voronoi's first kind:

\begin{proposition}\label{prop2.1}
Let $\overline E$ be an euclidean lattice of rank $n$, of Voronoi's
first kind,
and let $(\xi_1,\ldots,\,\xi_n)$ be a $\Z$-basis of $E$ adapted to
$\overline E$.

If the $($finite$)$ automorphism group $G:= {\rm Aut}\,\overline E$
of $\overline E$ acts
irreducibly on $E_\R$, then there exists $\lambda\in \R^*_+$ and a
subset $S$ of
$\bigl\{(i,j), \ 0\leq i < j\leq n\bigr\}$ such that, for any $x\in
E_\R$,
\[
\Vert x\Vert^2 = \lambda \sum_{(i,j)\in S} \,\bigl( \xi_i(x) -
\xi_j(x)\bigr)^2,
\]
where as above $\xi_0:= 0$. Moreover, the contragredient action of
$G$ on $\check E$
permutes transitively the pairs of vectors $±\,(\xi_i-\xi_j),\
(i,j) \in S$.
\end{proposition}

\proof
We may assume that $\overline E$ is $\Z^n$ and $\Vert \cdot\Vert^2$ a
quadratic form $\Psi$ on $\R^n$ in Voronoi's principal domain --- in
other words, $\xi_{0}=0,$ $\xi_{1}=X_{1},$\ldots,$\xi_{n}=X_{n}.$ This
form may be written
\[\Psi = \sum_{0\leq i< j\leq n} \,p_{ij}\,q_{ij},\]
where $(q_{ij})_{0\leq i< j\leq n}$ is the basis of the space
$S^2\,\R^{n\lor}$ of
quadratic form on $\R^n$ defined by
\[
q_{ij}:=\left\{
\begin{array}{lcr}
X^2_j & \mbox{if}& i=0\\
(X_i-X_j)^2 & \mbox{if} & i\geq 1,
\end{array}
\right.
\]
and where the Selling parameters $p_{ij}$ are non-negative. The
minimal face of the first Voronoi
decomposition of the cone of non-negative quadratic forms on $\R^n$
which contains $\Psi$ is the
simplicial cone
\[\sum_{(i,j)\in S} \,\R_+\, q_{ij}\]
where $S:= \bigl\{(i,j)\big\vert\, p_{ij} \not= 0\bigr\}$. Any
subgroup $G$ of $GL_n(\Z)$
preserving the form $\Psi$ also preserves this cone, and permutes its
extremal half-lines
$(\R_+\,q_{ij})_{(i,j)\in S}$ and consequently the forms
$(q_{ij})_{(i,j)\in
S},$ since the
action of $G$ preserves the integral structure of quadratic forms.

When the action of $G$ on $\R^n$ is irreducible, the form $\Psi$ is
--- up to scaling --- the
unique element of $S^2\,\R^{n\lor}$ which is $G$-invariant.
Consequently the
non-zero Selling
parameters 
are necessarily equal --- this
establishes the existence of $S$ and $\lambda$ --- and the action of
$G$ permutes transitively the quadratic forms
$q_{ij}:=(\xi_{i}-\xi_{j})^{2},$ $(i,j)\in S,$ or equivalently, the
pairs of vectors $±\,(\xi_i-\xi_j),\
(i,j) \in S$.
\qed

Observe that, with the notation of Proposition \ref{prop2.1}, the
linear forms
\[
\xi_i-\xi_j,\ (i,j)\in S
\]
are primitive vectors of $E^\lor$, and that their
squares $(\xi_i-\xi_j)^2,\ (i,j)\in S$, are rank 1 quadratic forms
which span a
$G$-invariant subspace of $S^2 \check E_\R$, containing the euclidean
form
$\Vert \cdot\Vert^2$ defining the euclidean structure of $\overline
E$.

Proposition \ref{prop2.1} and these observations may be used to show
that
some lattice with ``big'' automorphism groups are {\it not} of
Voronoi's first
kind.
For instance, we have:

\begin{proposition}\label{prop2.2}
For any integer $n\geq 4$, the root lattice $D_n$ is not of Voronoi's
first kind.
\end{proposition}

Surprisingly, this statement does not seem to appear in the
literature when $n\geq 5$.

{\small
\proof  We shall allow ourselves to leave  a few computational
details as exercises for the reader.

Recall that $D_n$ is defined by the lattice
\[
\bigl\{(x_i)_{1\leq i\leq n}\in \Z^n \big\vert \sum^n_{i=1}\,x_i \in
2\, \Z\bigr\}
\]
in $D_{n,\R} = \R^n$ equipped with the standard euclidean norm.

The dual lattice $D^\lor_n$ in $D^\lor_{n,\R} \simeq \R^n$ is easily
seen to be
$\Z^n \cup \bigl((\frac{1}{2})^{× n} + \Z^n\bigr)$. The
automorphism
group of $D_n$ contains (indeed, when $n\geq 5$, is equal to) the
semi-direct product
\[G:= \{± 1\}^n \rtimes \FS_n,\]
where $\{± 1\}^n$ acts diagonal, and $\FS_n$ by permutation of the
coordinates.
By considering the commutants of the action of $G$ on $D_{n,\R}$ and
on
$S^2 D^\lor_{n,\R}$, it is straightforward to check that the action
of $G$ on $D_{n,\R}$
is irreducible, and that there are precisely four non-zero
$G$-invariant subspaces
of $S^2 D^\lor_{n,\R}$ containing the standard euclidean form
$\sum\limits^n_{i=1} X^2_i$,
namely
\[
\begin{array}{lrcl}
& V_1 & := & \R\cdot \sum\limits^n_{i=1} X^2_i, \\
& V_2 & := & \bigoplus\limits^n_{i=1} \R\cdot  X^2_i, \\
& V_3 & := & \R\cdot \sum\limits^n_{i=1} X^2_i \oplus
\bigoplus\limits_{1\leq i< j\leq n} \R\cdot X_i X_j,\\
& V_4 & := & S^2 D^\lor_{n,\R}\,.
\end{array}
\]

Besides, by considering the vectors of minimal positive length in
$D_n$
--- namely the
vectors $(± e_i ± e_j),\ 1\leq i< j\leq n$, where $(e_1,\ldots ,\,
e_n)$ denotes
the canonical basis of $\R^n$ --- one sees that the euclidean lattice
$D_n$ is
indecomposable
\footnote{
Observe that, if an euclidean lattice $\overline E$ may be written as
an orthogonal
direct sum $\overline E_1 \oplus \overline E_2$, then any vector
$v\in E$ of minimal
positive length belongs to $E_1\cup E_2$.
Consequently, if the set $M$ of vectors of minimal positive length of
some euclidean
lattice $\overline E$ generates the $\R$-vector space $E_\R$ and
cannot be
partitioned
as $M=M_1 {\textstyle{\coprod}} M_2$ where any two vectors $e_1\in
M_1$ and $e_2 \in M_2$
are orthogonal, then $\overline E$ is indecomposable.}.

Let us assume that $D_n$ is of Voronoi's first kind.
To derive a contradiction, consider an adapted basis
$(\xi_1,\ldots,\, \xi_n)$ of $D^\lor_n$,
and apply Proposition \ref{prop2.1}, of which we now use the notation,
and the subsequent observations.
Consider in particular the $\R$-linear span $V$ of
$\bigl((\xi_i-\xi_j)^2\bigr)_{(i,j)\in S}$ in $S^2 D^\lor_{n,\R}$:
it is necessarily one of the spaces $V_i, \ 1\leq i\leq 4$ introduced
above.

Observe that $V$ cannot be $V_1$, which contains no quadratic form of
rank 1.

Assume now that $V$ is $V_2$. Then the action of $\{± 1\}^n \
(\subset G)$
on $V$ is trivial, and the quadratic forms $(\xi_i-\xi_j),\,(i,j) \in
S$, are fixed by this action.
This implies that each linear form $\xi_i-\xi_j,\,(i,j)\in S$, is
$± X_k$ for some $k\in \{1,\ldots,\, n\}$.
Lemma \ref{abs1.1} now shows that
\[
D^\lor_n  =  \sum_{(i,j)\in S} \Z\, (\xi_i-\xi_j)
  \subset  \bigoplus^n_{k=1}\, \Z \,X_k.
\]
This contradicts the fact that $D^\lor_n$ contains $\frac{1}{2}
\,\sum\limits^n_{k=1} X_k$.

It is straightforward to check that any quadratic form of rank 1 in
$V_3$ may be written
$\lambda \bigl(\sum\limits^n_{i=1} \varepsilon_k X_k\bigr)^2$ where
$\lambda \in \R^*$
and $(\varepsilon_k)_{1\leq k\leq n} \in \{± 1\}^n$. Consequently, if
$V=V_3$, then the
linear forms $\xi_i-\xi_j,\, (i,j)\in S$, may be written either
\begin{equation}
\label{gl2.2}
\sum^n_{k=1} \varepsilon_k X_k, \mbox{ where } (\varepsilon_k) \in
\{± 1\}^n \setminus \bigl\{ 1^{× n},\,(-1)^{× n}\bigr\}
\end{equation}
or
\begin{equation}
\label{gl2.3}
\frac{\varepsilon}{2}\,\sum^n_{k=1}\ X_k, \mbox{ where } \varepsilon
\in \{± 1\}.
\end{equation}
The sum of any two elements in $D^\lor_n$ of the form (\ref{gl2.2})
or
(\ref{gl2.3}) is never of the form (\ref{gl2.2}) or (\ref{gl2.3}).
This shows that, when $V=V_3$, there is no triple $(i,j,k), \,0\leq i<
j< k\leq n$,
such that $(i,j),\, (j,k)$ and $(i,k)$ belong to $S$.
In other words, the graph $\gamma(S)$ attached to $S$ has no cycle of
length 3.
An elementary counting argument shows that this contradicts the fact
that this
graph has $n+1$ vertices and
\[
\dim_\R\,V_3 = \frac{n(n+1)}{2} - (n-1)
\]
edges.

Finally, when $V=V_4$, then
\[
\vert S \vert = \dim\,V_4 = \frac{n(n+1)}{2},
\]
and consequently
\[
S= \bigl\{(i,j), \  0\leq i< j\leq n\bigr\}
\]
and, for every $x\in D_{n,\R},$
\[
\Vert x\Vert^2 = \lambda \sum_{0\leq i< j\leq n} \,\bigl(\xi_i(x) -
\xi_j(x)\bigr)^2.
\]
This implies that, up to scaling, the euclidean lattice $D_n$ is
isometric with
$A^*_n$ --- this is plainly wrong (compare the cardinality of their
sets of vectors
of minimal length, or of their automorphism groups).\qed
}

\subsection{The Voronoi cell of an euclidean lattice with strictly
obtuse superbase}\label{subsec.voronoi}

Let $n$ be a positive integer.
To any $n(n+1)/2$-tuple
\[
p=(p_{ij})_{0\leq i< j\leq n} \in (\R^*_+)^{\frac{n(n+1)}{2}},
\]
we attach an euclidean lattice $\overline V(p)$ of rank $n$ defined
as follows:
\[
\overline V(p) := (V, \langle \cdot , \cdot \rangle_p),
\]
where
\[
V:= \Z^{n+1}/\Z\cdot 1^{× (n+1)}
\]
and where $\langle \cdot,\cdot\rangle_p$ denotes the euclidean scalar
product on
$V_\R\simeq \R^{n+1}/\R\cdot 1^{× (n+1)}$ defined by
\[
\Big\langle\bigl[(x_i)_{0\leq i\leq n}\bigr],\bigl[(y_i)_{0\leq i\leq
n}\bigr]\Big\rangle_p
= \sum_{0\leq i< j\leq n} \,p_{ij} (x_i-x_j)\cdot(y_i-y_j),
\]
for any $(x_i)_{0\leq i\leq n}$ and $(y_i)_{0\leq i\leq n}$ in
$\R^{n+1}$
(where
$[\,\alpha]$ denotes the image in $V_\R$ of $\alpha \in \R^{n+1}$).

In the sequel, we shall often omit the subscript $p$ to simplify
notations.
In particular,
we shall write $v\cdot w$ (resp. $\Vert v\Vert^2$) instead of
$\langle v,w\rangle_p$ (resp. $\langle v,v\rangle_p$).

Let $(\varepsilon_0,\ldots,\, \varepsilon_n)$ be the canonical basis
of $\Z^{n+1}$.
Clearly
\[
\varepsilon_0 + \cdots + \varepsilon_n = 1^{× (n+1)}
\]
and
\[
(v_0,\ldots,\,v_n) :=
\bigl([\varepsilon_0],\ldots,\,[\varepsilon_n]\bigr)
\]
is a superbase of $\overline V$, which is strictly obtuse since
\begin{equation}
\label{gl3.1}
v_i\cdot v_j = - p_{ij} \ \mbox{ if } \ 0\leq i< j\leq n.
\end{equation}
It is convenient to define
\[p_{ij} := p_{ji} \ \mbox{ if } \ 0\leq j < i\leq n.\]
Then (\ref{gl3.1}) holds for any $(i,j) \in \{0,\ldots,\, n\}^2$ such
that $i \not= j$, and, for any $i\in \{0,\ldots,\, n\}$
\begin{equation}
\label{gl3.2}
\Vert v_i\Vert^2 = - v_i\cdot
\sum_{\genfrac{}{}{0pt}{}{0\leq j\leq n}{j\not= i}} \,v_j
= \sum_{\genfrac{}{}{0pt}{}{0\leq j\leq n}{j\not= i}}\, p_{ij}.
\end{equation}

More generally, for every $S\subset \{0,\ldots,\, n\}$ we let
\[
v_S:= \sum_{i\in S} \,v_i.
\]
Observe that $v_S = 0$ iff $S=\emptyset$ or $S=\{0,\ldots,\, n\}$.
Moreover, for any two elements $S_1$ and $S_2$ in
\[
\FP'\bigl(\{0,\ldots,\, n\}\bigr) := \FP \bigl(\{0,\ldots,\
n\}\bigr)\setminus \bigl\{\emptyset, \{0,\ldots,\,n\}\bigr\},
\]
we have
\[
v_{S_1} = v_{S_2} \ \mbox{ iff } \ S_1 =S_2,
\]
and
\[
v_{S_1} + v_{S_2} = 0 \ \mbox{ iff } \ \{0,\ldots,\, n\} = S_1
\,{\textstyle{\coprod}}\, S_2.
\]
A straightforward application of Selling's formula (\ref{gl1.2})
shows that the
$2^{n+1} - 2$ vectors $v_S, \ S\in \FP' \bigl(\{0,\ldots,\,
n\}\bigr)$, are
precisely the Voronoi vectors of $\overline V(p)$, and indeed are
strict Voronoi vectors
(\cite[Theorem 3]{conwaysloane92}).
For any such $S$, let
\begin{eqnarray*}
H_S & := & \Bigl\{ x \in V_\R \,\big\vert  \ \Vert x-v_S \Vert =
\Vert x\Vert \Bigr\} \\
& \, = & \bigl\{x \in V_\R\,\vert\, 2 v_S \cdot x = \Vert
v_S\Vert^2\bigr\}
\end{eqnarray*}
and let
\[
F_S := \mathcal{V}\bigl(\overline V(p)\bigr) \cap H_S
\]
be the corresponding face of $\mathcal{V}\bigl(\overline V(p)\bigr)$.
The $F_S, \ S\in \FP'\bigl(\{0,\ldots,\, n\}\bigr)$ are precisely the
$(n-1)$-dimensional
faces (also known as {\it facets}) of $\mathcal{V}\,\bigl(\overline
V(p)\bigr)$.

\begin{lemma}\label{abs3.1}
If the faces $F_S$ and $F_{S'}$ attached to two elements $S$ and $S'$
in
$\FP'\bigl(\{0,\ldots,\,n\}\bigr)$ are not disjoint, then $S\subset
S'$ or $S' \subset S.$
\end{lemma}

\proof
Define $I:= S\cap S',\ T:= S\setminus I$, and $T' := S'\setminus I$.
Then $S$ (resp. $S'$) is the disjoint union of $I$ and $T$ (resp. $I$
and $T'$).
Moreover, if $x$ is some element of
\[F_S \cap F_{S'} = H_S \cap H_{S'} \cap \mathcal{V}\bigl(\overline
V(p)\bigr),\]
then we have
\[\begin{array}{lllcl}
   & & 2\, v_S \cdot x & = & \Vert v_S\Vert^2, \\
&  & 2\, v_{S'} \cdot x & = & \Vert v_{S'} \Vert^2, \\
&  & 2 \,v_{S\cup S'} \cdot x & \leq & \Vert v_{S\cup S'}\Vert^2, \\
  \mbox{and}  && & & \\
  & & 2\, v_I \cdot x & \leq & \Vert v_I\Vert^2.
\end{array}\]
As
\[v_S + v_{S'} = v_{S\cup S'} + v_I,\]
this implies:
\[
\Vert v_S\Vert^2 + \Vert v_{S'}\Vert^2 \leq \Vert v_{S\cup S'} \Vert^2
+ \Vert v_I\Vert^2.\]
However, we have
\begin{eqnarray*}
\Vert v_{S\cup S'} \Vert^2 + \Vert v_I\Vert - \Vert v_S\Vert^2 -
\Vert v_{S'}\Vert^2 & = & \Vert v_T+ v_{T'} + v_I\Vert^2 + \Vert v_I
\Vert^2 \\
& & - \Vert v_T + v_I \Vert^2 - \Vert v_{T'} + v_I\Vert^2 \\
& = & 2\, \langle v_T, v_{T'} \rangle \\
& = & - 2 \sum\limits_{\genfrac{}{}{0pt}{}{i\in T}{j\in T'}} \,
p_{ij},
\end{eqnarray*}
and the last sum is negative when $T$ and $T'$ are not empty.
\qed

Let us consider the set $\mathcal{S}(n)$ of subsets of cardinality
$n$ of
$\FP'\bigl(\{0,\ldots ,\, n\}\bigr)$ which are totally ordered by
inclusion.
The group $\FS_{n+1}$ of permutations of $\{0,\ldots,\, n\}$ acts
naturally on $\mathcal{S}(n)$,
and the following lemma is straightforward:

\begin{lemma}\label{abs3.2}
The action of $\FS_{n+1}$ on $\mathcal{S}(n)$ is simply transitive.
In other words, the mapping
\[
\begin{array}{ccl}
\FS_{n+1} & \longrightarrow & \mathcal{S}(n) \\
\sigma & \longmapsto & \Bigl\{\bigl\{\sigma(1),\ldots,\,
\sigma(n)\bigr\},\
\bigl\{\sigma(2),\ldots,\,\sigma(n)\bigr\},\ldots,\,
\bigl\{\sigma(n)\bigr\}\Bigr\}
\end{array}
\]
is a bijection.
\end{lemma}

For any $A\in \mathcal{S}(n)$, the vectors $(v_S)_{S\in A}$ are
linearly independent,
and consequently the hyperplanes $(H_S)_{S\in A}$ have an unique
common point,
which we shall denote $s_A$.

\begin{lemma}\label{abs3.3}
If $A= \bigl\{\{i_1,\ldots,\, i_n\},\,\{i_2,\ldots,\, i_n\},
\ldots,\, \{i_n\}\bigr\}$,
then
\begin{equation}
\label{gl3.3}
\Vert \sum^n_{k=0}\, x_k v_{i_k} - s_A\Vert^2 - \Vert s_A\Vert^2 =
\sum_{0\leq k< \ell \leq n} \,p_{i_k i_\ell}\, \bigl[(x_\ell-x_k)^2 -
(x_\ell - x_k)\bigr].
\end{equation}
holds for any element $(x_0,\ldots,\, x_n) \in \R^{n+1}$.
\end{lemma}

\proof
The point $s_A$ is defined by the $n$ linear equations
\[
2 \,v_S \cdot s_A = \Vert v_S\Vert^2\,,\ S\in A,
\]
which may also be written:
\[
2\, v_{i_k} \cdot s_A= \Vert v_{i_k} + \cdots + v_{i_n} \Vert^2 -
\Vert v_{i_{k+1}} + \cdots + v_{i_n}\Vert^2, \ 1\leq k\leq n-1
\]
and
\[
2 \,v_{i_n}\cdot s_A = \Vert v_{i_n}\Vert^2.
\]
Using (\ref{gl3.1}) and (\ref{gl3.2}), these relations take finally
the form
\[
2 \,v_{i_k} \cdot s_A= \sum^n_{\ell=0} \, \varepsilon(\ell,k)\,
p_{i_k i_\ell},\ 1\leq k\leq n,
\]
where
\[
\varepsilon(\ell,k)  := \left\{
\begin{array}{lcr}
1  & \mbox{if}& \ell < k \\
0  & \mbox{if}& \ell=k \\
-1 & \mbox{if}& \ell > k.
\end{array}
\right.
\]
Consequently, for any $(x_0,\ldots,\, x_n) \in \R^{n+1},$
\begin{eqnarray*}
\Vert \sum^n_{k=0} x_k v_{i_k} - s_A\Vert^2-\Vert s_A\Vert^2 & =
&\Vert \sum^n_{k=0} x_k v_{i_k}\Vert^2 - 2 \sum^n_{k=0} x_k
v_{i_k}\cdot s_A \\
& = & \sum_{0\leq k< \ell\leq n}\, p_{i_ki_\ell}\,(x_k-x_\ell)^2 \\
& & -\sum_{0\leq k,\ell\leq n}\,\varepsilon(\ell,k)\, p_{i_k i_\ell}\,
x_k \\
& = & \sum_{0\leq k<\ell\leq n}\, p_{i_k i_\ell}\,
\bigl[(x_\ell-x_k)^2
+ x_k-x_\ell\bigr].
\end{eqnarray*}
\qed

When $(x_0,\ldots,\,x_n) \in \Z^{n+1}$, the right-hand side of
(\ref{gl3.3}) is non-negative.
This shows that, for any $v\in V$,
\begin{equation}
\label{gl3.4}
\Vert v- s_A\Vert^2 \geq \Vert s_A\Vert^2,
\end{equation}
or equivalently, that $s_A$ belongs to the Voronoi cell
$\mathcal{V}\,\bigl(V(p)\bigr)$ of $\overline V(p)$.

We are now in position to establish:

\begin{proposition}\label{prop3.4}
1) The mapping $(A\mapsto s_A)$ is a bijection from  $\mathcal{S}(n)$
onto
the set of vertices of $\mathcal{V}\,\bigl(\overline V(p)\bigr)$.

2) More generally, one defines an inclusion reversing bijection
from the set $\mathcal{O}(n)$ of non-empty subsets of
$\FP'\bigl(\{0,\ldots,\,n\}\bigr)$
which are totally ordered onto the set of faces of
$\mathcal{V}\,\bigl(\overline V(p)\bigr)$ by mapping $P\in
\mathcal{O}(n)$ to
$\Delta_P := \bigcap\limits_{S\in P}\, F_S$.
\end{proposition}

Observe that $\Delta_P$ is a face of dimension $n-\vert P\vert$ of
$\mathcal{V}\,\bigl(\overline V(p)\bigr)$.

When all the $p_{ij}$'s are equal to 1, the permutation group
$\FS_{n+1}$ acts
on $\overline V(p)$ by permutation of the vectors in the superbase
$(v_0,\ldots,\,v_n)$, and the vertices of
$\mathcal{V}\,\bigl(\overline V(p)\bigr)$
are the $(n+1)!$ points in the orbit under $\FS_{n+1}$ of
\[
\frac{1}{n+1}\,\Bigl[(-\frac{n}{2},\,-\frac{n}{2} +
1,\ldots,\,\frac{n}{2}-1,\,\frac{n}{2})\Bigr]
  = \frac{1}{n+1} \,\sum\limits^n_{i=0} (-\frac{n}{2}+i)\,v_i =
\frac{1}{n+1}\,\sum\limits^n_{i=0} i\cdot v_i.
\]
The Voronoi cell $\mathcal{V}\,\bigl(\overline V(p)\bigr)$ is then a
so-called {\it permutohedron}.
Actually in this case, $\overline V(p)$ is isometric (up to a
scaling) with $A^*_n$,
and this description of its Voronoi cell is classical (see for
instance
\cite[102 - 103]{voronoi09}, and 
\cite[Chapter 21, Theorem 7]{conwaysloane99}).

Proposition \ref{prop3.4} shows that the Voronoi cell of any lattice
of Voronoi's first kind
with positive Selling parameters still has the combinatorial type of
a permutohedron.
Actually, this is also a consequence of results in Voronoi's last
paper
(see \cite[104]{voronoi09}, notably p. 147). Let us finally point out
that,
when all the $p_{ij}$'s are 1, identity (\ref{gl3.3}) also appears in
this paper
(up to a permutation of the variables, it coincides with the
penultimate
equation in p. 140 of \cite[103]{voronoi09}).

\proof
We freely use basic facts concerning polytopes and their posets of
faces, as
described for instance in \cite[chapter 3]{grunbaum67}, and
\cite[chapter 2]{ziegler95}.

1) We have just shown in (\ref{gl3.4}) that, for any $A\in
\mathcal{S}(n)$,
the point $s_A$ belongs to $\mathcal{V}\,\bigl(\overline V(p)\bigr)$.
Moreover, it
is an extreme point of $\mathcal{V}\,\bigl(\overline V(p)\bigr)$,
since it is
the ``vertex'' of the intersection of $n$ ''half spaces''
\[
\bigcap_{S\in A} \bigl\{x \in V(p)_\R\,\vert \, 2 \,v_S\cdot x\leq
\Vert v_S\Vert^2\bigr\}
\]
which contains $\mathcal{V}\,\bigl(\overline V(p)\bigr)$.
(Observe that the vectors $(v_S)_{S\in A}$ are linearly independent.)

Conversely, any vertex $P$ of $\mathcal{V}\,\bigl(\overline
V(p)\bigr)$ is contained
in at least $n$ distinct facets $F_{S_1},\ldots,\, F_{S_n}$ of
$\mathcal{V}\,\bigl(\overline V(p)\bigr)$. Lemma \ref{abs3.1} shows
that
$A:= \{F_{S_1},\ldots,\, F_{S_n}\}$ is a totally ordered subset of
$\FS'\bigl(\{0,\ldots,\,n\}\bigr)$. Consequently, $P= s_A$.

2)
For any $P$ in $\mathcal{O}(n)$, the intersection
$\Delta_P:= \bigcap_{S\in P}\, F_S$ is not empty, since it  exists
some element $A$ of $\mathcal{S}(n)$ containing $P$ and therefore
$\Delta_P$ contains $s_A$.
Moreover, the vectors $(v_S)_{S\in P}$, orthogonal to the facets
$(F_S)_{S\in
P},$
are linearly independent, and therefore $\Delta_P$ is a face of
$\mathcal{V}\,\bigl(\overline V(p)\bigr)$ of dimension $n-\vert
P\vert$.

Conversely any face of dimension $n-p$ of
$\mathcal{V}\,\bigl(\overline V(p)\bigr)$ is
an intersection of $p$ distinct facets of
$\mathcal{V}\,\bigl(\overline V(p)\bigr)$ and
may therefore be written $\bigcap\limits_{S\in P}\,F_S$ where $P$ is
a subset of
$\FS'\bigl(\{0,\ldots,\,n\}\bigr)$ of cardinality $p$. Again Lemma
\ref{abs3.1} shows
that $P$ is totally ordered by inclusion.

This establishes that the map
\[\begin{array}{ccl}
\mathcal{O}(n) & \longrightarrow & \Bigl\{\mbox{faces of }
\mathcal{V}\bigl(\overline V(p)\bigr)\Bigr\} \\
P & \longmapsto & \Delta_P
\end{array}\]
is bijective. It is clearly inclusion reversing.
\qed
}

\end{appendix}

\bibliographystyle{alpha}

\bibliography{hvar-arXive}

\begin{thebibliography}{BBM82}

\bibitem[Bae34]{baer34}
R.~Baer.
\newblock {Erweiterung von Gruppen und ihren Isomorphismen.}
\newblock {\em Math. Z.}, 38:375--416, 1934.

\bibitem[Ban93]{banaszczyk93}
W.~Banaszczyk.
\newblock New bounds in some transference theorems in the geometry of numbers.
\newblock {\em Math. Ann.}, 296(4):625--635, 1993.

\bibitem[Ban95]{banaszczyk95}
W.~Banaszczyk.
\newblock Inequalities for convex bodies and polar reciprocal lattices in
  {${\bf R}\sp n$}.
\newblock {\em Discrete Comput. Geom.}, 13(2):217--231, 1995.

\bibitem[BBM82]{berthelotbreenmessing82}
P.~Berthelot, L.~Breen, and W.~Messing.
\newblock {\em Th\'eorie de {D}ieudonn\'e cristalline. {II}}, volume 930 of
  {\em Lecture Notes in Mathematics}.
\newblock Springer-Verlag, Berlin, 1982.

\bibitem[BGS94]{bostetal94}
J.-B. Bost, H.~Gillet, and C.~Soul{\'e}.
\newblock Heights of projective varieties and positive {G}reen forms.
\newblock {\em J. Amer. Math. Soc.}, 7(4):903--1027, 1994.

\bibitem[BK]{bostkuennemann2}
J.-B. Bost and K.~K{\"u}nnemann.
\newblock Hermitian vector bundles and extension groups on arithmetic schemes.
  {I}{I}.
\newblock In preparation.

\bibitem[Bor69]{borel69}
Armand Borel.
\newblock {\em {Introduction aux groupes arithm{\'e}tiques}}.
\newblock {Paris: Hermann \& Cie. 125 p. }, 1969.

\bibitem[Bor05]{borek05}
T.~Borek.
\newblock Successive minima and slopes of {H}ermitian vector bundles over
  number fields.
\newblock {\em J. Number Theory}, 113(2):380--388, 2005.

\bibitem[Bos96]{bost96a}
J.-B. Bost.
\newblock P\'eriodes et isogenies des vari\'et\'es ab\'eliennes sur les corps
  de nombres (d'apr\`es {D}. {M}asser et {G}. {W}\"ustholz).
\newblock {\em Ast\'erisque}, (237):Exp.\ No.\ 795, 4, 115--161, 1996.
\newblock S\'eminaire Bourbaki, Vol.\ 1994/95.

\bibitem[Cas67]{Casselsetal67}
{\em Algebraic number theory}.
\newblock Proceedings of an instructional conference organized by the London
  Mathematical Society. Edited by J. W. S. Cassels and A. Fr\"ohlich. Academic
  Press, London, 1967.

\bibitem[CLT01]{chamberttschinkel01}
A.~Chambert-Loir and Y.~Tschinkel.
\newblock Torseurs arithm\'etiques et espaces fibr\'es.
\newblock In {\em Rational points on algebraic varieties}, volume 199 of {\em
  Progr. Math.}, pages 37--70. Birkh\"auser, Basel, 2001.

\bibitem[Con00]{conrad00}
B.~Conrad.
\newblock {\em Grothendieck duality and base change}, volume 1750 of {\em
  Lecture Notes in Mathematics}.
\newblock Springer-Verlag, Berlin, 2000.

\bibitem[Cox51]{coxeter51}
H.~S.~M. Coxeter.
\newblock Extreme forms.
\newblock {\em Canadian J. Math.}, 3:391--441, 1951.

\bibitem[CS91]{conwaysloane91}
J.H. Conway and N.J.A. Sloane.
\newblock The cell structures of certain lattices.
\newblock In {\em Miscellanea mathematica, Festschr. H. G\"otze}, pages
  71--107. Springer-Verlag, Berlin, 1991.

\bibitem[CS92]{conwaysloane92}
J.~H. Conway and N.~J.~A. Sloane.
\newblock Low-dimensional lattices. {VI}. {V}orono\u\i\ reduction of
  three-dimensional lattices.
\newblock {\em Proc. Roy. Soc. London Ser. A}, 436(1896):55--68, 1992.

\bibitem[CS99]{conwaysloane99}
J.H. Conway and N.J.A. Sloane.
\newblock {\em {Sphere packings, lattices and groups. With additional
  contributions by E. Bannai, R. E. Borcherds, J. Leech, S. P. Norton, A. M.
  Odlyzko, R. A. Parker, L. Queen and B. B. Venkov. 3rd ed.}}
\newblock {Grundlehren der Mathematischen Wissenschaften. 290. New York, NY:
  Springer. lxxiv, 703 p.}, 1999.

\bibitem[Del70]{deligne70}
P.~Deligne.
\newblock {\em \'{E}quations diff\'erentielles \`a points singuliers
  r\'eguliers}.
\newblock Springer-Verlag, Berlin, 1970.
\newblock Lecture Notes in Mathematics, Vol. 163.

\bibitem[Del73]{deligne73}
P.~Deligne.
\newblock Expos\'e {X}{V}{I}{I} : Cohomologie \`a support propre.
\newblock In {\em Th\'eorie des topos et cohomologie \'etale des sch\'emas.
  {T}ome 3 -- SGA 4}, pages 250--461. Springer-Verlag, Berlin, 1973.
\newblock Lecture Notes in Mathematics, Vol. 305.

\bibitem[EM42]{eilenbergmaclane42}
S.~Eilenberg and S.~MacLane.
\newblock Group extensions and homology.
\newblock {\em Ann. Math.}, 43:757--831, 1942.

\bibitem[For38]{ford38}
L.R. Ford.
\newblock {Fractions.}
\newblock {\em Am. Math. Mon.}, 45:586--601, 1938.

\bibitem[GH78]{griffithsharris78}
P.~Griffiths and J.~Harris.
\newblock {\em Principles of algebraic geometry}.
\newblock Wiley-Interscience [John Wiley \& Sons], New York, 1978.
\newblock Pure and Applied Mathematics.

\bibitem[Gra84]{grayson84}
D.~R. Grayson.
\newblock Reduction theory using semistability.
\newblock {\em Comment. Math. Helv.}, 59(4):600--634, 1984.

\bibitem[Gra01]{graftieaux01}
P.~Graftieaux.
\newblock {Formal groups and the isogeny theorem.}
\newblock {\em Duke Math. J.}, 106(1):81--121, 2001.

\bibitem[Gri66]{griffiths66}
P.A. Griffiths.
\newblock {The extension problem in complex analysis. II: Embeddings with
  positive normal bundle}.
\newblock {\em Am. J. Math.}, 88:366--446, 1966.

\bibitem[Gri69]{griffiths69}
P.A. Griffiths.
\newblock {Hermitian differential geometry, Chern classes, and positive vector
  bundles}.
\newblock In {\em {Global Analysis, Papers in Honor of K. Kodaira 185-251 }}.
  1969.

\bibitem[Gro66]{grothendieck66}
A.~Grothendieck.
\newblock On the de {R}ham cohomology of algebraic varieties.
\newblock {\em Inst. Hautes \'Etudes Sci. Publ. Math.}, (29):95--103, 1966.

\bibitem[Gr{\"u}67]{grunbaum67}
B.~Gr{\"u}nbaum.
\newblock {\em Convex polytopes}.
\newblock With the cooperation of Victor Klee, M. A. Perles and G. C. Shephard.
  Pure and Applied Mathematics, Vol. 16. Interscience Publishers John Wiley \&
  Sons, Inc., New York, 1967.

\bibitem[GS90]{gilletsoule90}
H.~Gillet and C.~Soul{\'e}.
\newblock Arithmetic intersection theory.
\newblock {\em Inst. Hautes \'Etudes Sci. Publ. Math.}, 72:93--174, 1990.

\bibitem[GS91]{gilletsoule91}
H.~Gillet and C.~Soul{\'e}.
\newblock On the number of lattice points in convex symmetric bodies and their
  duals.
\newblock {\em Israel J. Math.}, 74(2-3):347--357, 1991.

\bibitem[GS92]{gilletsoule92}
H.~Gillet and C.~Soul{\'e}.
\newblock An arithmetic {R}iemann-{R}och theorem.
\newblock {\em Invent. Math.}, 110(3):473--543, 1992.

\bibitem[Har66]{hartshorne66}
R.~Hartshorne.
\newblock {\em Residues and duality}.
\newblock Lecture notes of a seminar on the work of A. Grothendieck, given at
  Harvard 1963/64. With an appendix by P. Deligne. Lecture Notes in
  Mathematics, No. 20. Springer-Verlag, Berlin, 1966.

\bibitem[Har77]{hartshorne77}
R.~Hartshorne.
\newblock {\em Algebraic geometry}.
\newblock Springer-Verlag, New York, 1977.
\newblock Graduate Texts in Mathematics, No. 52.

\bibitem[Ill71]{illusie71a}
L.~Illusie.
\newblock Expos\'e {I}{I} : Existence de r\'esolutions globales.
\newblock In {\em Th\'eorie des intersections et th\'eor\`emes de Riemann-Roch
  -- SGA 6}, pages 160--221. Springer-Verlag, Berlin, 1971.
\newblock Lecture Notes in Mathematics, Vol. 225.

\bibitem[Ill96]{illusie96}
L.~Illusie.
\newblock Frobenius et d\'eg\' en\'erescence de {H}odge.
\newblock In Jos{\'e} Bertin, Jean-Pierre Demailly, Luc Illusie, and Chris
  Peters, editors, {\em Introduction \`a la th\'eorie de {H}odge}, volume~3 of
  {\em Panoramas et Synth\`eses}. Soci\'et\'e Math\'ematique de France, Paris,
  1996.

\bibitem[Kit93]{kitaoka93}
Y.~Kitaoka.
\newblock {\em {Arithmetic of quadratic forms.}}
\newblock {Cambridge Tracts in Mathematics. 106. Cambridge: Cambridge
  University Press. x, 268 p.}, 1993.

\bibitem[Lan88]{lang88}
S.~Lang.
\newblock {\em Introduction to {A}rakelov theory}.
\newblock Springer-Verlag, New York, 1988.

\bibitem[LLS90]{lagariasetal90}
J.~C. Lagarias, H.~W. Lenstra, Jr., and C.-P. Schnorr.
\newblock Korkin-{Z}olotarev bases and successive minima of a lattice and its
  reciprocal lattice.
\newblock {\em Combinatorica}, 10(4):333--348, 1990.

\bibitem[Mar03]{martinet03}
J.~Martinet.
\newblock {\em Perfect lattices in {E}uclidean spaces}, volume 327 of {\em
  Grundlehren der Mathematischen Wissenschaften}.
\newblock Springer-Verlag, Berlin, 2003.

\bibitem[MH80]{Milnor73}
J.~Milnor and D.~Husemoller.
\newblock {\em Symmetric bilinear forms}, volume~73 of {\em Ergebnisse der
  Mathematik und ihrer Grenzgebiete}.
\newblock Springer-Verlag, Berlin-Heidelberg-New York, 1980.

\bibitem[ML95]{maclane95}
S.~Mac~Lane.
\newblock {\em Homology}.
\newblock Classics in Mathematics. Springer-Verlag, Berlin, 1995.

\bibitem[Moc96]{mochizuki96}
S.~Mochizuki.
\newblock A theory of ordinary {$p$}-adic curves.
\newblock {\em Publ. Res. Inst. Math. Sci.}, 32(6):957--1152, 1996.

\bibitem[Moc99]{mochizuki99p}
S.~Mochizuki.
\newblock The {H}odge-{A}rakelov theory of elliptic curves: {G}lobal
  discretization of local {H}odge theories.
\newblock RIMS Preprint Nos. 1255, 1256, (October 1999).

\bibitem[Neu99]{neukirch99}
J.~Neukirch.
\newblock {\em Algebraic number theory}, volume 322 of {\em Grundlehren der
  Mathematischen Wissenschaften}.
\newblock Springer-Verlag, Berlin, 1999.
\newblock Translated from the 1992 German original and with a note by Norbert
  Schappacher, With a foreword by G. Harder.

\bibitem[Rad64]{rademacher64}
H.~Rademacher.
\newblock {\em {Lectures on elementary number theory}}.
\newblock {A Blaisdell Book in the Pure and Applied Sciences. Introduction to
  Higher Mathematics. New York-Toronto-London: Blaisdell Publishing Company, a
  division of Ginn and Company. IX, 146 p. }, 1964.

\bibitem[Sel74]{selling74}
E.~Selling.
\newblock \"uber die bin\"aren und tern\"aren quadratischen formen.
\newblock {\em J. reine angew. Math.}, 77:143--229, 1874.

\bibitem[SGA03]{SGA1}
{\em Rev\^etements \'etales et groupe fondamental ({SGA} 1)}.
\newblock Documents Math\'ematiques (Paris), 3. Soci\'et\'e Math\'ematique de
  France, Paris, 2003.
\newblock S\'eminaire de g\'eom\'etrie alg\'ebrique du Bois Marie 1960--61.,
  Directed by A. Grothendieck, With two papers by M. Raynaud, Updated and
  annotated reprint of the 1971 original [Lecture Notes in Math., 224,
  Springer, Berlin;].

\bibitem[Sou97]{soule97}
C.~Soul{\'e}.
\newblock Hermitian vector bundles on arithmetic varieties.
\newblock In {\em Algebraic geometry---Santa Cruz 1995}, volume~62 of {\em
  Proc. Sympos. Pure Math.}, pages 383--419. Amer. Math. Soc., Providence, RI,
  1997.

\bibitem[Stu76]{stuhler76}
U.~Stuhler.
\newblock Eine {B}emerkung zur {R}eduktionstheorie quadratischer {F}ormen.
\newblock {\em Arch. Math. (Basel)}, 27(6):604--610, 1976.

\bibitem[Szp85]{szpiro85}
L.~Szpiro.
\newblock Degr\'es, intersections, hauteurs.
\newblock In {\em S\'eminaire sur les pinceaux arithm\'etiques: la conjecture
  de Mordell}, volume 127 of {\em Ast\'erisque}, pages 11--28. Soc. Math.
  France, Paris, 1985.

\bibitem[Tou72]{tougeron72}
J.-C. Tougeron.
\newblock {\em Id\'eaux de fonctions diff\'erentiables}.
\newblock Springer-Verlag, Berlin, 1972.
\newblock Ergebnisse der Mathematik und ihrer Grenzgebiete, Band 71.

\bibitem[TT90]{thomason90}
R.~W. Thomason and Th. Trobaugh.
\newblock Higher algebraic {$K$}-theory of schemes and of derived categories.
\newblock In {\em The Grothendieck Festschrift, Vol.\ III}, volume~88 of {\em
  Progr. Math.}, pages 247--435. Birkh\"auser Boston, Boston, MA, 1990.

\bibitem[vdW56]{waerden56}
B.L. van~der Waerden.
\newblock {Die Reduktionstheorie der positiven quadratischen Formen.}
\newblock {\em Acta Math.}, 96:265--309, 1956.

\bibitem[Vor08]{voronoi08a}
G.~Vorono\u\i.
\newblock Nouvelles applications des param\`etres continus \`a la th\'eorie des
  formes quadratiques. {I}. premier m\'emoire: Sur quelques propri\'et\'es des
  formes positives parfaites.
\newblock {\em J. reine angew. Math.}, 133:97--178, 1908.

\bibitem[Vor09]{voronoi09}
G.~Vorono\u\i.
\newblock Nouvelles applications des param\`etres continus \`a la th\'eorie des
  formes quadratiques. {III}. deuxi\`eme m\'emoire: Recherche sur les
  parall\'elo\`edres primitifs - seconde partie - domaines de formes
  quadratiques correspondant aux diff\'erents types de parall\'elo\`edres
  primitifs.
\newblock {\em J. reine angew. Math.}, 136:67--181, 1909.

\bibitem[Wei94]{weibel94}
C.~A. Weibel.
\newblock {\em An introduction to homological algebra}, volume~38 of {\em
  Cambridge Studies in Advanced Mathematics}.
\newblock Cambridge University Press, Cambridge, 1994.

\bibitem[Wit41]{witt41}
E.~Witt.
\newblock {Spiegelungsgruppen und Aufz\"ahlung halbeinfacher Liescher Ringe.}
\newblock {\em Abh. Math. Semin. Hansische Univ.}, 14:289--322, 1941.

\bibitem[Zie95]{ziegler95}
G.~Ziegler.
\newblock {\em Lectures on polytopes}, volume 152 of {\em Graduate Texts in
  Mathematics}.
\newblock Springer-Verlag, New York, 1995.

\end{thebibliography}

\end{document}